\documentclass[11pt, reqno]{amsart}
\usepackage{boxedminipage}

 \usepackage{amsmath}
    \usepackage{enumerate}
    \usepackage{amsfonts}
    \usepackage{amssymb}
    \usepackage{amsthm}
    \usepackage{bbm}
    \usepackage{xcolor}
    \usepackage{tabularx}
    \usepackage{graphicx}
    \usepackage{xr}
    \usepackage{xspace}
    \usepackage{thmtools}
    \usepackage{thm-restate}
    \usepackage{fourier}
    \usepackage[multiple]{footmisc}
    \usepackage{microtype}
    \usepackage[font=small,format=plain,labelfont=sf,up,textfont=sf,up]{caption}
    \usepackage[tmargin=1in,bmargin=1in,lmargin=1.24in,rmargin=1.24in]{geometry}
    \usepackage[export]{adjustbox}
    \usepackage{hyperref}
    
\definecolor{halfgray}{gray}{0.55}
\definecolor{webgreen}{rgb}{0,.5,0}
\definecolor{webbrown}{rgb}{.6,0,0}
\definecolor{Maroon}{cmyk}{0, 0.87, 0.68, 0.32}
\definecolor{RoyalBlue}{cmyk}{1, 0.50, 0, 0}
\definecolor{Black}{cmyk}{0, 0, 0, 0}
\definecolor{pinkish}{RGB}{255, 192, 203}

\hypersetup{%
    colorlinks=true, linktocpage=true, pdfstartpage=1, pdfstartview=FitV,%
    breaklinks=true, pdfpagemode=UseNone, pageanchor=true, pdfpagemode=UseOutlines,%
    plainpages=false, bookmarksnumbered, bookmarksopen=true, bookmarksopenlevel=1,%
    hypertexnames=true, pdfhighlight=/O,
    urlcolor=webbrown, linkcolor=RoyalBlue, citecolor=webgreen, 
}

\usepackage{todonotes}

\declaretheoremstyle[ 
    spaceabove=1em, spacebelow=1em,
    headfont=\scshape,
    notefont=\normalfont, notebraces={[}{]},
    bodyfont=\normalfont \itshape,
    postheadspace=0.5em,
    numbered=no,
                    ]{myStyle}

\theoremstyle{myStyle}

\theoremstyle{myStyle}
\declaretheorem[numberwithin=section]{theorem}
\declaretheorem[sibling=theorem, name=Lemma]{lemma}
\declaretheorem[sibling=theorem, name=Remark] {rk}
\declaretheorem[sibling=theorem, name=Definition]{defn}
\declaretheorem[sibling=theorem, name=Proposition]{prop}
\declaretheorem[sibling=theorem, name=Corollary]{corr}

\newtheorem{innercustomass}{Assumption}
\newenvironment{customass}[1]
  {\renewcommand\theinnercustomass{#1}\innercustomass}
  {\endinnercustomass}

\newtheorem{innercustomhypo}{Hypothesis}
\newenvironment{customhypo}[1]
    {\renewcommand\theinnercustomhypo{#1}\innercustomhypo}
    {\endinnercustomhypo}

\makeatletter
\def\@settitle{\begin{center}%
  \baselineskip14\p@\relax
  \bfseries
  \uppercasenonmath\@title
  \@title
  \ifx\@subtitle\@empty\else
     \\[1ex]\uppercasenonmath\@subtitle
     \footnotesize\mdseries\@subtitle
  \fi
  \end{center}%
}
\def\subtitle#1{\gdef\@subtitle{#1}}
\def\@subtitle{}
\makeatother

\newcommand{\bt}   {\begin{theorem}}
\newcommand{\et}   {\end  {theorem}}
\newcommand{\bl}   {\begin{lemma}}
\newcommand{\el}   {\end  {lemma}}
\newcommand{\bp}   {\begin{prop}}
\newcommand{\ep}   {\end  {prop}}
\newcommand{\bc}   {\begin{corr}}
\newcommand{\ec}   {\end  {corr}}
\newcommand{\bd}   {\begin{defn}}
\newcommand{\ed}   {\end  {defn}}
\newcommand{\ba}   {\begin{array}}
\newcommand{\ea}   {\end  {array}}
\newcommand{\bi}   {\begin{itemize}}
\newcommand{\ei}   {\end  {itemize}}
\newcommand{\be}   {\begin{enumerate}}
\newcommand{\ee}   {\end  {enumerate}}
\newcommand{\eq}[1]{\begin{equation}#1}
\newcommand{\en}   {\end{equation}}
\newcommand{\eqn}[1]    {\begin{equation} #1 \end{equation}}
\newcommand{\eqan}[1]   {\begin{align} #1 \end{align}}
\newcommand{\eqarray}   {\begin{eqnarray}}
\newcommand{\enarray}   {\end{eqnarray}}

\numberwithin{equation}{section}

    \newcommand{\vep}           {\varepsilon}

    \newcommand{\eps}           {\varepsilon}

    \renewcommand{\d}           {{\rm d}}
    \newcommand{\dd}            {\operatorname{d}^d \hspace{-2pt}}

    \newcommand{\sss}           { \scriptscriptstyle }

    \newcommand{\nn}            {\nonumber}
    \newcommand{\nnb}           {\nonumber \\}
    \newcommand{\shift}         {\!\!\!\!}

    \makeatletter
\newcommand\xconn[2][]{%
  \ext@arrow 9999{\longleftrightarrowfill@}{#1}{#2}}
\newcommand\longleftrightarrowfill@{%
  \arrowfill@\leftarrow\relbar\rightarrow}
\makeatother

\newcommand{\nndim}[1]{}
\renewcommand{\Re}{\mathrm{Re}}
\renewcommand{\Im}{\mathrm{Im}}

\DeclareMathAlphabet{\pazocal}{OMS}{zplm}{m}{n}
\newcommand{\Acal}  {\pazocal{A}}
\newcommand{\Bcal}  {\pazocal{B}}
\newcommand{\Ccal}  {\pazocal{C}}
\newcommand{\Dcal}  {\pazocal{D}}

\newcommand{\Fcal}  {\pazocal{F}}

\newcommand{\Kcal}  {\pazocal{K}}

\newcommand{\Ncal}  {\pazocal{N}}

\newcommand{\Scal}  {\pazocal{S}}

\newcommand{\Ucal}  {\pazocal{U}}

\newcommand{\Xcal}  {\pazocal{X}}

\newcommand{\B}     {\mathbb{B}}

\newcommand{\R}     {\mathbb{R}}
\newcommand{\N}     {\mathbb{N}}
\renewcommand{\P}   {\mathbb{P}}
\newcommand{\D}     {\mathbb{D}}
\newcommand{\E}     {\mathbb{E}}
\newcommand{\Z}     {\mathbb{Z}}
\newcommand{\1}     {\mathbbm{1}}

\newcommand{\T}     {\mathbb T}

\newcommand{\Rd}    {\mathbb{R}^d}
\newcommand{\Zd}    {\mathbb{Z}^d}
\newcommand{\edges} {\mathbb{B}}

\newcommand{\Td}    {{\mathbb{T}^d}}

\newcommand{\bb}    {\underline{b}}
\newcommand{\tb}    {\overline{b}}

\newcommand{\ulb}   {\underline{b}}
\newcommand{\olb}   {\overline{b}}

\newcommand{\uv}    {{\vec{u}}}

\newcommand{\bk}    {\mathsf{k}}
\newcommand{\bT}    {\mathsf{T}}
\newcommand{\boldt} {\mathsf{t}}
\newcommand{\bn}    {\mathsf{n}}

\newcommand{\tbe}   {\overline{e}}

\newcommand{\Piic}      {\mathbb{P}_{\mathsf{ \sss IIC}}}
\newcommand{\Eiic}      {\mathbb{E}_{\mathsf{ \sss IIC}}}
\newcommand{\Pp}        {\mathbb{P}_p}
\newcommand{\Ppc}       {\mathbb{P}_{p_c}}

\newcommand{\prob}      {{\mathbb P}}

\newcommand{\Ep}		{\mathbb{E}_p}
\newcommand{\expec}     {\mathbb{E}}

\newcommand{\indic}[2][]    {\1_{\{#2\}_{#1}}}
\newcommand{\indicwo}[1]    {\1_{{#1}}}

\newcommand{\ttau}		{\tilde{\tau}}

\newcommand{\taub}			{\tau^{\mathsf{b}}}
\newcommand{\ttaub}			{\tilde{\tau}^{\mathsf{b}}}

\newcommand{\htau}          {\hat{\tau}}
\newcommand{\htaub}          {\hat{\tau}^{\mathsf{b}}}

\newcommand{\indi}          {\mathbbm{1}}

\newcommand{\rhonT}[1]      {{\hat{\underline{\rm \rho}}}^{\sss (#1)}}
\newcommand{\taunT}[1]      {{\hat{\underline{\rm \tau}}}^{\sss (#1)}}

\newcommand{\Rho}           {\operatorname{P}}
\newcommand{\Tau}           {\operatorname{T}}
\renewcommand{\rho}         {\varrho}

\newcommand{\Phl}           {\Phi_{\rm  A}}

\newcommand{\Phr}           {\Phi_{\rm  B}}

\newcommand{\Dk}            {\hat{D}(k)}

\newcommand{\nin}       {\notin}
\newcommand{\ua}        {\nearrow}

\newcommand{\conn}      {\leftrightarrow}
\newcommand{\nconn}     {\nleftrightarrow}
\newcommand{\Conn}      {\Leftrightarrow}

\newcommand{\dbc}       {\Leftrightarrow}

\newcommand{\off}       {{\text{\rm{ off }}}}
\newcommand{\inn}       {{\text{\rm{ in }}}}

\DeclareMathOperator{\diam} {diam}

\newcommand{\tCcal}     {\tilde{\Ccal}}

\newcommand{\Piv}       {{\text{\rm Piv}}}

\newcommand{\am}        {{2\wedge  \alpha}}
\newcommand{\twa}       {{(2 \wedge \alpha)}}

\newcommand{\e}         {\mathrm{e}}

\newcommand{\tri}       {\triangle}
\newcommand{\btri}      {\bar{\triangle}}

\newcommand{\squz}      {\square_z}

\newcommand{\shortversion}[1]{}
\newcommand{\longversion}[1]{#1}

\graphicspath{{./images/}}

\begin{document}

\title{Backbone scaling limit of the {high}-{dimensional} IIC}
\longversion{\subtitle{Extended version}}
\author[M.\ Heydenreich]{Markus Heydenreich}
\address{Mathematisches Institut, Ludwig-Maximilians-Universit{\"a}t M{\"u}nchen, Theresienstr.~39, 80333~M{\"u}nchen, Germany.}
\email{m.heydenreich@lmu.de}
\author[R.\ van der Hofstad]{Remco van der Hofstad}
\address{Department of Mathematics and
	Computer Science, Eindhoven University of Technology, P.O.\ Box 513,
	5600 MB Eindhoven, The Netherlands.}
\email{rhofstad@win.tue.nl, w.j.t.hulshof@tue.nl}
\author[T.\ Hulshof]{Tim Hulshof}
\author[G.\ Miermont]{Gr\'egory Miermont}
\address{\'Ecole Normale Sup\'erieure de Lyon, Unit\'e de Math\'ematiques Pures et Appliqu\'ees, UMR CNRS 5669, 46, all\'ee d'Italie, 69364 Lyon Cedex 07, France.}
\email{gregory.miermont@ens-lyon.fr}
\date{\today}

\begin{abstract}
  We identify the scaling limit of the backbone of the {high}-di\-men\-sio\-nal incipient infinite cluster (IIC), both in the finite-range and the {long}-{range} setting. In the {finite}-{range} setting, this scaling limit is Brownian motion, in the long-range setting, it is a stable motion. The proof relies on a novel lace expansion that keeps track of the number of pivotal bonds.
\end{abstract}

\maketitle
\vspace{.6cm}
\noindent
{\small {\it MSC 2010.} 60K35, 60K37, 82B43.

}

\vspace{.6cm}
\hrule
\vspace{1em}

\longversion{
\tableofcontents
}

\section{Introduction}
\emph{Universality} --the emergence of the same macroscopic phe\-no\-mena from dissimilar microscopic laws-- and \emph{self-similarity} --the absence of identifying scales in macroscopic objects-- are two topics that have held the interest of mathematicians and physicists alike for many decades now. The ubiquity of these phenomena in real-life systems has inspired the close scrutiny of mathematical models that exhibit them as well. Critical percolation models, especially in two and high dimensions, stand out among these models as tractable instances with a remarkably rich structure. This has inspired long lines of research. One of the ultimate goals of this research is to determine the \emph{scaling limit} of critical percolation clusters. So far, this has yielded various new, deep insights (especially in two-dimensional models, where it was determined that the boundary scales to $\mathsf{SLE}_6$, and variants thereof \cite{Smir01,SmiWer01}). Yet in no setting is the scaling limit completely determined.
In this paper we take a step towards this goal by identifying the scaling limit of the \emph{backbone of large critical clusters and of the IIC} in high dimensions. (Roughly speaking, the IIC is a critical cluster conditioned to be infinite. The backbone of a cluster is the union of all self-avoiding paths between two marked vertices in the cluster; for the IIC, these are the vertex at the root and the vertex at infinity.)

We prove the following:
For critical nearest-neighbor percolation, when the dimension is sufficiently high, the backbone scales to a Brownian motion. This scaling limit also applies to sufficiently {spread}-{out} percolation in dimension $d \ge 15$, and to spread-out \emph{long-range} percolation with a sufficiently strong decay of the bond-rentention probabilities.

 Hence, our result demonstrates both universality, in that the underlying model does not affect the scaling limit much, and self-similarity, because the limit is a canonical self-similar object.
 Besides this, we also show that the backbone of long-range percolation in sufficiently high dimension but with weaker decay of bond-retention probabilities scales to a stable-motion whose index depends on the rate of decay, establishing another instance of universality and self-similarity in a percolation scaling limit.
 
 These results form the first scaling limit results for high-dimensional unoriented percolation paths. 
\medskip

\noindent
\subsubsection*{Our main innovations} To prove these scaling limits we derive a lace expansion for the two-point function (the probability of a connection event) with a {\em fixed number of pivotal bonds}, both for critical clusters, as well as for the incipient infinite cluster. These lace expansions are rather delicate, as we need to deal with various non-monotone events. This turns out to be substantially harder than the classical lace expansion, which involves only increasing events. We prove sharp asymptotics for such two-point functions, thus also deriving limits for the expected number of vertices in $\Zd$ that have precisely $n$ pivotal bonds along the percolation paths leading up to them. Then, by extending these results to convergence of finite-dimensional distributions and by proving tightness, we also prove convergence in path space. 
We do all this both in the critical percolation setting and in the incipient infinite cluster setting.

\longversion{
\subsubsection*{The abridged version of this paper.}
This is the extended version of the paper \cite{HeyHofHulMie17a}. The only significant difference between this version and \cite{HeyHofHulMie17a} is that this version contains two additional appendices with various standard technical estimates that are needed in the proofs, but that are not considered particularly interesting in their own right.}
\shortversion{
\subsubsection*{The extended version of this paper} 
Considering the length of the full proof, we have omitted some parts that are fairly standard in the lace expansion literature. These proofs, however, can be found in Appendices A and B of \cite{HeyHofHulMie17b}, the extended version of this paper.
}

\subsubsection*{Bond percolation on $\Zd$, a generalized setup.}
Our models are defined in terms of a (weight) function $D\colon \Zd \times \Zd \to [0,1]$ that satisfies $\sum_{v \in \Zd} D(u,v)=1$ for all $u \in 
\Zd$.  Let $p\ge0$ be a parameter chosen such that $p D(u,v) \le 1$ for all $\{u,v\} \in \Zd \times \Zd$. For arbitrary lattice sites $u,v\in\Zd$, we declare the bond $\{u,v\}$ {\em occupied} with probability $p D(u,v)$ and {\em vacant} otherwise.  The occupation statuses of the bonds are independent random variables. Mind that $p$ is {\em not} supposed to denote the probability of an event, but instead the expected degree of a vertex in the percolation model, and might exceed $1$.

Our results hold for a broad range of models. We state the precise assumptions we make on $D$ below in Assumptions \ref{ass:A}, \ref{ass:D} and \ref{ass:E}.
Before we state them, however, we first describe three ``standard'' models that satisfy these assumptions as examples to keep in mind.

The first ``standard'' model is also the best known model, namely the \emph{nearest-neighbor model}. Here $D(u,v)=1/(2d)$ for all $\{u,v\}$ with $|u-v|=1$, and $D(u,v)=0$ for all other bonds, where $| \, \cdot \, |$ denotes the Euclidean norm on $\R^d$. Note that this definition implies that each nearest-neighbor bond is occupied with probability $p/(2d)$, not with probability $p$.

The second ``standard'' model is a \emph{finite-range spread-out model,} where, for a given $L \in \N$, bonds of length up to $L$ are occupied with equal probability, and longer bonds are always vacant, 
i.e.,
\begin{equation}\label{e:defFinRangePerco}
    p D(u,v) = \frac{p}{(2 L +1)^d-1}\indi_{\{0 < \|u-v \|_{\infty} \le L\}},
\end{equation}
where $\| \,\cdot\, \|_\infty$ denotes the supremum norm.
The parameter $L$ serves to spread out the connections and is typically fixed at a large value.

The third ``standard'' model is a \emph{long-range spread-out model,} where the occupation probabilities decay as a power of the length of the bond. Indeed, for $\alpha \in (0,\infty)$, we define
\begin{equation}\label{e:defLongRangePerco}
    p D(u,v) = p \frac{\Ncal_L}{\max \{|u-v|/L , 1\}^{d+\alpha}},
\end{equation}
where $\Ncal_L$ is a normalizing constant. 

Our results apply in much greater generality than the three models described above: all that we require is that $D$ has the symmetries of $\Zd$ and that the Fourier transform of $D$ obeys certain bounds. To make this precise, we need some definitions. For $f\colon \Zd \to \R$, let $\hat{f}$ denote the Fourier transform of $f$, i.e.,
\begin{equation}\label{e:Fourierdef}
   	\hat{f}(k)=\sum_{x\in \Z^d}\e^{ik\cdot x}f(x)\, ,\qquad k\in \R^d\, .
\end{equation}
Note that $\hat{f}(k)$, when it exists, is a periodic function with fundamental domain $\Td := [-\pi,\pi)^d$.

We write $(a \wedge b) = \min \{a, b\}$. We further write $\hat f(k) \sim \hat g(k)$ if the ratio of $\hat f(k)/ \hat g(k) \to 1$ when $k \to 0$.

The results in this paper hold for models where $D$ satisfies the following three assumptions:
\begin{customass}{A}[Symmetry and asymptotics]\label{ass:A}
We assume that $D$ is invariant under the symmetries of $\Zd$, i.e., that $D(x, y)$ only depends on $x, y \in \Zd$ through $y - x$, so $D(x,y)= D(0,y-x)$, and that $D(x ,y) = D(x', y')$ for any $x', y' \in \Zd$ such that $x'- y'$ is equal to $x -y$ up to permutations and sign-changes of the coordinates. (Hence, we will frequently write $D(x) = D(0,x)$ when that is more convenient.)

We also assume that $D$ is normalized as
\begin{equation}\label{e:Ddef}
  	D(x) = \frac{h(x/L)}{\sum_{x\in \Zd} h(x/L)} \qquad \forall x \in \Zd,
\end{equation}
where $L \in (0,\infty)$ and $h$ is a nonnegative bounded function on $\mathbb{R}^d$ that is piecewise continuous, has the aforementioned symmetries, satisfies
\[
	\int\limits_{\mathbb{R}^d}  h(x) \dd x   =1,
\]
and that is either
\begin{enumerate}
	\item supported in $[-1, 1]^d$,
	\item exponentially decaying (i.e., there exist $C,\kappa >0$ such that $h(x)\le C \e^{- \kappa \|x \|_\infty}$),
	\item or that there exists $\alpha, c_1,c_2>0$ and $Q < \infty$ such that
		\begin{equation}\label{e:asymph}
   			c_1 | x|^{-d -\alpha} \le h(x) \le c_2 | x |^{-d -\alpha} \qquad \text{ for all } | x | \ge Q.
		\end{equation}
\end{enumerate}
We call a model that satisfies (i) or (ii) a \emph{finite-range model,} and a model that satisfies (iii) a \emph{long-range model.}

We consider $\alpha$ a parameter of the model, which we call the \emph{long-range parameter.} In such cases where $h$ has bounded support or decays exponentially we set $\alpha = \infty$.
We also consider $L$ a parameter of the model, which we call the \emph{spread-out parameter.}
\end{customass}

\begin{customass}{B}[Bounds on $\hat D$]\label{ass:D}
    Consider a percolation model that satisfies Assumption~\ref{ass:A} with long-range parameter $\alpha \in (0,\infty]$ and spread-out parameter $L$ (with $L=1$ for nearest-neighbor models). We assume that $D$ satisfies the following bounds:
    There exist constants $c_1, c_2>0$ such that
        \begin{align}
            \label{e:hatD1} 1- \hat D(k) &\ge c_1 L^\twa |k|^\twa   &\text{ if } \quad \|k\|_\infty \le L^{-1};\\
            \label{e:hatD2} 1- \hat D(k) &> c_2  &\text{ if }\quad \|k\|_\infty \ge L^{-1}.
        \end{align}
    Furthermore, there exists a constant $w$ with $0<w =O(L^\twa)$ such that, for $\vep >0$ sufficiently small,
        \begin{equation}\label{e:hatD4}
            1- \hat D(k) \le w |k|^\twa \qquad\qquad\qquad\qquad\;\text{ if }\;\quad|k| \le \vep.
        \end{equation}
\end{customass}
\begin{customass}{C}[Convergence of $\hat D$]\label{ass:E}
    Consider a percolation model that satisfies Assumption~\ref{ass:A} with long-range parameter $\alpha \in (0,\infty]$. We assume that there exists a constant $0 < v_\alpha < \infty$ such that, as $k \to 0$,
    \begin{equation} \label{hatD-asymp}
        1-\hat{D}(k)\sim
            \begin{cases}
                v_{\alpha}|k|^{(\am)} &\text{if }\alpha\neq  2,\\
                v_2|k|^2 \log(1/|k|) &\text{if }\alpha=2.
            \end{cases}
    \end{equation}
\end{customass}
Assumptions~\ref{ass:A} and \ref{ass:D} are common in the high-dimenional percolation literature. Assumption~\ref{ass:E} is not commonly needed, but it is a natural assumption to make for a scaling limit result. In fact, Assumption \ref{ass:E} is also needed to show that the scaling limit of a simple random walk with step distribution $D$ is Brownian motion or $\alpha$-stable motion. The same assumption is also made in \cite[(1.1)]{CheSak09} and \cite[Lemma 1.1]{Heyd11}. See \cite{Heyd11} for an in-depth discussion of the asymptotics in \eqref{hatD-asymp}.

We write $\P_p$ for the law of configurations of occupied bonds, and we write $\E_p$ for the corresponding expectation.
Given a configuration, we say that $x$ is connected to $y$, and write $x \conn y$, when there is a path of occupied bonds from $x$ to $y$
(or when $x = y$). Let $\B = \Zd \times \Zd$, so that $(\Zd, \B)$ is a complete graph.
We write $\Ccal(x)$ for the subgraph of the occupied percolation cluster that contains $x$ (and at times we abuse notation by also writing $\Ccal(x)$ for either just the sites or the occupied bonds of this subgraph).

We usually work at (or just below) the percolation threshold $p = p_c$ where $p_c$ is the \emph{critical value} of $p$, i.e., $p_c :=\inf\{p\, : \, \P_p(|\Ccal(0)|=\infty)>0\}$. In the parametrization of this paper, $p_c$ satisfies $p_c>1$, and $p_c$ tends to $1$ as either
$d \to \infty$ or $L \to \infty$, see \cite{HarSla90a,HeyHofSak08}.

Define the \emph{spatial variance} to be
\begin{equation} \label{def-sigma}
    \sigma^2 := \sum_{x \in \Z^d} |x|^2 D(x).
\end{equation}
We say that a model has finite variance when $D$ is such that $\sigma^2 <\infty$ and we say the model has infinite variance otherwise. Of course, $\sigma^2$ is finite for any finite-range model. The variance of a long-range model is finite when $\alpha>2$, but it is infinite when $\alpha \le 2$. Models with finite variance behave very differently from models with infinite variance.
For instance, the critical \emph{two-point function} $\Ppc(x \conn y)$, the probability that $x$ and $y$ are connected by a path of open edges at criticality, satisfies
\begin{equation}\label{e:twoptasymp}
	c |x-y|^{\twa - d} \le 	\Ppc(x \conn y) \le C |x-y|^{\twa - d}
\end{equation}
for some $C,c >0$.
This fact is proved for nearest-neighbor percolation with $d\ge 11$ in \cite{FitHof15b, Hara08}, for finite-range models with $d > 6$ and $L$ sufficiently large in \cite{HarHofSla03}, and for long-range models with $d > 3\twa$ and $L$ sufficiently large (and under some additional assumptions on $D$) in \cite{CheSak12}. In many of our results below the same behavior is also apparent.


\subsubsection*{The incipient infinite cluster.}
It is common knowledge for two-dimensional models and for models in sufficiently high dimension that there is no infinite cluster at the critical point. Nevertheless, we can think of the critical point as the point where the infinite cluster is at the verge of appearing. Thus, one might believe that we can construct infinite clusters at the critical point via reasonable conditioning and limiting schemes. Indeed, Kesten \cite{Kest86a} showed that conditioning the critical two-dimensional percolation measure on the event that the vertex at $0$ is connected to the boundary of a ball of radius $r$, when $r$ is taken to infinity, gives a well-defined probability measure.  The infinite cluster of $0$ of this measure is known as the \emph{incipient infinite cluster} (IIC). Kesten also showed that two other reasonable schemes for an infinite cluster at criticality yield \emph{the same measure} in two dimensions. Later, limit scheme constructions of the IIC were given for {high}-di\-men\-sio\-nal models as well:  The IIC for {spread}-{out} \emph{oriented} percolation above $4+1$ dimensions was constructed in \cite{HofHolSla02}, for {finite}-{range} {spread}-{out} percolation above 6 dimensions and nearest-neighbor percolation with $d$ large in \cite{HofJar04}, and in the general setting discussed above in~\cite{HeyHofHul14a}. For nearest-neighbor percolation this was extended to $d\ge 11$ in \cite{FitHof15b}.

The easiest to use construction of the IIC goes as follows. Define, for every cylinder event $F$ (i.e., every event that depends on the occupation status of finitely many bonds),
	\begin{equation}\label{def:Qiic}
    	\Piic(F) := \lim_{p\ua p_c}\frac{\sum_{x\in\Zd}\P_p(F\cap\{0\conn x\})}{\sum_{x\in\Zd}\P_p(0\conn x)}=\lim_{p\ua p_c}\frac{1}{\chi(p)}\sum_{x\in\Zd}\P_p(F\cap\{0\conn x\}),
	\end{equation}
where we define the \emph{susceptibility} $\chi(p):= \sum_x \P_p (0 \conn x)$ to be the expected size of a typical cluster. Because of the appearance of the factor $\chi(p)$ on the right-hand side of \eqref{def:Qiic}, we call this limit the {\em susceptibility limit}. It is proved in \cite{HeyHofHul14a} that this construction works in the generalized setting of high-dimensional percolation (i.e., if the dimension is sufficiently high and if Assumptions \ref{ass:A} and~\ref{ass:D} hold). Since the limiting measure $\Piic$ exists for all cylinder events, the Kolmogorov Extension Theorem implies that the measure also exists on their sigma algebra. It is also shown in \cite{HeyHofHul14a,HofJar04} that several related and natural constructions lead to the \emph{same} limit. This indicates that the IIC is a natural and robust object. 

The IIC in {high}-di\-men\-sio\-nal percolation has attracted considerable attention. For instance, it has been observed that random walk on the IIC is strongly subdiffusive. This phenomenon has been studied extensively in recent years (cf.\ \cite{BarJarKumSla08, HeyHofHul14b, KozNac09}).

Another aspect of the IIC that has been studied is its scaling limit. It is widely believed that the scaling limit of very large critical percolation clusters is super-Brownian motion (SBM). There is plenty of supporting evidence for this conjecture, much of it coming from studies of the IIC. Indeed, the asymptotics of the $r$-point functions of the oriented percolation IIC have been identified as those for the canonical measure of super-Brownian motion conditioned to survive forever (ICSBM) \cite{Hofs06a,HofSla03b}. The ICSBM measure was introduced by Evans \cite{Evan93}. It consists of a single infinite Brownian motion path (the immortal particle) together with super-Brownian motions branching off this path. The ICSBM can be viewed as an SBM conditioned to have infinite mass.

In this paper we take a step towards identifying the scaling limit of the IIC of {high}-di\-men\-sio\-nal percolation by identifying the scaling limit of a subgraph of the IIC that corresponds to the trace of the immortal particle of ICSBM.
Indeed, the IIC contains an (essentially) unique infinite path: The \emph{union} of all infinite self-avoiding paths in the IIC started from the origin is a well-defined topologically one-ended random object (cf.\ \cite{HofJar04}). Furthermore, it turns out that the \emph{intersection} of these paths yields an infinite, naturally ordered family of bonds, called {\em backbone pivotal bonds}.  
These objects should both play a similar role as the immortal particle for ICSBM.
We show that the scaling limit of the IIC backbone and the IIC backbone pivotals both yield a Brownian motion for finite-variance models, while they are a stable motion for infinite-variance models. This is consistent with the conjecture that the IIC is super-Brownian or super-$\alpha$-stable motion conditioned to survive forever, and it might bring its proof significantly forward.


\subsection{Main results}
\label{sec-res}
Given an event $A$ and a configuration $\omega \in A$, we say that a bond $b$ is \emph{pivotal} for $(A,\omega)$ if the configuration $\omega_b$, which is the same as $\omega$ \emph{except at $b$,} satisfies $\omega_b \notin A$. (We often simply say that $b$ is pivotal for $A$ and leave $\omega$ implicit.) We say that a bond $b$ is a \emph{backbone pivotal} if $b$ is pivotal for $\{0 \conn \infty\}$.

The backbone pivotal bonds are naturally ordered by their appearance as $(e_i)_{i=1}^{\infty}$, in such a way that every infinite self-avoiding path of occupied bonds started at the origin passes through $e_i$ before it passes through $e_{i+1}$. 
Moreover, the backbone-pivotal bonds can be naturally viewed as being \emph{directed} bonds $b=(x,y)$, where the direction is such that 
on the configuration $\omega_b$, 0 is connected to $x$ but not to $y$.
For a directed bond $b=(x,y)$, we write $\bb=x$ for its \emph{bottom}, and $\tb=y$ for its \emph{top}. Then we write
\begin{equation} \label{Sn-def}
    	S_n=\tbe_n
\end{equation}
for the lattice position of the top of the $n$th backbone pivotal bond $e_n$, with the convention that $S_0=\tbe_0=0$. 
We view the family $(S_n)_{n=0}^{\infty}$ as a discrete-time stochastic process. We will study its scaling limit in high dimensions, where geometry tends to trivialize. That geometry becomes trivial in high dimensions can be understood by noting that the displacement $S_n-S_{n-1}=\tbe_n-\tbe_{n-1}$ is the
displacement between two subsequent backbone pivotals, and, in high dimensions, these displacements should be only weakly dependent. Therefore, we expect that the scaling limit of $(S_n)_{n=0}^{\infty}$ is the same as the scaling limit of a random walk with independent and identically distributed steps. This suggests that the scaling limit is either a Brownian motion or a stable motion, depending on the number of existing (spatial) moments of
$\tbe_n-\tbe_{n-1}$. 

We define the scaling function as
\begin{equation}\label{fn-def}
    	f_{\alpha}(n)=
        \begin{cases}
            (v_{\alpha}n)^{-1/(2 \wedge \alpha)} &\text{if }\alpha\neq  2,\\
            (v_2n \log{n})^{-1/2} &\text{if }\alpha=2.
        \end{cases}
\end{equation}
Furthermore, we define the continuous-time stochastic process $X_n$ as
\begin{equation}\label{Xn-def}
    	X_n(t):=f_{\alpha}(n)\, S_{\lceil nt\rceil}, \qquad t\geq 0.
\end{equation}

Our results apply to the IIC but also to critical percolation. We need a few definitions to state the extension to critical percolation. In the context of critical percolation, whenever we say that a bond $b$ is \emph{pivotal}, we mean that it is pivotal for the event $\{x\conn y\}$ for some $x,y\in \Z^d$, i.e., when $x$ is connected to $y$ on the (possibly modified) configuration where
$b$ is made occupied, while $x$ is not connected to $y$ on the (possibly modified) configuration where $b$ is made vacant. For every $n\geq 1$, we define a probability measure on \emph{marked configurations,} i.e.,\ on the set of pairs $(\omega,x)$ (where $\omega$ is a percolation configuration and $x\in \Z^d$), by
\begin{equation}
   	 \E^*_{p_c,n}[F]=\frac{\E_{p_c}\Big[\sum_{x\in \Z^d}F(\cdot,x)
	 \indic{0\conn x \text{ with $n-1$ occupied pivotal bonds}} {}\Big]}
	 {\E_{p_c}\Big[\sum_{x\in \Z^d}\indic{0\conn x \text{ with $n-1$ occupied pivotal bonds}}{}
	 \Big]},
\end{equation}
for every non-negative measurable function $F=F(\omega,x)$.

 We write $x_*$ for the distinguished vertex under $\P^*_{p_c,n}$ (but note that $x_*$ is a random variable with respect to $\P^*_{p_c,n}$). Under this measure, we let $S'_0=0$, and $S'_1,S'_2,\ldots,S'_{n-1}$ be the top of the $i$th pivotal bond for the event $\{0\conn x_*\}$, and $S'_n=x_*$. This defines the random process $\big(S'_i\big)_{i=0}^n$.

Similarly to the IIC setting in \eqref{Xn-def}, let $Y_n$ be the rescaled version of $\big(S'_i\big)_{i=0}^n$ given by
\begin{equation}\label{eqDefY}
	Y_n(t):=f_{\alpha}(n)\, S'_{\lceil nt\rceil},
    	\qquad t\in[0,1].
\end{equation}

In the following theorem we let $(B^{\sss (2)}_t,t\geq 0)$ be a standard $d$-dimensional Brownian motion, and for $\alpha\in (0,2)$, we let
$(B^{\sss(\alpha)}_t,t\geq 0)$ be a symmetric stable process with characteristic function $\E[\exp(i \lambda B^{\sss(\alpha)}_1)]=\exp(-|\lambda|^\alpha)$ for $\lambda \in \R$. As a shorthand, we write $B^{\sss \twa}$ for these processes.

\begin{theorem}[Backbone scaling limit]
\label{thm-backbone-scal-limit}
Consider nearest-neighbor bond percolation on $\Zd$ with $d$ sufficiently large, or a percolation model that satisfies Assumptions~\ref{ass:A}, \ref{ass:D}, and \ref{ass:E} with the spread-out parameter $L$ sufficiently large and $d > 6\twa$.
Then there exists $K_\alpha>0$ (depending on $ D$) such that the following convergences in distribution hold as $n\to\infty$ in the space of right-continuous functions with left
limits, respectively $\D([0,\infty),\R^d)$ and $\D([0,1],\R^d)$, endowed with the Skorokhod $J_1$ topology:
\begin{equation}
	X_n
	\Rightarrow
	K_\alpha^{\frac1{(\am)}}\,B^{\sss(\am)}
	\qquad\text{ and }\qquad
	Y_n
	\Rightarrow
	K_\alpha^{\frac1{(\am)}}\,B^{\sss(\am)}.
\end{equation}
\end{theorem}

In the course of the proof we also derive a formula for $K_\alpha$, which can be found in \eqref{eqDefKalpha} below.

The condition of Theorem \ref{thm-backbone-scal-limit} is that $d$ or $L$ is ``sufficiently large'', but how large is sufficient? It turns out that we can associate a parameter to every distribution $D$. We call this the \emph{mean-field parameter} $\beta$. For nearest-neighbor models we let $\beta = O(d^{-1})$ and for spread-out models we let $\beta = O(L^{-d})$. The mean-field parameter $\beta$ can be made arbitrarily small by increasing $d$ or $L$. We use the lace expansion to establish bounds on certain quantities in terms of power series in $\beta$ \longversion{(see in particular Appendix~\ref{sec-LEcoefficients}).}\shortversion{(see in particular Appendix A of the extended version of this paper, \cite{HeyHofHulMie17b}).} A small $\beta$ ensures that these series converge. \nndim{To bring the dimension down to $d>14$ for nearest-neighbor percolation, we rely on the recent work of Fitzner with the second author \cite{FitHof15a, FitHof15b} that proves mean-field behavior for $d\geq 11$.}

The mean-field parameter $\beta$ also appears in the bound on the triangle diagram,
	\begin{equation}\label{eq-STC}
    	\sum_{x,y \in \Zd} \P_{p_c}(0 \conn x) \P_{p_c}(x \conn y) \P_{p_c}(y \conn 0) = 1+ O(\beta).
	\end{equation}
This bound is known as the \emph{strong triangle condition}, which holds for nearest-neighbor percolation when $d$ is sufficiently large \cite{HarSla90a}, and for models that satisfy Assumptions~\ref{ass:A} and \ref{ass:D} when $d> 3\twa$ \cite{HarSla90a,HeyHofSak08}. The \emph{triangle condition} merely states that the triangle diagram in \eqref{eq-STC} is {\em finite}, and is believed to be the necessary and sufficient condition for mean-field behavior. For nearest-neighbor percolation it is known to hold when $d \ge 11$ \cite{FitHof15b}. Moreover, it is believed that $d>3\twa$ is a necessary and sufficient condition for the triangle condition to hold. The value $3\twa$ is hence known as the \emph{upper critical dimension} for percolation. Most results in our paper are valid under the strong triangle condition with a sufficiently small $\beta$. In particular, the convergence of finite-dimensional distributions holds when $d>3(\am)$ for spread-out models. Only one aspect of the proof of Theorem \ref{thm-backbone-scal-limit}, namely tightness of the sequences $(X_n)_{n \ge 1}$ and $(Y_n)_{n \ge 1}$, requires the stronger condition that $d > 6 \twa$. 
We strongly believe that this constraint on the dimension is due to suboptimal bounds, and that Theorem~\ref{thm-backbone-scal-limit} holds when $d>3\twa$, but were unable to obtain tightness down to the upper critical dimension (see also Remark~\ref{rem-dimension-tightness} below).

Theorem \ref{thm-backbone-scal-limit} shows that the set of pivotal bonds of the IIC backbone is close to the image of a Brownian or $\alpha$-stable motion when properly renormalized. It is natural to ask whether the geometry of the entire backbone is well-captured by the pivotal bonds. Let $\Bcal_n$ be the set of vertices of $\Z^d$ that are in the backbone, and that are separated from the origin by at most $n$ backbone pivotal bonds. With this definition, the increasing union $\Bcal_\infty := \bigcup_{n \ge 0} \Bcal_n$ is the set of vertices of the backbone. For $n\geq 1$, we also let $\Scal_n=\Bcal_n\setminus \Bcal_{n-1}$, and let $\Scal_0=\Bcal_0$. The set $\Scal_n$ can be viewed as the vertex set of the subgraph that is induced by the $n$th ``sausage'' of the backbone. The sausages are the doubly-connected subgraphs of the backbone (in the sense that removing any one of its bonds cannot disconnect the graph). Note that $\Scal_n$ could be just a single vertex.

We let $\Kcal$ be the set of non-empty compact subsets of $\R^d$. The \emph{Hausdorff distance} $d_H$ between two subsets $A,A'\subseteq \R^d$ is given by
\begin{equation}
	d_H(A,A')=\sup_{x\in A}d(x,A')\vee \sup_{x'\in A'}d(x',A)\, ,
\end{equation}
where $d(x,A')=\inf_{x'\in A'}|x-x'|$. The space $(\Kcal,d_H)$ is a complete, locally compact metric space~\cite{BurBurIva01}.

In the remainder of the introduction, we will work under an extra hypothesis. We believe that this hypothesis is true in general, but we have only been able to prove it in certain cases ({finite}-{range} models are one such case). Under this hypothesis we can prove that the sausages are uniformly small in the scale $f_\alpha(n)$, so that compact subsets of the backbone are close -- in the Hausdorff sense -- to the set of pivotals. A common notion in the percolation literature is that for events $A,B$, the event $A\circ B$ is the event of disjoint occurrence of $A$ and $B$ (see \cite{Grim99} or Section~\ref{sec-BK-BKR} below for a definition and more details).

\begin{customhypo}{H}\label{hyp:H} There exists a finite constant $C>0$
such that for every $m\geq 1$,
	\begin{equation}
	\max_{|x|\leq m}\P_{p_c} \big(\exists y\in \Z^d\colon |y|>2m,\{0\conn y\}\circ\{y\conn x\}\big)
	\leq \frac{C}{m^{2(2\wedge \alpha)}}\, .
    \end{equation}
\end{customhypo}
While we do not have a proof that shows that Hypothesis \ref{hyp:H} holds under the most general set of
assumptions used in this paper, we prove it under slightly stronger assumptions:
\begin{prop}[Verification of Hypothesis \ref{hyp:H}]
\label{prop-assum-H}
Hypothesis \ref{hyp:H} holds for percolation models that satisfy Assumptions~\ref{ass:A} and \ref{ass:D} and that are
\begin{itemize}
	\item[(i)] {finite}-{range} and satisfy the strong triangle condition \eqref{eq-STC} with sufficiently small $\beta$, or
	\item[(ii)] long-range with $L$ sufficiently large and $d>4(2\wedge \alpha)$, or
	\item[(iii)] long-range with $L$ sufficiently large, $d>6$, $\alpha > 4$, and that satisfy \eqref{e:twoptasymp}.
\end{itemize}
\end{prop}
Note that Hypothesis~\ref{hyp:H} is thus essentially verified under the assumptions of Theorem~\ref{thm-backbone-scal-limit} (the only difference being that the sufficiency conditions for $\beta$ and $L$ may be more stringent).

\begin{theorem}[Hausdorff convergence of the IIC backbone]
\label{sec:main-results-1}
Assume Hypothesis \ref{hyp:H}.
Under the conditions of Theorem \ref{thm-backbone-scal-limit}, and
under the probability measure $\Piic$,
for every $T\geq 0$, the following convergence in
distribution holds in the space $(\Kcal,d_H)$:
\begin{equation}
	K_\alpha^{(\am)}\,f_\alpha(n)\,\Bcal_{\lfloor nT\rfloor}\Rightarrow
	\big\{B^{\sss(\am)}_t\colon 0\leq t\leq T \big\}^{\mathrm{cl}}\, ,
\end{equation}
where $A^{\mathrm{cl}}$ is the closure of $A\subseteq \R^d$.
\end{theorem}

Several variants of this result can be stated. For instance, we can view $(\Scal_{\lfloor nt\rfloor},t\geq 0)$ and $(\Bcal_{\lfloor nt\rfloor},t\geq 0)$ as stochastic processes in the Skorokhod space $\mathbb{D}([0,\infty),\Kcal)$ endowed with the inherited Skorokhod $J_1$ topology (see e.g.\ \cite{EthKur86}).
\begin{prop}\label{sec:main-results}
Assume Hypothesis \ref{hyp:H}.  Under the conditions of Theorem \ref{thm-backbone-scal-limit},
\begin{eqnarray*}
 	\left(K_\alpha^{(\am)}\,f_\alpha(n) \,\Scal_{\lfloor nt\rfloor},t\geq 0\right)
  	&\Rightarrow& \left(\{B^{\sss(\am)}_t\},t\geq 0\right)\, ,\\
 	\left(K_\alpha^{(\am)}\,f_\alpha(n) \,\Bcal_{\lfloor nt\rfloor},t\geq 0\right)
 	&\Rightarrow& \big(\{B^{\sss(\am)}_s:0\leq s\leq t\}^{\mathrm{cl}},t\geq 0\big)\, ,
\end{eqnarray*}
in distribution in $\mathbb{D}([0,\infty),\Kcal)$.
\end{prop}
This result is proved in Section \ref{sec-set-conv}. Theorem \ref{sec:main-results-1} follows as a corollary (see Section~\ref{sec:hausdorff} below).
An analogous result may be stated for Skorokhod convergence of the backbone under $\E^*_{p_c}$, but we omit this here.


\subsection{Further results}
\label{sec-fur-res}
In this section we state some results that are used in the proof of Theorem~\ref{thm-backbone-scal-limit} and that are interesting in their own right as well. The first such result is that the one-dimensional marginal distribution of the processes $X_n$ and $Y_n$ converge.

For $x\in \Z^d$, we define the IIC backbone two-point function by
\begin{equation} \label{rho-n-def}
    	\rho_n(x):=\Piic(S_n=x),
\end{equation}
i.e., $\rho_n$ is the probability mass function for the position of the top of the $n$th backbone-pivotal. We further study the two-point function with
a \emph{fixed} number of pivotal bonds,
\begin{equation} \label{taunx-def}
  	\tau^p_n(x):=\P_{p}(0\conn x \text{ with $n$ pivotal bonds}).
\end{equation}
For $\rho_n$ and $\tau^{p_c}_n$ we prove the following result:

\begin{theorem}[Weak convergence of end-to-end displacement]
\label{thm-endpoint}

Consider nearest-neighbor bond percolation on $\Zd$ with $d$ sufficiently large, or a percolation model that satisfies Assumptions~\ref{ass:A}, \ref{ass:D}, and \ref{ass:E} with the spread-out parameter $L$ sufficiently large and $d > 3\twa$. 

Let $k_n := f_\alpha (n) k$. 
There exist constants $K_{\alpha}, A\in(0,\infty)$ (depending on $ D$) such that the Fourier transforms of the IIC backbone two-point function $\rho_n$ and the two-point function with a fixed number of pivotals $\tau^{p_c}_n$ satisfy that 
\begin{equation}\label{EndpointDistribution}
	\hat \rho_n(k_n)\to\exp\big(-K_\alpha\,|k|^\twa\big), \qquad \hat \tau^{p_c}_n(k_n)
	\to A\exp\big(-K_\alpha\,|k|^\twa\big) \qquad\text{as $n\to\infty$}
\end{equation}
uniformly in $k$ on any compact subset of $\,\R^d$. 
\end{theorem}
The constants that appear in this theorem are defined in terms of the \emph{lace-expansion coefficients:} $K_{\alpha}$ is defined in \eqref{eqDefKalpha}, and $A$ in \eqref{eqDefA'} below.
The following corollary identifies $A$ as the expected number of vertices $x$ for which there are $n$ pivotals for $\{0\conn x\}$:

\begin{corr}[Expected number of vertices $n$ pivotals away from origin]
\label{corr:pivballboundary} 
Under the assumptions of Theorem~\ref{thm-endpoint}, as $n \rightarrow \infty$,
\begin{equation}
	 \sum_{x \in \Zd} \tau^{p_c}_n(x) =  A (1 + o(1)).
\end{equation}
\end{corr}
\proof Apply Theorem~\ref{thm-endpoint} with $k=0$. \qed
\medskip

The \emph{mean-$r$ displacement} along the IIC backbone is defined as the $r$th spatial moment of $\rho_n(x)$. We write $f\asymp g$ if there are uniform positive constants $c,C$ such that $cg\le f \le Cg$.
\begin{theorem}[Mean-$r$ displacement]\label{thm-xiR}
Under the assumptions of Theorem \ref{thm-endpoint},
for any $r<(\am)$,
\begin{equation}
\label{eqXiR}
    \left(\sum_{x\in\Zd}|x|^r\rho_n(x)\right)^{1/r}
    \asymp\left(\sum_{x\in\Zd}|x|^r\tau^{p_c}_n(x)\right)^{1/r}
    \asymp
        \begin{cases}
            n^{1/(\am)},\quad&\text{if $\alpha\neq2$},\\
            (n\,\log n)^{1/2}
            ,\quad&\text{if $\alpha=2$},
        \end{cases}
\end{equation}
as $n\to\infty$.
\end{theorem}
Chen and Sakai \cite{CheSak11} have shown that similar bounds as in \eqref{eqXiR}
hold for the mean-$r$ displacement of long-range self-avoiding walk and long-range oriented percolation for \emph{all} $r\in(0,\alpha)$. In their work, they also identify leading-order constants and give bounds on the error terms. The weaker statement (in terms of the constraint on $r$) in Theorem \ref{thm-xiR} lets us prove tightness of the sequences $X_n$ and $Y_n$, and thus suffices for our purposes.

\subsection{Discussion}
\label{sec-disc}
\paragraph{{\it Nearest-neighbor percolation}} 
In two recent works \cite{FitHof15a, FitHof15b}, Fitzner and the second author prove that mean-field behavior for critical percolation holds when $d\geq 11$. This extends the seminal result of Hara and Slade \cite{HarSla90a} that was previously known to work for $d\geq 19$ \cite{HarSla94}. (Hara, in private communication with the second author, later reported that their method even applies for $d\geq 15$.) 

It would be interesting to verify above which dimension the current scaling limit results apply for nearest-neighbor percolation. We believe that the result is valid whenever $d>6$, but a proof of this is far beyond our reach. Certainly, the scaling limit requires a rigorous proof of mean-field behavior, so at the moment Theorem \ref{thm-endpoint} cannot be verified unless $d\geq 11$. We doubt, however, that $d=11$ can be achieved without substantial changes to the proof, as our lace-expansion coefficients contain more terms than those in the classical lace expansion by Hara and Slade in \cite{HarSla90a} (roughly speaking, the number of terms in our expansion grows as $4^N$, whereas in theirs it grows as $2^N$). This means that we need a bound on the triangle diagram \eqref{eq-STC} that is substantially sharper than what is currently known. 

Besides this issue, our proof of Theorem \ref{thm-backbone-scal-limit} also uses some suboptimal estimates in the proof of tightness that limit the result to $d \ge 13$, although these estimates might be improved with fairly standard techniques, but some significant effort. 

\medskip

\paragraph{{\it Scaling limit of the IIC}}  It would be interesting to see if we can use our results to identify the scaling limit of the complete IIC. Indeed, the IIC can be viewed as a backbone ``decorated'' with critical clusters connected to it by a single bond, as is the case for the IIC on the tree. But on the tree the critical clusters attached to the backbone are \emph{independent and identically distributed}, whereas for the IIC on $\Z^d$, they are \emph{mutually avoiding}. 

Hara and Slade, in \cite{HarSla00a, HarSla00b} derived several scaling-limit-type results for critical percolation clusters conditioned on their size, and showed that these are the same as for the \emph{integrated super-Brownian excursion} (ISE), a tree-like measure-valued process introduced by Aldous \cite{Aldo93}. This suggests that the scaling limit of large critical high-dimensional percolation clusters is ISE. The picture is not complete, however, as there is no fixed dimension for which finite-dimensional distributions are established.

But perhaps these results of Hara and Slade can be combined with the scaling limit of the backbone to get the scaling limit for the full IIC. This approach may be easier for \emph{oriented percolation},
where our knowledge of the scaling limit of large critical clusters is more precise \cite{Hofs06a, HofSla03b}.
\medskip

\paragraph{{\it Random walk on the IIC}} The scaling limit of the IIC backbone is an important ingredient in the study of random walk on the {high}-di\-men\-sio\-nal incipient infinite cluster in \cite{HeyHofHul14b} (in particular, Assumption (S) therein). Indeed, Theorem \ref{thm-backbone-scal-limit} is used to estimate the number of backbone pivotals between the origin and the boundary of a large Euclidean ball. This in turn gets a lower bound on the effective resistance between the origin and the boundary of the ball. Effective resistances are often key quantities when studying random walk properties (cf.\ \cite{DoySne84}).
\medskip

\paragraph{{\it Scaling limit of random walk on the IIC}}
Recently, substantial progress was made on scaling limits of random walks on high-dimensional critical objects. Ben Arous, Cabezas, and Fribergh \cite{BenCabFri16a} consider random walk on the high-dimensional \emph{branching random walk} conditioned to have total population $n$. This is a natural candidate for the scaling limit of random walk on a critical high-dimensional percolation cluster with $n$ vertices, since branching random walk is often a good mean-field model for percolation. It is believed that such critical clusters of finite size have the same scaling limit as for critical branching random walk conditioned on the population size, namely ISE. In \cite{BenCabFri16a}, the authors prove that random walk on critical branching random walk conditioned on the population size being $n$ converges to  \emph{Brownian motion on Integrated Super-Brownian Excursion}, the latter object having been introduced by Croydon \cite{Croy09}. In \cite{BenCabFri16b}, the same authors extend such results to other critical objects under certain assumptions that have yet to be verified for critical high-dimensional percolation. Ben Arous and Fribergh \cite{BenFri16} and Fribergh \cite{Frib13} prove similar results for {\em biased} random walks.

It would be of great interest to extend the results of Ben Arous \emph{et al.}\ to scaling limits for random walks on the high-dimensional IIC. For this, we believe that our result is likely quite relevant, since the IIC has a single infinite path which is the backbone path. As a result, random walk can only escape to infinity along the backbone, while the clusters hanging off it act as {\em traps} for it, making the backbone scaling limit an important substructure. A first possible result using our results could be to establish that random walk restricted to the backbone has, apart from a constant rescaling of time, the same scaling limit as random walk on the random walk trace, which Croydon identifies in \cite{Croy08}.
\medskip

\paragraph{{\it Convergence to the canonical measure of super-Brownian motion (CSBM)}} The convergence of $\hat \tau^{p_c}_n(k_n)$ in Theorem \ref{thm-endpoint} can be viewed as the convergence of the two-point function at a fixed time. This could be extended to the convergence of all $r$-point functions, as was done for oriented percolation in \cite{HofSla03b} and for lattice trees in \cite{Holm08}. We define the $r$-point function as $\tau_{\vec{n}}^{p_c, (r)}(0,x_1,\dots,x_{r-1}) := \Ppc\big(\bigcap_{i=1}^{r-1} \{0\conn x_i\text{ with }n_i \text{ pivotals}\}\big)$. Together with the results of the second author and Holmes in \cite{HofHol12}, this would also imply the convergence of the {\em survival probability} $n\theta_n=n\prob_{p_c}(\exists x\text{ such that }0\conn x, |\Piv(0,x)|=n)$. In turn, Holmes and Perkins \cite{HolPer07} have shown that these would imply that the rescaled critical percolation cluster, with the number of pivotals interpreted as time, converges to CSBM. Recently, with Holmes and Perkins, the second author has identified a tightness criterion for convergence of critical spatial structures to CSBM in terms of finite-dimensional convergence \cite{HofHolPer16}. These questions constitute major challenges in high-dimensional percolation.
\medskip

\paragraph{{\it The intrinsic distance as a time variable}}  In this paper we consider the number of pivotals between (connected) vertices $x,y\in \Z^d$ as the time between them. Of course, an alternative and possibly more natural time variable could be the intrinsic distance along the percolation cluster. It would be of great interest to investigate whether $\tau_n'(x):=\prob_{p_c}\big(0\conn x \text{ with }d_{\Ccal(0)}(0,x)=n\big)$, where $d_{\Ccal(0)}(0,x)$ denotes the graph distance between $0$ and $x$ in $\Ccal(0)$, satisfies the same scaling behavior as $\tau^{p_c}_n(x)$ investigated here.


\subsection{Overview}
\label{sec-overview}
The rest of this paper is organized as follows.
In Section \ref{sec-overview} we give an outline of the proof of our main result, Theorem \ref{thm-backbone-scal-limit}.
In Section \ref{sec-le} we derive a lace expansion for $\rho_n(x)$ and $\tau^{p_c}_n(x)$ along the backbone pivotals.
In Section \ref{sec-bounds-lace-exp} we prove various diagrammatic estimates for the lace-expansion coefficients. These diagrammatic estimates are similar to the diagrammatic estimates for the classical lace expansion as first derived by Hara and Slade in \cite{HarSla90a}. Yet, they are also substantially different from these estimates, and require a detailed analysis. 

The bounds on the lace-expansion coefficients starting from the diagrammatic estimates are more or less standard, and are completed in \longversion{Appendices~\ref{sec-LEcoefficients} and \ref{sec-squares},}\shortversion{\cite[Appendices~A and~B]{HeyHofHulMie17b},} where we prove both spatial and temporal moment estimates for the lace-expansion coefficients. 
These estimates are formulated in several important propositions: The first, Proposition~\ref{prop:smallpi}, which is proved in \longversion{Appendix~\ref{sec-LEcoefficients},}\shortversion{\cite[Appendix~A]{HeyHofHulMie17b},} contains {\em weak bounds} and its proof relies only on bounds that are already known in the literature (in particular, from \cite{HarSla90a} and \cite{HeyHofSak08}). 
Using Proposition \ref{prop:smallpi} we prove so-called {\em infrared bounds} on generating functions of the form $\sum_n \hat{\tau}^{p_c}_n(k) z^n$ and $\sum_n \hat{\rho}_n(k) z^n$ for $z\in \mathbb{C}$. These bounds are formulated in Proposition \ref{prop:positivity}. They allow us to prove improved bounds on {\em spatial moments} and {\em temporal moments} of the lace-expansion coefficients, as stated in Propositions~\ref{prop-LEcoefficients2} and \ref{prop-LEcoefficients}, respectively. These bounds are proved in \longversion{Appendix~\ref{sec-squares}.}\shortversion{\cite[Appendix~B]{HeyHofHulMie17b}.} In Section \ref{sec-impr-bds-LE} we also prove Proposition~\ref{prop-LEcoefficients} subject to an important technical lemma, Lemma \ref{lem:Diagbd}, which gives a bound on the rate of divergence of a {\em weighted square diagram} as the weight vanishes. The improved bounds in Propositions~\ref{prop-LEcoefficients2} and~\ref{prop-LEcoefficients} are the key estimates needed in the proof of the sharp asymptotics of $\hat{\tau}^{p_c}_n(k)$ and $\hat{\rho}_n(k)$ in Theorem~\ref{thm-endpoint}. We carry out this proof in Section \ref{sec-pf-endpoint}. In Section \ref{sec-ErrorBound} we then use these moment estimates to prove Propositions~\ref{lemmaError} and \ref{propMain2}, which we combine to complete the proof of Theorem~\ref{thm-endpoint}.

As may have become clear from the above description, the order in which we prove our results is delicate but quite important, so let us elaborate on this point:

\begin{center}
\setlength{\fboxsep}{7pt}
\setlength{\fboxrule}{1pt}
\begin{boxedminipage}{15cm}
We {first} use classical results from \cite{HarSla90a} and \cite{HeyHofSak08} to prove {\em weak bounds} on the lace-expansion coefficients in Proposition \ref{prop:smallpi}. {Secondly,} we use these bounds to obtain {\em infrared bounds} on generating functions in Proposition \ref{prop:positivity}. These results together set the stage for an improved analysis that proves the sharp asymptotics in Theorem \ref{thm-endpoint}. Because the proof of Proposition~\ref{prop:smallpi} {\em only} relies on estimates on percolation quantities as proved in \cite{HarSla90a} and \cite{HeyHofSak08} and not on quantities that require knowledge of the backbone two-point functions or their lace expansions, we are certain that we avoid any kind of circular reasoning. Using the results from \cite{HarSla90a} and \cite{HeyHofSak08} as a starting point has the added benefit that we may forgo an analysis of the convergence of the lace expansion at $p_c$, which usually requires some complicated reasoning (e.g.\ the well-known ``bootstrap lemma'' of \cite{HarSla90a}).
To illustrate clearly that the argument is not circular, \longversion{we prove Proposition \ref{prop:smallpi} in Appendix~\ref{sec-LEcoefficients},}\shortversion{Proposition \ref{prop:smallpi} is proved in \cite[Appendix~A]{HeyHofHulMie17b},} which is self-contained and may be read immediately after Section~\ref{sec:LEprop}, so this will be established before we get to the improved bounds of Propositions~\ref{prop-LEcoefficients2} and \ref{prop-LEcoefficients}, which \longversion{we prove in Appendix~\ref{sec-squares}.}\shortversion{are proved in \cite[Appendix~B]{HeyHofHulMie17b}.}
\end{boxedminipage}
\end{center}

After this, we complete the proofs of the remaining results as follows.
In Section \ref{sec-meanr} we prove Theorem \ref{thm-xiR} on the mean-$r$ displacement.
In Section~\ref{sec-findim} we complete the proof of Theorem \ref{thm-backbone-scal-limit}, our main result, by proving Proposition~\ref{prop-FinDimConv} about convergence of finite-dimensional distributions, and Proposition~\ref{cor-tightness} about tightness.
In Section~\ref{sec-set-conv} we discuss convergence in path space and prove Theorem \ref{sec:main-results-1} on Hausdorff convergence, and Proposition \ref{prop-assum-H}, which verifies Hypothesis~\ref{hyp:H}.
\longversion{In Appendices \ref{sec-LEcoefficients} and \ref{sec-squares} we prove all the novel bounds on diagrams that are used in this paper.}
\shortversion{Appendices A and B of \cite{HeyHofHulMie17b} contain the proofs of the novel bounds on diagrams that are used in this paper.} For the most part these are quite standard computations.

\subsection{Notation}
We will henceforth write $\sum_x $ for $\sum_{x \in \Zd}$ when the sum is over vertices, and likewise $\sum_b$ for $\sum_{b \in \mathbb{B}}$ if the sum is over bonds, whenever it is clear from the context. We write $\sum_{n \ge i}$ for $\sum_{n =i}^\infty$. We write $C, c > 0$ for generic constants, which may change from line to line. We write $f(n) = O(g(n))$ if there exists a $C >0$ such that $|f(n)| \le C g(n)$ for all $n$ sufficiently large, and we write $f(n) = o(g(n))$ if $|f(n)| / g(n) \to 0$. Given non-negative $f(n), g(n)$ we also occasionally write $f(n) \gg g(n)$ or $g(n) \ll f(n)$ if $f(n) / g(n) \to \infty$.

\section{Outline of the proof of Theorem \ref{thm-backbone-scal-limit}}
\label{sec-overview}
In this section we give an overview of the proof of our main result, Theorem \ref{thm-backbone-scal-limit}. It is a classical result (see e.g.\ \cite[Theorem~13.1]{Bill99}) that $(P_n)_{n \ge 1}$ converges to $P$ in distribution in $\D([0,1],\R^d)$ if (1) the finite-dimensional distributions of $(P_n)_{n \ge 1}$ converge to those of $P$, and 
(2) $(P_n)_{n\ge 0}$ is tight on $\D([0,1],\R^d)$.
Our two main aims in this paper are thus to prove that these two properties hold for the backbone pivotal process of large critical clusters and the IIC, $Y_n$ and $X_n$.
To prove that $X_n$ converges in distribution in $\D([0,\infty),\R^d)$ it suffices to prove that the restriction of $X_n$ to the interval $[0,T]$ converges in distribution in $\D([0,T],\R^d)$ for every $T>0$. By a simple scaling argument this is equivalent to proving the case where $T=1$. Therefore, we only consider the restriction of the process $X_n$ to $[0,1]$ from now on.


\subsection{Convergence of finite-dimensional distributions}
By convergence of finite-dimensional distributions we mean that for every $N=1,2,3,\dots$, any $0<t_1<\cdots<t_N\le1$, and any bounded continuous function $g\colon\mathbb (\R^{d})^N\to\mathbb R$,
\begin{equation}\label{eqFinDimConv2}
    	\lim_{n\to\infty}\Eiic\big[g\big(X_n(t_1),\dots,X_n(t_N)\big)\big]
    	= \E\big[ g\big(B^{\sss (\am)}_{t_1},\dots,B^{\sss (\am)}_{t_N}\big)\big].
\end{equation}
If we have convergence of the characteristic functions, then convergence in distribution follows, so it suffices to consider functions $g$ of the form
\begin{equation}
    	g(x_1,\dots,x_N)=\exp\big( i\,\bk\cdot(x_1,\dots,x_N)\big),
\end{equation}
where $\bk=\big(k^{\sss (1)},\dots,k^{\sss (N)}\big)\in\R^{dN}$ and $x_i\in\Rd$, $i=1,\dots,N$.
The problem becomes easier still when we use the equivalent form
\begin{equation}
    	g(x_1,\dots,x_N)=\exp\big( i\,\bk\cdot(x_1,x_2-x_1,\dots,x_N-x_{N-1})\big).
\end{equation}
For $\bn=(n^{\sss (1)},\dots,n^{\sss (N)})\in\N^N$ with $n^{\sss (1)}<\dots<n^{\sss (N)}$,
we define
\begin{equation}
\label{eqDefCnN}
    	\rhonT{N}_{\bn}(\bk)=
    	\Eiic\Big[
      	\exp\Big(i\sum_{j=1}^Nk^{\sss (j)}\cdot\left(S_{n^{\sss (j)}}-S_{n^{\sss (j-1)}}\right)\Big)\Big]
\end{equation}
as the characteristic function of the increments of $(S_n)_{n\geq 0}$, with $n^{\sss (0)}=0$.
The quantity $\taunT{N}_{\bn}(\bk)$ is defined accordingly, with $S_n$ replaced by $S_n'$,
and $\Eiic$ in \eqref{eqFinDimConv2} replaced by $\E^*_{p_c,n}$.

\begin{prop}[Finite-dimensional distributions]\label{prop-FinDimConv}
Let $N$ be a positive integer,
$k^{\sss (1)},\dots,$ $k^{\sss (N)}\in\R^d$,
$0=t_0 <t_1 <\dots<t_N  \le 1$.
Write
\begin{equation}\label{e:knvecdef}
    	\bk_n
    	=\big(k_n^{\sss (1)},\dots,k_n^{\sss (N)}\big)
    	=f_\alpha(n)\,\big(k^{\sss (1)},\dots,k^{\sss (N)}\big),
    	\qquad
    	n\boldt=\big(\lfloor nt_1 \rfloor,\dots,
    	\lfloor nt_N \rfloor\big).
\end{equation}
Under the conditions of Theorem \ref{thm-endpoint},
\begin{align}
   	 \lim_{n\to\infty}\rhonT{N}_{n\boldt}(\bk_n)
    	=\lim_{n\to\infty}\frac{\taunT{N}_{n\boldt}(\bk_n)}{\hat\tau_{n}(0)}
    	=\exp\Big(-K_\alpha\,\sum_{j=1}^N|k^{\sss (j)}|^\twa \;(t_j -t_{j-1})\Big).
\end{align}
\end{prop}

We conclude that the finite-dimensional distributions of the {finite}-{range} and long-range IIC backbone converge to those of Brownian motion or of an $\alpha$-stable L\'evy motion. This also proves that Brownian or $\alpha$-stable motion is the only possible scaling limit for the backbone process.


\subsection{Tightness}
 Using the moment estimates in Theorem \ref{thm-xiR} we can prove tightness under the stronger condition $d>6(\alpha\wedge 2)$.
We prove the following statement in Section \ref{sec-tightness}:

\begin{prop}\label{cor-tightness}
Under the conditions of Theorem \ref{thm-backbone-scal-limit} the sequences $(X_n)_{n\ge 0}$ in \eqref{Xn-def} and $(Y_n)_{n\ge 0}$ in \eqref{eqDefY} are tight in $\mathbb D\big([0,1],\R^d\big)$.
\end{prop}

\subsection{Proof of Theorem~\ref{thm-backbone-scal-limit} subject to Propositions~\ref{prop-FinDimConv} and \ref{cor-tightness}}
Since $B^{\sss \twa}$ is a standard $d$-dimen\-sional Browian motion when $\alpha \ge 2$ and a standard symmetric $d$-dimensional $\alpha$-stable motion when $\alpha \in (0,2)$, the finite-dimensional distributions of $B^{\sss \twa}$ are characterized by
\begin{equation}
	\E\Big[\exp\big( i \mathsf{k} \cdot (B^{\sss \twa}_{t_1},\dots,B^{\sss \twa}_{t_N } - B^{\sss \twa}_{t_{N-1} })\big)\Big] = \exp\Big(- \sum_{j=1}^N |k^{\sss (j)}|^\twa (t_j - t_{j-1})\Big).
\end{equation}
Thus, by Proposition~\ref{prop-FinDimConv}, the finite-dimensional distributions of $X_n$ and $Y_n$ converge to those of $K_\alpha^{\frac{1}{\twa}} B^{\sss \twa}$. Moreover, by Proposition~\ref{cor-tightness}, the sequences $(X_n)_{n \ge 0}$ and $(Y_n)_{n \ge 0}$ are tight. Therefore, by \cite[Theorem~13.1]{Bill99} the claimed convergence holds. \qed
\medskip

\subsection{Lace expansion}
The key to the proofs of Propositions~\ref{prop-FinDimConv} and \ref{cor-tightness} is that we develop a lace expansion for the two-point functions ${\rho}_{n}(x)$ and ${\tau}^p_{n}(x)$. Let us now explain in more detail what these expansions entail. 

Define the \emph{two-point function}
\begin{equation}\label{tau-x}
    	\tau^p(x) := \P_p(0\conn x).
\end{equation}
Hara and Slade's original percolation lace expansion \cite{HarSla90a} gives an expansion for $\tau^p (x)$: they show that there exists a function $\pi^p(x) : \Zd \to [0,1]$ such that
\begin{equation}\label{HarSlaLE}
	\tau^p(x) = \pi^p(x) + (\pi^p * p D * \tau^p)(x),
\end{equation}
where $f\ast g$ denotes the convolution between two summable functions from $\Z^d$ to $\R$.
Moreover, the \emph{lace expansion coefficient} $\pi^p(x)$ is expressible as a series expansion through repeated applications of the inclusion-exclusion principle. When the strong triangle condition \eqref{eq-STC} is satisfied for sufficiently small $\beta$, this expansion can be shown to be convergent, which allows one to give bounds on $\pi^p(x)$ that are similar to Feynman diagrams.
So far, almost all results for high-dimensional percolation have been proved with this lace expansion.

Our proofs use a lace expansion of the form
\begin{equation}
	\label{eqExpansionX}
    	{\rho}_{n+1}(x) = \psi_{n+1}(x)+ \sum_{m=0}^{n} \left({\pi}_m 
	\ast p_cD\ast {\rho}_{n-m}\right)(x)
\end{equation}
for certain lace expansion coefficients $\pi_n(x)$ and $\psi_n(x)$. We derive \eqref{eqExpansionX} in Section \ref{sec-le}. We also derive a lace expansion for the critical two-point function with
a fixed number of pivotals in Section \ref{sec-le}, i.e., for $\tau_n^p$ as defined in \eqref{taunx-def} with $p \le p_c$. This expansion reads
\begin{equation}
\label{eqExpansionTauX}
    	\tau^p_{n+1}(x) = \pi^p_{n+1}(x)+ \sum_{m=0}^{n} \left({\pi}^p_m 
	\ast p D\ast {\tau}^p_{n-m}\right)(x).
\end{equation}
The coefficients $\pi^p_m(x)$ are the same as the ones appearing in
\eqref{eqExpansionX} at $p =p_c$, i.e., $\pi_m$ in \eqref{eqExpansionX} equals $\pi_m^{p_c}$ in \eqref{eqExpansionTauX}. Moreover, as we discuss in Remark~\ref{rem-HS-exp} below, there is a simple expression for $\pi^p$ in terms of $\pi^p_m$.

We can use the expansion in \eqref{eqExpansionX} to prove Theorem \ref{thm-endpoint}.
If we multiply \eqref{eqExpansionX} and \eqref{eqExpansionTauX} by $z^{n+1}$ ($z\in\mathbb C$) and sum over $n\ge0$, then we get
\begin{equation}\label{eqExpansionX2}
	\Rho_z(x)=\Psi_z(x)+(zp_c D\ast \Rho_z\ast\Pi^{p_c}_z)(x), \qquad \Tau^p_z(x),
	=\Pi^p_z(x)+(zp D\ast \Tau^p_z\ast\Pi^p_z)(x)
\end{equation}
where we define
\begin{equation}\label{Gz-def}
	\Rho_z(x):=\sum_{n \ge 0} \rho_n(x)z^n, \qquad 
	\Tau^p_z(x):=\sum_{n \ge 0} \tau^{p}_n(x)z^n,
\end{equation}
and
\begin{equation}\label{eqDefPi}
	\Pi^p_z(x):=\sum_{n \ge 0} \pi^p_n(x)z^n, \qquad 
    	\Psi_z(x):=\sum_{n \ge 0} \psi_n(x)z^n.
\end{equation}
The generating functions $\Rho_z(x)$ and $\Tau^p_z(x)$ are power-series in $z$. Since $\sum_x \rho_n(x)=1$ for all $n\geq 0$,
	\begin{equation}
	\label{eqRho1}
    	\sum_x \Rho_z(x)=\frac{1}{1-z}.
	\end{equation}
It follows that the radius of convergence $z_c$ of the power-series $\sum_x \Rho_z(x)$ equals $z_c=1$. Furthermore, we can apply the Fourier transform on \eqref{eqExpansionX} to identify
\begin{equation}\label{eqRhoExp}
    	\hat{\Rho}_z(k)=\frac{\hat\Psi_z(k)}{1-zp_c \Dk\hat\Pi^{p_c}_z(k)},
\end{equation}
and for $\Tau^p_z$,
	\begin{equation}
	\label{eqTauExp}
    	\hat{\Tau}^p_z(k)=\frac{\hat{\Pi}^p_z(k)}{1-z p  \Dk\hat\Pi^p_z(k)}.
	\end{equation}
When we compare \eqref{eqRhoExp} and \eqref{eqTauExp} for $k=0$ we see that $z_c=1$ is also the radius of convergence of $\sum_x \Tau^{p_c}_z(x)$, if $\hat\Psi_z(0) = \sum_x \Psi_z(x)$ is uniformly bounded in $|z|\leq 1$. This condition turns out to be a simple consequence of bounds that we prove later (see Proposition \ref{prop-LEcoefficients}), but until the end of this outline we simply assume this boundedness.

Note that if we set $z=1$ and use \eqref{eqDefPi} and \eqref{eqRho1}, we can identify the critical value $p_c$ in terms of the new lace expansion:
\begin{equation}\label{eqPcPi}
	1=p_c\hat\Pi_1^{p_c}(0)=p_c\sum_{x}\sum_{n \ge 0}\pi^{p_c}_{n}(x).
\end{equation}
(Here we also rely on Remark \ref{rem-HS-exp} below, which explains that the lace-expansion coefficient $\Pi_1^{p_c}$ is also the coefficient for the classical Hara-Slade lace expansion \eqref{HarSlaLE}, i.e., $\Pi_1^{p_c} \equiv \pi^{p_c}$.)

\subsection{Overview of the proof of Theorem \ref{thm-endpoint}} The lace-expansion equations in \eqref{eqExpansionX2} form the starting point of our analysis. In order to prove results on $\hat{\tau}_n(k_n)$ and $\hat{\rho}_n(k_n)$ as stated in Theorem~\ref{thm-endpoint}, we prove asymptotics for $\hat{\Tau}_z(k)$ and $\hat{\Rho}_z(k)$ for $z\in\mathbb C$ with $|z|\leq 1$ and $z$ close to 1. Using a Tauberian Theorem, we can convert these results into results on the coefficients of the generating functions of  $\hat{\Tau}_z(k)$ and $\hat{\Rho}_z(k)$, which are $\hat{\tau}_n(k)$ and $\hat{\rho}_n(k)$, including explicit error estimates. The Tauberian Theorem requires estimates on the generating functions $\hat{\Tau}_z(k)$ and $\hat{\Rho}_z(k)$ with complex~$z$.

\section{A lace expansion for the backbone two-point function}\label{sec-le}

In this section we state and derive the lace expansion for the backbone two-point function. Recall \eqref{taunx-def}, and write $\tau^p_m(x,y)=\tau^p_m(y-x)$. For functions $f,g\colon (\Z^d)^2\to \R$, we define their convolution $f*g$ by
\begin{equation}
	(f*g)(x,y)=\sum_{z} f(x,z)g(z,y).
\end{equation}

The expansion for $\tau^p_n$ derived in this section is valid for any percolation parameter $p$ such that $\Pp(0 \conn \infty) =0$, i.e., when the cluster of $0$ (or any other vertex in $\Zd$) is almost surely finite. The results in \cite{HarSla90a,HeyHofSak08} imply that $\Ppc (0 \conn \infty)=0$ and that the expansion converges for our percolation models (see also Remark \ref{rem-HS-exp} below). When deriving the expansion for the backbone two-point function $\rho_n$ in Section \ref{sec-lace-exp-rho} we need to work with subcritical $p$ first, and only then take the limit $p\nearrow p_c$ or the IIC limit. We therefore write the percolation parameter $p$ as a superscript, e.g.\ $\pi^p_n(x)$. When it is clear from context we will often omit the superscript $p$.

For ease of notation we define 
\begin{equation}
	J^p(x,y) := p D(x,y),
\end{equation}
i.e., $J^p(x,y)$ is the probability that the bond $\{x,y\}$ is occupied.

The lace expansion for $\tau_n^p(x)$ is formulated in the following proposition:
\begin{prop}[Lace expansion]
\label{prop-lace-exp}
For every $x,y\in \Z^d$, any $p$ such that $\Pp(0 \conn \infty) =0$, and $n\geq 0$, there exists a function $(x,y,m,p)\mapsto \pi^p_m(x,y)$ such that
	\begin{equation}
	\tau^p_{n+1}(x,y) = \sum_{m=0}^n (\pi^p_m * J^p * \tau^p_{n-m})(x,y) + \pi^p_{n+1}(x,y).
	\end{equation}
The lace-expansion coefficients $\pi^p_m$ are defined in \eqref{pimNdef} and \eqref{pimdef}. 
\end{prop}
\medskip

The remainder of this section is organized as follows. In Section \ref{sec-fac-lem} we state the main technical tool in the derivation of the lace expansion, the so-called {\it Factorization Lemma}.  In Section \ref{sec-der-lace-exp} we derive the first step of the lace expansion. In Section \ref{sec-compl-lace-exp} we complete the derivation of the lace expansion. In Section~\ref{sec-ext-lace-exp-KJK} we prove a special expansion around interior points, which we need later in Section \ref{sec-tightness}. In Section~\ref{sec-lace-exp-rho} we extend the lace expansion to the backbone two-point function $\rho_n$. Finally, in Section~\ref{sec-non-nega-bds-LE} we make some preparations for bounding the lace-expansion coefficients, and we give an outline of what remains of the proof.
\medskip

We write $\{x\Conn  y\}$ when there exist two bond-disjoint
paths of occupied bonds that connect $x$ to~$y$,
and we adopt the convention that $\{x\Conn  x\}$ is the full probability space.
For $x,y\in {\mathbb Z}^d$, we write $\Piv(x,y)$ for the \emph{ordered list of occupied and pivotal} (directed) bonds, i.e.,
\begin{equation}\label{eqPivDef}
	\Piv(x,y):=
	\begin{cases}
		(b_1,b_2,\dots,b_m)
		\quad&\text{if $x\Conn  \bb_1, \tb_1\Conn  \bb_2, \ldots, \tb_m\Conn  y$, and}\\
		&\qquad\quad \text{$b_1$,\dots,$b_m$ are occupied and pivotal for $x\conn y$};\\
		(\;\;)\quad&\text{if $x\Conn  y$;}\\
		\varnothing\quad&\text{if $x$ and $y$ are not connected}.
	\end{cases}
\end{equation}
Note the difference between ``$(\;\;)$'' ($x$ and $y$ are connected, but without any pivotal bonds) and ``$\varnothing$'' ($x$ and $y$ are not connected). 
With this definition, we partition $\tau_{n+1}$ according to the positions of the $n+1$ pivotal bonds as
\begin{equation}
    \label{tau-m-def}
    \tau_{n+1}(x,y)=\sum_{b_1, \ldots, b_{n+1}} \Pp \big(\Piv(x,y)=(b_1, \ldots, b_{n+1})\big),
    \qquad n\ge 0,
\end{equation}
and $\tau_0(x,y)=\Pp \big(\Piv(x,y)=(\;\;)\big)$. We let $|\Piv(x,y)|$ denote the number of elements in the list $\Piv(x,y)$ (which, by convention, we set to be 0 in the second and third case in \eqref{eqPivDef}). 
Equation~\eqref{tau-m-def} is the starting point of our analysis.


\subsection{The Factorization Lemma}
\label{sec-fac-lem}
Before we can perform the lace expansion, we need to introduce some notation and recall a useful lemma.
\begin{defn}
\label{def-onin2}
\begin{enumerate}

	\item[(i)] Given a (deterministic or random) set of vertices $A$ and a bond configuration 
	$\omega$, we define $\omega_A$, the \emph{restriction of $\omega$ to $A$}, to be
	\begin{equation}
    		\omega_A(\{x,y\}) =
   		\left\{
   		\begin{array}{ll}
  	 		\omega(\{x,y\})   &\text{if  }x,y \in A,\\
   			0           &\text{otherwise},
  	 	\end{array}
   		\right.
	\end{equation}
	for every bond $\{x,y\}$. In other words, we get $\omega_A$ from $\omega$ by making
	every bond that does not have both endpoints in $A$ vacant.

	\item[(ii)] Given a (deterministic or random) set of vertices $A$ and an event $E$, we say that
	 \emph{$E$ occurs in $A$}, and write $\{E\inn A\}$, if $\omega_A\in E$. In other words, 
	 $\{E\inn A\}$ means that $E$ occurs on the (possibly modified) configuration in which every 
	 bond that does not have both endpoints in $A$ is made vacant.  We adopt the convention 
	 that $\{x\conn x\inn A\}$ occurs if and only if $x\in A$. We further say that $E$ 
	 \emph{occurs off} $A$, and write $\{E\off A\}$, when $E$ occurs in $A^c$.

	\item[(iii)] Given a bond configuration and $x \in \Z^d$, we define $\Ccal(x)$ to be the set of
	 vertices to which $x$ is connected, i.e., $\Ccal(x)=\{y \in \Z^d\colon x\conn y\}$. Given a bond
	 configuration and a bond $b$, we define $\tilde{\Ccal}^{b}(x)$ to be the set of vertices 
	 $y \in \Ccal(x)$ to which $x$ is connected in the (possibly modified) configuration in which 
	 $b$ is made vacant. 
	 \end{enumerate}
\end{defn}

In terms of the above definition,
\begin{multline}
    	\big\{b_1\text{ is the first occupied and pivotal bond for }\{x\conn y\} \big\} \\
    	= \big\{x\Conn  \bb_1 \inn \tilde{\Ccal}^{b_1}(x) \big\}
    	\cap \big\{b_1 \text{ occupied}\big\}\cap \big\{\tb_1\conn y\off \tilde{\Ccal}^{b_1}(x) \big\}.
\end{multline}
Similarly, we get the following crucial identity:
\begin{multline}
    	\big\{\Piv(x,y)=(b_1, \ldots, b_{n+1})\big\}\\
    	=\big\{x\Conn  \bb_1 \inn \tilde{\Ccal}^{b_1}(x)\big\}
    	\cap \big\{b_1 \text{ occupied}\big\}\cap \big\{\Piv(\tb_1,y)=(b_2, \ldots, b_{n+1})\off \tilde{\Ccal}^{b_1}(x)\big\}.
\end{multline}
Hence, we can rewrite
\begin{equation}
    \label{tau-m-rewr}
   	\tau_{n+1}(x,y) =\sum_{b_1, \ldots, b_{n+1}} \Pp \big(
    	\{x\Conn  \bb_1 \inn \tilde{\Ccal}^{b_1}(x)\} 
    	\cap \{b_1 \text{ occupied}\}\cap \{\Piv(\tb_1,y)=(b_2, \ldots, b_{n+1})
	\off \tilde{\Ccal}^{b_1}(x)\}\big).
\end{equation}

We next investigate the probabilities on the right-hand side of \eqref{tau-m-rewr}. A useful tool in this analysis is the
\emph{Factorization Lemma.} This lemma is the workhorse of our expansion method.
For a proof, we refer to \cite[Lemma 2.2]{HofHolSla07b}.

\begin{lemma}[Factorization Lemma]
    \label{lem-onin}
    For any $p$ such that $\Pp(0 \conn \infty)=0$, any bond $b = \{\ulb, \olb\} \in\edges$,
    vertex $x$ and events $E$, $F$,
    \begin{align}
    \label{lemcut(i)}
        \shift\Pp \left(E \inn \tilde{\Ccal}^{b}(x),
        F \off \tilde{\Ccal}^{b}(x)\right)
        &=
        \E_{\sss (0)}\left( \indic{E \inn
        \tilde{\Ccal}^{b}_{\sss 0}(x)}
        \E_{\sss (1)}\big[\indi_{\{F\off\tilde{\Ccal}^{b}_{\sss 0}(x)\}}\big]\right).
     \end{align}
     Moreover, when $E\subseteq \{\ulb \in \tilde{\Ccal}^{b}(x),
     \olb \nin \tilde{\Ccal}^{b}(x)\}$, the event on
     the left-hand side of \eqref{lemcut(i)} is independent of the
     occupation status of $b$.
\end{lemma}
\medskip

In the nested expectation on the right-hand side of \eqref{lemcut(i)}, the set $\tilde{\Ccal}_{0}^{b}(x)$ is {\it random} with respect to the outer expectation, but {\it deterministic} with respect to the inner expectation. As is common in the lace-expansion literature for percolation, we have added a subscript ``0'' to $\tilde{\Ccal}^{b}(x)$ and for the same reason we have also added subscripts ``(0)'' and ``(1)'' to the expectations on the right-hand side of \eqref{lemcut(i)} to emphasize this difference. The inner expectation on the right-hand side is with respect to a second, independent percolation model on a second lattice. The second model is thus dependent on the first model via the set $\tilde{\Ccal}^{b}_{0}(x)$.


\subsection{Derivation of the lace expansion: the first step}
\label{sec-der-lace-exp}
Now that we have stated the preliminaries for the lace-expansion derivation, we take the first step towards proving it. By the Factorization Lemma \ref{lem-onin} we may rewrite \eqref{tau-m-rewr} as
\begin{equation}
    \label{tau-m-rewr-b}
    \tau_{n+1}(x,y)=\sum_{b_1, \ldots, b_{n+1}} J (b_1)\, \E_{\sss (0)}\Big[
    \indic{x\Conn  \bb_1 \inn \tilde{\Ccal}^{b_1}_0(x)}
    \P_{\sss (1)}\big(\Piv(\tb_1,y)=(b_2, \ldots, b_{n+1})\off \tilde{\Ccal}^{b_1}_0(x)\big)\Big].
\end{equation}
Here and henceforth we abbreviate $J(b)=J(\bb,\tb)$.  
We can replace the event
$\{x\Conn  \bb_1 \inn \tilde{\Ccal}^{b_1}_0(x)\}$ by the event
$\{x\Conn  \bb_1\}$, since if $\{x\Conn  \bb_1\}$ occurs but
$\{x\Conn  \bb_1 \inn \tilde{\Ccal}^{b_1}_0(x)\}$ does not, then
$\tb_1\in \tilde{\Ccal}^{b_1}_0(x)$. But if this is the case, then
\begin{equation}
    	\P_{\sss (1)}\big(\Piv(\tb_1,y)=(b_2, \ldots, b_m)\off \tilde{\Ccal}^{b_1}_0(x)\big)=0.
\end{equation}
Therefore,
\begin{equation}
	\label{tau-m-rewr-c}
    	\tau_{n+1}(x,y)=\sum_{b_1, \ldots, b_{n+1}} J (b_1) \E_{\sss (0)}\Big[\indic{x\Conn  \bb_1}
	\P_{\sss (1)}\big(\Piv(\tb_1,y)=(b_2, \ldots, b_{n+1})\off \tilde{\Ccal}^{b_1}_0(x)\big)\Big].
\end{equation}
When $n=0$, the situation simplifies somewhat, because then $\Piv(\tb_1,y)=(\;\;)$, so that we can replace $\{\Piv(\tb_1,y)=(b_2, \ldots, b_{n+1})\}$ by $\{\Piv(\tb_1,y)=(\;\;)\}=\{\tb_1\Conn  y\}$.
\medskip

To continue the expansion we introduce some more notation:
\begin{defn}[Pivotals off $A$]
\label{def-off-A}
\color{white}.\color{black}

\begin{enumerate}

\item[(i)] Given a (deterministic) set of vertices $A$, and two vertices $v,y\in \Z^d$, we let $\Piv^{A}(v,y)$ be the collection of pivotal bonds for the event $\{v\conn y\off A\}$.

\item[(ii)] Given a (deterministic) set of vertices $A$, and two vertices $v,y\in \Z^d$, we let 
\begin{equation}
	\tau_{n}^{p,A} (v,y):=\Pp \big(|\Piv^A(v,y)|=n \big)
\end{equation}
(and write $\tau_n^A = \tau_n^{p,A}$ when it is clear what $p$ is).
\end{enumerate}
\end{defn}
\medskip

We can rewrite \eqref{tau-m-rewr-c} in terms of Definition \ref{def-off-A} as
\begin{equation}
	\label{tau-m-rewr-d}
    	\tau_{n+1}(x,y)=\sum_{b_1} J (b_1) \,\E_{\sss (0)} \big[
	\indic{x\Conn \bb_1}\tau_n^{ \tilde{\Ccal}^{b_1}_0(x)}(\olb_1,y)\big], \qquad n\ge0.
\end{equation}
For $m \ge 0$ we define
\begin{equation}\label{e:pinzerodef}
	\pi_m^{p,\sss(0)}(x,y):= \delta_{0,m}\Pp (x\Conn y),
\end{equation}
where $\delta_{0,m}$ is a Kronecker delta, and rewrite \eqref{tau-m-rewr-d} as
\begin{equation}\label{eqRewriteTau}
		\tau_{n+1}(x,y)  =  \sum_{m=0}^n (\pi^{ \sss(0)}_m * J * \tau_{n-m})(x) - 
		\sum_{b_1} J (b_1)\Ep \big[\indi_{\{x \Conn \bb_1\}} 
		\big(\tau_{n}(\tb_1,y) - 
		\tau_{n}^{ \tCcal^{b_1}(x)} (\tb_1,y)\big)\big].
\end{equation}
(The sum over $m$ is obviously not necessary but will turn out to be convenient below.)

\begin{figure}
\includegraphics[width = .8\textwidth]{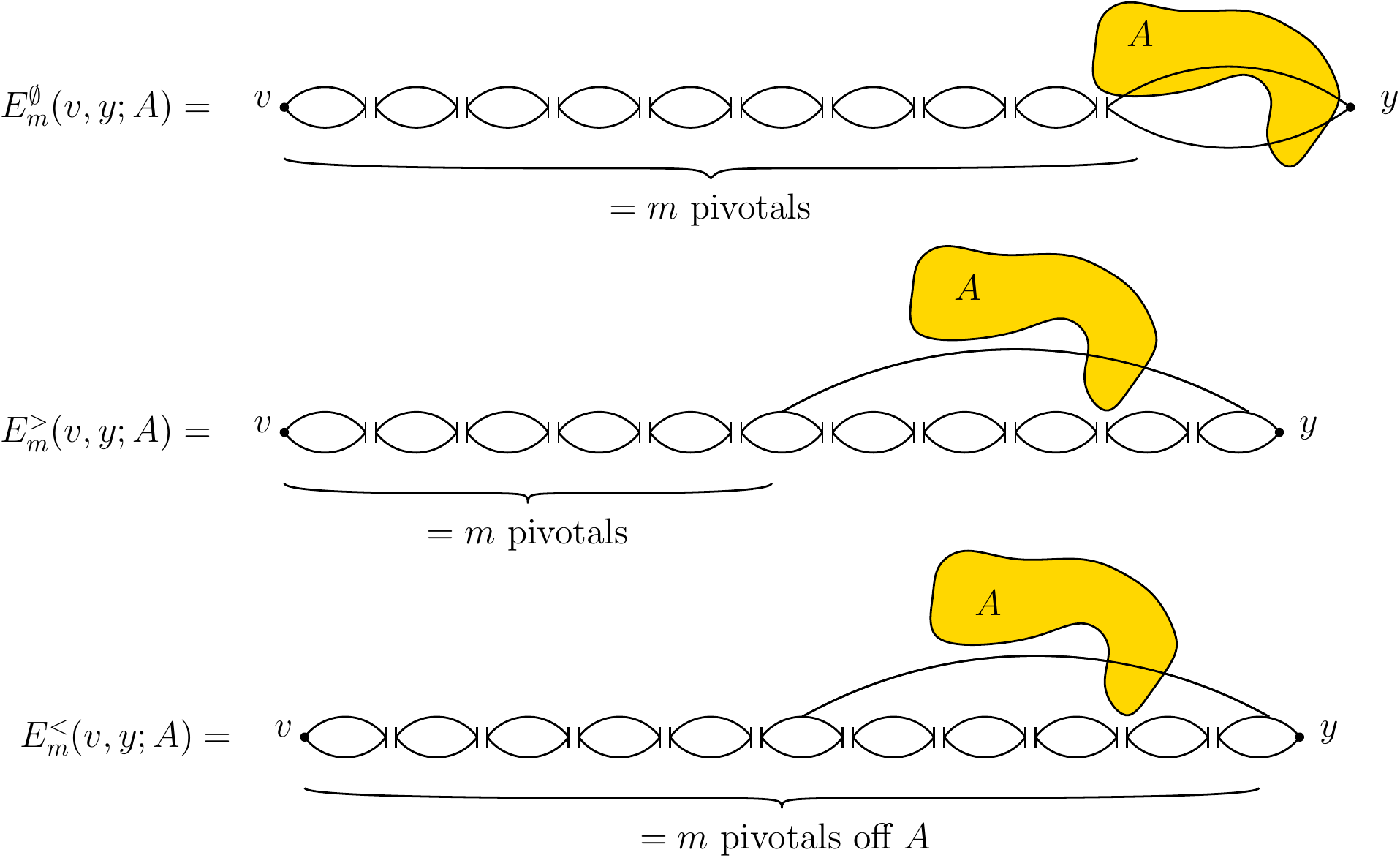}
\caption{\label{fig:Em3} A sketch of $E_m^{\sss \bullet}$ for $\bullet \in \{\varnothing,>,<\}$. The connection between $v$ and $y$ is drawn as a ``string of sausages'', with the pivotal bonds represented with two vertical dashes like $| \, |$, and the double connections between two consecutive pivotal bonds drawn as a ``sausage''.}	
\end{figure}

Since $|\Piv^A(x,y)| \ge |\Piv(x,y)|$ holds on the event $\{x\conn y\off A\}$, we can write, for all $A \subseteq \Zd$,
\begin{equation} \label{eqExpansionTauN}
	\begin{split}
		\tau_n(v,y) - \tau^{A}_n(v,y) & = \Pp \big(|\Piv(v,y)|=n, \Piv^A(v,y) 
		= \emptyset \big) + \Pp \big(|\Piv(v,y)|=n, |\Piv^A(v,y)| > n \big) \\
		& \quad - \Pp \big(|\Piv(v,y)| < n, |\Piv^A(v,y)| = n \big).
	\end{split}
\end{equation}
Now define for integers $m\ge0$ the events
\begin{align}
	\label{E-nothing-def}
	E_{m}^{\sss\varnothing}(v,y;A) := & \big\{|\Piv(v,y)|=m \big\} \cap \big\{ \Piv^A(v,y) = \emptyset \big\} \\
	& \quad \cap \big\{\nexists b \in \Piv(v,y) \text{ s.t. } |\Piv(v,\bb)| \neq |\Piv^A(v,\bb)| \big\},\nn\\
	\label{E->-def}
	E_{m}^{\sss>}(v,y;A) := & \big\{|\Piv(v,y)|=m \big\} \cap \big\{|\Piv^A(v,y)|>m \big\} \\
	& \quad \cap \big\{\nexists b \in \Piv(v,y) \text{ s.t. } |\Piv(v,\bb)| \neq |\Piv^A(v,\bb)| \big\},\nn\\
	\label{E-<-def}
	E_m^{\sss<}(v,y;A) := &  \big\{|\Piv(v,y)|<m \big\} \cap \big\{|\Piv^A(v,y)|=m \big\} \\
	& \quad \cap \big\{\nexists b \in \Piv(v,y) \text{ s.t. } |\Piv(v,\bb)| \neq |\Piv^A(v,\bb)| \big\}.\nn
\end{align}
Here we note that $E_0^{\sss<}(v,y;A) = \emptyset$. See Figure \ref{fig:Em3} for a sketch of these three events. Further, for all $0 \le m<  n$ and $\bullet \in \{\varnothing ,>, <\}$,
\begin{equation}
	F^{\bullet}_{m,n}(x,u,v,y;A) :=  E^{\bullet}_{m}(x,u;A) \cap \big\{(u,v) \in \Piv(x,y) \big\} 
	\cap \big\{|\Piv(v,y)| = n-m-1 \big\}.
\end{equation}

The following lemma gives a partition of the event of interest in terms of the above events:
\begin{lemma}[Partition along cutting bonds]
\label{lem-partition}
For all $n \ge 1$ integer, $v,y \in \Zd$ and $A \subseteq \Zd$,\\
(a)
\begin{equation} 
	\big\{|\Piv(v,y;A)| = n, \Piv^A(v,y) = \emptyset\big\} 
	= E^{\sss\varnothing}_{n}(v,y;A)  \cup \Bigg( \bigcup_{b} \bigcup_{m =0}^{ n-1} 
	F^{\sss\varnothing}_{m,n}(v,\bb,\tb,y;A)\Bigg),
\end{equation}
(b)
\begin{equation} 	
	\big\{|\Piv(v,y)|=n, |\Piv^A(v,y)| > n\big\} =  E^{\sss>}_{n}(v,y;A) \cup 
	\Bigg( \bigcup_{b} \bigcup_{m =0}^{n-1} F^{\sss>}_{m,n}(v,\bb,\tb,y;A)\Bigg),
\end{equation}
(c)
\begin{equation} 
	\big\{|\Piv(v,y)| < n, |\Piv^A(v,y)| = n\big\} = E^{\sss<}_{n}(v,y;A)  \cup 
	\Bigg(\bigcup_{b}  \bigcup_{m =0}^{ n-1}  F^{\sss<}_{m,n}(v,\bb,\tb,y;A)\Bigg),
\end{equation}
where each of the unions is in fact a \emph{disjoint} union.
\end{lemma}
\medskip

The proof of Lemma \ref{lem-partition} follows by partitioning on the first bond $b$ and $m$ for which the event $F^{\bullet}_{m,n}(v,\bb,\tb,y;A)$ holds, or the absence thereof. We leave the details to the reader.
\bigskip

It follows that
\begin{align}
	\tau_n(v,y) - \tau^{A}_n(v,y)  = \, & \Pp (E_n^{\sss\varnothing}(v,y;A)) + \sum_{b}\sum_{m=0}^{n-1} \Pp (F_{m,n}^{\sss\varnothing}(v,\bb,\tb,y;A))\\
	& + \Pp (E_n^{\sss>}(v,y;A)) + \sum_{b}\sum_{m = 0}^{n-1} \Pp (F_{m,n}^{\sss>}(v,\bb,\tb,y;A))\nn\\
	& - \Pp (E_n^{\sss<}(v,y;A)) - \sum_{b}\sum_{m = 0}^{n-1} \Pp (F_{m,n}^{\sss<}(v,\bb,\tb,y;A)).\nn
\end{align}

We continue by investigating the probabilities $\Pp(F_{m,n}^{\bullet}(v,\bb,\tb,y;A))$:
\begin{lemma}[Cutting bond lemma] 
\label{lem-cutting-bond}
For $\bullet \in \{\varnothing, >, <\}$, $0\le m < n$, $A\subset \Zd$, $v,y\in\Z^d$, $b\in\edges$, and all $p$ such that $\Pp(0 \conn\infty)=0$,
\begin{equation}
	\Ep \big[\indi_{F^{\bullet}_{m,n}(v,\ulb,\olb,y;A)}\big] 
	=  J^p (b)\Ep \big[\indi_{E_{m}^{\bullet}(v,\ulb;A)} \tau_{n-m-1}^{p,\tCcal^{b}(v)} (\olb,y)\big].
\end{equation}
\end{lemma}

\proof This follows directly from the Factorization Lemma \ref{lem-onin}, together with the fact that
\begin{equation}\label{e:offexample}
	\indic{E_{m}^{\bullet}(v,\ulb;A)\text{ in }\tCcal^{b}(v)}\tau_{n-m-1}^{p,\tCcal^{b}(v)} (\olb,y)
	=\indicwo{E_{m}^{\bullet}(v,\ulb;A)}\tau_{n-m-1}^{p, \tCcal^{b}(v)} (\olb,y),
\end{equation}
the proof of which is similar to that given right below \eqref{tau-m-rewr-b}.
\qed
\medskip

For $m\ge0$ define
\begin{equation}
	\label{Um-def}
	\Ucal_m(v,y; A) := \indi_{E^{\sss\varnothing}_{m}(v,y; A)}+
	\indi_{E^{\sss>}_{m}(v,y; A)}-\indi_{E^{\sss<}_{m}(v,y; A)}.
\end{equation}
Applying this definition and Lemma \ref{lem-cutting-bond} to \eqref{eqExpansionTauN}, we obtain that
\begin{equation}
	\label{first-exp-taun}
	\tau_n(v,y) - \tau^A_n(v,y)  = \Ep[\Ucal_{n}(v,y;A)]  
	+   \sum_{m=0}^{n-1}\sum_{b} J (b)  \Ep \big[\Ucal_m (v,\bb;A) 
	\tau_{n-m-1}^{\tCcal^{b}(v)} (\tb,y) \big].
\end{equation}

Now define
\begin{equation}
	\pi^{p,\sss(1)}_{0}(x,y) := 0 \quad\text{and}\quad
	\pi^{p,\sss(1)}_{m}(x,y) := \sum_{b_1} J^p (b_1) \E_{\sss (0)}\Big[\indi_{\{x \Conn \bb_1\}} 
	\E_{\sss (1)}\big[\Ucal_{m-1}(\tb,y;\tCcal_0^{b_1}(x))\big]\Big], \quad m\ge1,
\end{equation}
where we write the subscript ``$0$'' in $\tCcal_0^{b_1}(x)$ to indicate that this set is random with respect to $\E_{\sss (0)}$, but fixed with respect to $\E_{\sss (1)}$.
We can then write
\begin{equation}
	\label{e-conv1}
		\tau_{n+1}(x,y) =  (\pi^{\sss(0)} * J * \tau_{n})(x,y) - \pi^{\sss(1)}_{n+1}(x,y)
		  - \sum_{m=1}^{n-1}  (\pi^{\sss(1)}_{m} * J * \tau_{n-m})(x,y)  - R_{n+1}^{\sss(1)}(x,y),
\end{equation}
where $R_{n+1}^{\sss(1)}(x,y)$ is defined such that it contains all the remaining terms, i.e.,
\begin{align}\label{eqRi1}
 	R_{n+1}^{\sss(1)}(x,y) := \,& \sum_{m=1}^{n-1} \sum_{b_1} J(b_1) \sum_{b_2} 
	J (b_2)  \E_{\sss (0)} \Big[\indi_{\{x \Conn \bb_1\}} \E_{\sss (1)} 
	\big[\Ucal_{m-1}(\tb_1,\bb_2;\tCcal_0^{b_1}(x)) \\
	& \qquad \times  (\tau_{n-m}(\tb_2,y) - \tau_{n-m}^{\tCcal^{b_2}_1(\tb_1)}(\tb_2,y))\big]\Big].\nn
\end{align}

We continue by expanding $R_{n+1}^{\sss(1)}(x,y)$. To this end, define
\begin{align}
	\pi^{p,\sss(2)}_{m}(x,y) := \, & \sum_{m_1=1}^{m-1}\sum_{b_1}J^p (b_1) \sum_{b_2} J^p (b_2) 
	\E_{\sss (0)} \Big[\indi_{\{x \Conn \bb_1\}} \E_{\sss (1)}
	\big[\Ucal_{m_1-1}(\olb_1,\ulb_2;\tCcal_0^{b_1}(x))\nn \\
	& \qquad  \times  \E_{\sss (2)} [\Ucal_{m-m_1-1}(\olb_2,y; \tCcal_1^{b_2}(\ulb_1))]\big]\Big]
\end{align}
for $m\ge1$, and let $\pi^{\sss(2)}_{0}(x,y) \equiv 0$. 
We can then insert \eqref{first-exp-taun} into \eqref{eqRi1} and obtain 
\begin{equation}\label{e:Rn1}
	R_{n+1}^{\sss(1)}(x,y) =  \pi^{\sss(2)}_{n+1}(x,y) 
	+ \sum_{m=0}^{n-1} (\pi^{\sss(2)}_m * J * \tau_{n-m})(x,y) 
	+ R_{n+1}^{\sss(2)}(x,y),
\end{equation}
where $R_{n+1}^{\sss(2)}(x,y)$ is defined as 
\begin{align}
	\nonumber
	R_{n+1}^{\sss(2)}(x,y) := \, 
	& \sum_{b_1, b_2, b_3} \Bigg[\prod_{i=1}^3 J (b_i) \Bigg]  \sum_{m_1=1}^{n-1}\sum_{m_2=1}^{n-m_1} \E_{\sss (0)}\Bigg[\indi_{\{x \Conn \bb_1\}} 
	 \E_{\sss (1)} \Big[\Ucal_{m_1}(\tb_1,\bb_2;\tCcal_0^{b_1}(x))\\
	&   {}\times  
	\E_{\sss (2)} \big[\Ucal_{m_2}(\tb_2,\bb_3; \tCcal_1^{b_2}(\tb_1))
	\;(\tau_{n-m_1-m_2}(\tb_3,y)-\tau_{n-m_1-m_2}^{\tCcal^{b_3}_2(\tb_2)}(\tb_3,y))\big]\Big]\Bigg].
	 \label{e:Rn2}
\end{align}
Now we can write
\begin{align}
	\tau_{n+1}(x,y) =\, & \sum_{m=0}^{n} (\pi^{\sss(0)}_m * J * \tau_{n})(x,y) 
	- \pi^{\sss(1)}_{n+1}(x,y)  - \sum_{m =0}^{ n} (\pi_m^{\sss(1)} * J * \tau_{n-m})(x,y) \\
	&+  \pi_{n+1}^{\sss(2)}(x,y) + \sum_{m = 0}^{n}  (\pi^{\sss(2)}_{m} * J * \tau_{n-m})(x,y) 
	+ R_{n+1}^{\sss(2)}(x,y).\nn
\end{align}
This gives the second step of the expansion.
 

\subsection{Completing the expansion}
\label{sec-compl-lace-exp}

Define
\begin{align}
	\label{pimNdef}
	\pi_m^{p,\sss(N)}(x,y) := \sum_{b_1, \dots, b_N} & \left[\prod_{i=1}^N J^p (b_i) \right]
	 \sum_{\substack{m_1, \dots, m_N\geq0 \\ m_1 + \dotsm + m_N = m-N}} 
	\E_{\sss (0)} \Big[\indi_{\{x \Conn \ulb_1\}}  
	\E_{\sss (1)}\big [\Ucal_{m_1}(\olb_1, \ulb_2; \tCcal^{b_1}_0(x))\nnb
	& \times \E_{\sss (2)} \big[\Ucal_{m_2}(\olb_2,\ulb_3; \tCcal_1^{b_2}(\olb_1) )\, \dotsm \, 
	\E_{\sss (N)} \big[\Ucal_{m_N} (\olb_N, y ; \tCcal^{b_N}_{N-1}(\olb_{N-1}))\big] \dotsm \big]\big]\Big].
\end{align}
In terms of this notation the expansion of $\tau_{n+1}$ after $M$ iterations can be written as
\begin{equation} 
	\tau_{n+1}(x,y) = \sum_{m=0}^{n} \sum_{N=0}^M (-1)^N (\pi_m^{\sss(N)} * J * \tau_{n-m})
	(x,y) + \sum_{N=0}^M (-1)^N \pi_{n+1}^{\sss(N)}(x,y) + R_{n+1}^{\sss(M)}(x,y),
\end{equation}
where $R_{n+1}^{\sss (M)}(x,y)$ is the obvious inclusion-exclusion remainder term (similar to $R_{n+1}^{\sss (1)}$ in \eqref{e:Rn1} and $R_{n+1}^{\sss (2)}$ in \eqref{e:Rn2}). Note that $R_{n+1}^{\sss(M)}(x,y)= 0$ when $M>n+1$, since every expansion step reduces the number of remaining pivotals by at least one, and we have only $n+1$ pivotals to begin with.
Define 
\begin{equation}
	\label{pimdef}
	\pi_m^p (x,y) := \sum_{N \ge 0} (-1)^N \pi_m^{p,\sss(N)}(x,y),
\end{equation}
where again we note that $\pi_m^{p,\sss(N)}(x,y)\equiv 0$ when $N>m$, so the sum is in fact finite.
In terms of this notation we obtain
 \begin{equation}
	\label{taun-exp}
 	\tau^p_{n+1}(x,y) = \sum_{m=0}^{n} (\pi^p_m * J^p * \tau^p_{n-m})(x,y) + \pi^p_{n+1}(x,y),
\end{equation}
for all $p$ such that $\Pp(0 \conn \infty)=0$. This completes the derivation of the lace expansion in Proposition \ref{prop-lace-exp}.
\qed

\begin{rk}[Relation to the classical Hara-Slade lace expansion]
\label{rem-HS-exp}
To obtain an expansion for the percolation two-point function $\tau^{p}(x,y) := \Pp(x \conn y)$ for $p$ such that $\Pp(0 \conn \infty)=0$, we only need to sum \eqref{taun-exp} over $n\geq -1$, since $\sum_{n \ge -1} \tau_{n+1}^p(x,y) = \tau^p(x,y)$. This yields
 \begin{equation}
	\label{tau-exp}
 	\tau^{p}(x,y) =\Pi^{p}_1 (x,y) + (\Pi^{p}_1 * J^p * \tau^{p})(x,y).
\end{equation}
where $\Pi^{p}_1 (x,y) :=\sum_{m \ge 0} \pi^p_m(x,y)$. 
Since, by \eqref{Um-def} and \eqref{E-nothing-def}--\eqref{E-<-def},
\begin{equation}
	\sum_{m \ge 0} \Ucal_{m}(a,b; A)=\indicwo{E^{\sss \varnothing}(a,b;A)}\, , 
\end{equation}
with
\begin{equation}
	E^{\sss \varnothing}(a,b;A):=\bigcup_{m \ge 0} E^{\sss \varnothing}_m(a,b;A),
\end{equation}	
we obtain that $\sum_{m \ge 0} \pi^{p, \sss(N)}_m(x,y)$ agrees with \eqref{pimNdef}, except for the fact that every factor $\Ucal_{m}(a,b; A)$ is replaced by 
$\indicwo{E^{\sss \varnothing}(a,b;A)}$,  
because the other two terms cancel each other: 
\begin{equation}
	\bigcup_{m \ge 0} E^{\sss <}_m(a,b;A) = \bigcup_{m \ge 0} E^{\sss >}_m(a,b;A), 
\end{equation}
for all $a,b \in \Zd$ and $A \subset \Zd$.
This retrieves the classical Hara-Slade expansion in \cite{HarSla90a}, and thus \eqref{tau-exp} and \eqref{HarSlaLE} are identical. In particular, this implies that $\Pi_1^{p,\sss(N)}(x)=\sum_{m \ge 0} \pi^{p, \sss(N)}_m(0,x)$ is the {\em classical} lace-expansion coefficient. We can thus use bounds derived elsewhere on these lace-expansion coefficients (and other facts about them), when we are not interested in the dependence on $m$ of $\pi_m^{p, \sss(N)}$.
\end{rk}


\subsection{Extension: lace expansion around an interior point}
\label{sec-ext-lace-exp-KJK}
In this section we extend the above lace-expansion analysis to the case where the set of pivotal bonds is fixed. This extension will be useful below in the proofs of the convergence of finite-dimensional distributions and tightness, and it is also a useful ingredient in the expansion of the backbone two-point function $\rho_n$ derived in the next subsection. 

Recall the definition of $\tau_n$ in \eqref{tau-m-def}. Given a sequence of bonds $(b_i)_{i \ge 1}$, for each $n,m \in \N$, define the vectors $\vec{b}_{[n,m]}:=(b_n, \ldots, b_m)$ (with $b_{[n,m]} := \varnothing$ if $n > m$) and $\vec{b}_{[n]} := \vec{b}_{[1,n]}$. We define the two-point function with fixed pivotals as
\begin{equation}
	\label{tau-n-def-pivs}
	\tau^p_n(x,\vec{b}_{[n]}, y):= \Pp \big(\Piv(x,y)=\vec{b}_{[n]}\big),
	\qquad n\ge1,
\end{equation}
so that \eqref{tau-m-def} can be rewritten as
\begin{equation}
	\label{tau-n-def-rep}
	\tau^p_n(x,y)=\sum_{b_1, \ldots, b_n} \tau^p_n(x,\vec{b}_{[n]}, y),
	\qquad n\ge1.
\end{equation}
To work with fixed pivotals will be useful, for example when dealing with the finite-dimensional distributions $\taunT{r}_{\bn}(\bk)$ as introduced below \eqref{eqDefCnN}.
Such finite-dimensional distributions can be obtained by fixing $\tb_{n_i}=x_i$. 

We start by setting up some notation. Analogous to \eqref{E-nothing-def}--\eqref{E-<-def}, define
\begin{align}
	\label{E-nothing-def-pivs}
	E_{m}^{\sss\varnothing}(x,\vec{b}_{[m]},y;A) := 
	& \{\Piv(x,y)=\vec{b}_{[m]}\} \cap E_{m}^{\sss\varnothing}(x,y;A),\\
	\label{E->-def-pivs}
	E_{m}^{\sss>}(x,\vec{b}_{[m]},y;A) := 
	& \{\Piv(x,y)=\vec{b}_{[m]}\} \cap E_{m}^{\sss>}(x,y;A),\\
	\label{E-<-def-pivs}
	E_m^{\sss<}(x,\vec{b}_{[m]},y;A) := 
	&  \{\Piv^A(x,y)=\vec{b}_{[m]}\} \cap E_{m}^{\sss<}(x,y;A).
\end{align}
Analogously to \eqref{Um-def}, we define
\begin{equation}
	\label{Um-def-pivs}
	\Ucal_m(x,\vec{b}_{[m]},y; A) := \indi_{E^{\sss\varnothing}_{m}(x,\vec{b}_{[m]},y; A)}
	+\indi_{E^{\sss>}_{m}(x,\vec{b}_{[m]},y; A)}-\indi_{E^{\sss<}_{m}(x,\vec{b}_{[m]},y; A)}\, .
\end{equation}
Following the same steps as leading to \eqref{first-exp-taun}, we can write
	\eqan{
	\label{first-exp-taun-pivs}
	&\tau_n(x,\vec{b}_{[n]},y) - \tau^A_n(x,\vec{b}_{[n]},y)\\
	&\qquad= \Ep \big[\Ucal_{n}(x,\vec{b}_{[n]},y;A)\big]  +   \sum_{m \leq n} J (b_{m})  
	\Ep \big[\Ucal_{m-1} (x,\vec{b}_{[m-1]}, \ulb_m ;A) 
	\tau_{n-m}^{\tCcal^{b_{m}}(x)} (\bar{b}_{m},\vec{b}_{[m+1,n]},y)\big].\nn
	}
Iteration of this equation as before leads to
 \begin{equation}
	\label{taun-exp-pivs}
 	\tau_{n}(x,\vec{b}_{[n]}, y) = \sum_{m=1}^{n} \pi_{m-1}(x,\vec{b}_{[m-1]},\bb_{m}) J (b_{m})
	 \tau_{n-m}(\bar{b}_{m},\vec{b}_{[m+1,n]},y) + \pi_{n}(x,\vec{b}_{[n]}, y),
\end{equation}
with
\begin{equation}
	\label{pimdef-pivs}
	\pi_{m}(x,\vec{b}_{[m]}, y):=\sum_{N \ge 0} (-1)^N \pi_{m}^{\sss(N)}(x,\vec{b}_{[m]}, y)	
\end{equation}
and $ \pi_{m}^{\sss(N)}(x,\vec{b}_{[m]}, y)$ the analogue of  $\pi_{m}^{\sss(N)}(x, y)$ in \eqref{pimNdef} where each occurrence of $\Ucal_m(v,y; A)$ is replaced with $\Ucal_m(v,\vec{b}_{[m]},y; A)$ for an appropriate 
$m$ and $\vec{b}_{[m]}$. Further, writing
	\begin{equation}
	s_i=1+\sum_{j=1}^{i} m_j
	\end{equation}
for the total number of pivotals allocated to the expectations $\E_{\sss (0)}$ up to $\E_{\sss (i)}$,
\begin{equation}
\begin{split}
	\label{pimNdef-pivs}
	\pi_m^{p,\sss(N)}(x,\vec{b}_{[m]}, y) := & \sum_{\substack{m_1, \dots, m_N\colon \\ m_1 
	+ \dotsm + m_N = m}} 
	\left[\prod_{i=1}^N J^p (b_{s_i}) \right] 
	\E_{\sss (0)} \Big[\indi_{\{x \Conn \ulb_1\}}  \E_{\sss (1)} \big[\Ucal_{m_1-1}(\olb_{1}, 
	\vec{b}_{[2,s_1-1]}, \ulb_{s_1}; \tCcal^{b_{1}}_0(0))\\
	& \qquad \times \E_{\sss (2)} \big[\Ucal_{m_2-1}(\olb_{s_1},\vec{b}_{[s_1+1,s_2-1]},
	\ulb_{s_2}; \tCcal_1^{b_{s_1}}(\olb_{1}) ) \\
	&\qquad \times \dotsm \, \times \E_{\sss (N)}\big[\Ucal_{m_N-1} (\olb_{s_{N-1}}, 
	\vec{b}_{[s_{N-1}+1,s_N-1]}, y ; \tCcal^{b_{s_N}}_{N-1}(\olb_{s_{N-1}}))\big] 
	\dotsm \big]\big]\Big].
\end{split}
\end{equation}

\begin{figure}
	\includegraphics[width = \textwidth]{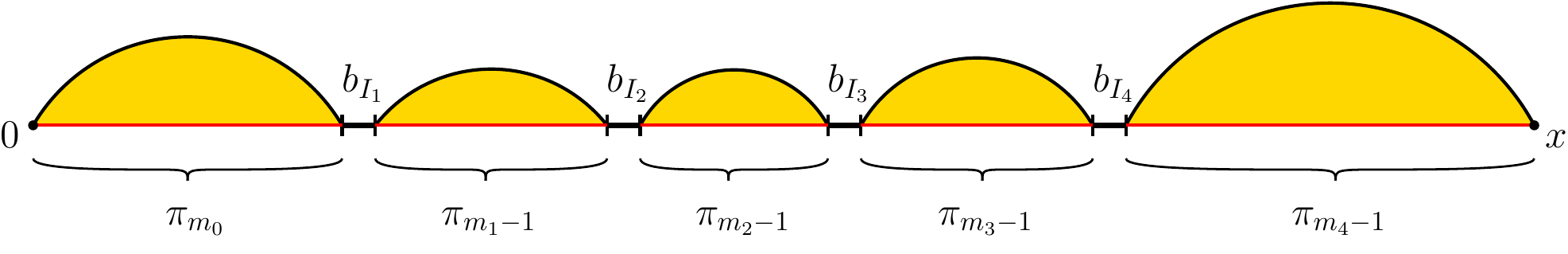}
	\caption{A schematic drawing of a term in \eqref{taun-exp-pivs-indef} with $M=4$. The red line represents the backbone pivotals $\bar{b}_{ [n]}$, the black dashes represent the pivotals $b_{I_1},\dots b_{I_4}$. The yellow filled arcs represent the factors $\pi$. \label{fig:pifixed}}
\end{figure}

Iterating the expansion in \eqref{taun-exp-pivs} $n$ times, we obtain
	 \begin{equation}
	\label{taun-exp-pivs-indef}
 	\tau_{n}(x,\vec{b}_{[n]}, y) = \sum_{M=0}^{n} \sum_{\substack{m_0, \dots, m_M\colon \\ m_0 + \dotsm + m_M = n}}
	\pi_{m_0}(x,\vec{b}_{[m_0-1]},\bb_{m_0}) \prod_{i=1}^M J(b_{I_{i-1}}) \pi_{m_{i}-1}(\bar{b}_{I_{i-1}},\vec{b}_{[I_{i-1}+1,I_{i}-1]},\bb_{I_{i}}),
	\end{equation}
where $I_i=\sum_{j=0}^i m_i$ for $i\in\{0,\ldots, M\}$ and with the convention that $\bb_{I_{M}}=y$ and the empty product, arising when $M=0$, equals 1. 
For fixed $M$ the above sum can be seen as yielding a partition of the interval $[0,n]$ into the $M+1$ disjoint intervals $[I_{i-1}-1, I_{i}]$ (with the convention that $I_{-1}=1$).
See Figure~\ref{fig:pifixed}.
This perspective will be useful in studying the $r$-point functions. It will also be a useful in the derivation of the lace expansion for the backbone two-point function
$\rho_n(y)$, which we derive next. 


\subsection{Lace expansion for the backbone two-point function}
\label{sec-lace-exp-rho}

We have derived a lace expansion for $\tau_n(x,y)$, and shall now extend this to an expansion for $\rho_n(x)=\Piic(S_n=x)$, proving \eqref{eqExpansionX}. 
The IIC backbone two-point function $\rho_n(x)=\Piic(S_n=x)$ is well-defined: it is proved in \cite{HeyHofHul14a} that the IIC measure $\Piic$ exists under the conditions of Theorem~\ref{thm-endpoint}.  
We use the construction of $\Piic$ in \cite{HeyHofHul14a} to derive the lace expansion for $\rho_n(x)$.

The main result of this section is formulated in Proposition \ref{prop-lace-exp-2} below. For the purpose of the expansion we write $\rho_n(x,y):=\rho_n(y-x)$.
\begin{prop}[Lace expansion for IIC-connectivity]\label{prop-lace-exp-2}
Assume that the conditions of Theorem~\ref{thm-endpoint} hold. Then, for all $x,y\in\Z^d$ and $n\ge0$, there exist functions $(x,y,m)\mapsto \pi_m(x,y)$ and $(x,y,m)\mapsto \psi_m(x,y)$ such that
\begin{equation}
	\label{psiexpansion}
	{\rho}_{n+1}(x,y) = \psi_{n+1}(x,y)+ \sum_{m=0}^{n} \left(\pi_m^{p_c} \ast J^{p_c} \ast {\rho}_{n-m}\right)(x,y), 
\end{equation}
where $\psi_m$ is defined in \eqref{psimdef} and $\pi^{p_c}_m$ is the same as in Proposition \ref{prop-lace-exp}. 
\end{prop}

Note that the conditions of Proposition~\ref{prop-lace-exp-2} are stronger than those of Proposition~\ref{prop-lace-exp}. This is because $\rho_n$ is defined with respect to $\Piic$, which requires that we have a firmer understanding of the behavior of the percolation model at $p_c$.

For $n\leq m$ write
\begin{equation}
	\label{pnm-def}
	\pi_{n,m}^p(y, w):=\sum_{\substack{b_1,\dots,b_m:\\ \tb_n=y}} \pi_{m}^p(0,\vec{b}_{[m]}, w)
\end{equation}
and similarly
\begin{equation}
	\label{pnmN-def}
	\pi_{n,m}^{p,{\sss (N)}}(y, w):=\sum_{\substack{b_1,\dots,b_m:\\ \tb_n=y}} \pi_{m}^{p, {\sss (N)}}(0,\vec{b}_{[m]}, w).
\end{equation}

The proof of the lace expansion for $\rho_n$ uses the following lemma:
\begin{lemma}[Boundedness and left-continuity at $p_c$ for lace-expansion coefficients]\label{lem:rhon-exp-lemma}
Under the assumptions of Theorem~\ref{thm-endpoint}:

(a) There exists a $C > 0$ such that uniformly in $n\ge 0$,
\begin{equation}
	\label{pinm-ass}
	\sup_{p\leq p_c} \sum_{w,y } \sum_{m\geq n} |\pi_{n,m}^p(y, w)|<C.
\end{equation}

(b) For all $y\in \Z^d$ and $n\ge 0$,
\begin{equation}
	\label{pinm-ass-2}
	\lim_{p\nearrow p_c}\sum_{l\geq n} \sum_{u}\pi_{n,l}^p(y,u)
	=\sum_{l\geq n} \sum_{u}\pi_{n,l}^{p_c}(y,u).
\end{equation}
\end{lemma}
This kind of lemma is fairly standard in the lace-expansion literature.
\longversion{We prove it in Appendix~\ref{app:cont}.}\shortversion{It is proved in \cite[Appendix~A.4]{HeyHofHulMie17b}.}

\proof[Proof of Proposition~\ref{prop-lace-exp-2} subject to Lemma~\ref{lem:rhon-exp-lemma}]

We extend our notion of ordered pivotal bonds in \eqref{eqPivDef} by writing $\Piv(x,\infty)$ for the ordered list of pivotal bonds for the event 
\begin{equation}
	\{x\conn\infty\}:=\bigcap_{r\in\N}\big\{x\conn\big(\Z^d\setminus [-r,r]^d\big) \big\}.
\end{equation}
For $x=0$, this is $\Piic$-a.s.\ an infinite list, and we define $[\Piv(x,\infty)]_n$ to be the projection onto the \emph{first} $n$ entries of $\Piv(x,\infty)$, together with the convention that $[\Piv(x,\infty)]_n =\emptyset$ whenever $|\Piv(x,\infty)|<n$ (this is a null-event under $\Piic$ and all $n \ge 1$). 

This new notion allows us to rewrite \eqref{rho-n-def} as
\begin{equation}
	\rho_n(y)=\sum_{\substack{b_1,\dots,b_n\colon \\ \tb_n=y}}
	\Piic\big( [\Piv(0,\infty)]_n =(b_1,\dots,b_n)\big).
\end{equation}

The starting point of the expansion is the \emph{Backbone Limit Reversal Lemma} \cite[Lemma~3.1]{HeyHofHul14a} stating that for any $n\in\N$ and bonds $b_1,\dots,b_n$, 
\begin{align}
	\Piic\big( (b_1,\dots,b_n)=[\Piv(0,\infty)]_n \big)
	= \lim_{p\nearrow p_c}\frac1{\chi(p)}\sum_{w}\prob_p
	\big( [\Piv(0,w)]_n =(b_1,\dots,b_n) \big).
\end{align}
We now sum over $b_1, \dots, b_{n}$ with $\olb_n=y$ fixed to obtain
\begin{align}\label{e:rhonlimit}
	\rho_n(y)=
	\lim_{p\nearrow p_c}
	\sum_{m\geq n}
	\sum_{\substack{b_1,\dots,b_m:\\ \tb_n=y}}
	\frac{1}{\chi(p)}\sum_{w}
	\tau^p_m (0,\vec{b}_{[m]},w).
\end{align}
We can write \eqref{taun-exp-pivs} as
 \begin{equation}
	\label{taun-exp-pivs-rep}
 	\tau_{m}^p(0,\vec{b}_{[m]}, w) = \sum_{l=0}^{m-1} \pi_l^p(0,\vec{b}_{[l-1]},\bb_{l}) 
	J^p (b_{l}) \tau_{m-l}^p(\bar{b}_{l},\vec{b}_{[l+1,m]},w) + \pi_{m}^p(0,\vec{b}_{[m]}, w).
\end{equation}

We now perform the sum over $b_1,\dots,b_m$ such that $\tb_n=y$ for $n\le m$. We consider the following three contributions to \eqref{taun-exp-pivs-rep} separately:
\begin{enumerate}
	\item the contribution due to $\pi_{m}^p(0,\vec{b}_{[m]}, w)$,
	
	\item the contributions due to the sum over $0 \le l\leq n-1$, and
	
	\item the contributions due to the sum over $n \le l \le m-1$.
\end{enumerate}
Summed over $\vec{b}_{[m]}$, the contribution (i) simply equals $\pi_{n,m}^p(y, w)$.
We rewrite contribution (ii) as
\begin{equation}
	\sum_{(u,v)} \pi_l^p(u)J^p (u,v) \tau_{n-l, m-l}^p(v,y,w),
\end{equation}
where
\begin{equation}
	\tau_{s, t}^p(v,y,w):=\sum_{\substack{b_1,\dots,b_t\colon\\ \tb_s=y}} \P_p\big(\Piv(v,w)=(b_1, \ldots, b_t)\big),
\end{equation}
and we rewrite contribution (iii) as
\begin{equation}
	\sum_{(u,v)} \pi_{n,l}^p(y,u)J^p(u,v) \tau_{m-l-1}^p(v,w),
\end{equation}
where $\pi_{n,l}(y,u)$ is defined in \eqref{pnm-def}. 
We can thus write \eqref{taun-exp-pivs-rep} as
\begin{align}
	\label{rhon-lim-LE}
	\rho_n(y)&=
	\lim_{p\nearrow p_c}
	\frac{1}{\chi(p)}\sum_{w}\sum_{m\geq n}
	\Big[\sum_{(u,v)} \sum_{l= 0}^{n-1} \pi_l^p(u)J^p (u,v) \tau_{n-l-1, m-l-1}^p(v,y,w)\\
	&\qquad \qquad \qquad\qquad \qquad\qquad+\sum_{(u,v)} \sum_{l\geq n}\pi_{n,l}^p(y,u) J^p (u,v) 
	\tau_{m-l-1}^p(v,w)+\pi_{n,m}^p(y, w)\Big].\nn
\end{align}
We now take the limit as $p\nearrow p_c$. 
By Lemma~\ref{lem:rhon-exp-lemma}(a) the final term in \eqref{rhon-lim-LE} vanishes in this limit, since $\chi(p)\rightarrow \infty$ as $p\nearrow p_c$. 

For the second sum in \eqref{rhon-lim-LE}, we note that by Lemma~\ref{lem:rhon-exp-lemma}(a) we may interchange the infinite sums over $l$ and $m$. Therefore, for every $l\geq n$ fixed, 
\begin{equation}
	\frac{1}{\chi(p)}\sum_{m\geq n}\sum_{w} \tau_{m-l-1}^p(v,w)=\frac{1}{\chi(p)}\sum_{m\geq l+1}\sum_{w} \tau_{m-l-1}^p(v,w)=1.
\end{equation}
By Lemma~\ref{lem:rhon-exp-lemma}(b) and the Dominated Convergence Theorem the second sum in \eqref{rhon-lim-LE} converges to
\begin{equation}
	\label{psimdef}
	 \psi_n(y):= p_c\sum_{l\geq n} \sum_{u}\pi^{p_c}_{n,l}(y,u).
\end{equation}
Finally, for fixed $v,y\in \Z^d$ and $l\leq n-1$, by \eqref{e:rhonlimit},
\begin{multline}
	\lim_{p\nearrow p_c} \sum_{m\geq n}
	\frac{1}{\chi(p)}\sum_{w} \tau_{n-l-1, m-l-1}^p(v,y,w) \\
	= \lim_{p \nearrow p_c} \frac{1}{\chi(p)} 
	\sum_{w} \Pp \big([\Piv(v,w)]_{n-l-1} 
	= (b_1, \dots, b_{n-l-1}), v \conn w \big) 
	=\rho_{n-l-1}(v,y).
\end{multline}
Note that for every $y\in \Z^d$ and $n\in \N$,
\begin{equation}
	\sum_{m\geq n} \frac{1}{\chi(p)}\sum_{w} \tau_{n-l-1, m-l-1}^p(v,y,w)\leq 1,
\end{equation}
and by Lemma~\ref{lem:rhon-exp-lemma}(a) and using that $\pi_l^p(u)=\pi_{0,l}^{p}(0,u)$,
\begin{equation}
	\sup_{p\leq p_c} \sum_{l \ge 0} |\pi_l^p(u)| J^p(u,v)
	\leq p_c \sup_{p\leq p_c} \sum_{w,y} \sum_{m \ge 0} |\pi_{0,m}^p(y, w)|<\infty,
\end{equation}
Furthermore, by \eqref{pinm-ass-2} and again using that $\pi_l^p(u)=\pi_{0,l}^{p}(0,u)$,
\begin{equation}
	\label{pinm-ass-2-cons}
	\lim_{p\nearrow p_c}\sum_{l \ge 0} \sum_{u}\pi_l^p(u)
	=\sum_{l\geq 0} \sum_{u}\pi_{0,l}^{p_c}(0,u) 
	= \sum_{l\geq 0} \sum_{u}\pi^{p_c}_l (u).
\end{equation}
Thus, by dominated convergence, this proves that the first sum in \eqref{rhon-lim-LE} converges to
\begin{equation}
	\sum_{(u,v)} \sum_{l\leq n-1} \pi^{p_c}_l(u) J^{p_c} (u,v) \rho_{n-l-1}(v,y).
\end{equation}
This completes the proof of Proposition~\ref{prop-lace-exp-2}.\qed 
\medskip

For future use, we extend the results to a version of $\psi_m$ with fixed pivotals, as in \eqref{taun-exp-pivs-indef}. 
We assume that an extension of Lemma \ref{lem:rhon-exp-lemma} holds, as formulated in Lemma~\ref{lem:rhon-exp-lemma-2} below.
In its statement, we write, recalling the definition of $\pi_{m}^p(y, \vec{b}_{[m]}, w)$ in \eqref{pimdef-pivs},
	\eqn{
	\label{pimdef-pivs-2}
	\pi_{n,m}^p(y, \vec{b}_{[n]}, w)=\sum_{b_{n+1}, \ldots, b_m} \pi_{m}^p(y, \vec{b}_{[m]}, w).
	}

\begin{lemma}[Boundedness and left-continuity at $p_c$ for lace-expansion coefficients continued]
\label{lem:rhon-exp-lemma-2}
Under the assumptions of Theorem~\ref{thm-endpoint}:

(a) Uniformly in $n\ge 0$ and $\vec{b}_{[n]} \in \edges^n$,
\begin{equation}
	\label{pinm-ass-rep}
	\sup_{p\leq p_c} \sum_{w,y } \sum_{m\geq n} |\pi_{n,m}^p(y, \vec{b}_{[n]}, w)|<\infty.
\end{equation}

(b) For all $y\in \Z^d$, $n\ge 0$ and $\vec{b}_{[n]} \in \edges^n$,
\begin{equation}
	\label{pinm-ass-2-rep}
	\lim_{p\nearrow p_c}\sum_{l\geq n} \sum_{u}\pi_{n,l}^p(y, \vec{b}_{[n]},u)
	=\sum_{l\geq n} \sum_{u}\pi_{n,l}^{p_c}(y,\vec{b}_{[n]}, u).
\end{equation}
\end{lemma}
\longversion{We prove Lemma~\ref{lem:rhon-exp-lemma-2} in Appendix~\ref{app:cont}.}
\shortversion{Lemma~\ref{lem:rhon-exp-lemma-2} is proved in \cite[Appendix~A.4]{HeyHofHulMie17b}.}
\medskip

\begin{figure}
	\includegraphics[width = \textwidth]{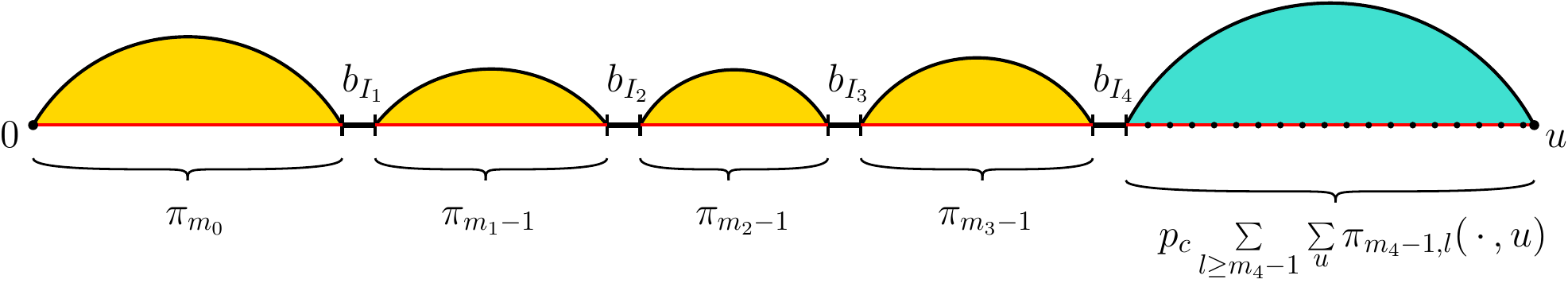}
	\caption{A schematic drawing of a term in \eqref{rhon-exp-pivs-indef} with $M=4$. The turquoise filled arc represents the factor $p_c\sum_{l\geq m_{4}-1} \sum_{u}\pi_{m_{4}-1,l}(\, \cdot \, ,u)$. \label{fig-intervals-a}}
\end{figure}

We omit further details of the extension of the lace expansion to $\rho_{n}(x,\vec{b}_{[n]}, y)$. It states that 
	\begin{align}
	\label{rhon-exp-pivs-indef}
 	\rho_{n}(x,\vec{b}_{[n]}, y) &= \sum_{M \ge 0} \sum_{\substack{m_0, \dots, m_M\colon \\ m_0 + \dotsm + m_M = n}}
	\pi_{m_0}(x,\vec{b}_{[m_0-1]},\bb_{m_0}) \prod_{i=1}^{M-1} J^{p_c} (b_{I_{i-1}}) \pi_{m_{i}-1}(\bar{b}_{I_{i-1}},\vec{b}_{[I_{i-1}+1,I_{i}-1]},\bb_{I_{i}})\nn\\
	&\qquad\qquad \times J^{p_c} (b_{I_{M-1}}) \psi_{m_{M}-1}(\bar{b}_{I_{M-1}},\vec{b}_{[I_{M-1}+1,I_{M}-1]},\bb_{I_{M}}),
	\end{align}
which is identical to \eqref{taun-exp-pivs-indef} except for the last term, which is $\psi_{m_{M}-1}$ rather than $\pi_{m_{M}-1}$. Given the representation in \eqref{psimdef}, this means that the last $\pi_{m_{M}-1}(\,\cdot\,)$ factor is replaced by 
	\begin{equation}
	p_c\sum_{l\geq m_{M}-1} \sum_{u}\pi_{m_{M}-1,l}(\,\cdot\,,u).
	\end{equation}
This can be pictorially expressed by saying that in the expansion of $\rho_{n}(x,\vec{b}_{[n]}, y)$ compared to $\tau_{n}(x,\vec{b}_{[n]}, y)$, there is a $\pi$-factor that ``reaches over the boundary point'' $n$. 
See Figure \ref{fig-intervals-a} for this representation. 
The rewrite in \eqref{rhon-exp-pivs-indef} will be useful to study the (backbone) finite-dimensional distributions $\rhonT{N}_{\bn}(\bk)$ defined in \eqref{eqDefCnN}. 

\subsection{Non-negative upper bounds on the lace-expansion coefficients}
\label{sec-non-nega-bds-LE}
It turns out to be useful below to modify the definition of $\Ucal_m(v,y; A)$ in \eqref{Um-def} by changing the minus into a plus and setting  
\begin{equation}
	\label{UmBar-def}
	\bar\Ucal_m(v,y; A) := \indi_{E^{\sss\varnothing}_{m}(v,y; A)}
	+\indi_{E^{\sss>}_{m}(v,y; A)}+\indi_{E^{\sss<}_{m}(v,y; A)}.
\end{equation}
Define $\bar\pi_m^{\sss(N)}$ and $\bar\pi_m$ as in \eqref{pimNdef} and \eqref{pimdef} with $\Ucal_m$  replaced by $\bar\Ucal_m$ and $\bar\psi_m^{\sss(N)}$ and $\bar\psi_m$ by replacing $\pi_{n,m}$ in \eqref{psimdef} by $\bar{\pi}^{\sss(N)}_{n,m}$ and $\bar{\pi}_{n,m}$, respectively. 
Now, obviously, $\bar\pi_m^{\sss(N)}$ and $\bar\psi_m^{\sss(N)}$ are non-negative.

Moreover, by  \eqref{psimdef}, for any $\vep\ge0$ and all $N \ge 0$,
\begin{equation}
	\label{psi-pi-rep}
	\sum_{m \ge 1} \sum_x m^{\vep}|\psi_m^{\sss(N)}(y)|
	\leq \sum_{m \ge 1} \sum_x m^{\vep}\bar{\psi}_m^{\sss(N)}(y)\le
	\sum_{m \ge 1} \sum_{n\geq m} \sum_{x,y} m^{\vep}\bar{\pi}_{m,n}^{\sss(N)}(x,y)
	\leq \sum_{n \ge 1} \sum_x n^{1+\vep}\bar\pi_n^{\sss(N)}(x).
\end{equation}
This will be helpful in concluding bounds on $\psi_m$ from those derived for $\pi_{m,n}$.
\medskip

In the remainder of this article we will mostly work with $\bar \pi_m$, $\bar \psi_m$, and related non-negative quantities.
\medskip

\noindent
\paragraph{\bf How we proceed.} We have derived the relevant lace expansions and proved \eqref{eqExpansionX} and \eqref{eqExpansionTauX}. The lace-expansion coefficients that arise in \eqref{eqExpansionX} and \eqref{eqExpansionTauX} are defined in terms of the rather complicated events $E_{m}^{\bullet}(v,y;A)$ in \eqref{E-nothing-def}--\eqref{E-<-def}. These events are more involved than the events used in the classical lace expansion, in which the restriction to a fixed number of pivotals does not appear. 
In the next section we derive bounds on the probability of the events $E_{m}^{\bullet}(v,y;A)$ in terms of so-called \emph{lace-expansion diagrams,} for which it will turn out to be helpful to sum out over $m$. In this analysis we will also need weighted sums, where the probability of $E_{m}^{\bullet}(v,y;A)$ is multiplied by $m$ or even $m(m-1)z^m$ before performing the sum. It is here that we need to carefully deal with the non-monotonicity issues of the events $E_{m}^{\bullet}(v,y;A)$.


\section{Bounds on the lace-expansion coefficients in terms of diagrams}
\label{sec-bounds-lace-exp}
In this section we bound the probabilities of the complicated events appearing in the definition of $\pi_m(x)$ in \eqref{pimNdef} and \eqref{pimdef}. More precisely, we will establish diagrammatic bounds on 
	\eqn{
	\label{pi-singly}
	\sum_x \sum_{m \ge 0} \bar{\pi}_m(x),
	\qquad 
	\sum_x  \sum_{m \ge 1} m\bar{\pi}_m(x),
	\qquad
	\sum_x \sum_{m \ge 0} |x|^\delta \bar{\pi}_m(x),
	}
that are phrased only in terms of the two-point function $\tau^{p_c}(x)=\prob_{p_c}(0\conn x)$ (as discussed in more detail in Section \ref{sec-overview}). These bounds are established in \eqref{piexpansion} and Lemmas~\ref{lem:singleLE} and \ref{lem:singleLE-spat} below. Such diagrammatic estimates can thus be bounded using the results in \cite{HarSla90a,HeyHofSak08} only. Furthermore, we also prove diagrammatic bounds on
	\eqn{
	\label{pi-doubly}
	\sum_x \sum_{m \ge 2} m(m-1) z^m \bar{\pi}_m(x)
	}
for $z\in [0,1]$, which also involve $\hat{\Tau}_z(k)$. See in particular Lemma \ref{lem:doubleLE}. We prove these bounds by bounding
	\eqn{
	\sum_{m \ge 0} \prob(E_{m}^{\sss \bullet}(v,y;A))
	\qquad
	\text{and}
	\qquad
	\sum_{m \ge 1} m \prob(E_{m}^{\sss \bullet}(v,y;A))
	}
in terms of two-point functions $\tau^{p}(x,y)=\prob_p(x\conn y)$ using the BK-inequality, while 
	\eqn{
	\sum_{m \ge 2} m(m-1)z^m \prob(E_{m}^{\sss \bullet}(v,y;A))
	}
is bounded in terms of two-point functions $\tau^{p}(x,y)$ and a single $\Tau_z(x,y)$. In the latter case we cannot avoid considering non-monotone events. To bound the probability of these we use the BKR-inequality.

Before starting with this, let use describe the BK and BKR inequalities.


\subsection{Preliminaries: the BK and BKR-inequality}
\label{sec-BK-BKR}
In this section we discuss a tool of the trade in high-di\-men\-sio\-nal percolation: the BK-inequality and its extension, the BKR-inequality.
The {\em van den Berg-Kesten or BK-inequality} \cite{BerKes85} gives an inequality that allows us to bound probabilities of complicated events involving many disjoint connections in terms of products of two-point functions. The most general version is proved by Reimer \cite{Reim00} and allows for events that are possibly non-monotone, and is called the van den Berg-Kesten-Reimer (BKR) inequality (see also~\cite{BorChaRan99}).
Recall from Definition~\ref{def-onin2}(ii) that given an event $E$ and a set of vertices $A$, we say an event $E$ ``occurs in $A$'' if and only if it occurs independently of the status of the set of bonds that do not have both endpoints in $A$. This event is easily generalized to any other set of bonds $B \subseteq \mathbb{B}$. For two events $E,F$, we let $E \circ F$ denote the event
\begin{equation}
	\label{BKR-def}
	E \circ F = \big\{ \omega \colon \exists B_1, B_2 \subseteq \mathbb{B}\text{ with }B_1 \cap B_2 
	= \varnothing, \omega \in \{E \text{ in } B_1
	\} \cap \{F \text{ in } B_2 \} \big\}.
\end{equation}
We refer to the (possibly random) sets of bonds $B_1, B_2$ as {\em witnesses} for the events $E$ and $F$, respectively. The BKR-inequality states that for any $E$ and $F$ depending on finitely many edges,
\begin{equation}
	\label{BKR-ineq}
	\prob_p(E\circ F)\leq \prob_p(E)\prob_p(F).
\end{equation}

A truncation argument to a finite set of bonds then allows one to prove inequalities of the form
	\begin{equation}
	\label{BK-conn}
    	\prob_p\big(\{x\conn y\}\circ\{u\conn v\}\big)\leq
    	\prob_p(x\conn y)\prob_p(u\conn v),
	\qquad u,v,x,y\in\Zd
	\end{equation}
as well as versions with more than two connections. These will be crucial in bounding the lace-expansion coefficients.

The hard part in applying the BKR-inequality is finding good sets of witnesses. One class of events for which witness sets are relatively easy to establish are the \emph{increasing events:} an event $A$ is increasing if $\omega \in A$ implies $\omega' \in A$ for all $\omega'$ such that $\omega(b) \le \omega'(b)$ for all $b \in \mathbb{B}$. For instance, for $\{x \conn y\} \circ \{u \conn v\}$ we can simply choose as witness sets any pair of paths of occupied bonds in $\omega$ from $x$ to $y$ and from $u$ to $v$ that have no bonds in common. When the BKR-inequality is applied to the disjoint occurrence of increasing events, it is often referred to as the BK-inequality, for historical reasons. 

The BKR-inequality \eqref{BKR-ineq} is also not difficult to apply when we have {\em precisely one} non-increasing event, as the decreasing part of the non-increasing event can be witnessed of by all closed bonds. To give an example of how the BKR-inequality in \eqref{BKR-ineq} can be applied, we use it to prove a bound on $\rho_n$ in terms of $\tau_n$:

\begin{lemma} 
\label{lem-rhon-taun}
For every $x\in \Z^d$ and $n\geq 0$,
\begin{equation}
	\rho_{n+1}(x)\leq (J^{p_c} \ast \tau^{p_c}_n)(x).
\end{equation}
\end{lemma}

\proof By the Backbone Limit Reversal Lemma \cite[Lemma 3.1]{HeyHofHul14a}, for any $n\in\N$ and bonds $b_1,\dots,b_{n+1}$, 
	\begin{align}
	\Piic\big( (b_1,\dots,b_{n+1})=[\Piv(0,\infty)]_{n+1} \big)
	= \lim_{p\nearrow p_c}\frac1{\chi(p)}\sum_{w}
	\Pp \big( [\Piv(0,w)]_{n+1} =(b_1,\dots,b_{n+1}) \big).
	\end{align}
Therefore,
\begin{equation}
	\rho_{n+1}(y)=\lim_{p\nearrow p_c}
	\sum_{\substack{b_1,\dots,b_{n+1}:\\ \tb_{n+1}=y}}
	\frac{1}{\chi(p)}\sum_{w}
	\Pp \big([\Piv(0,w)]_{n+1} =(b_1,\dots,b_{n+1})\big).
\end{equation}
We claim that
\begin{equation}
	\big\{[\Piv(0,w)]_{n+1} =(b_1,\dots,b_{n+1})\big\}
	\subseteq
	\big\{\Piv(0,\bb_{n+1})=(b_1,\dots,b_n) \big\}\circ \big\{b_{n+1}\text{ occ.}\big\}\circ \big\{\tb_{n+1}\conn w \big\}.
\end{equation}
The witnesses for the three events on the right-hand side are
\begin{itemize}
	\item[$\rhd$] the bonds (open and closed) attached to $\tilde{\Ccal}^{b_{n+1}}(0)$ 
	except for $b_{n+1}$ for the first,
	
	\item[$\rhd$] $b_{n+1}$ for the second, and 
	
	\item[$\rhd$] an occupied path from $\tb_{n+1}\conn w$ avoiding the bonds in 
	$\tilde{\Ccal}^{b_{n+1}}(0)$ and $b_{n+1}$ for the third.
\end{itemize}
An application of the BKR-inequality \eqref{BKR-ineq} for the first $\circ$, followed by an application of the BK-inequality for the second implies the claim.
\qed


\subsection{Bounds on simple diagrams}
\label{sec-simple-diagrams}
Recall from Definition \ref{def-onin2} that $\tCcal^{b}(x)$ is the cluster of $x$ with the modification that the bond $b$ has been closed. We write $\tCcal^{b}(x)_i$ to denote this cluster on the space associated with $\E_{\sss (i)}$, the $i$-th level in the nested expectation that was introduced to derive the lace expansion and to define the $\pi_m^{\sss(N)}$ diagrams, see e.g.\ \eqref{pimNdef}. To simplify this notation we will write
\begin{equation}
	\tCcal_{i-1} := \tCcal^{b_i} (\olb_{i-1})_{i-1}.
\end{equation}
We also write $\E_{\sss \otimes [0,n]}$ to denote the nested expectation $\E_{\sss (0)}[ \dotsm \E_{\sss (n)}[ \dotsm] \dotsm]$, and we write a subscript $i$ on an event to indicate that the event is associated with the $i$-th level of this nested expectation.

By the definition of $\pi^{\sss(N)}_m(x)$ in \eqref{pimNdef} and the modification \eqref{UmBar-def} that gives $\bar\pi^{\sss (N)}_m(x)$, for $z\in [0,1]$,
\begin{align}
	\label{pimN-bd-1}
	\sum_{m \ge 0} \bar{\pi}_m^{\sss(N)}(x) z^m \le \, & \sum_{m \ge 0} 
	\sum_{b_1,\dots,b_{\sss N}} \Bigg[\prod_{i=1}^N z J (b_i)\Bigg] \\
	& \qquad \times \sum_{\substack{m_1,\dots,m_{\sss N} \ge0: \\ m_1 + \dotsm 
	+ m_{\sss N} = m-N}} \E_{\sss\otimes[0,N]} \Bigg[\indi_{\{0 \Conn \ulb_1\}_0} 
	\prod_{i=1}^N \sum_{\bullet \in \{\varnothing, > , <\}} z^{m_i} 
	\indi_{E^{\bullet}_{m_i} (\olb_i, \ulb_{i+1}; \tCcal_{i-1})_i}\Bigg]\nn\\
	= \,& \sum_{b_1,\dots,b_{\sss N}} \Bigg[\prod_{i=1}^N z J (b_i)\Bigg] 
	\E_{\sss\otimes[0,N]} \Bigg[\indi_{\{0 \Conn \ulb_1\}_0} 
	\prod_{i=1}^N \sum_{\bullet \in \{\varnothing, > , <\}} 
	\sum_{m_i \ge 0} z^{m_i} \indi_{E^{\bullet}_{m_i} (\olb_i, \ulb_{i+1}; \tCcal_{i-1})_i}\Bigg].\nn
\end{align}
In this section we discuss how these coefficients can be bounded. As we will see, we need to give bounds on $\prob_p(E^{\bullet}_{m} (v,y;A))$ and related quantities.
It turns out that since $E^{\bullet}_{m} (v,y;A)$ includes knowledge about the precise number of pivotal bonds, such estimates are quite hard. 
Therefore, instead, we aim to reduce them to estimates in which we sum out over $m$, such as
\begin{equation}
	\label{simple-prob-bds}
	\sum_{m\geq 0}\prob_p(E^{\bullet}_{m} (v,y;A)) 
	\qquad
	\text{and}
	\qquad
	\sum_{m\geq 0}m\prob_p(E^{\bullet}_{m} (v,y;A)).
\end{equation}
These bounds will be set up in such a way that they can be iterated inside the nested expectations. For this, we will have to identify the effect that the bounds on expectations involving $E^{\bullet}_{m_i} (\olb_i, \ulb_{i+1}; \tCcal_{i-1})_i$ have on those involving $E^{\bullet}_{m_{i-1}} (\olb_{i-1}, \ulb_{i}; \tCcal_{i-2})_{i-1}$. As is common in the lace-expansion literature, we achieve such bounds by constructing lace-expansion diagrams: complicated products of (mostly) two-point functions, derived by repeated application of the BK-inequality, that can most easily be represented in the form of Feynman diagrams. This is the key step in the bounds below. It turns out that we need to bound two types of diagrams: (a) those that are finite even when $z=1$, and (b) those that blow up when $z\nearrow 1$. The latter are more tricky, and require us to determine how fast the blow-up is. It is here that we will have to rely on the BKR-inequality, as the event $\{|\Piv(a,b)|=m\}$ is neither increasing nor decreasing.

The rest of this section is structured as follows: We start in Sections \ref{sec-unweighted-diagrams}, \ref{sec-unweighted-lace expansion}, and \ref{sec-singly-weighted-lace expansion} with bounds the diagrams that do not diverge at $z=1$. Then, in Section \ref{sec-doubly-weighted-diagrams} we deal with the more involved divergent diagrams. Our aim in both cases is to provide bounds that are ``plug-and-play'', and that make it easy to obtain insight in how the various bounds are related to one another. This will be achieved by introducing {\em constructions} that describe how to adapt a bound on, say, $\sum_{m\geq 0}\prob_p(E^{\bullet}_{m} (v,y;A))$ so that it may serve as a bound on $\sum_{m\geq 0}m\prob_p(E^{\bullet}_{m} (v,y;A))$. 


\subsection{Bounds on unweighted- and singly-weighted diagrams}
\label{sec-unweighted-diagrams}
In this section we determine several elementary bounds necessary for our bounds on
	\begin{equation}
	\label{singly-weighted-diagrams}
	\sum_x \sum_{m \ge 0} \bar\pi_m^{\sss(N)}(x),
	\qquad
	\sum_x  \sum_{m \ge 1}  m \bar\pi_m^{\sss(N)}(x),
	\qquad
	\sum_x \sum_{m \ge 0}  |x|^{\delta} \bar\pi_m^{\sss(N)}(x).
	\end{equation}
We call the first of these diagrams {\em unweighted,} and the other two {\em singly-weighted diagrams} (for reasons that will become apparent below).
We will derive these bounds from bounds on \eqref{simple-prob-bds}, using the following lemma:

\begin{lemma}[Bounds on $E_m^\bullet$]
\label{lem-bd-E}
For any $p$ and any $A\subseteq \Z^d$ and $\bullet \in \{\emptyset,<,>\}$,
\begin{equation}
	\label{bd-E}
	\sum_{m\geq 0} \prob_p(E^{\bullet}_{m} (v,y;A))\leq \sum_{a\in A} \sum_t\tau(t-v)\tau(a-t)\tau(y-a)\tau(y-t).
\end{equation}
Furthermore,
\begin{align}
	\label{bd-E-line}
	\sum_{m\geq 0} \prob_p \big(E^{\bullet}_{m} (v,y;A)\cap\{v\conn s\}\big) 
	&\leq \sum_{a\in A} \sum_{t,w} \tau(s-w)
	\big[\tau(w-v)\tau(t-w)\tau(a-t)\tau(y-a)\tau(y-t)\nn\\
	&\hskip1.5cm +\tau(t-v)\tau(a-w)\tau(w-t)\tau(y-a)\tau(y-t)\\
	&\hskip1.5cm+\tau(t-v)\tau(a-t)\tau(y-w)\tau(w-a)\tau(y-t)\nn\\
	&\hskip1.5cm+\tau(t-v)\tau(a-t)\tau(y-a)\tau(w-t)\tau(y-w)\big].\nn
\end{align}
\end{lemma}
\medskip

\proof Define 
\begin{equation}
	E^{\bullet}(v,y;A):=\bigcup_{m\geq 0} E^{\bullet}_{m} (v,y;A).
\end{equation}
Since this is a disjoint union, we obtain
\begin{equation}
	\label{e:Ebulletsum}
	\sum_{m\geq 0} \prob_p(E^{\bullet}_{m} (v,y;A))=\prob_p(E^{\bullet}(v,y;A)).
\end{equation}
By taking the union over $m \ge 0$, we effectively forget the pivotal structure of each of the three events $E^{\bullet}$, $\bullet \in \{\varnothing,>,<\}$, but we retain the information that there exists some bond $b = (s,t)$ such that $b$ is the last bond that is contained in both $\Piv(v,y)$ and $\Piv^A(v,y)$, and that the event occurs that there are disjoint paths between (1) $v$ and $t$, (2) $t$ and $y$, (3) $t$ and $a$ for some $a \in A$, and (4) $a$ and $y$. (Figure~\ref{fig:Em3} shows this at a glance.)
Therefore, 
\begin{equation}
	\label{E-bd-circ}
	E^{\bullet}(v,y;A)\subseteq \bigcup_{a\in A}\bigcup_{t} \{v\conn t\}\circ 
	\{t\conn a\}\circ \{a\conn y\}\circ \{t\conn y\}.
\end{equation}
The claim in \eqref{bd-E} now follows by the union bound, followed by repeated application of the BK-inequality \eqref{BKR-ineq}. The claim in \eqref{bd-E-line} can be concluded similarly, by noting that the path $v\conn s$ needs to branch off one of the four disjoint connections in \eqref{E-bd-circ}.
\qed
\bigskip

We see that already the bound in \eqref{bd-E-line} is quite involved, so we introduce some notation that will simplify the analysis of such bounds: 
We introduce the notion of a \emph{diagram}, which in our setting is a formula that is the product of two-point functions $\tau$, pivotal generating functions $z \Tau_z$, and convolutions of $\tau$'s, $z\Tau_z$'s and $J$ functions (e.g.\ $(z \Tau_z * J * \tau * J * \tau)(u-v)$), where possibly some spatial coordinates are summed out). Moreover, in a diagram we will allow some of the two-point functions to carry the label $\mathsf{b}$, so we write $\taub(u-v)$ (say). This label is a bookkeeping tool for the constructions defined below, but does not change the value of the function, i.e., $\taub(\, \cdot \,) \equiv \tau(\, \cdot \,)$.
An important instance of a diagram is
\begin{equation}
	\label{diag-E-def}
	F(v,y,a):= \sum_t  \taub(t-v)\tau(a-t)\tau(y-a)\taub(y-t).
\end{equation}
We extend our terminology by also referring to a sum of diagrams as a diagram.
We refer to any of the factors $\tau$ in a diagram as a {\em line}, and we refer to a labelled two-point function $\taub$ as a \emph{backbone line.} The motivation for this terminology comes from the fact that such formulas have an unambiguous graphical representation that has the form of a diagram, which is illustrated in Figure \ref{fig:diagFvya}. Strictly speaking we never need this correspondence, and in this section it plays no role other than illustration, \longversion{but for reasons of brevity we will use it extensively in the proofs in Appendix \ref{sec-LEcoefficients}, where we also elaborate on this graphical representation of diagrams.}
\shortversion{but it is heavily used in the extended version of this paper, \cite{HeyHofHulMie17b}, where it provides us with a clear, compact, and unambiguous notational device. The precise rules of the graphical representation can be found in \cite[Appendix~A]{HeyHofHulMie17b}.} 

\begin{figure}[tb]
	\includegraphics[width =.85\textwidth]{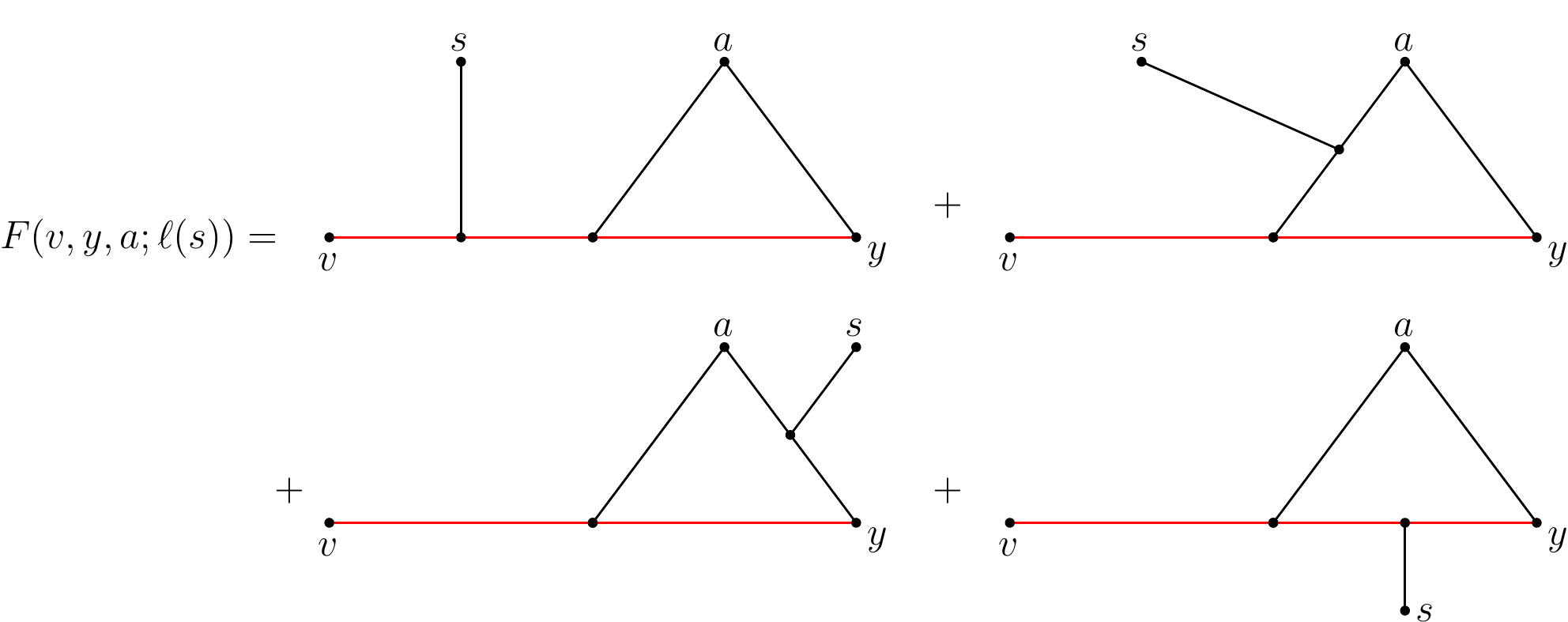}
	\caption{\label{fig:diagFvyaells} A diagrammatic representation of $F(v,y,a; \ell(s))$. The unlabelled vertices represent summation over the variables $t$ and $w$.}
\end{figure}

We now explain some simple rules by which we can express complicated diagrams in terms of simpler ones:

\begin{figure}[t]
	\includegraphics[width = .5\textwidth]{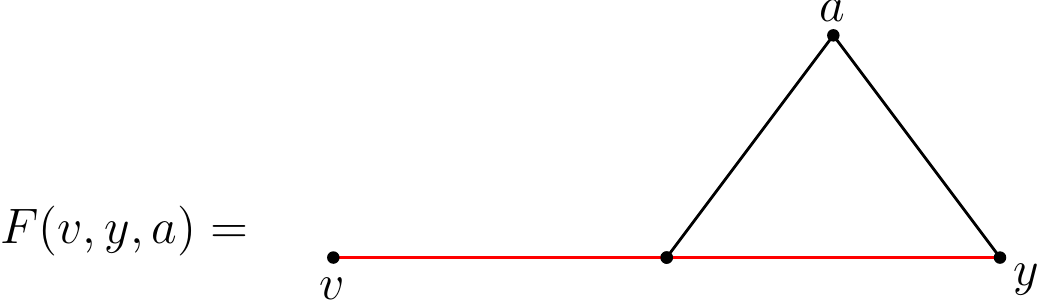}
	\caption{\label{fig:diagFvya} A diagrammatic representation of $F(v,y,a)$. The labelled vertices represent the arguments of $F(v,y,a)$. The unlabelled vertex represents the summation variable $t$. The red lines represent the backbone lines $\taub(t-v)$ and $\taub(y-t)$, and the black lines represent $\tau(a-t)$ and $\tau(y-a)$.}
\end{figure}
\pagebreak[2]
\begin{defn}[Constructions 1 and $\ell$]
\label{def-con1}
\color{white}.\color{black}
\begin{itemize}

    \item[(a)]
     Given a diagram and $w\in \Z^d$, \emph{Construction~$1(w)$} is the operation in which a new vertex $w$ is inserted in one of the lines of the diagram, followed by a summation over all possible lines in which the new vertex can be inserted. Explicitly, this means that a two-point function is replaced by a product of two two-point functions of the same type and label, so $\tau(v-u)$ (say) is replaced by $\tau(v-w) \tau(w-u)$,  and $\taub(v-u)$ (say) is replaced by $\taub(v-w)\taub(w-u)$, and this is done for all lines in the diagram, followed by a sum over possible lines to which the construction may be applied.

    \item[(b)] Given a diagram and $w\in \Z^d$, \emph{Construction~$1$} is the
    operation in which a new vertex $w$ is inserted in one of the lines of the diagram followed by a sum over $w$ and over all possible lines. Explicitly, this means that a two-point function is replaced by a convolution of two two-point functions of the same type, so $\tau(v-u)$ (say) is replaced by $(\tau * \tau)(v-u)$, and $\taub(v-u)$ is replaced by $(\taub * \taub)(v-u)$, followed by a sum over all possible lines to which the construction may be applied.

    \item[(c)] Given a diagram and $s\in \Z^d$, \emph{Construction~$\ell(s)$} is the operation in which Construction~$1(w)$ is performed for $w \in \Zd$, followed by a multiplication with $\tau(s-w)$ and a sum over $w$. Explicitly, this means that a two-point function $\tau(v-u)$ (say) is
    replaced by $\sum_w \tau(v-w)\tau)(w-u)\tau(s-w)$.

    \item[(d)] Given a diagram $G (\,\cdot\,)$, we write $G(\,\cdot\, ; \ell(s))$ for the result of applying Construction~$\ell(s)$ to $G(\,\cdot\,)$.
\end{itemize}
\end{defn}
See Figure \ref{fig:diagFvyaells} for a diagrammatic representation of $F(v,y,a; \ell(s))$.

\medskip

In terms of these bounds and \eqref{diag-E-def} we can compactly summarize \eqref{bd-E} as
\begin{equation}
	\label{bd-E-rep}
	\sum_{m\geq 0} \prob_p(E^{\bullet}_{m} (v,y;A))\leq 
	\sum_{a \in A}  F(v,y,a)
\end{equation}
and \eqref{bd-E-line} as
\begin{equation}
	\label{bd-E-line-rep}
	\sum_{m\geq 0} \prob_p \big(E^{\bullet}_{m} (v,y;A)\cap\{v\conn s\}\big)\leq 
	\sum_{a \in A}  F(v,y,a; \ell(s)).
\end{equation}


\subsection{Diagrammatic bounds on unweighted lace-expansion coefficients}
\label{sec-unweighted-lace expansion}

In this section we use the bounds in the previous section to give diagrammatic bounds on unweighted- and singly-weighted diagrams as in \eqref{singly-weighted-diagrams}.
We aim to repeatedly use \eqref{bd-E-rep} and \eqref{bd-E-line-rep}. Let us describe how this is achieved. Write $\ulb_{\sss N}=x$. Using \eqref{bd-E-rep}, we bound the innermost expectation in \eqref{pimN-bd-1} by
\begin{equation}
	\label{bd-E-rep2}
	\sum_{m_{\sss N}\geq 0} \prob_{\sss (N)} \big(E^{\bullet}_{m_{\sss N}} 
	(\olb_{\sss N-1}, \ulb_{\sss N}; \tCcal_{\sss N-1})_{\sss N} \big)
	\leq 3\sum_{z_{\sss N}}  F(\olb_{\sss N-1}, \ulb_{\sss N},z_{\sss N})
	\indic{z_{\sss N}\in \tCcal_{\sss N-1}}.
\end{equation}
We now see why we needed \eqref{bd-E-line} in the first place. For $N\geq 2$, we insert \eqref{bd-E-rep2} into  \eqref{pimN-bd-1}, and investigate the contribution due to the preceding expectation, $\E_{\sss (N-1)}$,
\begin{equation}
	\sum_{m_{\sss N-1}\geq 0} \prob_{\sss (N-1)}\big(E^{\bullet}_{m_{\sss N-1}} (\olb_{\sss N-2}, \ulb_{\sss N-1}; \tCcal_{\sss N-2})_{\sss N-1} \cap \{ z_{\sss N}\in \tCcal_{\sss N-1}\}\big).
\end{equation}
Now we use \eqref{bd-E-line-rep} to bound
\begin{equation}
	\label{bd-E-line-rep2}
	\sum_{m_{\sss N-1}\geq 0} \prob_{\sss (N-1)} \big(E^{\bullet}_{m_{\sss N-1}} (\olb_{\sss N-2}, \ulb_{\sss N-1}; \tCcal_{\sss N-2})_{\sss N-1}, z_{\sss N}\in \tCcal_{\sss N-1}\big)\leq 
	\sum_{z_{\sss N-1}}  F \big(\olb_{\sss N-2}, \ulb_{\sss N-1},z_{\sss N-1}; \ell(z_{\sss N})\big)\indic{z_{\sss N-1}\in \tCcal_{\sss N-2}}.
\end{equation}
We repeat this procedure for each subsequent nested expectation in \eqref{pimN-bd-1}, until we come to the outermost expectation, $\E_{\sss (0)}[\indi_{\{0 \dbc \bb_1\}} \indi_{\{z_{\sss 1} \in \tCcal_{\sss 0}\}}]$. To bound this expectation we define the diagram
\begin{equation}
	F_{0}(0,x):=\tau(x)^2.
\end{equation}
Using the standard arguments involving the union bound and the BK-inequality we obtain the bound
\begin{equation}
	\label{bd-E-line-rep3}
	\prob_{\sss (0)} \big(\{0\dbc \bb_1\}_{\sss 0}, z_{\sss 1}\in \tCcal_{\sss 0}\big)\leq 
	 F_{0}(0,\bb_1;\ell(z_{\sss 1})).
\end{equation}

We now put all of these terms together.
We can bound \eqref{pimN-bd-1} for $z=1$ from above by
\begin{equation}
\begin{split}
	\label{pimN-bd-2}
	\sum_{m \ge 0} \bar \pi_m^{\sss(N)}(x)\le \, 
	& 3^N \sum_{b_1,\dots,b_{\sss N}} \sum_{z_{\sss 1}, \ldots, z_{\sss N}}
	\Bigg[\prod_{i=1}^N z J(b_i)\Bigg]F_0(0,\bb_1;\ell(z_{\sss 1}))\\
	&\qquad \times\Bigg[\prod_{i=1}^{N-1} 
	F(\olb_{\sss i-1}, \ulb_{\sss i},z_{\sss i}; \ell(z_{\sss i+1}))\Bigg]
	F(\olb_{\sss N-1}, \ulb_{\sss N},z_{\sss N}).
\end{split}
\end{equation}
This bound is an important step, but we can still improve on it by reorganizing our notation.

\begin{figure}
\begin{center}
\includegraphics[width=.8\textwidth]{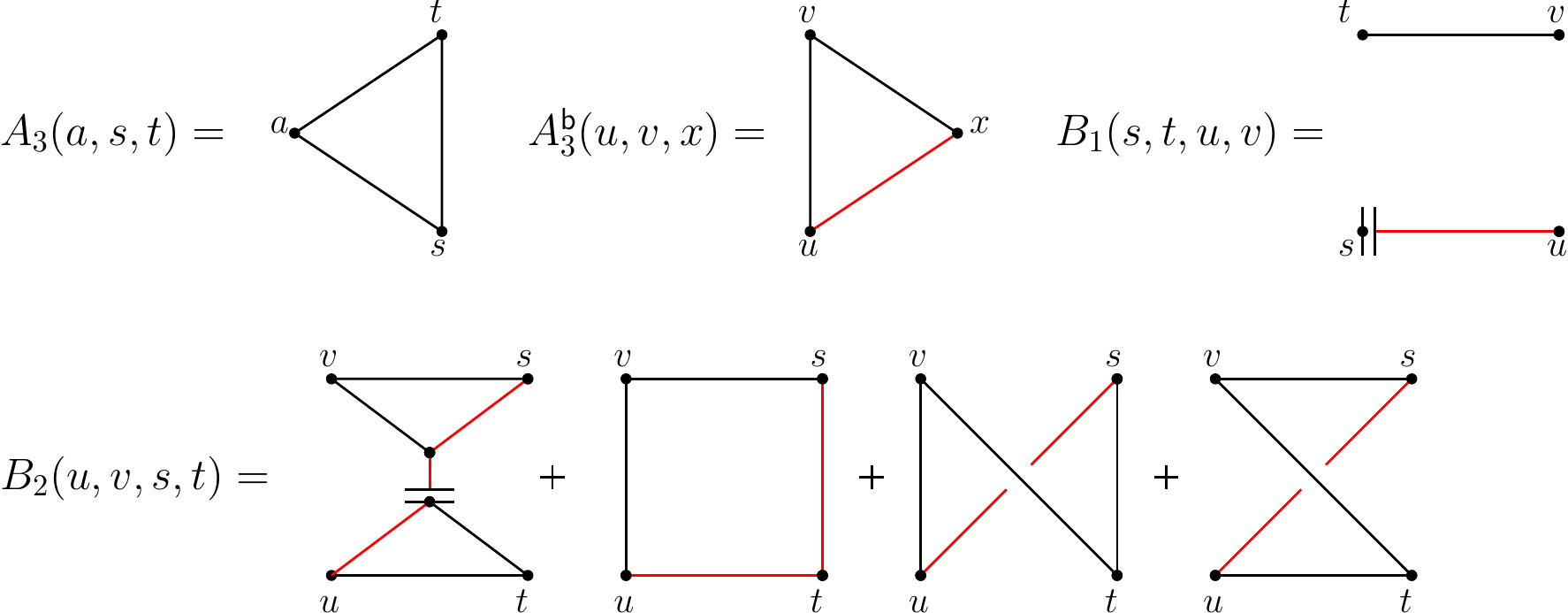}
\caption{Diagrammatic representation of $A_3(s,u,v)$, $A_3^\mathsf{b}(s,u,v)$, $B_1(s,t,u,v)$, and $B_2(u,v,s,t)$. The backbone lines are colored red.}
\label{figA3B1B21}
\end{center}
\end{figure}

To do this we introduce the following functions that are similar to those in other lace expansions (see e.g.\ \cite[pp.\ 109--110]{Slad06}):
Define	
\begin{equation}
	\label{eqDefTauTT}
	 \ttau(x): =(J*\tau)(x).
\end{equation}
and similarly $\ttaub(x) := (J * \taub)(x)$, and
\begin{eqnarray}
	\label{eqA3def} 
	A_3(a,s,t)
	&:=& \tau(s-a)\,\tau(t-a)\,\tau(s-t),\\
	A_3^\mathsf{b}(u,v,x)
	&:=& \taub(x-u)\,\tau(x-v)\,\tau(v-u),\\
	B_1(s,t,u,v)&:=&	\ttaub(u-s)\,\tau(v-t)\\
	B_2^{\sss(0)}(u,v,s,t) &:=&\tau(t-u)\,\tau(s-v)\sum_{a,b}\taub(a-u)\,\ttaub(b-a)\,\tau(t-a)\,
	\tau(b-v)\,\taub(s-b),\\
	B_2^{\sss(1)}(u,v,s,t)&:=&{}\taub(t-u)\,\tau(s-v)\,\tau(v-u)\,\taub(s-t),\\
	B_2^{\sss(2)}(u,v,s,t)&:=&{}\tau(t-v)\,\taub(s-u)\,\tau(v-u)\,\tau(s-t),\\
	B_2^{\sss(3)}(u,v,s,t)&:=&{}\tau(t-u)\,\tau(s-v)\,\tau(v-t)\,\taub(s-u),\\
    	\label{eqB2def} 
	B_2 (u,v,s,t) &:=& \sum_{i=0}^3 B_2^{\sss(i)}(u,v,s,t).
\end{eqnarray}
In Figure \ref{figA3B1B21} we give the graphical representations of these diagrams.

In terms of  \eqref{eqA3def}--\eqref{eqB2def}, and renaming $\bb_i=u_i, \tb_i=v_i$, we can rewrite
\begin{align}
    	F_0(0,u_{\sss 0};\ell(z_{\sss 1})) &\leq
    	\sum_{w_{\sss 0}} A_3(0,u_{\sss 0},w_{\sss 0})\tau(z_{\sss 1}- w_{\sss 0}),
    	\label{F0bd1}\\
    	\sum_{v_{\sss N-1}} J(u_{\sss N-1},v_{\sss N-1}) F(v_{\sss N-1}, x, z_{\sss N})&\leq
    	\sum_{t_{\sss N}} \frac{B_1(u_{\sss N-1},w_{\sss N-1},t_{\sss N},z_{\sss N})}
    	{\tau(z_{\sss N}-w_{\sss N-1})}
    	A_3^\mathsf{b}(x,t_{\sss N},z_{\sss N})\label{F0bd2}
\end{align}
(where we note that the right-hand side of \eqref{F0bd2} actually does not depend on $w_{\sss N-1}$).
Similarly,
\begin{align}
    	&\sum_{v_{i-1}}J(u_{\sss i-1},v_{\sss i-1})
	F(v_{\sss i-1}, u_{\sss i},z_{\sss i}; \ell(z_{\sss i+1}))\nonumber\\
    	&\qquad \qquad \leq
    	\sum_{w_i, t_i} \frac{B_1( u_{\sss i-1},w_{\sss i-1},  t_i, z_i)}{\tau(z_i- w_{\sss i-1})}
    	B_2(t_i,z_i, u_i, w_{i})
    	\tau(z_{\sss i+1}- w_i)
    	\label{F1bd}
\end{align}
(where again the right-hand side of \eqref{F1bd} does not depend on $w_{\sss i-1}$).

Using \eqref{eqA3def}--\eqref{eqB2def}, we can rewrite the right-hand side in \eqref{pimN-bd-2} as 
\begin{align}
	\label{piexpansion}
	\sum_{m \ge 0}
	\bar \pi_m^{\sss (N)} (x)
	&\le 3^N
	\sum_{s_1,\dots,s_{\sss N}}
	\sum_{t_1,\dots,t_{\sss N}}
	\sum_{u_1,\dots,u_{\sss N}}
	\sum_{v_1,\dots,v_{\sss N}}
	A_3(0,s_1,t_1)\\
	 &\qquad \times \Bigg[\prod_{j=1}^{N-1}B_1(s_j,t_j,u_j,v_j)\,B_2(u_j,v_j,s_{j+1},t_{j+1})\Bigg]\nn\\
	&\qquad\qquad \times B_1(s_{\sss N},t_{\sss N},u_{\sss N},v_{\sss N})\,A_3^\mathsf{b}(u_{\sss N},v_{\sss N},x) =: \Phi^{\sss (N)}(x).\nn
\end{align}
Here we note that the product of the factors $\tau (z_i-w_{\sss i-1})$ is telescoping, and thus cancels out. Again, note the similarity with the usual lace-expansion diagrams for percolation, e.g.\ \cite[(10.53)]{Slad06}.

Since we'll encounter below several large diagrams related to $\Phi^{\sss (N)}$, we introduce the following abbreviated notation where we omit the arguments of the diagrams:
\begin{equation}\label{e:pimnbdshort}
	\sum_{m \ge 0}
 	\Phi^{\sss (N)}(x) = 3^N \sum_{(\Zd)^{4N}} A_3 \Bigg[ \prod_{j=1}^{N-1} 
	B_1 B_2 \Bigg] B_1 A_3^\mathsf{b}(x). 
\end{equation}
Moreover, for completeness we define the diagram
\begin{equation}\label{e:phizerodef}
	\Phi^{\sss (0)}(x) := \tau(x)^2,
\end{equation}
so that by \eqref{e:pinzerodef} and the BK-inequality, $\sum_{m \ge 0}  \bar \pi_m^{\sss (0)}(x) \le \Phi^{\sss (0)}(x)$.

\begin{figure}
\begin{center}
\includegraphics[width= .9\textwidth]{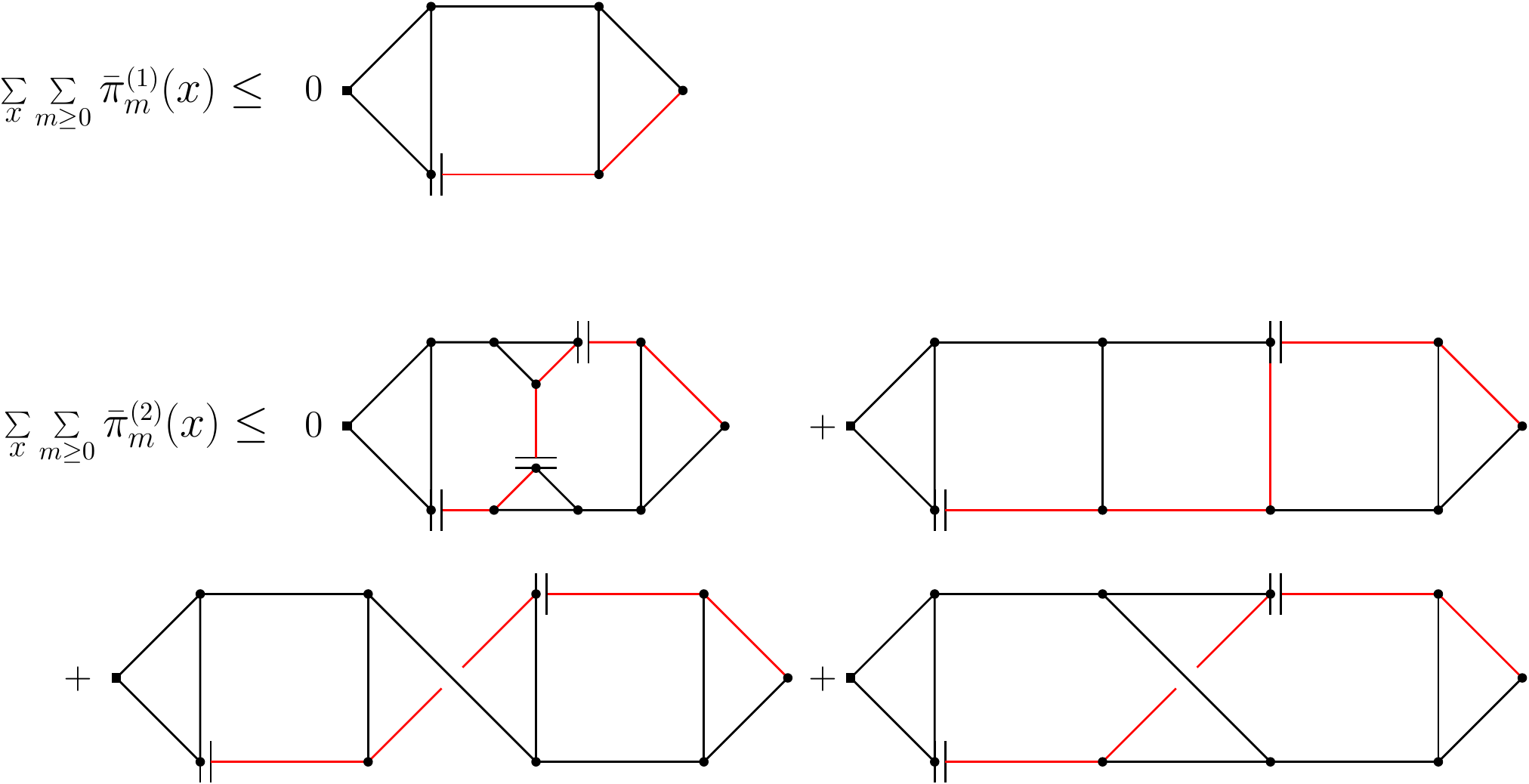}
\caption{The diagrams $\Phi^{\sss (1)}(x)$ and $\Phi^{\sss (2)}(x)$ bounding $\sum_x \sum_{m \ge 0}
	\bar \pi_m^{\sss (1)} (x)$ and $ \sum_x \sum_{m \ge 0}
	\bar \pi_m^{\sss (2)} (x)$ (omitting the factors $3$ and $3^2$ on the right-hand side).}
\label{figPiN12}
\end{center}
\end{figure}

The resulting diagrams for $N=1$ and $N=2$ are shown in Figure \ref{figPiN12}.


\subsection{Bounds on singly-weighted diagrams and lace-expansion coefficients}
\label{sec-singly-weighted-lace expansion}

We now state two weighted versions of the above bounds. One, in which the weight of the coefficient $\bar \pi_m^{\sss (N)}$ is given by $m$, and another in which the weight is given by $|x|^{\delta}$ for some small $\delta >0$. We also give a similar bound for the quantity in Lemma~\ref{lem:rhon-exp-lemma}(a).

\begin{figure}
	\includegraphics[width =.85\textwidth]{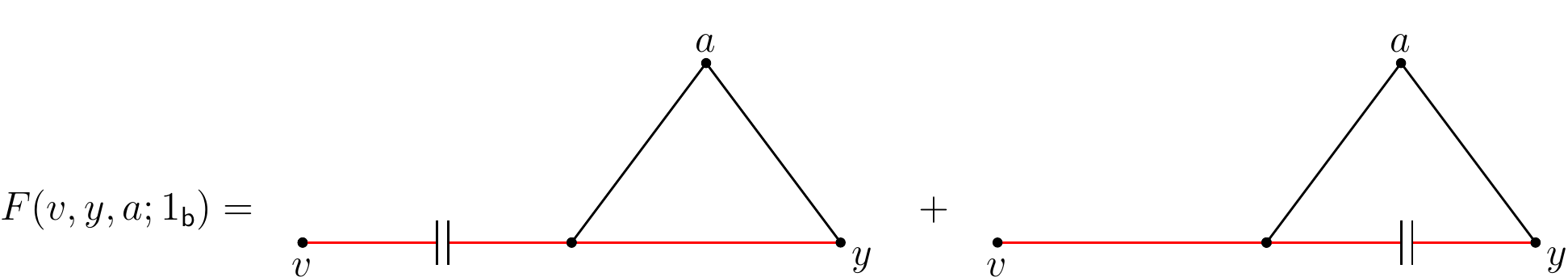}
	\caption{\label{fig:diagFvyaoneb} A diagrammatic representation of $F(v,y,a; 1_\mathsf{b})$. The unlabelled vertex corresponds to the summation variable $t$. The red line that is broken by the double vertical dashes $| \, |$ represents the function $(\taub * J * \taub)(t-v)$ in the first term, and $(\taub * J * \taub)(y-t)$ in the second term.}
\end{figure}

We start with the $m$-weighted diagrams. For this, we define a version of Construction 1 that is only applied to backbone lines, and that inserts a bond rather than a vertex: 

\begin{defn}[Construction $1_\mathsf{b}$]
\label{def-con1b}
\color{white}.\color{black}
\begin{itemize}
    
    \item[(a)] Given a diagram and a bond $(s,t)\in \Z^d\times \Z^d$, \emph{Construction~$1_\mathsf{b}(s,t)$} is the operation in which a new bond $(s,t)$ is inserted in one of the backbone lines of the diagram and is multiplied by $J(s,t)$. Explicitly, this means that a backbone two-point function $\taub(v-u)$ (say) is replaced by $\taub(v-s)J(s,t) \taub(t-u)$. 
    
    \item[(b)] Given a diagram with a certain collection of backbone lines, and $(s,t)\in \Z^d\times \Z^d$, \emph{Construction~$1_\mathsf{b}$} is the operation in which a new bond $(s,t)$ is inserted in one of the backbone lines of the diagram and is multiplied by $J(s,t)$, followed by a sum over $(s,t)$. Explicitly, this means that a backbone two-point function $\taub(v-u)$ (say) is replaced by $(\taub * J * \taub)(v-u)$.
    
    \item[(c)] Given a diagram $G (\,\cdot\,)$, we write $G(\,\cdot\,; 1_\mathsf{b})$ for the result of applying Construction~$1_\mathsf{b}$ to $G (\,\cdot\,)$.
\end{itemize}
\end{defn}
See Figure \ref{fig:diagFvyaoneb} for a diagrammatic representation of $F(v,y,a; 1_\mathsf{b})$.
\medskip

Using Definition \ref{def-con1b}, we can prove the following bounds on the weighted $E$ probabilities:
\begin{lemma}[Bounds on weighted $E$]
\label{lem-bd-weighted-E}
For any $p$, any $A\subseteq \Z^d$ and $\bullet \in \{\emptyset,<,>\}$,
\begin{equation}
	\label{bd-E-m}
	\sum_{m\geq 1} m \prob_p(E^{\bullet}_{m} (v,y;A))\leq \sum_{a\in A} F(v,y,a; 1_\mathsf{b}).
\end{equation}
Furthermore,
\begin{equation}
	\label{bd-E-line-m}
	\sum_{m\geq 1} m \prob_p \big(E^{\bullet}_{m} (v,y;A) \cap \{v \conn s\} \big)
	\leq \sum_{a\in A} F(v,y,a; \ell(s), 1_\mathsf{b}).
\end{equation}
\end{lemma}
The constructions $\ell(s)$ and $1_\mathsf{b}$ commute, but some of the constructions introduced below do not. Hence, we now fix the notational convention that given a diagram $G( \, \cdot \,)$ and constructions $\mathsf{a}_1,\dots,\mathsf{a}_n$, we write $G(\, \cdot \, ; \mathsf{a}_1,\dots, \mathsf{a}_n)$ for the diagram where the constructions are applied \emph{in the given order} (read from left to right).

\proof We start with $\bullet \in \{\varnothing,>\}$. 
The event $E^{\bullet}_m(v,y;A)$ with $\bullet \in \{\varnothing,>\}$ implies that $|\Piv(v,y)| = m$ occurs. Observe the following equivalence: 
\begin{equation}\label{e:msumispivsum}
	|\Piv(v,y)| = m \qquad \text{ if and only if } \qquad \sum_{b} \indic{b \in \Piv(v,y)} = m.
\end{equation}
 So as an alternative to summing over $m$, we may sum over all pivotal bonds for the connection between $v$ and $y$. We use this equivalence to rewrite
\begin{equation}
	 m\prob_p(E^{\bullet}_{m} (v,y;A))=\sum_{b} \prob_p \big(E^{\bullet}_m(v,y;A)\cap \{b\in \Piv(v,y)\} \big).
\end{equation}
Now we can sum over $m$, using \eqref{e:Ebulletsum}, to write
\begin{equation}
		 \sum_{m \ge 1} m\prob_p(E^{\bullet}_{m} (v,y;A))
		 =\sum_{b} \prob_p \big(E^{\bullet} (v,y;A)\cap \{b\in \Piv(v,y)\}\big).
\end{equation}
Having summed over $m \ge 0$, we have lost the information about the number of pivotals, but we have retained the structure of the diagram. By the same argument as in the proof of Lemma \ref{lem-bd-E} we can conclude that if $E^{\bullet}(v,y;A)\cap \{b\in \Piv(v,y)\}$ occurs, then there must exist disjoint connections (1) from $v$ to $t$ (the end-point of the last pivotal in $\Piv(v,y)$), (2) from $t$ to some point $a \in A$, (3) from $t$ to $y$, and (4) from $a \in A$ to $y$. Moreover, we know that the path from $v$ to $t$ passes through the bond $b$, and that this bond must be occupied. So applying the union bound and the BK-inequality we bound
\begin{equation}
\begin{split}
	\sum_{m\geq 1} m\prob_p(E^{\bullet}_{m} (v,y;A)) 
	&\le \sum_{a \in A} \sum_{t} (\taub * J * \taub)(t-v) \tau(a-t) \taub(y-t) \tau(y-a)\\
	& \le \sum_{a \in A} F(v,y,a; 1_{\mathsf{b}}).
\end{split}
\end{equation}
(Note that the middle quantity corresponds to the first term in $\sum_{a \in A} F(v,y,a;1_{\mathsf{b}})$ in Figure \ref{fig:diagFvyaoneb}.)

Similarly, if $\bullet = <$, then each configuration that has $|\Piv^A(v,y)|=m$ can be summed over in the same way, using the $\Piv^A(v,y)$ analogue of \eqref{e:msumispivsum}:
\begin{equation}
	\sum_{m\geq 1} m\prob_p(E^{\sss<}_{m} (v,y;A))=\sum_{b} 
	\Pp \big(E^{\sss<}(v,y;A)\cap \{b\in \Piv^A(v,y) \} \big).
\end{equation}
Now proceed as before. This time, however, we need to distinguish two cases: either $b \in \Piv(v,y)$, or $b \in \Piv^A(v,y) \setminus \Piv(v,y)$. In the former case, the bond is inserted on the line from $v$ to $t$, while in the latter case, the bond is inserted on the line from $t$ to $y$. (See Figure~\ref{fig:Em3}.) Hence we get the bound that involves both terms of $F(v,y,a;1_{\mathsf{b}})$ in this case.

 The bound on \eqref{bd-E-line-m} is similar, and we omit the proof.
\qed
\bigskip

Note that in \eqref{bd-E-line-m}, the order in which we perform Constructions $\ell(s)$ and $1_\mathsf{b}$ is irrelevant, due to the fact that we update the labels of backbone lines after the application of Construction $1_\mathsf{b}$. As a result, we may immediately apply the construction $1_\mathsf{b}$ to the diagram $\Phi^{\sss (N)}$ for an upper bound on the lace-expansion coefficients weighted by $m$:

\begin{lemma}[Singly-weighted LE coefficients]
\label{lem:singleLE} 
For $N \ge 1$,
\begin{equation}
	\sum_{m \ge 1} m  \bar \pi_m^{\sss (N)}(x)  \le \Phi^{\sss (N)}(x; 1_\mathsf{b}).
\end{equation}
\end{lemma}

The right-hand side can be written more explicitly as
\begin{equation}
\begin{split}
	\Phi^{\sss (N)}(x; 1_\mathsf{b})
	& = 3^N \sum_{k=1}^{N-1}  \sum_{(\Zd)^{4N}} A_3 \Bigg[\prod_{j=1}^{k-1} B_1 B_2\Bigg]
	\big(B_1(1_\mathsf{b}) B_2 + B_1 B_2(1_\mathsf{b})\big) \Bigg[\prod_{j'=k+1}^{N-1} 
	B_1 B_2\Bigg] B_1 A_3^\mathsf{b}(x) \\
	& \quad +   \sum_{(\Zd)^{4N}} A_3 \Bigg[\prod_{j=1}^{N-1}B_1 B_2 \Bigg] B_1(1_\mathsf{b}) 
	A_3^\mathsf{b}(x) + \sum_{(\Zd)^{4N}} A_3 
	\Bigg[\prod_{j=1}^{N-1}B_1 B_2 \Bigg] B_1 A_3^\mathsf{b}(x; 1_\mathsf{b}).
\end{split}
\end{equation}
\medskip

Now recall the definition of $\pi_{n,m}^p(y,w)$ in \eqref{pnm-def}. A first step towards proving Lemma~\ref{lem:rhon-exp-lemma}(a) is to apply diagrammatic estimates. Following the same reasoning as above, but now using construction $1_{\mathsf{b}}(b_n)$ to fix the position of the $n$th pivotal bond, we can bound
\begin{equation}
		\sup_{p \le p_c} \sum_{w,y} \sum_{m \ge n} |\pi_{n,m}^{p, {\sss (N)}} (y,w)| 
		 \le \sup_{p \le p_c} \sum_{w,y} \sum_{m \ge n} \bar \pi_{n,m}^{p,{\sss (N)}}(y,w) 
		 \le \sum_{w,y} \sum_{v : (v,y) \in \edges} \Phi^{\sss (N)}(w ; 1_{\mathsf{b}}(v,y)).
\end{equation}
Incorporating the sum over $v$ and $y$ in the diagram yields:
\begin{lemma}[LE coefficients with a fixed pivotal bond]\label{lem:fixedbond}
\begin{equation}
	\sup_{p \le p_c} \sum_{w,y} \sum_{m \ge n} |\pi_{n,m}^{p, {\sss (N)}}(y,w)| \le
	 \sum_w \Phi^{\sss (N)}(w ; 1_{\mathsf{b}}).
\end{equation}
\end{lemma}
\smallskip
Recalling the definition of $\psi_m (x)$ from \eqref{psimdef} it also immediately follows that:
\begin{lemma}[Diagrammatic bound on $\psi^{\sss (N)}_m$]
\label{lem:psidiagram}
\begin{equation}
	\sum_{m \ge 0} \psi_m^{\sss (N)}(x) \le p_c \sum_{v,w} \Phi^{\sss (N)}(w ; 1_{\mathsf{b}}(v,x))
\end{equation}
for all $N \ge 0$.
\end{lemma}

The bound on $\sum_{m \ge 0} |x|^{\delta} \bar \pi_m^{\sss(N)}(x)$ is comparable. We use a similar (but simpler) approach as in the proof of  \cite[Proposition 2.5]{HeyHofHul14b}, to which we refer for the full details of the bound. Here we shall only give the important steps and withhold several details.

Note that in the diagram $\Phi^{\sss (N)}(x)$ on the right-hand side of \eqref{piexpansion} there exists a (non-unique) path of lines that connects 0 to $x$. For instance, the path that starts along the top and swaps from top to bottom every time a factor $B_2^{\sss(2)}$ is present. Write $x_1,\ldots, x_M$ for the displacements along the this path, where $M=2N+1$. Then we can bound
\begin{equation}\label{e:xdeltabd}
	|x|^{\delta}\leq M^{\delta}\big(|x_1|^{\delta}+\cdots
	+ |x_M|^{\delta}\big).
\end{equation}

The resulting bound is very similar to \eqref{piexpansion}, differing in that it is the sum over all distinct diagrams where the lines $\tau(x_i)$ in \eqref{piexpansion} are replaced by $|x_i|^{\delta}\tau(x_i)$. 
We again use a construction to formalize these modifications:
\begin{defn}[Construction $1^{x,y}_\delta$]
\label{def-con1delta}
 Given a diagram $G(\, \cdot \, , x,y)$ containing a simple path of $M$ lines connecting $x$ to $y$, \emph{Construction~$1^{x,y}_\delta$} is the operation in which a simple path of lines connecting $x$ to $y$ is chosen according to some predetermined but arbitrary rule, and where one of the lines from the path, $\tau(v-u)$ (say) is multiplied with $M^\delta |v-u|^{\delta}$. This is summed over all $M$ possible choices of lines in the path from $x$ to $y$. We write $G(\, \cdot \, , x,y; 1^{x,y}_\delta)$ for the result of applying Construction $1^{x,y}_\delta$ to $G(\, \cdot \, , x,y)$. We write $1_{\delta}^x$ for the construction $1_{\delta}^{0,x}$.
\end{defn}
\medskip

The above reasoning, inequality \eqref{e:xdeltabd} in particular, proves the following:
\begin{lemma}[Spatially weighted LE coefficients]\label{lem:singleLE-spat} For $N \ge 0$,
\begin{equation}
	\sum_{m \ge 0} |x|^{\delta} \bar \pi_m^{\sss(N)}(x) \le \Phi^{\sss (N)}(x ; 1^x_\delta).
\end{equation}
\end{lemma}


\subsection{Bounds on doubly-weighted diagrams and lace-expansion coefficients}
\label{sec-doubly-weighted-diagrams}

We continue with bounds on doubly-weighted diagrams, i.e., diagrams in which the factor of $m$ in \eqref{bd-E-m}--\eqref{bd-E-line-m} is replaced by $m(m-1)$ or by $|x|^\delta m$. These bounds, especially the ones involving $m(m-1)$, are much more subtle, and we will extend them to diagrams that include $z$-dependence for some $z\in[0,1]$. Before we state the bounds, we again introduce a relevant construction:
\pagebreak[2]
\begin{defn}[Construction $1^z$]
\label{def-con1z}
\color{white}.\color{black}
\begin{itemize}
    \item[(a)] Given a diagram with a certain collection of backbone lines and a bond $(s,t)\in \Z^d\times \Z^d$, \emph{Construction~$1^z(s,t)$} is the operation in which a backbone line is replaced with a pivotal generating function $z \Tau_z$ as defined in \eqref{Gz-def} that goes from one end to $s$, is multiplied with $J(s,t)$, and from $t$ connects to the other end with an ordinary line $\tau$. Explicitly, this means that a backbone line $\taub(v-u)$ (say) is replaced by $z \Tau_z(v-s)J(s,t) \tau(t-u)$.
    
    \item[(b)] Given a diagram with a certain collection of backbone lines, \emph{Construction~$1^z$} is the operation in which we perform Construction~$1^z(s,t)$ followed by a sum over all $(s,t) \in \Zd \times \Zd$. Explicitly, this means that a backbone line $\taub(v-u)$ (say) is replaced by $(z \Tau_z * J * \tau)(v-u)$.
    
    \item[(c)] Given a diagram $G (\,\cdot\,)$, we write $G(\,\cdot\,; 1^z)$ for the result of applying Construction~$1^z$ to $G (\,\cdot\,)$.
\end{itemize}
\end{defn}

\begin{figure}
	\includegraphics[width =.85\textwidth]{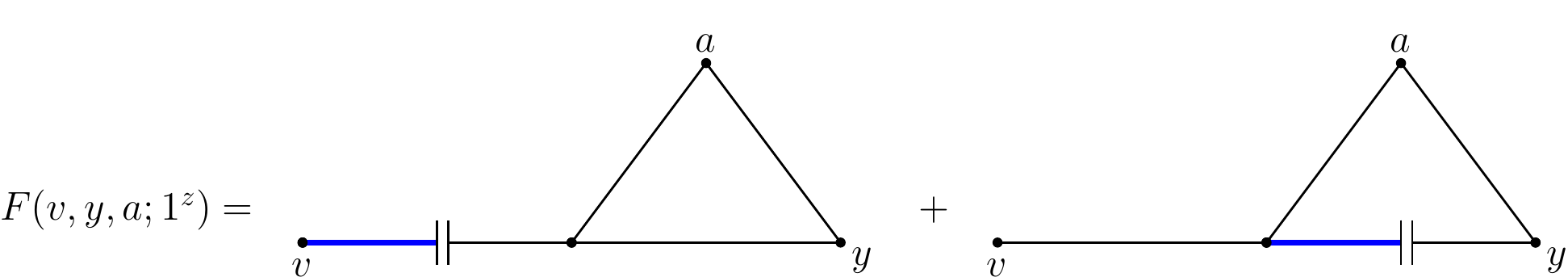}
	\caption{\label{fig:diagFvyaonez} A diagrammatic representation of $F(v,y,a; 1^z)$. The thick blue line represents the pivotal generating function $z \Tau_z$. The vertical dashes represent $J$. The unlabelled vertex corresponds to the summation variable $t$. }
\end{figure}

See Figure \ref{fig:diagFvyaonez} for a diagrammatic representation of $F(v,y,a; 1^z)$.
\medskip
 
 We will only apply Construction $1^z$ at most once, and only after all Constructions $1_\mathsf{b}$ have been performed. Therefore, we do not need to keep track of the backbone lines in Construction $1^z$. In terms of Construction $1^z$, the main bound on the doubly-weighted probabilities $\prob_p(E^{\bullet}_{m} (v,y;A))$ is given by the following proposition:

\begin{prop}[Bounds on doubly-weighted $E$]
\label{prop-bd-doubly-weighted-E}
For every $p\leq p_c$, $z\in[0,1]$, $A\subseteq \Z^d$ and $\bullet \in \{\varnothing,<,>\}$,
\begin{equation}
	\label{bd-E-m2}
	\sum_{m\geq 2} \tfrac{1}{2}m(m-1)z^m \prob_p(E^{\bullet}_{m} (v,y;A))
	\leq \sum_{a\in A} F(v,y,a; 1_\mathsf{b}, 1^z).
\end{equation}
Furthermore,
\begin{equation}
	\label{bd-E-line-m2}
	\sum_{m\geq 2} m(m-1)z^m \prob_p \big(E^{\bullet}_{m} (v,y;A)\cap\{v\conn s\} \big)
	\leq \sum_{a\in A}  F(v,y,a; \ell(s), 1_\mathsf{b}, 1^z).
\end{equation}
\end{prop}

We prove Proposition \ref{prop-bd-doubly-weighted-E} below. We start with a lemma that allows us to replace a $z$-dependent {\em restricted} two-point function $\tau^A$ (see Definition \ref{def-off-A})  by a similar {\em unrestricted} two-point function $\tau$. It is here that we rely on the fact that $z\in[0,1]$:

\begin{lemma}[Pivotals off $A$]
\label{lem-piv-off}
For every $z\in[0,1]$, $p \in [0, \|D\|_{\infty}^{-1}]$, and $A\subseteq \Z^d$,
\begin{equation}
	\label{piv-off}
	\sum_{m\geq 0} z^m \tau_m^A(u,v)\leq \sum_{m\geq 0} z^m \tau_m(u,v).
\end{equation}
\end{lemma}
\medskip

\proof We write the left-hand side as
\begin{equation}
	\sum_{m\geq 0} z^m \tau_m^A(u,v)=\expec_p\big[z^{|\Piv^A(u,v)|} \indic{u\conn v \off A}\big].
\end{equation}
Observe that if $u$ and $v$ are connected $\off$ $A$, then 
$|\Piv^A(u,v)|\geq |\Piv(u,v)|$. Thus, since $z\in[0,1]$, we can bound
\begin{equation}
	\sum_{m\geq 0} z^m \tau_m^A(u,v)
	\leq \expec_p\big[z^{|\Piv(u,v)|} \indic{u\conn v \off A}\big]
	\leq \expec_p\big[z^{|\Piv(u,v)|} \indic{u\conn v}\big]= \sum_{m\geq 0} z^m \tau_m(u,v). \qed
\end{equation}
\medskip

\proof[Proof of Proposition \ref{prop-bd-doubly-weighted-E}]
The proofs of \eqref{bd-E-m2} and \eqref{bd-E-line-m2} are very similar, but the proof of \eqref{bd-E-line-m2} is longer to write down, so we shall only give the proof of \eqref{bd-E-m2}.

We start with \eqref{bd-E-m2} for $\bullet \in \{\varnothing,>\}$. By a similar argument as in \eqref{e:msumispivsum}, we can write the factor $m(m-1)$ as a sum over two distinct pivotal bonds, to obtain
\begin{equation}
	\sum_{m\geq 2} \tfrac{1}{2}m(m-1)z^m \prob_p(E^{\bullet}_{m} (v,y;A))
	=\tfrac{1}{2}\sum_{b_1\neq b_2} \sum_{m\geq 2} z^m 
	\Pp \big(E^{\bullet}_m(v,y;A)\cap \{b_1,b_2\in \Piv(v,y)\} \big).
\end{equation}
If $b_1, b_2 \in \Piv(v,y)$, then either $b_1 \in \Piv(v, \ulb_2)$ or $b_2 \in \Piv(v, \ulb_1)$ (the indices do not imply an order). Both cases are symmetrical when $b_1$ and $b_2$ are summed over all bonds, so
\begin{align}
	\label{doubly-1}
	&\sum_{m\geq 2} \tfrac{1}{2}m(m-1)z^m \prob_p(E^{\bullet}_{m} (v,y;A))\\
	&\qquad\leq \sum_{b_1, b_2} \sum_{m_1\geq 0} z^{m_1 + 1} 
	\prob_p \big(E^{\bullet}(v,y;A)\cap \{b_1 \in \Piv(v,\ulb_2)\} \cap \{b_2\in \Piv(v,y)\}
	\cap \{|\Piv(v,\bb_1)|=m_1\} \big),\nn
\end{align}
where the inequality arises because we replace $z^m$ with $z^{m_1 +1}$, and clearly, $m_1 +1 \le m$ whenever $b_1 \in \Piv(v,y)$, $|\Piv(v,\bb_1)|=m_1$, and $|\Piv(v,y)|=m$ all occur.

\begin{figure}
	\includegraphics[width=.9\textwidth]{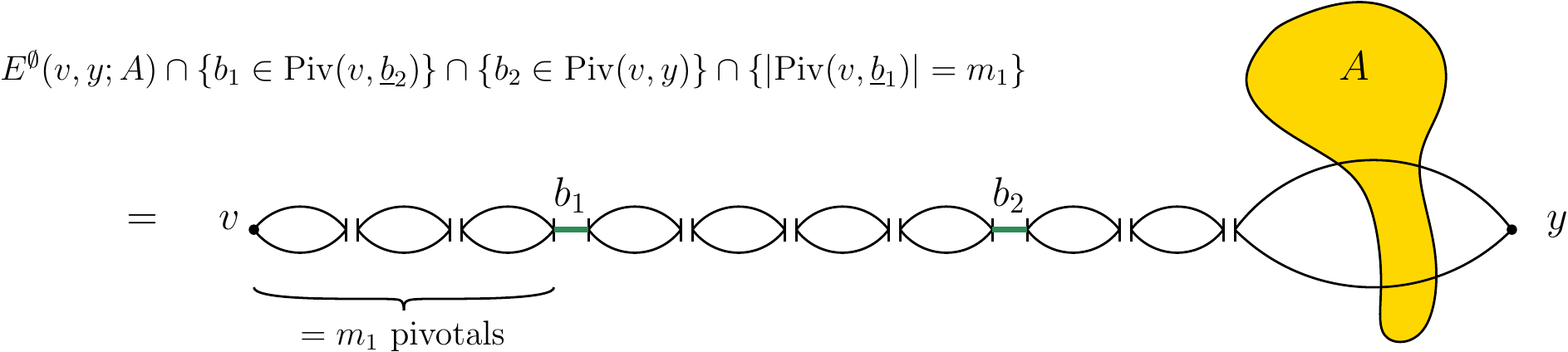}
	\caption{\label{fig:m2E} A sketch of the event on the left-hand side of \eqref{E-bd-m2-circ}.}
\end{figure}

Now consider only the case $\bullet=\varnothing$. 
We claim that
\begin{align}
	\label{E-bd-m2-circ}
	&E^{\sss\varnothing}(v,y;A)\cap \{b_1 \in \Piv(v,\ulb_2)\} \cap \{b_2\in \Piv(v,y)\} 
	\cap \{|\Piv(v,\bb_1)|=m_1\}\\
	&\qquad \subseteq
	\bigcup_{a\in A}\bigcup_{t} \{|\Piv(v,\bb_1)|=m_1\}
	\circ \{b_1\text{ occ.}\}\circ \{\tb_1\conn\bb_2\}\circ \{b_2\text{ occ.}\}\circ \{\tb_2\conn t\}\nn\\
	&\qquad\qquad\qquad\circ 
	\{t\conn a\}\circ \{a\conn y\}\circ \{t\conn y\}.\nn
\end{align}
The proof of this claim goes along the same lines as the proof of Lemmas \ref{lem-bd-E} and \ref{lem-bd-weighted-E}, so we omit it here. See Figure \ref{fig:m2E} for a sketch. Also compare \eqref{E-bd-m2-circ} with \eqref{E-bd-circ} and note that the event $\{v\conn t\}$ in \eqref{E-bd-circ} is simply replaced here by the event on the middle line of \eqref{E-bd-m2-circ} (excluding the unions over $a\in A$ and $t$).

Let us continue by bounding the probability of the right-hand side of \eqref{E-bd-m2-circ}. We follow the ideas in the proof of Lemma \ref{lem-rhon-taun}. We see that $\{|\Piv(v,\bb_1)|=m_1\}$ is the only non-increasing event on the right-hand side of \eqref{E-bd-m2-circ}. To apply the BKR-inequality, we need to identify sets of witnesses. For the (increasing) path events, we take the paths of open bonds that realize the connections as witness sets. The witness set for $\{b_i$ occupied$\}$ is simply the bond $b_i$. The witness set for $\{|\Piv(v,\bb_1)|=m_1\}$ consists of all the open bonds in all paths between $v$ and $\bb_1$, together with {\em all} the closed bonds in the graph. By definition, these sets of witnesses are all disjoint. Indeed, if they would not, then a path from $v\conn \bb_1$ would intersect any of the other paths that are enforced by the other events on the right-hand side of \eqref{E-bd-m2-circ}, which would make some of the bonds in $\Piv(v,\bb_1)$ not pivotal for $v\conn y$, which is a contradiction. This means that we can apply the BKR-inequality \eqref{BKR-ineq}, to obtain
\begin{equation}
	\begin{split}
	&\sum_{m\geq 2} \tfrac{1}{2}m(m-1)z^m \prob_p(E^{\sss\varnothing}_{m} (v,y;A))\\
	&\qquad\leq \sum_{b_1, b_2} \sum_{a\in A} \sum_t \sum_{m_1\geq 0} z^{m_1+1} 
	\tau_{m_1}(\bb_1-v)
	J(b_1)\tau(\bb_2-\tb_1)J(b_2)\tau(\bb_2-v)\tau(a-t)\tau(y-a)\tau(y-t)\\
	&\qquad= \sum_{a\in A} \sum_t (z \Tau_z * J * \tau * J * \tau)(t-v)\tau(a-t)\tau(y-a)\tau(y-t)\\
	&\qquad \le \sum_{a\in A} \sum_tF(v,t,y,a; 1_\mathsf{b}, 1^z).
	\end{split}
\end{equation}
\medskip

By \eqref{doubly-1}, the proof of \eqref{bd-E-m2} for $\bullet=>$ is analogous.
\medskip

\begin{figure}
	\includegraphics[width=.9\textwidth]{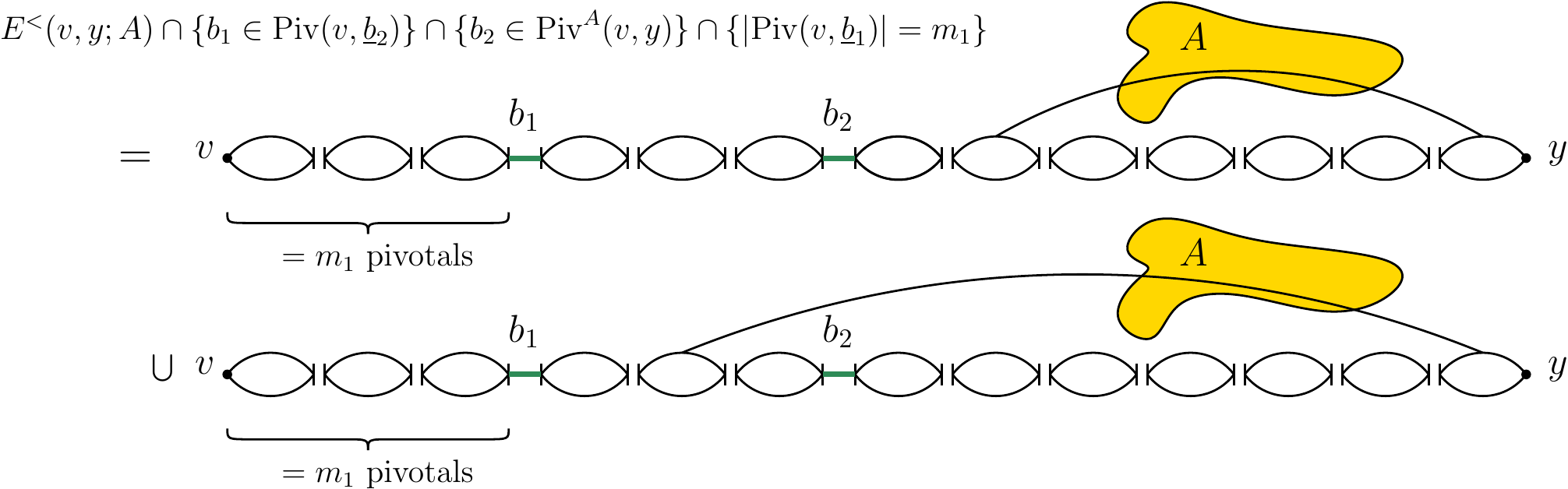}
	\caption{\label{fig:m2E2} A sketch of the event on the left-hand side of \eqref{E-bd-m2-circ-2}.}
\end{figure}

To prove \eqref{bd-E-m2} for $\bullet=<$, we start by noting that by an identity similar to \eqref{e:msumispivsum},
\begin{equation}
	\sum_{m\geq 2} \tfrac{1}{2}m(m-1)z^m \prob_p(E^{\sss <}_{m} (v,y;A))=\tfrac{1}{2}\sum_{b_1\neq b_2} \sum_{m\geq 2} z^m \prob_p \big(E^{\sss <}_m(v,y;A)
	\cap \{b_1,b_2\in \Piv^A(v,y)\} \big).
\end{equation}
We obtain
\begin{align}
	\label{doubly-2}
	&\sum_{m\geq 2} \tfrac{1}{2}m(m-1)z^m \prob_p(E^{\sss <}_{m} (v,y;A))\\
	&\qquad\leq \sum_{b_1 \neq b_2} \sum_{m_1\geq 0} z^{m_1 +1} \prob_p \big(E^{\sss <}(v,y;A)\cap \{b_1\in \Piv^A(v,\ulb_2)\} \cap \{b_2 \in \Piv^A(v,y)\} \cap \{|\Piv^A(v,\bb_1)|=m_1\} \big).\nn
\end{align}
There are two cases, depending on whether $b_1\in \Piv(v,y)$ or not. When $b_1\in \Piv(v,y)$ (i.e., $b_1$ is both pivotal and pivotal off $A$), then
\begin{equation}\label{E-bd-m2-circ-2}
	\begin{split}
	E^{\sss <}(v,y;A)\cap & \{b_1\in \Piv(v,\ulb_2)\}\cap \{b_2\in \Piv^A(v,y)\}
	\cap \{|\Piv(v,\bb_1)|=m_1\}\\
	&\subseteq
	\Bigg(\bigcup_{a\in A}\bigcup_{t} \{|\Piv(v,\bb_1)|=m_1\}
	\circ \{b_1\text{ occ.}\}\circ \{\tb_1\conn\bb_2\}\circ \{b_2\text{ occ.}\}\circ \{\tb_2\conn t\}\\
	&\qquad\qquad\circ 
	\{t\conn a\}\circ \{a\conn y\}\circ \{t\conn y\}\Bigg)\\
	&\qquad \cup \Bigg(\bigcup_{a\in A}\bigcup_{t} \{|\Piv(v,\bb_1)|=m_1\}
	\circ \{b_1\text{ occ.}\}\circ \{\tb_1\conn t\}\\
	&\qquad\qquad\circ 
	\{t\conn a\}\circ \{a\conn y\}\circ \{t\conn \bb_2\}\circ \{b_2\text{ occ.}\}
	\circ \{\tb_2\conn t\}\Bigg).
	\end{split}
\end{equation}
Here the first term corresponds to the case where also $b_2\in 	\Piv(v,y),$ and the second to the case where $b_2\in \Piv^A(v,y)\setminus \Piv(v,y)$.
See Figure \ref{fig:m2E2} for a sketch. 
	
Taking the probability of the right-hand side of \eqref{E-bd-m2-circ-2} we can again apply the BKR-inequality. We thus obtain the bound
\begin{multline}\label{e:firstboundless}
	\sum_{b_1, b_2} \sum_{m_1 \ge 0} z^{m_1 +1} \Pp \big( E^{\sss <}(v,y;A)\cap 
	 \{b_1\in \Piv(v,\ulb_2)\}\cap \{b_2\in \Piv^A(v,y)\}\cap \{|\Piv(v,\bb_1)|=m_1\} \big)\\
	 \le \sum_{a\in A} \sum_t (z \Tau_z * J * \tau * J * \tau)(t-v)\tau(a-t)\tau(y-a)\tau(y-t) \\
	  + \sum_{a\in A} \sum_t (z \Tau_z * J * \tau)(t-v)\tau(a-t)\tau(y-a)(\tau* J * \tau)(y-t).
	\end{multline}
Note that both terms on the right-hand side are contributions to $\sum_{a \in A} F(v,y,a;1_\mathsf{b}, 1^z)$.

\begin{figure}[tb]
	\includegraphics[width=.9\textwidth]{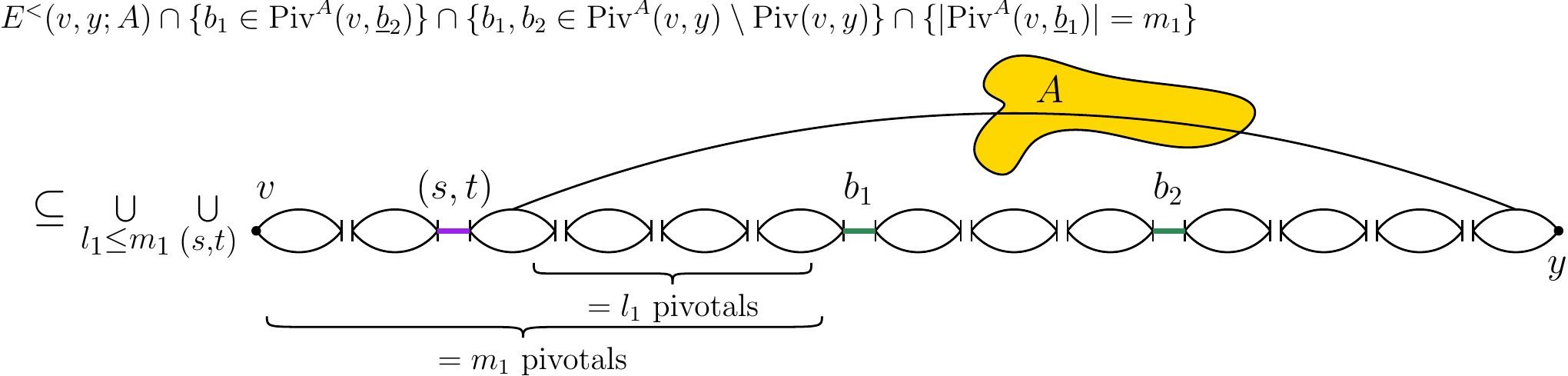}
	\caption{\label{fig:m2E3} A sketch of the inclusion in \eqref{E-bd-m2-cap2}.}
\end{figure}

When $b_1\in \Piv^A(v,y)\setminus \Piv(v,y)$, necessarily also $b_2\in \Piv^A(v,y)\setminus \Piv(v,y)$. Then, let $b$ be the last element of $\Piv(v,y)\cap \Piv^A(v,y)$, so that necessarily $\{\tb\dbc y\}$ occurs (when such a $b$ does not exist we take $\tb=v$). Write $t=\tb$. 
Then,
\begin{multline} \label{E-bd-m2-cap2}
	E^{\sss <}(v,y;A)\cap \{b_1 \in \Piv^A(v,\ulb_2)\} \cap \{b_1,b_2\in 
	\Piv^A(v,y)\setminus \Piv(v,y)\}\cap \{|\Piv^A(v,\bb_1)|=m_1\}\\
	 \subseteq \bigcup_{l_1\leq m_1}\bigcup_{(s,t)} \{b_1\in 
	\Piv^A (v,\ulb_2)\}\cap \{b_2\in \Piv^A(v,y)\}\cap \{|\Piv^A(v,\bb_1)|=m_1\}\\
	\cap \{(s,t)\in \Piv(v,y)\} \cap \{t\dbc y\}
	\cap \{|\Piv^A(t,\bb_1)|=l_1\}.
\end{multline}
See Figure \ref{fig:m2E3} for a sketch of this inclusion. 
Now we write the event in terms of disjointly occurring events as 
\begin{align}
	&E^{\sss <}(v,y;A)\cap \{b_1,b_2\in \Piv^A(v,y)\setminus \Piv(v,y)\}
	\cap \{t\dbc y\}\cap \{(s,t)\in \Piv(v,y)\}
	\cap \{|\Piv^A(t,\bb_1)|=l_1\}\\
	&\quad \subseteq
	\bigcup_{a\in A}\bigcup_{t } \{v\conn t\}\circ \{|\Piv^A(t,\bb_1)|=l_1\}
	\circ \{b_1\text{ occ.}\}\circ \{\tb_1\conn\bb_2\}\circ \{b_2\text{ occ.}\}
	\circ \{\tb_2\conn y\}\circ \{t\conn a\}\circ\{a\conn y\}
	.\nn
\end{align}
Therefore, using the BKR-inequality (where the disjointness is verified by finding witness sets as in \eqref{E-bd-m2-circ} above), and using that $l_1 \leq m_1$, and Lemma~\ref{lem-piv-off}, we obtain
\begin{equation}
	\begin{split}
		\sum_{b_1, b_2} \sum_{m_1 \ge 0} z^{m_1+1} \Pp \big( & E^{\sss <}(v,y;A) 
		\cap \{b_1 \in \Piv^A(v,y) \setminus \Piv(v,y)\} \cap \{|\Piv^A(v, \ulb_1)| = m_1\}\big)\\
		& \le \sum_{a \in A} \sum_t \tau(t-v) (z \Tau_z * J * \tau * J * \tau)(y-t) \tau(a-t) \tau(y-a).
	\end{split}
\end{equation}
Combining this bound with \eqref{e:firstboundless} proves the claim for $\bullet = <$, and thus completes the proof of \eqref{bd-E-m2}. \qed
\medskip

We proceed by bounding the doubly-weighted lace-expansion coefficients. An important example of such a doubly-weighted lace-expansion coefficient is
	\begin{equation}
	\label{doubly-weighted-diagrams}
	\sum_{x} \sum_{m \ge 1} m^2 \bar \pi_m^{\sss(N)}(x) z^m.
	\end{equation}
This coefficient appears when bounding $\sum_{m \ge 1} \sum_x m^{1+\vep} \bar \pi_m^{\sss(N)}(x)$ for $\vep>0$ sufficiently small. We use the result in Proposition \ref{prop-bd-doubly-weighted-E}, and we start by reducing the estimate to this setting. In \eqref{pimNdef} the coefficient $\pi_m^{\sss(N)}(x)$ is given explicitly as a sum over the positions of $b_1,\dots, b_n$, with $m=N+m_1+\cdots+m_{\sss N}$.
We apply \eqref{e:msumispivsum} to sum over the pivotals to get a factor $m^2$, using the bound 
\begin{equation}\label{eqMSquareBound}
	m^2\leq N(N^2 + m_1^2+\cdots+m_{\sss N}^2)\leq N^3 + Nm +2N(\tfrac{1}{2}m_1(m_1-1)+\cdots+\tfrac{1}{2}m_{\sss N}(m_{\sss N}-1)).
\end{equation}
The first and second term were already bounded in terms of diagrams in Sections \ref{sec-unweighted-lace expansion} and \ref{sec-singly-weighted-lace expansion}, respectively. To bound the third term, we apply Proposition \ref{prop-bd-doubly-weighted-E}. The result of this is that we start with the bound on $\sum_{x\in \Z^d} \bar \pi^{\sss(N)}_m(x)$ in \eqref{piexpansion}, and apply Constructions $1_\mathsf{b}$ and $1^z$ to two lines that originate from one nested expectation. To formalize this we introduce some notation.
Define
\begin{multline}\label{eqphildef}
    	\Phl^{\sss (i)} (s_i, t_i) := \sum_{s_1, \dotso, s_{i-1}}\sum_{t_1, \dotso, t_{i-1}}
	\sum_{u_1, \dotso, u_{i-1}}\sum_{v_1, \dotso, v_{i-1}} A_3 (0, s_1, t_1)\\
     	\times \Bigg[\prod_{k=1}^{i-1} B_1 (s_k, t_k, u_k, v_k) B_2 (u_k, v_k, s_{k+1}, t_{k+1}) \Bigg],
\end{multline}
\begin{multline}
    	\Phr^{\sss (N-j)} (u_{N-j+1}, v_{N-j+1}, x) \\
        := \sum_{s_{N-j+2}, \dotso, s_N}\sum_{t_{N-j+2}, \dotso, t_N}\sum_{u_{N-j+2}, \dotso, u_N}
        \sum_{v_{N-j+2}, \dotso, v_N}
     	\Bigg[\prod_{k=N-j+1}^{N-1}B_2(u_k, v_k, s_{k+1}, t_{k+1})\\
     	\times B_1 (s_{k+1}, t_{k+1}, u_{k+1}, v_{k+1})\Bigg] A_3^{\mathsf{b}} (u_N, v_N, x),
\end{multline}
\begin{equation}\label{eqomegadef}
	\Omega(s,t,u,v) := \sum_{a,b,c,d} B_1 (s, t, a, b) B_2(a, b, c, d) B_1 (c, d, u, v),
\end{equation}
and 
\begin{equation}\label{equpsilondef}
	\Upsilon(s,t,x) := \sum_{u,v} B_1(s,t,u,v)A_3^\mathsf{b}(u,v,x).
\end{equation}
So we can write $\Phi^{\sss (N)}$ for $N \ge 1$ as 
\begin{equation}\label{phiexp1}
    	\Phi^{\sss (N)}(x) = \sum_{s,t} \Phl^{\sss (N)}(s,t) \Upsilon(s,t,x)
\end{equation}
and for any $1 \le i \le N-1$ as
\begin{equation}\label{phiexp2}
	\Phi^{\sss (N)}(x) =  \sum_{s,t,u,v} \Phl^{\sss (i)}(s,t) \Omega(s,t,u,v) \Phr^{\sss (N-i-1)}(u,v,x).
\end{equation}

Plugging \eqref{eqMSquareBound} and the definitions \eqref{eqphildef}--\eqref{phiexp2} into \eqref{piexpansion} readily obtains our next lemma:
\begin{lemma}[Doubly weighted LE coefficients]
\label{lem:doubleLE}
		 For all $N \ge 1$ and $z\in[0,1]$,
\begin{equation}
	\begin{split}
		\sum_{m \ge 1} m^2 \bar \pi_m^{\sss(N)}(x) z^m & \le N^3 \Phi^{\sss (N)}(x) + N \Phi^{\sss (N)}(x; 1_\mathsf{b})\\
		& \quad
		 + 2N \sum_{i=1}^{N-1} \sum_{s,t,u,v} \Phl^{\sss (i)}(s,t) 
		 \Omega(s,t,u,v; 1_\mathsf{b}, 1^z) \Phr^{\sss (N-i-1)}(u,v,x)\\
		& \quad + \sum_{s,t} \Phl^{\sss (N)}(s,t) \Upsilon(s,t,x; 1_\mathsf{b}, 1^z).
	\end{split}
\end{equation}	
\end{lemma}
This completes the bound on the doubly-weighted diagram in \eqref{doubly-weighted-diagrams}.
\medskip

Our proof requires a bound on another doubly-weighted diagram, namely, the diagram with mixed weight $|x|^\delta m$. The bound here is simpler than the one for weights $m(m-1)$, because for sufficiently small $\delta$ these diagrams will be convergent also when $z=1$. For this we may simply apply the arguments of Section~\ref{sec-singly-weighted-lace expansion}, together with first an application of Construction $1_{\mathsf{b}}$ and then an application of Construction $1_{\delta}^x$ to obtain the following bound:
\begin{lemma}[Mixed weighted LE coefficients]
\label{lem:timespaceLE} 
For $N \ge 1$, $\delta > 0$,
\begin{equation}
	\sum_{m \ge 1} |x|^\delta m  \bar \pi_m^{\sss (N)}(x)  \le \Phi^{\sss (N)}(x; 1_\mathsf{b}, 1_\delta^x).
\end{equation}
\end{lemma}
\smallskip

Finally, combining Lemma~\ref{lem:singleLE} with the arguments that led to Lemmas~\ref{lem:fixedbond} and \ref{lem:psidiagram}, respectively, we may conclude the following two lemmas:
\begin{lemma}[Doubly weighted $\pi_{m,n}$ coefficients]
\label{lem:pidouble}
For $N \ge 0$, $\delta >0$, uniformly in $n \ge 1$,
\begin{equation}
	\sum_{x,y}  \sum_{m \ge n} (|x|^\delta + |y|^\delta) \bar \pi^{\sss (N)}_{n,m}(y,x) \le  \sum_{v,x,y} \Phi^{\sss (N)}(x; 1_{\mathsf{b}}(v,y), 1^x_\delta) + \sum_{v,x,y} \Phi^{\sss (N)}(x; 1_{\mathsf{b}}(v,y),  1^y_\delta)
\end{equation}
\end{lemma}
\begin{lemma}[Spatially weighted $\psi_m$ coefficients]
\label{lem:psispatial}
For $N \ge 0$, $\delta >0$,
\begin{equation}
	\sum_x \sum_{m \ge 0} |x|^\delta \bar \psi^{\sss (N)}_m(x) \le p_c \sum_{v,w,x} \Phi^{\sss (N)}(w; 1_{\mathsf{b}}(v,x), 1^x_\delta).
\end{equation}
\end{lemma}
\medskip

\noindent
\paragraph{\bf How we proceed.} In the next two sections we will state bounds on the lace-expansion coefficients that follow from the diagrammatic estimates derived in this section. This goes in three steps. The first step, in Section~\ref{sec:LEprop}, is to establish some {\em weak bounds}, as given in Proposition~\ref{prop:smallpi}. These bounds follow from bounds on the two-point function $\tau_p(x,y)$ derived in \cite{HarSla90a,HeyHofSak08}, and are proved in \longversion{Appendix~\ref{sec-LEcoefficients}}\shortversion{\cite[Appendix~A]{HeyHofHulMie17b}} using standard techniques. In the second step, also in Section~\ref{sec:LEprop}, we use these weak bounds to prove \emph{infrared bounds} on $\hat{\Tau}_z(k)$ in Proposition~\ref{prop:positivity}. The third step, in Section~\ref{sec-impr-bds-LE}, is to use the weak and infrared bounds to establish (what we call) {\em improved bounds} on spatial fractional derivatives of the lace-expansion coefficients in Proposition~\ref{prop-LEcoefficients2} and on temporal fractional derivatives in Proposition~\ref{prop-LEcoefficients}. These improved bounds are proved in Section~\ref{sec-impr-bds-LE} and \longversion{Appendix~\ref{sec-squares},}\shortversion{\cite[Appendix~B]{HeyHofHulMie17b},} again using fairly standard techniques. These propositions then allow us to complete the proofs of all our main results.


\section{Weak bounds and infrared bounds}
\label{sec:LEprop}
In this section we state the weak bounds and state and prove the infrared bounds (subject to the proof of the weak bounds).

\begin{prop}[Weak bounds on the lace-expansion coefficients]\label{prop:smallpi}
Under the assumptions of Theorem \ref{thm-endpoint} there exists $\tilde c >0$ such that for all $p \le p_c$,
\begin{eqnarray}   \label{eqsmallpi1} 
    	\sum_{x} \sum_{m \ge 0}  \bar \pi^p_m (x) 
    	&\le& 1 + \tilde c \beta^{1/6},\\
        \label{eqsmallpi2} \sum_{x} \sum_{m \ge 1}  m \bar  \pi^p_m (x) 
        &\le& \tilde c \beta^{1/6},\\
\end{eqnarray}
\end{prop}
The proof of this proposition uses the bounds derived in Section~\ref{sec-bounds-lace-exp} and is given in \longversion{Appendix~\ref{sec-LEcoefficients}.}\shortversion{\cite[Appendix~A]{HeyHofHulMie17b}.} The bound on \eqref{eqsmallpi1} that is proved is actually $1 + O(\beta^{1/2})$. Moreover, the bounds in Proposition \ref{prop:smallpi} can be improved to $1+O(\beta)$ and $O(\beta)$, respectively, but this requires significantly more effort, and we do not need such strong bounds. 

Let us emphasize here that the proof of Proposition \ref{prop:smallpi} relies on known results, in particular, on the infrared bound $\hat{\tau}^{p}(k)\leq C/[1-\hat{D}(k)]$ for $p\le p_c$ proved in \cite{HarSla90a,HeyHofSak08}. This means that we do not need to use a ``bootstrap argument'' to first establish convergence of $\sum_{x} \sum_{m \ge 0}  \bar \pi^p_m (x)$ without reference to the two-point function asymptotics, as is common in many lace-expansion analyses (e.g.\ in \cite{HarSla90a,HeyHofSak08}). 
\longversion{The reader may now wish to read Appendix~\ref{sec-LEcoefficients} and verify that this is indeed the only ingredient to the proof of Proposition \ref{prop:smallpi}.}

\begin{center}
\setlength{\fboxsep}{7pt}
\setlength{\fboxrule}{1pt}
\begin{boxedminipage}{15cm}
From here until the end of this paper \longversion{(excluding Appendix~\ref{sec-LEcoefficients})} we will only need to consider $p=p_c$. We will therefore omit most superscripts $p_c$ from here onward.
\end{boxedminipage}
\end{center}

We now use the weak bounds from Proposition~\ref{prop:smallpi} to establish a number of bounds that we need in many of the subsequent proofs. In particular, we establish bounds on the divergence of $\hat \Pi_z(k)$ in $z$ and $k$. We also establish bounds on the divergence of $\hat \Tau_z(k)$ in $k$ and $z$ that are useful in the settings where $z$ is close to $1$ and where $k$ is close to $0$, respectively. These settings are analogous to certain ``low energy approximations'' in Quantum Field Theory, and we therefore call them ``infrared bounds'', as is common in QFT.

A number of our results is obtained by analyzing a generating function in terms of an auxiliary variable $z$. So far, it was sufficient to consider $z\in[0,1]$ only. From here onwards, we shall use complex values of $z$ in order to apply a Tauberian theorem in the proofs of Propositions \ref{lemmaError} and~\ref{lemmaError-Tau} below.

\begin{prop}[Infrared bounds]
\label{prop:positivity} 
Under the assumptions of Theorem \ref{thm-endpoint}:
\begin{enumerate}
 	\item There exists $\tilde c >0$ such that uniformly in $k \in \R^d$ and $z\in\mathbb C$ with $0<|z|\le1$,
	\begin{itemize}
		\item[(i.a)]
		\begin{equation} \label{eq-PiPi1}
       			 \left\lvert \hat \Pi_1 (k) - \hat \Pi_z (k) \right\rvert \le  \tilde c \beta^{1/6} |1-z|, 
		\end{equation}
		
		\item[(i.b)]
		\begin{equation} \label{eq-PiPi2}
       			 {\lvert \hat \Pi_z(k) - 1\rvert \le \tilde c \beta^{1/6}},
		\end{equation}
		
		\item[(i.c)]
		\begin{equation} \label{eq-PiPi3}
			{\lvert \hat \Pi_1(0) - \hat \Pi_1(k)\rvert \le \tilde c \beta^{1/6}}[1-\hat{D}(k)].
		\end{equation}
	\end{itemize}
	\item There exists constants $C, C'>0$ such that uniformly in $k \in \R^d$ and $z\in\mathbb C$ with $|z|<1$,
	\begin{itemize}
	\item[(ii.a)]
	\begin{equation}\label{eq-T-bd-pos1}
        		|\hat \Tau_z (k)| \le \frac{C}{|1-z|} ,
	\end{equation}

	\item[(ii.b)] for $k \neq 0$,
	\begin{equation}\label{eq-T-bd-pos2}
	 	|\hat \Tau_z(k)|  \le \frac{C'}{1-\hat D(k)}.
	\end{equation}
	\end{itemize}
\end{enumerate}
\end{prop}
\proof[Proof of Proposition \ref{prop:positivity} subject to Proposition \ref{prop:smallpi}]
(i.a)
We bound
\begin{equation}\label{eq-Pi1-Piz}
        \left\lvert \hat \Pi_1 (k) - \hat \Pi_z (k) \right\rvert 
        = \left\lvert \sum_{x} \sum_{m \ge 0} (1-z^m) \pi_m(x) \e^{ik\cdot x} \right\rvert 
        \le \sum_{x} \sum_{m \ge 0} |1-z^m|\, \bar \pi_m (x),
\end{equation}
where the absolute summability required for the equality is due to Proposition~\ref{prop:smallpi}.
Using $|1-z^m| = |1-z|\left|\sum_{\ell=0}^{m-1}z^\ell\right| \le m \left|1-z\right|$ for $|z|\le1$ and $m\ge0$, and using Proposition \ref{prop:smallpi} again, it follows that
\begin{equation}
    	\sum_{x} \sum_{m \ge 1} |1-z^m|\, \bar \pi_m (x) 
	\le |1-z| \sum_{x} \sum_{m \ge 1} m \bar  \pi_m (x)  
	\le \tilde c \beta^{1/6}|1-z|.
\end{equation}
\smallskip

(i.b) Since $\pi_0 (0)= 1$,  
\begin{equation}
        \left\lvert \hat \Pi_z (k) - 1\right\rvert 
        = \left\lvert \sum_{x} \sum_{m \ge 0} (1-\delta_{0,x} \delta_{0,m}) 
        \pi_m (x) z^m \e^{ik \cdot x} \right\rvert 
        \le  \sum_{x} \sum_{m \ge 0} (1-\delta_{0,x} \delta_{0,m})
        \bar \pi_m (x) . 
\end{equation}
Using again $\pi_0 (0) = 1$ shows that this last expression is equal to $\sum_{x } \sum_{m\ge 0} \bar \pi_m (x) -1$, which is bounded above by $\tilde c \beta^{1/6}$ by  \eqref{eqsmallpi1}. 
\medskip

(i.c) This follows directly from Remark~\ref{rem-HS-exp} and \cite[Proposition 2.5]{HeyHofSak08}.
\medskip

(ii.a) 
Write $c_i$ ($i\in\N$) for positive constants that depend neither on $k$ nor on $z$. 
We use \eqref{eqTauExp}
to write 
\begin{equation}\label{eqPos0}
	\lvert \hat\Tau_z(k)\rvert^{-1}=\frac1{|\hat\Pi_z(k)|}\left|1-zp_c\Dk\hat\Pi_z(k)\right|
	\ge c_1\left|1-zp_c\Dk\hat\Pi_z(k)\right|
\end{equation}
by \eqref{eq-PiPi2} for $\beta$ small enough. Hence, 
\begin{equation}\label{eqPos1}
	c_1^{-1}\lvert \hat\Tau_z(k)\rvert^{-1}
	\ge \left|1-p_c\Dk\hat\Pi_1(k) + (1-z)p_c\Dk\hat\Pi_1(k)
	+zp_c\Dk\left(\hat\Pi_1(k)-\hat\Pi_z(k)\right)\right|. 
\end{equation} 
By \eqref{eqPcPi} and \eqref{eq-PiPi1}, we further bound the first two terms on the right-hand side by
\begin{equation}\label{eqPionekbd}
   \begin{split}
	\left|1-p_c\Dk\hat\Pi_1(k)\right|
	&=\left|1-\Dk+p_c \Dk\left(\hat\Pi_1(0)-\hat\Pi_1(k)\right)\right| \\
	&\ge 1-\Dk - p_c |\Dk|\,\tilde c \beta^{1/6} |1-\hat{D}(k)|.
   \end{split}
\end{equation}
We use the bound $|p_c-1| \le c_2\beta$ (cf.\ \cite{HarSla90a, HeyHofSak08}), \eqref{eq-PiPi1}, and \eqref{eq-PiPi2} to obtain for some $c_3>0$,
\begin{align}
	c_1^{-1}\lvert \hat\Tau_z(k)\rvert^{-1}
	&\ge \left|1-p_c\Dk\hat\Pi_1(k)\right| -
	\left|(1-z)p_c\Dk\hat\Pi_1(k) + zp_c\Dk\left(\hat\Pi_1(k)-\hat\Pi_z(k)\right)\right|
	\nnb
	&=[1-\Dk]\left(1-c_3\beta^{1/6}\right)-\left(1+c_3\beta^{1/6}\right)|1-z|\,|\Dk|. 
	\label{eqPos2}
\end{align}

We fix some $\eps\in(0,1)$ and choose $\delta>0$ such that 
\begin{equation}
	\frac{1-\eps}{3+2c_3\beta^{1/6}}>\delta,
\end{equation}
which is equivalent to 
\begin{equation}
	(1-\delta)- (1+c_3\beta^{1/6})2\delta>\eps.
\end{equation}
We consider two cases: first we derive the bound for $\{k\in\R^d:|\Dk|\le \delta\}$ and later consider the remaining domain. 
For $|\Dk|\le \delta$, we get from \eqref{eqPos1}--\eqref{eqPos2} that 
\begin{equation}
	\label{eqPos3}
	c_1^{-1}\lvert \hat\Tau_z(k)\rvert^{-1}
	\ge (1-\delta)-(1+c_3\beta^{1/6})2\delta>\eps,
\end{equation}
so that  
\begin{equation}
	\label{eqPos4}
	\lvert \hat\Tau_z(k)\rvert^{-1}
	\ge \eps\,c_1
	\ge \eps\,c_1\,\tfrac12\,|1-z|
\end{equation} 
for such $k$ and $|z|<1$, where we use that $\frac12\,|1-z|\leq 1$ when $|z|<1$. 

Finally, consider $k$ for which $|\Dk|> \delta$. 
Then, using $|z|^2=\Re(z)^2+\Im(z)^2$ for any $z\in\mathbb C$, the absolute value of the first three terms on the right-hand side of \eqref{eqPos1} can be written as
\begin{multline}\label{eqPos5}
	\Big|\left(1-p_c\Dk\hat\Pi_1(k)\right)+\left((1-z)p_c\Dk\hat\Pi_1(k)\right)\Big| \\
	=\sqrt{\left(1-p_c\Dk\hat\Pi_1(k)+(1-\Re(z))p_c\Dk\hat\Pi_1(k)\right)^2
	+\left(\Im(z)p_c\Dk\hat\Pi_1(k)\right)^2}.
\end{multline}
If $|\Im(z)|>|1-\Re(z)|$ then $|\Im(z)|^2 > \frac12 |1-z|^2$, so
\begin{equation}
	\left(\Im(z)p_c\Dk\hat\Pi_1(k)\right)^2
	\ge \tfrac12|1-z|^2\left(p_c\hat\Pi_1(k)\Dk\right)^2.
\end{equation}
On the other hand, if $|\Im(z)|\le|1-\Re(z)|$, then we use that $(1-\Re(z))^2\ge \frac12 \lvert1-z\rvert^2$ and that $1-p_c\Dk\hat\Pi_1(k)\ge0$ (the latter by \eqref{eqTauExp}, \eqref{eq-PiPi2}, and the fact that $\hat \Tau_1(k)= \hat\tau(k) \ge0$ by \cite[Lemma 3.3]{AizNew84}) and $1-\Re(z)>0$ to see that the term in the square root in \eqref{eqPos5} is bounded below by 
\begin{equation}
	\left((1-\Re(z))p_c\Dk\hat\Pi_1(k)\right)^2
	\ge \frac12|1-z|^2\left(p_c\hat\Pi_1(k)\Dk\right)^2.
\end{equation}
Thus, for any $z\in\mathbb C$ with $|z|<1$, 
\eqref{eqPos5} is bounded below by 
\begin{equation}
	\frac1{\sqrt{2}}|1-z|p_c\lvert\hat\Pi_1(k)\rvert\lvert\Dk\rvert.
\end{equation}
We plug this estimate into \eqref{eqPos1} for the lower bound 
\begin{align}
	c_1^{-1}\lvert \hat\Tau_z(k)\rvert^{-1}
	&\ge\Big|1-p_c\Dk\hat\Pi_1(k) + (1-z)p_c\Dk\hat\Pi_1(k)\Big|
	-\left|zp_c\Dk\left(\hat\Pi_1(k)-\hat\Pi_z(k)\right)\right|\nnb
	&\ge\frac1{\sqrt{2}}|1-z|p_c\lvert\hat\Pi_1(k)\rvert\lvert\Dk\rvert
	-|z|p_c\lvert\Dk\rvert\left|\hat\Pi_1(k)-\hat\Pi_z(k)\right|\nnb
	&\ge\frac{1-c_4\beta^{1/2}}{\sqrt{2}}|1-z|p_c\lvert\Dk\rvert
	-c_4\beta^{1/6}|1-z|p_c\lvert\Dk\rvert,
\end{align}
where the final inequality is due to \eqref{eq-PiPi1}. 
Thus 
\begin{equation}\label{eqPos6}
	c_1^{-1}\lvert \hat\Tau_z(k)\rvert^{-1}
	\ge \frac{1-(1+\sqrt2)c_4\beta^{1/6}}{\sqrt{2}}\left(1-c_2\beta\right)\delta\lvert1-z\rvert,
\end{equation}
and this is positive if $\beta$ is small enough. 
The bounds \eqref{eqPos4} and \eqref{eqPos6} together establish \eqref{eq-T-bd-pos1} for \emph{all} $k\in \Rd$. 

\medskip
(ii.b) Finally, we prove  $|\hat \Tau_z (k)| \le C' [1-\hat D(k)]^{-1}$ for $k \neq 0$. 
The first bound in \eqref{eq-T-bd-pos1} readily gives us the result as long as $1-\Dk<2|1-z|$. It therefore suffices to consider $1-\Dk\ge 2|1-z|$. 
We start by recalling from \eqref{eqTauExp} that
\begin{equation}
   	 \left|\hat \Tau_z (k)\right|^{-1} = \left\lvert {\hat \Pi_z (k)}\right\rvert^{-1}
	 \left\lvert{1-z p_c \hat D(k) \hat \Pi_z (k)}\right\rvert. 
\end{equation}
The factor $\lvert {\hat \Pi_z (k)}\rvert^{-1} $ is arbitrarily close to 1 for $\beta$ sufficiently small by \eqref{eq-PiPi2}, thus (for such $\beta$) there exist $c_1>0$ such that 
\begin{align}
    	\left|\hat \Tau_z (k)\right|^{-1} 
    	&\ge c_1\left\lvert{1-\Dk+\Dk\left(1-z+z p_c(\hat \Pi_1 (0) - \hat \Pi_z (k))\right)}\right\rvert\nnb
    	&\ge c_1\left(1-\Dk-|\Dk|\left\lvert1-z+z p_c\left(\hat \Pi_1 (0) - \hat\Pi_1(k)+ \hat\Pi_1(k)
	- \hat \Pi_z (k)\right)\right\rvert\right),
    \label{eq TzkBd1}
\end{align}
where in the first line we have used that $p_c\hat \Pi_1 (0)=1$. 

We now bound $|\Dk|\le1$, $|1-z|\le \frac12 (1-\Dk)$ by our assumption, and $|zp_c|\le p_c$. We further bound
\begin{equation}
	\left|\hat \Pi_1 (0) - \hat\Pi_1(k)\right|
	\le 2 \tilde c\beta^{1/6}[1-\hat D(k)]
\end{equation}
by \eqref{eq-PiPi3}, and 
\begin{equation}
	\left|\hat\Pi_1(k)- \hat \Pi_z (k)\right|
	\le \tilde c\beta^{1/6}|1-z|
	\le \tilde c\beta^{1/6}\tfrac12 [1-\Dk]
\end{equation}
by \eqref{eq-PiPi1} and again our assumption.  
When substituting these into \eqref{eq TzkBd1}, we obtain 
\begin{equation}
	\left|\hat \Tau_z (k)\right|^{-1} 
	\ge c_1\left(1-\left(\tfrac12+p_c\left(2+\tfrac12\right)\tilde c\beta^{1/6}\right)\right)[1-\Dk].
\end{equation}
The factor that multiplies $[1-\Dk]$ is positive for $\beta$ sufficiently small. 
This completes the proof of \eqref{eq-T-bd-pos2}. \qed

\section{Improved bounds on the lace-expansion coefficients}
\label{sec-impr-bds-LE}
In this section we state the improved bounds on the lace-expansion coefficients that allow us to prove sharp asymptotics of two-point functions as in Theorem \ref{thm-endpoint}. The following proposition deals with \emph{spatial fractional derivatives:}
\begin{prop}[Bounds on spatial fractional derivatives]\label{prop-LEcoefficients2} 
Under the assumptions of Theorem~\ref{thm-endpoint}, there exist $\delta_i >0$ for $i=1,2,3$ such that
\begin{enumerate}
    	\item 
    	\begin{equation}\label{eqSpatDerivative} 
    		\sum_{x} \lvert x \rvert^{(\am)+\delta_1} \lvert  \Pi_1 (x) \rvert < \infty,
	\end{equation}
    
    	\item 
	\begin{equation} \label{eqSD2} 
		\sum_{x} \sum_{m \ge 0} \lvert x \rvert^{\delta_2} m \bar \pi_m (x) < \infty,
	\end{equation}
    
    	\item \begin{equation} \label{eqSD3} 
		\sum_{x} \sum_{m \ge 0} \lvert x \rvert^{\delta_3} \bar  \psi_m (x)  < \infty.
	\end{equation}
    
    	\item  There exists $C,\delta_4>0$ such that, uniformly in $n\ge 1$,
        \begin{equation} \label{eqSD4} 
        		\sum_{x,y} \sum_{m=0}^n (|x|^{\delta_4} + |y|^{\delta_4})\lvert \pi_{m,n}(y,x) \rvert < C.
        \end{equation}
        
        \item There exists a $C'>0$ such that for all $k \in \R^d$,
        \begin{equation}\label{eqNewPiKBd}
        	\sum_x \sum_{m \ge 0} [1-\cos(k \cdot x)] \bar \pi_m(x) \le C' [1-\hat D(k)].
				\end{equation}
\end{enumerate}
\end{prop}
Proposition~\ref{prop-LEcoefficients2}(i) is a direct consequence of \cite[Proposition~2.5 and Remark~2.6]{HeyHofHul14a} combined with Remark~\ref{rem-HS-exp} above. \longversion{We prove the other parts of this proposition in Appendix~\ref{sec:spatfrac}.}\shortversion{The remaining parts of this proposition are proved in \cite[Appendix~B.2]{HeyHofHulMie17b}.} (Note that even though $|\hat \Pi_1(0) - \hat \Pi_1 (k)| = \sum_x \sum_{m \ge 0} [1-\cos(k \cdot x)] |\pi_m(x)|$, Proposition~\ref{prop-LEcoefficients2}(v) does not follow directly from Proposition~\ref{prop:positivity}(1.c) because it may be the case that $\bar \pi_m(x)  \ge |\pi_m(x)|$.) 
\medskip

The upcoming proposition contains the most involved bounds of this section. 
It gives bounds on \emph{temporal fractional derivatives} (with $m$ playing the role of ``time''). It is here that the infrared bounds in Proposition \ref{prop:positivity} will be crucial.
\begin{prop}[Bounds on temporal fractional derivatives]\label{prop-LEcoefficients}
Under the assumptions of Theorem~\ref{thm-endpoint} there exists $\tilde c >0$ and $\eps\in(0,d-3(\am) \wedge 1)$ such that
\begin{enumerate}
\item 
	\begin{equation}\label{eqTempDerivative}
	\sum_{x}\sum_{m \ge 1} m^{1+\vep}\bar{\pi}_m(x) \le \tilde c \beta^{1/6};
	\end{equation}

\item 
	\begin{equation}
	\label{psi-bd}
	\sum_{x}\sum_{m \ge 1} m^{\vep} \bar{\psi}_m (x) \le  \tilde c \beta^{1/6}.
	\end{equation}
\end{enumerate}
\end{prop}

The proof of this proposition uses the full power of our new lace expansion, since only this lace expansion makes the $m$-dependence explicit. The proof of \eqref{eqTempDerivative}, given below, is based on the following lemma that combines the \emph{temporal fractional derivatives} as described in \cite[Section~6.3]{MadSla93} and the diagrammatic bounds developed in Section~\ref{sec-bounds-lace-exp}.

Define the \emph{square diagram} $\squz$ as
\begin{equation}\label{e:squzdef}
        \squz   := \sup_{a,b,c \in \Td}  \int\limits_{\Td} \frac{\d^d k}{(2 \pi^d)} 
        |\hat{D} (k+a)| \htau(k+a) \htau(k+b) \htau(k+c) z |\hat \Tau_z (k)|.
\end{equation}
\begin{lemma}[Diagrammatic bounds]\label{lem:Diagbd}
Under the assumptions of Theorem \ref{thm-endpoint} there exists a constant $\beta>0$ such that for $N \ge 0$, $C_1, C_2 > 0$ and $z\in[0,1]$,
\begin{equation}\label{eqdiagbd}
	\sum_{x} \sum_{m \ge 1} m^2 \bar{\pi}_{m}^{\sss (N)}(x) z^m
        	\le C_1 N^3 \left(C_2 \beta^{1/6}\right)^{(N-2)\vee 0}   \squz.
\end{equation}
\end{lemma}
Lemma~\ref{lem:Diagbd} is proved in \longversion{Appendix~\ref{sec-squares}.}\shortversion{\cite[Appendix~B]{HeyHofHulMie17b}.}
Here is a very brief outline of the proof: 
Using Constructions~$1_{\mathsf{b}}$ and~$1^z$ defined in Section~\ref{sec-bounds-lace-exp} we ``distribute'' the factor $m^2$ over the $\bar\pi_m^{\sss (N)}$ diagrams by adding vertices to backbone lines and changing one two-point function $\tau$ to pivotal generating functions $\Tau_z$. 
\longversion{Then, in Appendix~\ref{sec-squares}, we bound the resulting diagrams in terms of triangle diagrams and a single square diagram (weighted with the factor $z$ originating from the $\Tau_z$ factor), using the bounds already derived in Section~\ref{sec:LEprop} and Appendix~\ref{sec-LEcoefficients}.}
\shortversion{In \cite[Appendix~B]{HeyHofHulMie17b}, the resulting diagrams are then bounded in terms of triangle diagrams and a single square diagram (weighted with the factor $z$ originating from the $\Tau_z$ factor), using the bounds derived in Section~\ref{sec:LEprop} and \cite[Appendix~A]{HeyHofHulMie17b}.}
The final step of the proof is then to bound the triangle diagrams using the strong triangle condition to obtain the factors of $\beta$.

\proof[Proof of Proposition \ref{prop-LEcoefficients} subject to Lemma \ref{lem:Diagbd}.]
We start with the bound on the left-hand side of \eqref{eqTempDerivative}.
For $\vep \in (0,1)$ we have the identity (cf.\ \cite[(6.3.5)]{MadSla93})
\begin{equation}
    	m^{\vep} = \frac{m}{(1-\vep) \Gamma(1-\vep)} 
	\int\limits_{0}^{\infty} \d \lambda\, \e^{-m \lambda^{1/(1-\vep)}},
\end{equation}
which implies, using Fubini's Theorem along the way, that
\begin{equation}
    	\sum_{N \ge 0} \sum_{x} \sum_{m \ge 1} m^{1+\vep} \bar{\pi}_{m}^{\sss (N)} (x) 
	= \frac{1}{(1-\vep)\Gamma(1-\vep)}\sum_{N \ge 0}  
	\int\limits_{0}^{\infty} \d \lambda\,\sum_{x} \sum_{m \ge 1} m^2 \bar{\pi}_{m}^{\sss (N)} (x) \e^{-m \lambda^{1/(1-\vep)}}.
\end{equation}
Applying Lemma \ref{lem:Diagbd} with $z_{\lambda} = \e^{-\lambda^{1/(1-\vep)}}$ to the right-hand side gives
\begin{equation}\label{square2int}
    	C \sum_{N \ge 0}  \int\limits_{0}^{\infty} \d \lambda\, \sum_{x} \sum_{m \ge 1}
	m^2 \bar{\pi}_{m}^{\sss (N)} (x) \e^{-m \lambda^{1/(1-\vep)}} 
	\le  \sum_{N \ge 0} C_1 N^3 \left(C_2 \beta^{1/6}\right)^{(N-2)\vee 0} 
	\int\limits_{0}^{\infty} \d \lambda\,\square_{z_\lambda}.
\end{equation}
Applying the bound $\htau(k) \le C/[1-\hat{D}(k)]$ (cf.\ \eqref{eq-T-bd-pos2} or \cite{HarSla90a, HeyHofSak08}), we obtain
\begin{equation}
	\label{integralsplit}
    	\int\limits_{0}^{\infty} \d \lambda\,\square_{z_\lambda} \le
	C \sup_{a,b,c \in \Td} \int\limits_{\Td} \frac{\d^d k}{(2 \pi^d)} 
	\frac{|\hat D(k+a)|}{1-\hat{D}(k+a)}\frac{1}{1-\hat D(k+b)} \frac{1}{1-\hat D(k+c)} 
	\int\limits_{0}^{\infty} \d \lambda\, z_{\lambda} |\hat \Tau_{z_\lambda} (k)|.
\end{equation}

The aim is now to show that the integral over $\lambda$ is small compared to $[1-\hat{D}(k)]^{-1}$, so that the integral over $k$ is effectively the same as the integral over a triangle diagram.
For any $\vep \in (0,1)$ we can use both upper bounds in Proposition \ref{prop:positivity}(ii) to bound the integral over $\lambda$: substitute $s = \lambda^{1/(1-\vep)}$, then
\begin{equation}
    \begin{split}
        \int\limits_{0}^{\infty} \d \lambda\, z_\lambda \hat{\Tau}_{z_\lambda}(k)
        & = \int\limits_{0}^{\infty} \d s\,(1-\vep) s^{-\vep} \e^{-s} | \hat{\Tau}_{\e^{-s}} (k)| \\
        & \le c_2 (1-\vep) \int\limits_{0}^{1-\hat{D}(k)} \d s\, \frac{s^{-\vep}}{1 - \hat{D}(k)} 
        + c_1 (1-\vep) \int\limits_{1-\hat{D}(k)}^{\infty} \d s\, \frac{s^{-\vep} \e^{-s}}{1- \e^{-s}} \\
        & \le c_2 [1-\hat{D}(k)]^{-\vep} + c_1 (1-\vep) \int\limits_{1 -\hat{D}(k)}^{1} \d s 
        \left(s^{-1-\vep} + \frac{s^{-\vep}}{2} + o(1)\right) 
        + c_2 (1-\vep) \int\limits_{1}^{\infty} \d s\, \e^{-s} \\
        & \le c_2 [1-\hat{D}(k)]^{-\vep} + \frac{c_1 (1-\vep)}{\vep} [1-\hat{D}(k)]^{-\vep} 
        + O(1) \le C_\vep [1-\hat{D}(k)]^{-\vep},
    \end{split}
\end{equation}
where $C_{\vep}$ is a constant that depends only on $\vep$. When we apply this bound to \eqref{integralsplit}  we get 
\begin{equation}
	\label{integralsplit2}
    	\int\limits_{0}^{\infty} \d \lambda\,\square_{z_\lambda} \le
	C \sup_{a,b,c \in \Td} \int\limits_{\Td} \frac{\d^d k}{(2 \pi^d)} 
	\frac{|\hat D(k+a)|}{1-\hat{D}(k+a)}\frac{1}{1-\hat D(k+b)} 
	\frac{1}{1-\hat D(k+c)} \frac{1}{[1-\hat D(k)]^{\vep}}.
\end{equation}
By an application of H\"older's inequality (see \cite[Lemma 7.3]{HeyHofHul14a}), the facts that $|\hat D(k)| \le 1$ (by Parseval's Theorem),  that $[1-\hat D(k)]^{-1}$ is real and non-negative by Assumption \ref{ass:D}, that the Fourier transform is periodic with period $2\pi$ along each axis by Assumption~\ref{ass:A}, and by \cite[Proposition 2.2]{HeyHofSak08} this can be bounded by
\begin{equation}
	\left(\,\, \int\limits_{\Td} \frac{\d^d k}{(2 \pi^d)} 
	\frac{\hat D(k)^{2}}{[1-\hat{D}(k)]^{3+\vep}} \right)^{\frac{1}{3 + \vep}} 
	\left( \,\,\int\limits_{\Td} \frac{\d^d k}{(2 \pi^d)} 
	\frac{1}{[1-\hat{D}(k)]^{3+\vep}} \right)^{\frac{2 + \vep}{3 + \vep}} \le C \beta^{1/4},
\end{equation}
where the final inequality holds when $\beta$ is small enough, $\vep < (d - 3\twa) \wedge 1$, and $d>3(\am)$ by \eqref{hatD-asymp}.
Insterting this bound into \eqref{square2int} completes the proof of \eqref{eqTempDerivative} in Proposition \ref{prop-LEcoefficients}. 
\medskip

Recall \eqref{psiexpansion}, \eqref{psimdef} and \eqref{psi-pi-rep}, from which the claim in \eqref{psi-bd} in Proposition \ref{prop-LEcoefficients} follows.
\qed
\medskip

\noindent
\paragraph{\bf 
Where we stand now.} 
The proof of Proposition \ref{prop:smallpi} as well as the proofs of Lemmas \ref{lem:rhon-exp-lemma} and \ref{lem:rhon-exp-lemma-2} are in \longversion{Appendix~\ref{sec-LEcoefficients},}\shortversion{\cite[Appendix~A]{HeyHofHulMie17b},} and only rely on the infrared bound $0 \le \hat{\tau}^{p_c}(k)\leq C/[1-\hat{D}(k)]$ as proved in e.g.\ \cite{HarSla90a, HeyHofSak08}. 
Furthermore, \longversion{Appendix~\ref{sec-squares}}\shortversion{\cite[Appendix~B]{HeyHofHulMie17b}} contains the proofs of Proposition~\ref{prop-LEcoefficients2}(ii)--(iv) and Lemma~\ref{lem:Diagbd}. 
We thus validated the lace expansions of $\tau_{n}(x)$ (Proposition~\ref{prop-lace-exp}) and $\rho_n(x)$ (Proposition~\ref{prop-lace-exp-2}) and verified all the bounds from Sections~\ref{sec:LEprop} and \ref{sec-impr-bds-LE}, which will be used frequently in the remainder of this paper. We now continue by proving our main results, starting with the proof of Theorem~\ref{thm-endpoint} in the next section.

\section{Scaling of end-to-end displacement: proof of Theorem \ref{thm-endpoint}}
\label{sec-pf-endpoint}
In this section we use Propositions~\ref{prop-LEcoefficients2} and \ref{prop-LEcoefficients} to prove Theorem~\ref{thm-endpoint}. We start with Theorem \ref{thm-endpoint} for $\rho_n$, which is slightly simpler.

\subsubsection*{Analyzing the expansion for $\Rho_z$.}\label{secRho}
We proceed by proving Theorem \ref{thm-endpoint}, as well as some consequences of these propositions that we will state below.
Our argument is roughly the same as the arguments of \cite{CheSak09} and \cite[Section 2.1]{Heyd11}. The main difference is that we have to deal with the term $\hat\Psi_z(k)$ in the numerator of \eqref{eqProof1}, which complicates the analysis significantly.

Recall from \eqref{eqRhoExp} that
\begin{equation}\label{eqProof1}
    	\hat\Rho_z(k)=\frac{\hat\Psi_z(k)}{1-zp_c \Dk\hat\Pi_z(k)} 
	=\frac{\hat\Psi_z(0)}{1-zp_c \Dk\hat\Pi_z(k)}
	-\frac{\hat\Psi_z(0)-\hat\Psi_z(k)}{1-zp_c \Dk\hat\Pi_z(k)}.
\end{equation}
We first rewrite the denominator on the right-hand side of \eqref{eqProof1} as
\begin{equation}\label{eqProof1b}
    	{1-zp_c \Dk\hat\Pi_z(k)} =(1-z)\hat A(k)+\hat B(k)+\hat E_z(k),
\end{equation}
with
\begin{eqnarray}
  \label{eqProof3}
  	\hat A(k) &:=& p_c\Dk\left(\hat\Pi_1(k)+\partial_z\hat\Pi_z(k)|_{z=1}\right),\\
  \label{eqProof4}
	\hat B(k) &:=&  [1-\Dk]+p_c\Dk\left(\hat\Pi_1(0)-\hat\Pi_1(k)\right),\\
  \label{eqProof5}
  	\hat E_z(k) &:=& zp_c\Dk\left(\hat\Pi_1(k)-\hat\Pi_z(k)\right)-(1-z)\,p_c\Dk\;\partial_z\hat\Pi_z(k)|_{z=1}.
\end{eqnarray}
We then write
\begin{equation}\label{eqProof7}
  	\hat\Rho_z(k) = \frac{\hat\Psi_z(0)}{(1-z)\hat A(k)+\hat B(k)} -\hat \Theta_z(k),
\end{equation}
where
\begin{equation}\label{eqProof8}
    	\hat \Theta_z(k) =\frac{\hat\Psi_z(0)\,\hat E_z(k)}{\left((1-z)\hat A(k)
	+\hat B(k)\right)\left(1-zp_c \Dk\hat\Pi_z(k)\right)} 
	+\frac{\hat\Psi_z(0)-\hat\Psi_z(k)}{1-zp_c \Dk\hat\Pi_z(k)}.
\end{equation}
For the first (and main) term in (\ref{eqProof7}) we write
\begin{equation}\label{eqProof9}
   	 \frac{\hat\Psi_z(0)}{(1-z)\hat A(k)+\hat B(k)} 
	 =\frac{\hat\Psi_z(0)}{\hat A(k)+\hat B(k)}
	 \sum_{n \ge 0} z^n\left(\frac{\hat A(k)}{\hat A(k)+\hat B(k)}\right)^n.
\end{equation}
The geometric sum converges whenever $|z|<(\hat A(k)+\hat B(k))/\hat A(k)$ which approaches 1 as $|k|\to0$. For $|z|<1$, we can write $\hat \Theta_z(k)$ as a power series in $z$ as well, i.e.,
\begin{equation}\label{eqProof10}
    	\hat \Theta_z(k)=\sum_{n \ge 0} \hat \theta_n(k)\,z^n.
\end{equation}

\begin{prop}[Lace expansion error bound for $\Rho_z$]
\label{lemmaError}
Under the assumptions of Theorem~\ref{thm-endpoint} there exist $0<\eps<1$ and $N=N(\vep)$ such that for $n \ge N$, 
\begin{equation}
	|\hat \theta_n(k)|\le O(n^{-\eps} + |k|^\vep \log n) \qquad \text{ for all $k\in\R^d$}.
\end{equation}
\end{prop}

The proof of this proposition relies on certain decay rates of $\hat E_z(k)$ and is deferred to Section \ref{sec-ErrorBound}.

Recall that $\hat\Rho_z(k)=\sum_{n \ge 0} \hat \rho_n(k)z^n$ and $\hat \Psi_z(k) = \sum_{n \ge 0} \hat \psi_n (k) z^n$, and recall the definition of $f_\alpha(n)$ from \eqref{fn-def}. Proposition \ref{lemmaError} implies that there exist $0<\eps<1$ and $N=N(\vep)\ge 0$ such that for all $n\ge N$ and $k\in\Rd$ that satisfy $k_n=f_\alpha (n) k$,
\begin{equation}\label{eqRhoKn}
    	\hat \rho_n(k_n) =\frac{1}{\hat A(k_n)
	+\hat B(k_n)}\sum_{m=0}^n \hat\psi_m(0)
	\left(\frac{\hat A(k_n)}{\hat A(k_n)+\hat B(k_n)}\right)^{n-m}
    	+O\left(n^{-\eps} + |k_n|^\vep \log n\right).
\end{equation}

To determine the values of these coefficients, we proceed to study the small-$k$ behavior of the lace-expansion coefficients $\hat \Psi_1(k)$ and $\hat \Pi_1(k)$, as well as the decay properties of $\hat{\psi}_m(0)$ for $m$ large. That is the content of the upcoming propositions.
\begin{prop}[Small $k$ behaviour of $\hat \Pi_1$]\label{propMain2}
Under the assumptions of Theorem \ref{thm-endpoint},
\begin{equation}\label{eqMain2}
	   \lim_{|k|\to 0}\frac{\hat\Pi_{1}(0)-\hat\Pi_{1}(k)}{1-\Dk} =
	   \begin{cases}
	            1,\quad&\text{if $\alpha\le2$};\\
	            1+ (2dv_\alpha)^{-1}\sum_{x}|x|^2\,\Pi_{1}(x),\quad&\text{if $\alpha>2$,}
	   \end{cases}
\end{equation}
where $v_\alpha$ is the constant from Assumption~\ref{ass:E}.
\end{prop}
The proof of Proposition \ref{propMain2} relies on Assumptions~\ref{ass:A}, \ref{ass:D}, and \ref{ass:E}, and on Proposition~\ref{prop-LEcoefficients2}. It is identical to the proof of \cite[Proposition 2.3]{Heyd11}, so we omit it here (but mind that in \cite{Heyd11} the term $1$ has been explicitly excluded from $\hat\Pi_z(k)$ in the expansion).
\medskip

By Proposition \ref{propMain2},
\begin{equation}\label{eqProof13}
	\lim_{|k|\to0}\frac{\hat B(k)}{1-\hat D(k)}=
        \begin{cases}
            1+p_c, \quad&\text{if $\alpha\le2$};\\
            1+p_c(2dv_\alpha)^{-1}\sum_{x}|x|^2\,\Pi_{1}(x), \quad&\text{if $\alpha>2$.}
        \end{cases}
\end{equation}
Moreover, by \eqref{eqDefPi} and \eqref{eqPcPi},
\begin{equation}
	\lim_{|k| \to 0} \hat A(k) = 1 + p_c \sum_x \sum_{n \ge 0} n \pi_n(x).
\end{equation}

Observe furthermore that $k_n$ has been chosen such that by Assumption~\ref{ass:E},
\begin{equation}\label{[1-D]-n-asymp}
    	\lim_{n\rightarrow \infty} n[1-\hat D(k_n)]=|k|^{(\am)},\qquad k\in\R^d.
\end{equation}
If a sequence $h_n$ converges to a limit $h$, then the sequence $(1-h_n/n)^n$ converges to $\e^{-h}$. We apply this to \eqref{eqRhoKn} with
\begin{equation}\label{eqProof14}
	h_n=\frac{n\,\hat B(k_n)}{\hat A(k_n)+\hat B(k_n)} 
	=n[1-\hat D(k_n)]\;\frac{\hat B(k_n)}{1-\hat D(k_n)}\;\frac{1}{\hat A(k_n)+\hat B(k_n)}
	\quad\stackrel{n\to\infty}{\longrightarrow} \quad K_\alpha\,|k|^{(\am)},
\end{equation}
where
\begin{equation}\label{eqDefKalpha}
    	K_{\alpha} :=\lim_{n\to\infty}\left(\frac{\hat B(k_n)}{1-\hat D(k_n)}\cdot\frac1{\hat A(k_n)
    	+\hat B(k_n)}\right)
	=
	\begin{cases}\displaystyle
            	\frac{1+p_c}{1+p_c\sum_{x}\sum_{n \ge 0}n\pi_{n}(x)},
            	\hskip.5em&\text{if $\alpha\le2$};\\  \\
            	\displaystyle
            	\frac{1+p_c(2dv_\alpha)^{-1}\sum_{x}\sum_{n \ge 0}|x|^2\,\pi_{n}(x)}
            	{1+p_c\sum_{x}\sum_{n \ge 0}\,n\pi_{n}(x)},
            	\hskip.5em&\text{if $\alpha>2$.}
        \end{cases}
\end{equation}
This gives an explicit description of the constant $K_\alpha$ that was introduced in Theorem \ref{thm-endpoint}. 
Moreover, $0<K_\alpha < \infty$ by Propositions~\ref{prop:smallpi} and \ref{prop-LEcoefficients2}. 

We now split \eqref{eqRhoKn} as
\begin{equation}\label{eqRhoKn-b}
    	\hat \rho_n(k_n) =\left(\frac{\hat A(k_n)}{\hat A(k_n)
	+\hat B(k_n)}\right)^{n} \frac{\hat\Psi_1(0)}{\hat A(k_n)+\hat B(k_n)}
	+\hat{\xi}_n(k_n)
    	+O\left(n^{-\eps} + |k_n|^\vep \log n\right),
\end{equation}	
where, with $\hat{a}(k_n):=\hat A(k_n)/(\hat A(k_n)+\hat B(k_n))$,
\begin{equation}\label{eqRhoKn-c}
	\hat{\xi}_n(k_n)=\frac{1}{\hat A(k_n)+\hat B(k_n)}\sum_{m=0}^n \hat\psi_m(0)
	\Big[\hat{a}(k_n)^{n-m}-\hat{a}(k_n)^{n}\Big]
	-\frac{\hat{a}(k_n)^{n}}{\hat A(k_n)+\hat B(k_n)}\sum_{m \ge n+1} \hat\psi_m(0).
\end{equation}
Noting that $|\hat{a}(k_n)|\leq 1$, we can bound the second contribution in \eqref{eqRhoKn-c} by
\begin{equation}
	\frac{\hat{a}(k_n)^{n}}{\hat A(k_n)+\hat B(k_n)}\sum_{m \ge n+1} |\hat\psi_m(0)|
	\leq \frac{1}{\hat A(k_n)+\hat B(k_n)}n^{-\vep} \sum_{m \ge 0} m^{\vep} |\hat\psi_m(0)|
	=O(n^{-\vep})
\end{equation}
by Proposition \ref{prop-LEcoefficients}(ii). For the first contribution in \eqref{eqRhoKn-c}, we use that for $a,\vep \in [0,1]$, and $0\leq m\leq n$,
\begin{equation}
	a^{n-m}-a^n
	 = a^{n-m}(1-a^m)^{1-\eps}\left(\frac{1-a^m}{1-a}\right)^\eps(1-a)^\eps
	\leq 1 \cdot \left(\sum_{k=0}^{m-1}a^k\right)^\eps(1-a)^\eps
	\leq m^\vep(1-a)^{\vep}.
\end{equation}
Therefore, the first contribution in \eqref{eqRhoKn-c} is bounded by
\begin{equation}
	\frac{1}{\hat A(k_n)+\hat B(k_n)}\sum_{m=0}^n \big(1-\hat{a}(k_n)\big)^{\vep} m^{\vep}|\hat\psi_m(0)|
	=O(n^{-\vep}),
\end{equation}
since \eqref{eqProof14} implies that $1-\hat{a}(k_n)=O(n^{-1})$. We conclude that
\begin{equation}\label{eqRhoKn-d}
    	\hat \rho_n(k_n) =\left(\frac{\hat A(k_n)}{\hat A(k_n)
	+\hat B(k_n)}\right)^{n} \frac{\hat\Psi_1(0)}{\hat A(k_n)+\hat B(k_n)}
    	+O\left(n^{-\eps} + |k_n|^\vep \log n\right).
\end{equation}		
	
When we apply \eqref{eqProof14} to \eqref{eqRhoKn-d}, we arrive at
\begin{equation}\label{eqProof15}
	\lim_{n\to\infty}\hat \rho_n(k_n) 
	=\frac{\hat\Psi_1(0)}{\hat A(0)}\exp\!\left(-K_\alpha\,|k|^\am\right).
\end{equation}
Finally, if we set $k=0$ in \eqref{eqProof15} and we use $\hat\rho_n(0)=1$, then we get $\hat\Psi_1(0)=\hat A(0)$, and so we have identified the first limit in \eqref{EndpointDistribution} subject to Propositions \ref{lemmaError} and \ref{propMain2}. Next we consider the limit behaviour of $\hat\tau_n$.


\subsubsection*{Analyzing the expansion for $\Tau_z$.}\label{secTau} 
We next extend the analysis to $\hat\Tau_z(k)$. We start from
	\begin{equation}
	\label{eqProof2a}
    	\hat\Tau_z(k)^{-1} =  \frac1{\hat\Pi_z(k)} -zp_c\Dk,
	\end{equation}
which follows from \eqref{eqTauExp}.

Recall that $p_c\hat\Pi_1(0)=p_c\hat\Pi_1^{p_c}(0)=1$ by \eqref{eqPcPi}. Like we did for $\hat\Rho_z(k)$, we write 
\begin{equation}
	\hat\Tau(k)^{-1}=(1-z)\hat A'(k)+\hat B'(k)+\hat E'_z(k),
\end{equation}
where
\begin{eqnarray}
  \label{eqProof3a}
  	\hat A'(k) &:=&	-\partial_z\Bigg(\frac1{\hat\Pi_z(k)}-zp_c\Dk\Bigg)\Bigg|_{z=1}
	=\frac{\partial_z\hat\Pi_z(k)|_{z=1}}{\hat\Pi_1(k)^2}+p_c\Dk,\\
  \label{eqProof4a}
	\hat B'(k) &:=&
	\frac1{\hat\Pi_1(k)}-\frac1{\hat\Pi_1(0)}+p_c[1-\Dk],\\
  	\label{eqProof5a}
  	\hat E'_z(k) &:=&
	\left(\frac1{\hat\Pi_z(k)}-\frac1{\hat\Pi_1(k)}\right)
	+(1-z)\partial_z\frac1{\hat\Pi_z(k)}\Bigg|_{z=1}.
\end{eqnarray}
Now we repeat the analysis as for $\hat{\Rho}_z (k)$ 
but with $\hat{\Psi}_z(k)$ replaced by 1,
so that $\hat{\xi}_n(k_n)$ in \eqref{eqRhoKn-b} vanishes. Define $\hat \Theta'_z(k)$ and $\hat \theta'_n(k)$ analogously to \eqref{eqProof8} and \eqref{eqProof10}. Then we obtain the following bound on the error term:

\begin{prop}[Lace expansion error bound for $\hat{\Tau}_z$]
\label{lemmaError-Tau}
Under the assumptions of Theorem~\ref{thm-endpoint}, there exist $0<\eps<1$ and $N=N(\vep)$ such that, for $n \ge N$, 
\begin{equation}
	|\hat \theta_n'(k)|\le O(n^{-\eps} + |k|^\vep \log n) \qquad \text{ for all $k\in\R^d$}.
\end{equation}
\end{prop}

The proof of Proposition \ref{lemmaError-Tau}  follows the same steps as that of Proposition \ref{lemmaError}, and is performed in Section \ref{sec-ErrorBound}.
In particular, the bounds in Section \ref{sec-ErrorBound} and the proof of Lemma \ref{lem:errorhelp} require only minor modifications to yield the correct bounds, if we replace $\hat \Psi$ by $1$, $\hat A$ by $\hat A'$, $\hat B$ by $\hat B'$, and $\hat E_z$ by $\hat E'_z$.
It follows that if we use \eqref{eqPcPi} and if we reason as in \eqref{eqProof9}--\eqref{eqProof15},
we find that 
\begin{equation}
	\lim_{n\to\infty}\hat\tau_n(k_n)=A\exp\big(-K_\alpha|k|^\twa\big),
\end{equation}
with
\begin{equation}\label{eqDefA'}
	A=\lim_{|k|\to0}\frac1{\hat A'(k)+\hat B'(k)}
	=\frac1{\hat A'(0)}
	=\frac1{p_c^2\sum_{x}\sum_{n \ge 0} n \,\pi_{n}(x)  + p_c},
\end{equation}
and
\begin{equation}\label{eqDefK'}
	K_\alpha
	=\lim_{|k|\to0}\left(\frac{\hat B'(k)}{1-\hat D(k)}\cdot\frac1{\hat A'(k)+\hat B'(k)}\right)
	=\begin{cases}\displaystyle
            \frac{1+ p_c}{1+ p_c\sum_{x}\sum_{n\ge 0} n \,\pi_{n}(x)},
            \quad&\text{if $\alpha\le2$};\\  \\
            \displaystyle
            \frac{1+ p_c(2dv_\alpha)^{-1}\sum_{x}\sum_{n\ge 0}|x|^2\,\pi_{n}(x)}
            {1+ p_c\sum_{x}\sum_{n\ge 0}\,n\,\pi_{n}(x)},
            \quad&\text{if $\alpha>2$,}
        \end{cases}
\end{equation}
as we set out to prove.
The limit in \eqref{eqDefK'} is equal to the limit in \eqref{eqDefKalpha}. This completes the proof of Theorem~\ref{thm-endpoint} subject to Propositions~\ref{lemmaError} and \ref{lemmaError-Tau}, which we prove in the upcoming section. \qed


\section{Error bounds for the lace expansion: proof of Propositions \ref{lemmaError} and \ref{lemmaError-Tau}}
\label{sec-ErrorBound}
The proofs in this section involve analyzing several variables introduced in Section~\ref{sec-pf-endpoint}, which are defined as power series in the complex variable $z$. Hence derivatives with respect to $z$, henceforth denoted by $\partial_z$, can be interpreted as term-by-term derivatives. We crucially rely on the infrared bounds in Proposition~\ref{prop:positivity}, and on the improved bounds on the lace-expansion coefficients in Propositions~\ref{prop-LEcoefficients2} and \ref{prop-LEcoefficients}.

We start by stating a lemma that contains the main bounds used in this proof:

\begin{lemma}\label{lem:errorhelp}
There exists $c_1,c_2,c_3,c_4>0$ such that for $\vep$ as in Proposition \ref{prop-LEcoefficients}, all $k\in \R^d$ and $z\in\mathbb C$ with $|z|<1$,
    \begin{eqnarray}
        \label{e:help1} \lvert \hat E_z (k) \rvert &\le&  c_1 \beta^{1/6} |1-z|^{1+\vep},\\
        \label{e:help1a} \lvert \partial_z\hat E_z (k) \rvert &\le&  c_2 \beta^{1/6} |1-z|^{\vep},\\
        \label{e:help2} \lvert 1 - z p_c \hat D (k) \hat \Pi_z (k)\rvert &\ge&  c_3\,|1-z|,\\
        \label{e:help3} \lvert (1-z)\hat A(k) + \hat B (k)\rvert &\ge& c_4 |1-z|.
    \end{eqnarray}
\end{lemma}
We prove this lemma at the end of the section.

\proof[Proof of Proposition \ref{lemmaError} subject to Lemma \ref{lem:errorhelp}]
In this proof we extract bounds on the coefficients of a certain power series from its divergence properties near the critical radius. For this we use the following tool:

Let $f(z)= \sum_{n \ge 0} a_n z^n$ have radius of convergence $1$. Derbez and Slade, in \cite[Lemma~3.2]{DerSla98}, prove the following Tauberian theorem: If there exist $C,b\ge1$ such that for $z\in \mathbb C$, $|z|<1$, the bound $\lvert f(z) |\le C |1-z|^{-b}$ holds, then it follows that $| a_n | \le C' n^{b-1}$ if $b>1$ and $|a_n| \le C' \log n$ if $b=1$ for a uniform constant $C'>0$. 
Consequently, if $\lvert \partial_z f(z) |\le C |1-z|^{-b}$ for all $|z|<1$ and $b>1$, then $| a_n | \le C' n^{b-2}$. 

We apply this reasoning to the power series
\begin{equation}\label{e:priordiscussion}
    	\hat\Theta_z(k)=\sum_{n \ge 0} \hat\theta_n(k)z^n,
\end{equation}
which has radius of convergence 1. 
Using the decomposition in \eqref{eqProof8}, we write
\begin{equation}
    	\hat \Theta_z(k)
    	=  \hat \Psi_z(0)\hat\Theta^{\sss (1)}_{z}(k) + \hat\Theta_z^{\sss (2)}(k),
\end{equation}
with
\begin{align}
	\hat\Theta^{\sss(1)}_{z}(k)&:=
	\frac{\hat E_z(k)}{\left((1-z)\hat A(k)+\hat B(k)\right)
	\left(1-zp_c \Dk\hat\Pi_z(k)\right)},\nnb
	\hat\Theta_z^{\sss (2)}(k)&:=
	\frac{\hat\Psi_z(0)-\hat\Psi_z(k)}{1-zp_c \Dk\hat\Pi_z(k)}.
\label{eqThetaDer}
\end{align}

A key ingredient in the proof of Proposition~\ref{lemmaError} is to show that there exist $0<\eps<1$ and $C>0$ such that
\begin{equation}\label{eqThetaDerBound1}
	|\partial_z\hat\Theta_z^{\sss (1)}(k)|\le C \lvert1-z\rvert^{-(2-\eps)} ,
	\quad \text{ and }\quad	
	|\hat\Theta_z^{\sss (2)}(k)|\le C |k|^{\vep} \lvert1-z\rvert^{-1}
\end{equation}
uniformly in $|z|<1$ and $k$. 
We then apply \cite[Lemma~3.2]{DerSla98} with $b=1$ on $|\hat\Theta_z^{\sss (2)}(k)|$ to obtain a bound of $O( |k|^{\vep} \log{n})$ on this contribution to $\hat{\theta}_n(k)$. 
We proceed now to prove Proposition~\ref{lemmaError} under the assumption that \eqref{eqThetaDerBound1} holds, and then establish \eqref{eqThetaDerBound1} to complete the proof.
\medskip

The terms due to $\hat \Psi_z(0)\hat\Theta^{\sss (1)}_{z}(k)$ are somewhat involved. For it, we use the Tauberian theorem with $b=2-\vep$ on $|\partial_z\hat\Theta_z^{\sss (1)}(k)|$, to obtain
	\eqn{
	\label{theta-1-bd}
	|\hat\theta^{\sss (1)}_{n}(k)|\leq O(n^{-\eps}).
	}
For a power-series $f(z)=\sum_{n\geq 0} a_n z^n$, write $[f(z)]_n :=a_n$. We use the convolution equation to obtain
	\eqn{
	\label{conv-inversion-theta-Psi}
	\big[\hat \Psi_z(0) \hat\Theta^{\sss (1)}_{z}(k)\big]_n = \sum_{m=0}^n \hat \psi_m(0) \hat\theta^{\sss (1)}_{n-m}(k).
	}
We claim that 
\begin{equation}
	\big[\hat \Psi_z(0) \hat\Theta^{\sss (1)}_{z}(k)\big]_n= O(n^{-\eps}).
\end{equation}
For this, we consider separately the cases where 
$m\geq \lceil n/2 \rceil$ or $m\leq \lfloor n/2 \rfloor$ in the sum on the right-hand side of \eqref{conv-inversion-theta-Psi}. The contribution due to $m\geq \lceil n/2 \rceil$ can be bounded for some sufficiently small $\vep$ using \eqref{eqThetaDerBound1} and Proposition~\ref{prop-LEcoefficients}(ii) by
	\eqn{
	 n^{-\vep} \sum_{m= \lceil n/2 \rceil}^n m^{\vep} |\hat \psi_m(0)| |\hat\theta^{\sss (1)}_{n-m}(k)|
	 \leq C n^{-\vep} \sum_{m\geq 0} m^{\vep} |\hat \psi_m(0)| =O(n^{-\vep}),
	 } 
while the contribution due to $m\leq \lfloor n/2 \rfloor$ can be bounded using \eqref{theta-1-bd} and Proposition~\ref{prop-LEcoefficients}(ii) by
	\eqn{
	\sum_{m=0}^{\lfloor n/2 \rfloor} |\hat \psi_m(0)| |\hat\theta^{\sss (1)}_{n-m}(k)|\leq O(n^{-\vep}) \sum_{m\ge 0} |\hat \psi_m(0)| =O(n^{-\vep}),
	}
as required. This establishes Proposition~\ref{lemmaError} subject to \eqref{eqThetaDerBound1}. We complete the proof by establishing the bounds in \eqref{eqThetaDerBound1}.
\medskip

We start by bounding $|\partial_z\hat\Theta^{\sss (1)}_{z}(k)|$. A straightforward calculation yields 
\begin{align}
	\partial_z\hat\Theta^{\sss (1)}_{z}(k)
	=& \frac{\partial_z\hat E_z(k)}{\left((1-z)\hat A(k)
	+\hat B(k)\right)\left(1-zp_c \Dk\hat\Pi_z(k)\right)}
	+ \frac{\hat E_z(k)\hat A(k)}{\left((1-z)\hat A(k)
	+\hat B(k)\right)^2\left(1-zp_c \Dk\hat\Pi_z(k)\right)}\nnb
	&{}+\frac{\hat E_z(k)\left[p_c\Dk\hat\Pi_z(k)
	+zp_c\Dk\partial_z \hat \Pi_z(k)\right]}{\left((1-z)\hat A(k)
	+\hat B(k)\right)\left(1-zp_c \Dk\hat\Pi_z(k)\right)^2}.
\end{align}
Proposition \ref{prop-LEcoefficients} implies
\begin{equation}
    	\left|p_c\Dk\hat\Pi_z(k)+zp_c\Dk\partial_z\hat \Pi_z(k)\right|
    	\le p_c\sum_{x }\sum_{n \ge 1} \left(\lvert z\rvert^{n} 
    	+ n \lvert z\rvert^{n}\right) \bar \pi_n(x)
    	\le C
\end{equation}
for some $C>0$. 
Using this and all four bounds of Lemma \ref{lem:errorhelp} thus establishes 
\begin{equation}\label{eqThetaProof7}
        \lvert \partial_z\hat \Theta_z^{\sss (1)}(k) \rvert 
        \le \frac{C}{|1-z|^{2-\vep}},
\end{equation}
as required.

It remains to show the bound on $\lvert \hat \Theta_z^{\sss (2)}(k) \rvert$. 
When $\eps\in(0,2\wedge\delta_3]$, Proposition \ref{prop-LEcoefficients2}(iii) implies
\begin{align}\label{eqThetaProof8}
	|\hat\Psi_z(0)-\hat\Psi_z(k)|
	\le&\sum_{x}\sum_{n \ge 1}[1-\cos(k\cdot x)]\,|\psi_n(x)|
	\le\sum_{x}\sum_{n \ge 1}(k\cdot x)^\eps\,\bar \psi_n(x)\nnb
	\le&|k|^\eps\sum_{x}\sum_{n \ge 1}|x|^\eps \bar \psi_n(x)
	\le C\,|k|^\eps.
\end{align}
Together with the lower bound  \eqref{e:help2} on the denominator of \eqref{eqThetaDer}, this shows
\begin{equation}
    	\lvert \hat \Theta_z^{\sss (2)}(k) \rvert \le C\lvert k \rvert^\vep \lvert1-z\rvert^{-1},
\end{equation}
as required.
This concludes the proof of the bounds in \eqref{eqThetaDerBound1}, and thus completes the proof of Proposition~\ref{lemmaError}.
\qed
\smallskip

\proof[Proof of Proposition \ref{lemmaError-Tau} subject to Lemma \ref{lem:errorhelp}] We next extend the above analysis to $\hat{\theta}_n'(k)$. This is easier, since $\hat \Theta_z'(k)$ is simpler than $\hat \Theta_z(k)$, namely
	\eqn{
	\hat \Theta_z'(k) =\frac{\hat E_z'(k)}{\left((1-z)\hat A'(k)
	+\hat B'(k)\right)\left(1-zp_c \Dk\hat\Pi_z(k)\right)},
	}
where $\hat A'(k)$, $\hat B'(k)$ and $\hat E_z'(k)$ are defined in \eqref{eqProof3a}, \eqref{eqProof4a} and \eqref{eqProof5a}, respectively. 
Since the bounds in Lemma \ref{lem:errorhelp} continue to hold also for $\hat A'(k)$, $\hat B'(k)$ and $\hat E_z'(k)$, the proof of Proposition~\ref{lemmaError}
can be followed straightforwardly. In fact, this proof is easier, since the difficult factor $\hat{\Psi}_z(0)$ is absent. We omit further details.
\qed

\proof[Proof of Lemma \ref{lem:errorhelp}]
We start with \eqref{e:help1}. Recall the definition of $\hat E_z(k)$ in \eqref{eqProof5}.
We bound
\begin{equation}\label{e:Ezbd1}
    \begin{split}
    	\lvert \hat E_z (k) \rvert &= \left\lvert p_c \hat D (k) \left(\hat \Pi_1 (k) - \hat \Pi_z (k)\right) 
	- (1-z) p_c \hat D (k) \partial_z  \hat \Pi_z(k) \vert_{z=1} 
	- (1-z) p_c \hat D(k) \left(\hat \Pi_1 (k) - \hat \Pi_z (k)\right)\right\rvert\\
    	& \le \lvert1-z\rvert p_c |\hat D (k)| \left\lvert \frac{\hat \Pi_1 (k) 
	- \hat \Pi_z (k)}{1-z} - \partial_z \hat \Pi_z (k) \vert_{z=1} \right\rvert 
	+ \lvert1-z\rvert p_c |\hat D (k)| \left\lvert \hat \Pi_1(k) - \hat \Pi_z (k)\right\rvert\\
    	&\le \lvert1-z\rvert p_c |\hat D (k)| \left\lvert \frac{\hat \Pi_1 (k) - \hat \Pi_z (k)}{1-z} 
	- \partial_z \hat \Pi_z (k) \vert_{z=1} \right\rvert + \tilde c \beta^{1/6} \lvert1-z\rvert^2,
    \end{split}
\end{equation}
where the final inequality follows from Proposition \ref{prop:positivity}(i.a).
To bound the remaining term on the right-hand side of \eqref{e:Ezbd1}, we choose $\eps$ such that Proposition \ref{prop-LEcoefficients} holds.
Then, we use that for any $|z|<1$, $\eps\in(0,1]$ and integers $n\ge0$,
\begin{equation}
 	1-z^n
	=\left(1-z^n\right)^{1-\eps}\left(\frac{1-z^n}{1-z}\right)^\eps\left(1-z\right)^\eps
	=\left(1-z^n\right)^{1-\eps}  \left(\sum_{\ell=0}^{n-1}z^\ell\right)^\eps\left(1-z\right)^\eps,
\end{equation}
so that 
\begin{equation}\label{eq1-zBd}
	|1-z^n| \le 2 n^\eps\left|1-z\right|^\eps.
\end{equation}
and therefore,
\begin{equation}
   	 \sum_{\ell=0}^{n-1}\lvert1-z^\ell\rvert 
   	 \le 2\sum_{\ell=0}^{n-1} \ell^\vep \lvert1-z\rvert^\vep  
    	\le 2n(n-1)^\vep\lvert1-z\rvert^\vep.
\end{equation}
Using Proposition \ref{prop-LEcoefficients} and the above bound, for some sufficiently small $\vep$,
\begin{multline}
	\left|\partial_z\hat\Pi_z(k)|_{z=1}
	-\frac{\hat\Pi_1(k)-\hat\Pi_z(k)}{1-z}\right|
	= \left|\sum_{x}\sum_{n \ge 1} n\,\pi_n(x) \e^{ik\cdot x}
	-\sum_{x}\sum_{n \ge 1} \left(\frac{1-z^n}{1-z}\right)  \pi_n(x) \e^{ik\cdot x}\right|\\
	= \left|\sum_{x}\sum_{n \ge 1} \left(\sum_{\ell=0}^{n-1}(1-z^\ell)\right)  
	\pi_n(x) \e^{ik\cdot x}\right|
	\le 2\lvert1-z\rvert^\eps\sum_{x}\sum_{n \ge 1} n(n-1)^\eps\,\bar \pi_n(x)
	\le 2\tilde c \beta^{1/6}\lvert1-z\rvert^\eps.
	\label{eqThetaProof1}
\end{multline}
Combining \eqref{e:Ezbd1} and \eqref{eqThetaProof1} completes the proof of \eqref{e:help1}.

\proof[Proof of \eqref{e:help1a}] 
Differentiating \eqref{eqProof5} with respect to $z$ and rearranging the terms yields 
\begin{equation}\label{e:Ezbd1a}
    \begin{split}
    	\partial_z\hat E_z (k)
    	=& p_c\Dk \left((1-z)\left(\frac{\hat \Pi_1 (k) - \hat \Pi_z (k)}{1-z} 
	- \partial_z \hat \Pi_z (k) \vert_{z=1}\right)
    	+2(1-z)\left(\partial_z \hat \Pi_z (k) \vert_{z=1}\right)\right)\nnb
    	&{}+p_c\Dk z\left(\partial_z \hat \Pi_z (k) \vert_{z=1}-\partial_z \hat \Pi_z (k)\right).
    \end{split}
\end{equation}
For the second line we estimate
\begin{equation}
	\begin{split}
		\left|\partial_z \hat \Pi_z (k) \vert_{z=1}-\partial_z \hat \Pi_z (k)\right|
		&=\left| \sum_{x}\sum_{n \ge 1} (1-z^{n-1})n\pi_n(x) \e^{ik\cdot x}\right| \\
		& \le 2\lvert1-z\rvert^\eps\sum_{x }\sum_{n \ge 1} n^{1+\eps} \bar \pi_n(x) 
		\le 2\tilde c \beta^{1/6} |1-z|^\vep,
	\end{split}
\end{equation}
by \eqref{eq1-zBd}.  
Consequently, using \eqref{eqThetaProof1} and  Proposition \ref{prop-LEcoefficients}, 
\begin{equation}
\begin{split}
	\left|\partial_z\hat E_z (k)\right|
	&\le 2\tilde c\beta^{1/6}\lvert1-z\rvert^{1+\eps}+\tilde c\beta^{1/6}\lvert1-z\rvert
	+2\tilde c\beta^{1/6}\lvert1-z\rvert^{\eps} \\
	&\le (8\tilde c+2\tilde c+2\tilde c)\beta^{1/6}\lvert1-z\rvert^{\eps},
\end{split}
\end{equation}
as claimed.

\proof[Proof of \eqref{e:help2}] By \eqref{eqProof2a} and Proposition \ref{prop:positivity}(i) and (ii), there exists a $c' >0$ such that
\begin{equation}\label{eqLemii}
    	\lvert 1 - z p_c \hat D (k) \hat \Pi_z (k) \rvert 
	= \left\lvert \frac{\hat \Pi_z (k)}{\hat \Tau_z (k)} \right\rvert \ge c_3 \lvert1-z\rvert.
\end{equation}

\proof[Proof of \eqref{e:help3}] By \eqref{eqTauExp} and \eqref{eqProof1b} we can bound
\begin{equation}
    	\lvert (1-z) \hat A (k) + \hat B (k) \rvert 
	\ge \left\lvert \frac{\hat \Pi_z (k)}{\hat \Tau_z (k)}\right\rvert - \lvert\hat E_z (k) \rvert 
    	\ge c_3\lvert1-z\rvert - c_1 \beta^{1/6} \lvert1-z\rvert^{1+\vep} 
    	\ge c_4 \lvert1-z\rvert,
\end{equation}
when $\beta$ is small enough. 
For the second inequality we used \eqref{e:help1} and \eqref{e:help2}.
\qed


\section{The mean-$r$ displacement: Proof of Theorem \ref{thm-xiR}}
\label{sec-meanr}

We follow the proof of \cite[Theorem 1.4]{Heyd11}.
Write $x_1$ for the first coordinate of the vector $x\in\Zd$.
Since $|x_1|^r\le|x|^r\le d^{r/2}\sum_{i=1}^d|x_i|^r$, and by Assumption~\ref{ass:A} the model is invariant under rotation by $\pi/2$, it is sufficient to prove the existence of constants $c,C>0$ such that
\begin{equation}\label{eqProofXiR0}
    	c\,f_\alpha(n)^{-r}
    	\stackrel{(a)}{\le}\sum_{x}|x_1|^r\rho_{n+1}(x)
    	\stackrel{(b)}{\le} \sum_{x}|x_1|^r(J \ast \tau_n)(x)
    	\stackrel{(c)}{\le} C f_\alpha(n)^{-r},
\end{equation}
where $f_\alpha(n)$ is the scaling function defined in \eqref{fn-def}.

Inequality (b) follows directly from Lemma \ref{lem-rhon-taun}.

Now we prove inequality (a). Write $u_n := (f_\alpha(n),0,\dots,0)\in\Rd$.
Theorem \ref{thm-endpoint} implies
\begin{equation}
	\lim_{n\to\infty}1-\hat\rho_n(u_n)=1-\e^{-K_\alpha}>0.
\end{equation}
Hence, for $n$ sufficiently large,
\begin{equation}
	0<\frac12 (1-\e^{-K_\alpha})\le 1-\hat\rho_n(u_n)
	= \sum_{x}\big[1-\cos(f_\alpha(n)\, x_1)\big]\rho_n(x)
	\le\sum_{x}f_\alpha(n)^r|x_1|^r\rho_n(x),
\end{equation}
where the last bound follows from $1-\cos(t)\le|t|^r$ for any $r\le(\am)$.
This implies inequality (a).

Inequality (c) in \eqref{eqProofXiR0} requires a little more work.
We start by observing that for any $r\in(0,2)$ there exists $c_r\in(0,\infty)$ such that
\begin{equation}\label{eqIntId}
	t^r=c_r\int_0^\infty\frac{1-\cos(ut)}{u^{1+r}}\,\d u
\end{equation}
for all $t>0$.
With this in hand, we write $\uv=(u,0,\dots,0)$ for $u>0$ and obtain
\begin{multline}
	\sum_{x}|x_1|^r(J \ast \tau_n)(x)
	= c_r\sum_{x}\int_0^\infty
	\frac{\d u}{u^{1+r}}\, [1-\cos(\uv\cdot x)]\,(J \ast \tau_n)(x)\\
	\le p_c c_r\sum_{x,y}\int_0^\infty
	\frac{\d u}{u^{1+r}}\, 2\big([1-\cos(\uv\cdot y)]+[1-\cos(\uv\cdot (x-y))]\big)\,D(y)\, \tau_n(x-y).
\end{multline}
The split $1-\cos(u+v)\leq 2[1-\cos(u)]+2[1-\cos(v)]$ that we used is proved in e.g. \cite[Lemma~2.13]{FitHof15a}.
The term involving $[1-\cos(\uv\cdot y)]$ is of smaller order.
Indeed, using \eqref{e:hatD4}, 
\begin{equation}
\begin{split}
	\sum_{y}\int_0^\infty \frac{\d u}{u^{1+r}}\,[1-\cos(\uv\cdot y)]D(y)
	&=\int_0^\infty\frac{\d u}{u^{1+r}}\,[1-\hat D(\uv)] \\
	&\le \int_0^\varepsilon\frac{\d u}{u^{1+r}}\,[1-\hat D(\uv)]
	+\int_\varepsilon^\infty\frac{2\,\d u}{u^{1+r}}
	\le C
\end{split}
\end{equation}
for a suitably chosen constant $C$. Furthermore, by Corollary~\ref{corr:pivballboundary}, for fixed $y \in \Zd$,
\begin{equation}
	\sum_x \tau_n (x-y) = O(1).
\end{equation}
Since $\sum_y D(y) =1$, by translation invariance we are thus left to prove 
\begin{equation}\label{eqRemainingBound}
	\sum_{x}\int_0^\infty \frac{\d u}{u^{1+r}}\,[1-\cos(\uv\cdot x)]\, \tau_n(x)
	\le C f_\alpha(n)^{-r}.
\end{equation}
To this end, consider for $z\in\mathbb C, |z|<1$ the generating function
\begin{equation}\label{eqProofXiR1}
	H_{z,r}:=\sum_{x}\sum_{n \ge 0} |x_1|^r \tau_n(x)z^n
	=c_r\sum_{x}\sum_{n \ge 0} \int_0^\infty
        \frac{\d u}{u^{1+r}}\,
        [1-\cos(\uv\cdot x)]\,\tau_n(x)\,z^n,
\end{equation}
where the last identity uses \eqref{eqIntId}.
To prove \eqref{eqRemainingBound} it suffices to show that for any $r\in(0,(\am))$ there exists a constant $C>0$ such that for any $z\in\mathbb C$ with $|z|<1$, 
\begin{equation}\label{eqHzr-bd}
	|H_{z,r}|\le C|1-z|^{-1-r/(\am)}
\end{equation}
for $\alpha\neq2$ and $|H_{z,r}|\le C|1-z|^{-1-r/2}\log|1-z|^{-1/2}$ for $\alpha=2$, because if these bounds hold we can apply the Tauberian theorem \cite[Lemma 3.2]{DerSla98} (cf.\ the discussion above \eqref{e:priordiscussion}) to obtain $\sum_{x}|x_1|^r\tau_n(x) \le C n^{r/(\am)}$ when $r < (\am)$ and to obtain $\sum_{x}|x_1|^r\tau_n(x) \le C n \log n$ when $r=(\am)$.

We will now prove \eqref{eqHzr-bd}. First consider the case $\alpha\neq2$.
We split the integral as
\begin{equation}
\label{eqProofXiR2}
	|H_{z,r}| \le c_r\int_0^{|1-z|^{1/(\am)}}
        \hspace{-.5em}
        \Big|\hat \Tau_z(0)-\hat \Tau_z(\uv)\Big| \frac{\d u}{u^{1+r}}\;
        + c_r\int_{|1-z|^{1/(\am)}}^\infty
        \left(|\hat\Tau_z(0)|+|\hat\Tau_z(\uv)|\right)
        \frac{\d u}{u^{1+r}}.
\end{equation}

Applying Proposition \ref{prop:positivity}(ii) to the second integral in \eqref{eqProofXiR2} gives
\begin{equation}\label{eqProofXiR1a}
\begin{split}
	\int_{|1-z|^{1/(\am)}}^\infty
        \left(|\hat\Tau_z(0)|+|\hat\Tau_z(\uv)|\right)
        \frac{\d u}{u^{1+r}}\;
	& \le
    	C|1-z|^{-1}\int_{|1-z|^{1/(\am)}}^\infty
        \frac{\d u}{u^{1+r}} \\
        & =
        \frac{C(\am)}r \,|1-z|^{-1-r/(\am)},
\end{split}
\end{equation}
as required.
For the first integral in \eqref{eqProofXiR2} we first rewrite the integrand using \eqref{eqProof2a}, and then we use Proposition \ref{prop:positivity}(ii) to bound
\begin{equation}\label{eqProofXiR3}
\begin{split}
	\left|\hat \Tau_z(0)-\hat \Tau_z(\uv)\right|
	&= \left|\hat\Tau_z(0)\,\hat\Tau_z(\uv)
	\left(\hat\Tau_z(\uv)^{-1}-\hat\Tau_z(0)^{-1}\right)\right| \\
	& \le \frac{C}{|1-z|^2} \left|\frac{\hat\Pi_z(0)-\hat\Pi_z(\uv)}{\hat\Pi_z(0)\;\hat\Pi_z(\uv)}
	+zp_c[1-\hat D(\uv)]\right|.
\end{split}
\end{equation}
Observe that it follows from Proposition~\ref{prop-LEcoefficients2}(v) that
\begin{equation}
    	\left|\hat\Pi_z(0)-\hat\Pi_z(\uv)\right|
    	\le \sum_{x}\sum_{m \ge 0} [1-\cos(\uv\cdot x)]\,|z|^n \bar \pi_m(x)
    	\le C\, [1-\hat D (\uv)]
\end{equation}
for some constant $C>0$ that does not depend on $u$ and $z$. We apply this bound, Assumption~\ref{ass:E}, and Proposition~\ref{prop:positivity}(i.b) to the right-hand side of \eqref{eqProofXiR3} to obtain $C|1-z|^{-2}|u|^{(\am)}$ as an upper bound.
Hence, using $r<(\am)$, the first integral in \eqref{eqProofXiR2} is bounded above by
\begin{equation}
\label{eqProofXiR1b}
\begin{split}
	\int_0^{|1-z|^{1/(\am)}}
        \Big|\hat \Tau_z(0)-\hat \Tau_z(\uv)\Big|
	\frac{\d u}{u^{1+r}}
	&\le C|1-z|^{-2} \int_0^{|1-z|^{1/(\am)}}
	\frac{u^{(\am)}}{u^{1+r}}\,\d u \\
	&\le C\,\frac{|1-z|^{-1-r/(\am)}}{(\am)-r}.
\end{split}
\end{equation}
We combine \eqref{eqProofXiR2}, \eqref{eqProofXiR1a} and \eqref{eqProofXiR1b}
to get the desired bound \eqref{eqHzr-bd}. This finishes the argument for the case $\alpha\neq2$.

To prove Theorem \ref{thm-xiR} for the case $\alpha=2$, we have to take the logarithmic corrections into account.
The way to do this has been demonstrated in \cite[Theorem 1.4]{Heyd11}, and we omit the proof here.
\qed


\section{Convergence as a stochastic process: Proof of Theorem \ref{thm-backbone-scal-limit}}\label{sec-findim}

This section is devoted to the proofs of Propositions \ref{prop-FinDimConv} and \ref{cor-tightness}, thus completing the proof of Theorem \ref{thm-backbone-scal-limit}.


\subsection{Finite-dimensional distributions}
We start with the proof of Proposition \ref{prop-FinDimConv}. Instead of proving it directly we consider the following generalized version:
\begin{prop}[Finite-dimensional distributions, generalized version]\label{prop-FinDimConv2}
Let $N$ be a positive integer, $k^{\sss (1)},\dots,$ $k^{\sss (N)}\in\R^d$, $0=t^{\sss (0)}<t^{\sss (1)}<\dots<t^{\sss (N)}<1$, and 
\begin{equation}
	g\in G := \big\{(g_n)\in\R^\N\mid 0\le g_n\le \log n/n \text{ for all }n\in\N\big\}.
\end{equation}
We write
\begin{eqnarray}
    	\bk_n &:=& \big(k_n^{\sss (1)},\dots,k_n^{\sss (N)}\big) 
	=f_\alpha(n)\,\big(k^{\sss (1)},\dots,k^{\sss (N)}\big),\\
    	n\bT &:=& \big(\lfloor nt^{\sss (1)}\rfloor,\dots,\lfloor nt^{\sss (N-1)}\rfloor,\lfloor nT\rfloor\big),
\end{eqnarray}
with $T:=t^{\sss (N)}(1-g_n)$.
Under the conditions of Theorem \ref{thm-endpoint},
\begin{align}\label{eqFinDimConv}
    	\lim_{n\to\infty}\rhonT{N}_{n\bT}(\bk_n)
    	=\frac1A\lim_{n\to\infty}\taunT{N}_{n\bT}(\bk_n)
    	=\exp\Big(-K_\alpha\,\sum_{j=1}^N|k^{\sss (j)}|^{(\am)}\;(t^{\sss (j)}-t^{\sss (j-1)})\Big)
\end{align}
holds uniformly for $g\in G$, where $A$ is the same constant as in Theorem \ref{thm-endpoint} and $\rhonT{N}_{n\bT}$ and $\taunT{N}_{n\bT}$ are as defined in \eqref{eqDefCnN}.
\end{prop}

The proof is carried out by induction on $N$. We use the sequence $g$ to ensure that we can advance the induction, because $g$ gives us a little flexibility in our choice of the location of the end-point.

\proof[Proof of Proposition \ref{prop-FinDimConv2}]
The proof that we present here follows a well-known approach, see for example the proofs of \cite[Theorem 6.6.2]{MadSla93} or \cite[Theorem 1.6]{Heyd11}. We will first give the proof for $\taunT{N}_{n\bT}$, and then discuss the necessary changes for $\rhonT{N}_{n\bT}$. For convenience, we write $nt^{\sss (j)}$ and $ nT$ instead of $\lfloor nt^{\sss (j)}\rfloor$ and $\lfloor nT\rfloor$.

We start the induction by applying Theorem \ref{thm-endpoint}: since $\taunT{1}_{n\bT}(\bk_n)=\hat \tau_{nT}(k^{\sss(1)}_n)$, we can replace $n$ by $nT$ in \eqref{eqProof14} and the claim follows.

The identity in \eqref{taun-exp-pivs-indef} along with \eqref{tau-n-def-rep} provides a decomposition of $\tau_n(x)$ into $M+1$ $\pi$-coefficients separated by bonds $b_1,\dots,b_M$. 
We now partially resum \eqref{taun-exp-pivs-indef} and sum over the intermediate pivotal bonds (using \eqref{tau-n-def-rep} again) so that we only keep the $\pi$-coefficient containing the $(nt^{\sss (N-1)})$-th pivotal bond (which stretches from the $I_1$-th to the $I_2$-th pivotal bond) and obtain 
\begin{align}
	\label{eqFinDimProof2}
	\tau_n(x)
	=&\sum_{\substack{I=[I_1,I_2]\\0< I_1\le nt^{\sss (N-1)} \le I_2< nT}}
	  \sum_{b_{I_1},b_{I_2+1}}
	   J (b_{I_1}) J (b_{I_2+1}) 
	 \tau_{I_1-1}(\bb_{I_1})\;\pi_{|I|}(\bb_{I_2+1}-\tb_{I_1})\;\tau_{n-I_2-1}(x-\tb_{I_2+1})\nnb
	 &{}\qquad+\sum_{\substack{I=[0,I_2]\\nt^{\sss (N-1)} \le I_2< nT}}
	  \sum_{b_{I_2+1}}
	  J (b_{I_2+1})
	 \pi_{|I|}(\bb_{I_2+1})\;\tau_{n-I_2-1}(x-\tb_{I_2+1}) \nnb
	 &{}\qquad +\sum_{\substack{I=[I_1,n]\\0< I_1\le nt^{\sss (N-1)}}}
	 \sum_{b_{I_1}}
	  J (b_{I_1})
	 \tau_{I_1-1}(\bb_{I_1})\;\pi_{|I|}(x-\tb_{I_1}),
\end{align}
where the first line contains the main contribution, and the second and third line deal with the cases where the $(nt^{\sss (N-1)})$-th pivotal bond lies in the very first or the very last $\pi$-coefficient (which turn out to be negligible contributions). 

We can similarly rewrite the characteristic function for increments, $\hat{\underline{\rm \tau}}^{\sss (N)}_{n\bT}(\bk_n)$, as
\begin{align}
	\label{eqFinDimProof2a}
    	\hat{\underline{\rm \tau}}^{\sss (N)}_{n\bT}(\bk_n)
    	=&\sum_{\substack{I=[I_1,I_2]\\0< I_1\le nt^{\sss (N-1)} \le I_2< nT}}
	\hat J(k_n^{\sss (N-1)}) \hat J(k_n^{\sss (N)}) \; \taunT{N-1}_{(nt^{\sss(1)},\dots,nt^{\sss(N-2)},I_1-1)}
        \big(k_n^{\sss (1)},\dots,k_n^{\sss (N-1)}\big)\nnb
        &\qquad \qquad\qquad\times 
        \hat{\underline\pi}_{(nt^{\sss(N-1)}-I_1,\,|I|)}(k_n^{\sss (N-1)},k_n^{\sss (N)})
        \,\hat \tau_{nT-I_2-1}(k_n^{\sss (N)}) + \hat{\mathcal E}^{\sss (N)}_{n\bT}(\bk_n) 
\end{align}
where, like \eqref{pimdef} and recalling \eqref{pimdef-pivs}--\eqref{pimNdef-pivs}, we define, for $0\le m\le n$ and $k_1, k_2\in\Rd$,
\begin{equation}\label{eqPsiExp2}
    	\hat{\underline\pi}_{(m,n)}(k_1,k_2)
    	:=  \sum_{b_1,\ldots, b_n}\sum_{\bb_{n+1}} 
    	\exp\big\{i(k_1\cdot \bb_m+k_2\cdot(\bb_n-\bb_m))\big\}\;
    	\pi_n(0,\vec{b}_{[n]},\bb_{n+1})
\end{equation}
and $\hat{\mathcal E}^{\sss (N)}_{n\bT}(\bk_n)$ is the contribution from those summands where $[I_1,I_2]\nin [nt^{\sss (N-2)}+1,nT-1]$, which includes the terms due to the second and third line in \eqref{eqFinDimProof2}. 
Observe that $\hat{\underline\pi}_{(m,n)}(0,0) = \hat \pi_n (0)$.

Using \eqref{eqFinDimProof2a} we can split $\hat{\underline{\rm \tau}}^{\sss (N)}_{n\bT}(\bk_n)$ into the contribution of short intervals $I$ and long intervals. We write ${{\underline{\rm \sigma}}}^{\sss (N)}_{n\bT}(\bk_n)$ for the contribution to $\hat{\underline{\rm \tau}}^{\sss (N)}_{n\bT}(\bk_n)$ that comes from short intervals, with length $|I|=I_2-I_1 \le \log n$, and we write ${{\underline{\rm \lambda}}}^{\sss (N)}_{n\bT}(\bk_n)$ for the contributions from the long intervals, with length $|I|>\log n$.

We start by showing that ${\underline{\rm \lambda}}^{\sss (N)}_{n\bT}(\bk_n)$ is negligible.
It follows from Theorem \ref{thm-endpoint} that $\hat\tau_m(k)\le\hat\tau_m(0)\le C$ for some $C>0$. Applying this bound (and using \eqref{eqFinDimProof2a}),
\begin{equation}\label{eqFinDimProof9}
	{\underline{\rm \lambda}}^{\sss (N)}_{n\bT}(\bk_n)
	\le {\underline{\rm \lambda}}^{\sss (N)}_{n\bT}(\underline{0})
    	= \sum_{\substack{I\ni nt^{\sss (N-1)}\\|I|> \log n}}
	p_c^2\hat\tau_{I_1-1}(0)\,|\hat\pi_{|I|}(0)|\,\hat\tau_{nT-I_2-1}(0)
	\le C^2\sum_{m \ge \log n+1} (m+1)\,|\hat\pi_{m}(0)|.
\end{equation}
We get the factor $m+1$ here because there are precisely $m+1$ ways to choose the interval $I\ni nt^{\sss (N-1)}$ under the restriction $|I|=m$. Since the right-hand side of \eqref{eqFinDimProof9} is finite for all $n\ge1$ by Proposition~\ref{prop-LEcoefficients}, it vanishes as $n\to\infty$. By the same argument, $\hat{\mathcal E}^{\sss (N)}_{n\bT}(\bk_n)$ vanishes, because $n(t^{\sss (N-1)}-t^{\sss (N-2)}) \gg \log n$ and $n(T-t^{\sss (N-1)})\gg\log n$ when $n$ is large.

We now establish the asymptotics for ${\underline{\rm \sigma}}^{\sss (N)}_{n\bT}(\bk_n)$. Assume that $n$ is large enough that both $(nt^{\sss (N-1)}-nt^{\sss (N-2)}) \ge \log n$ and $(nt^{\sss (N)}-nt^{\sss (N-1)})\ge \log n$. The induction hypothesis is that
\begin{equation}\label{eqFinDimProof4}
    	\taunT{N-1}_{(nt^{\sss(1)},\dots,nt^{\sss(N-2)},I_1-1)}
	\big(k_n^{\sss (1)},\dots,k_n^{\sss (N-1)}\big)\\
    	=A\, \exp\Big(-K_\alpha\,\sum_{j=1}^{N-1}|k^{\sss (j)}|^{(\am)}\;(t^{\sss (j)}-t^{\sss (j-1)})\Big)
	+E_1(I),
\end{equation}
where $E_1(I)$ is an error term that converges to $0$ as $n\to\infty$ uniformly in $|I|\le \log n$.

A slight generalization of the case $N=1$ shows that
\begin{equation}\label{eqFinDimProof5}
    	\hat \tau_{nT-I_2-1}(k_n^{\sss (N)})
    	=A\, \exp\Big(-K_\alpha\,|k^{\sss (N)}|^{(\am)}\;(t^{\sss (N)}-t^{\sss (N-1)})\Big)+E_2(I),
\end{equation}
where $E_2(I)$ is an error term that is due to Theorem \ref{thm-endpoint}. Note that $E_2(I)$ converges to $0$ as $n\to\infty$ uniformly in $|I|\le \log n$. Hence,
\begin{multline}
	\label{eqFinDimProof6}
    	{\underline{\rm \sigma}}^{\sss (N)}_{n\bT}(\bk_n)
    	=\left(A^2\,\exp\Big(-K_\alpha\,\sum_{j=1}^{N}|k^{\sss (j)}|^{(\am)}
	\;(t^{\sss (j)}-t^{\sss (j-1)})\Big)+E_3\right)\\
    	\times \sum_{\substack{I\ni \,nt^{\sss (N-1)}\\|I|\le \log n}}
	\hat J(k_n^{\sss (N-1)}) \hat J(k_n^{\sss (N})
        \hat{\underline\pi}_{(nt^{\sss(N-1)}-I_1,\,|I|)}(k_n^{\sss (N-1)},k_n^{\sss (N)}),
\end{multline}
where $E_3$ is the error term that comes from $E_1$ and $E_2$. Note that $E_3$ is uniform in the sequences $g$ that satisfy $g_n\le \log n/n$.

The proof is complete when we show that the second line in \eqref{eqFinDimProof6} converges to $1/A$. Using that $\hat J = p_c \hat D$, we begin by writing
\begin{multline}\label{eqFinDimProof7}
	\sum_{\substack{I\ni \,nt^{\sss (N-1)}\\|I|\le \log n}}
	\hat J (k_n^{\sss (N-1)}) \hat J(k_n^{\sss (N})
        \hat{\underline\pi}_{(nt^{\sss(N-1)}-I_1,\,|I|)}(k_n^{\sss (N-1)},k_n^{\sss (N)})\\
	=\sum_{\substack{I\ni \,nt^{\sss (N-1)}\\|I|\le \log n}}p_c^2\hat\pi_{|I|}(0)
		-\sum_{\substack{I\ni \,nt^{\sss (N-1)}\\|I|\le \log n}}
		p_c^2\left(\hat\pi_{|I|}(0)-
		\hat{\underline\pi}_{(nt^{\sss(N-1)}-I_1,\,|I|)}(k_n^{\sss (N-1)},k_n^{\sss (N)})
		\right)\\
	   -\sum_{\substack{I\ni \,nt^{\sss (N-1)}\\|I|\le \log n}}p_c^2
		\left[1-\hat D(k_n^{\sss (N-1)})\right] \left[1- \hat D(k_n^{\sss (N})\right]
		\hat{\underline\pi}_{(nt^{\sss(N-1)}-I_1,\,|I|)}(k_n^{\sss (N-1)},k_n^{\sss (N)}).
\end{multline}
It follows from \eqref{eqDefA'} and the fact that by Proposition \ref{prop-LEcoefficients}, $\sum_{m=n}^\infty (m+1)|\hat \pi_m (0)| \to 0$ as $n \to \infty$, that the first term converges to $1/A$. Indeed,
\begin{equation}\label{eqFinDimProof10}
	\sum_{\substack{I\ni \,nt^{\sss (N-1)}\\|I|\le \log n}}p_c^2\hat\pi_{|I|}(0)
	=\sum_{I\ni \,nt^{\sss (N-1)}}p_c^2\hat\pi_{|I|}(0)
	-\sum_{\substack{I\ni \,nt^{\sss (N-1)}\\|I|> \log n}}p_c^2\hat\pi_{|I|}(0)
	\quad\stackrel{n\to\infty}\longrightarrow\quad
	p_c^2\sum_{m\ge0}(m+1)\hat\pi_{m}(0)
	=1/A,
\end{equation}
where the factor $m+1$ arises because there are $m+1$ intervals of length $|I|=m$ that contain the point $nt^{\sss (N-1)}$.

Now we show that the second term on the right-hand side of \eqref{eqFinDimProof7} converges to $0$. Recall the definition of $\pi_{m,n}(y,x)$ in \eqref{pnm-def}. We use the spatial symmetry of the model to replace the exponential factor in \eqref{eqPsiExp2} by a cosine. We also use $\left|1-\cos(a)\cos(b)\right|\le 2|a|^\delta + 2|b|^{\delta}$ for all $\delta \in [0,2]$ to obtain the bound
\begin{equation}\label{eqFinDimProof13}
    	\left|\hat\pi_{|I|}(0)- \hat{\underline\pi}_{(nt^{(N-1)}-I_1,|I|)}(k_1,k_2) \right| 
	\le 2\sum_{y,x} (|k_1 \cdot y|^\delta + |(k_2-k_1) \cdot x|^\delta)  \lvert \pi_{nt^{(N-1)}-I_1,|I|}(y,x) \rvert.
\end{equation}

Again, there are $m+1$ of intervals of length $|I|=m$, hence, using \eqref{eqFinDimProof13}, uniformly in $k^{\sss (N-1)},k^{\sss (N)} \in \R^d$,
\begin{multline}
	\label{eqFinDimProof8}
	\sum_{\substack{I\ni \,nt^{\sss (N-1)}\\|I|\le \log n}}
		p_c^2\left|\hat{\underline\pi}_{(nt^{\sss(N-1)}-I_1,\,|I|)}(k_n^{\sss (N-1)},k_n^{\sss (N)})
		-\hat\pi_{|I|}(0)\right|\\
	\le C \sum_{x,y} \sum_{m=0}^{\log n} \sum_{m_1=0}^{m}  (|f_\alpha(n)\,x|^{\delta} 
	+ |f_\alpha(n)\,y|^{\delta})\lvert \pi_{m_1,m}(y,x)\rvert \\
	\le C f_\alpha(n)^\delta\,(\log n+1) \max_{m \le \log n} \sum_{x,y} \sum_{m_1=0}^{m} (|x|^{\delta}
	+ |y|^\delta) \lvert \pi_{m_1,m}(y,x) \rvert,
\end{multline}
and this converges to 0 as $n\to\infty$ when $\delta$ is sufficiently small since the sum is bounded uniformly in $n$ by Proposition \ref{prop-LEcoefficients2}(iv). It now follows from $\lim_{n\to\infty}n\,[1-\hat D(k_n)]=|k|^{(\am)}$ (cf.\ \eqref{[1-D]-n-asymp}), and \eqref{eqFinDimProof8} that the second line on the right-hand side of \eqref{eqFinDimProof7} vanishes as $n\to\infty$. This completes the proof that the second line in \eqref{eqFinDimProof6} converges to $1/A$, and thus we have completed the advancement of the induction. This completes the proof of the asymptotics of $\hat{\underline{\rm \tau}}^{\sss (N)}_{n\bT}(\bk_n)$.
\medskip

With the result for $\taunT{N}_{n\bT}$ in hand, we can derive the statement for $\rhonT{N}_{n\bT}$. Now we start with the identity \eqref{rhon-exp-pivs-indef} rather than \eqref{tau-n-def-rep}, and again do a partial resummation to obtain 
\begin{align}
    	\hat{\underline{\rm \rho}}^{\sss (N)}_{n\bT}(\bk_n)
    	=&\sum_{\substack{I=[I_1,I_2]\\0\le I_1\le nt^{\sss (N-1)} \le I_2\le nT}}
	\hat J (k_n^{\sss (N-1)})\hat J(k_n^{\sss (N)})
        \taunT{N-1}_{(nt^{\sss(1)},\dots,nt^{\sss(N-2)},I_1-1)}
        \big(k_n^{\sss (1)},\dots,k_n^{\sss (N-1)}\big)\nnb
        &\qquad\qquad\times
        \hat{\underline\pi}_{(nt^{\sss(N-1)}-I_1,\,|I|)}(k_n^{\sss (N-1)},k_n^{\sss (N)})
        \,\hat \rho_{nT-I_2-1}(k_n^{\sss (N)})+\hat{\mathcal F}^{\sss (N)}_{n\bT}(\bk_n), 
	\label{eqFinDimProof2b}
\end{align}
where we recall the definition of $\psi_n$ in \eqref{psimdef} and, similar to \eqref{eqPsiExp2}, let
	\begin{equation}\label{eqPsiExp3}
    	\hat{\underline\psi}_{(m,n)}(k_1,k_2)
    	=  \sum_{l\ge n}\sum_{b_1,\ldots, b_l}\sum_{\bb_{l+1}\in\Zd} 
    	\exp\big( i(k_1\cdot \bb_m+k_2\cdot(\bb_l-\bb_m))\big)\;
    	\pi_l(0,\vec{b}_{l},\bb_{l+1})
	\end{equation} 
for $	0\le m\le n$ and $k_1, k_2\in\Rd$ and let $\hat{\mathcal F}^{\sss (N)}_{n\bT}(\bk_n)$ be the contribution due to the ``long'' $\pi$-coefficients.
However, these $\pi$'s are ``too long'' to contribute to the final sum, and by the calculation in \eqref{eqFinDimProof9} they vanish in the limit as $n\to\infty$, see also \eqref{pinm-ass}. 

We are therefore left to analyse the first term on the right-hand side of  \eqref{eqFinDimProof2b}. 
The term involving $\hat{\underline{\tau}}$ gives rise to a factor
\begin{equation}
    	A\, \exp\big(-K_\alpha\,\sum_{j=1}^{N-1}|k^{\sss (j)}|^{(\am)}
	\;(t^{\sss (j)}-t^{\sss (j-1)})\big) + \tilde E,
\end{equation}
where $\tilde E$ is an error term, similar to \eqref{eqFinDimProof4}.
Likewise, the term involving $\hat{\underline\pi}$ gives rise to the factor $1/A$, by \eqref{eqFinDimProof10}.
Finally, by Theorem \ref{thm-endpoint},
\begin{equation}
	\hat \rho_{nT-I_2-1}(k_n^{\sss (N)})
	\to \exp\Big(-K_\alpha\,|k^{\sss (N)}|^{(\am)}\;(t^{\sss (N)}-t^{\sss (N-1)})\Big)
	\qquad\text{as $n\to\infty$,}
\end{equation}
so that the statement for $\rhonT{N}_{n\bT}$ follows. \qed


\subsection{Tightness}
\label{sec-tightness}
In this section we prove tightness of $X_n$ and $Y_n$. We start with tightness of $Y_n$.

Tightness of $Y_n$ follows from the bound
\begin{equation}
    	\label{tightness-condition}
    	\E^\ast_{p_c,n}\!\left[ |Y_n(t_2)-Y_n(t_1)|^r
    	|Y_n(t_3)-Y_n(t_2)|^r\right]
    	\leq C|t_3-t_1|^a,
\end{equation}
for some $r>0$, $a>1$, $C>0$, as shown in \cite[Theorem 13.5]{Bill99}(where (13.13) in that reference is replaced by the stronger moment condition (13.14)).

To establish \eqref{tightness-condition} we again rely on the lace expansion as in \eqref{taun-exp-pivs-indef}. Denote $n_i=\lceil nt_i\rceil$. 
Then, 
\begin{equation}
    \label{tightness-condition-rep}
    \E^\ast_{p_c,n}\!\left[ |Y_n(t_2)-Y_n(t_1)|^r
    |Y_n(t_3)-Y_n(t_2)|^r\right]
    =n^{-2r/(\am)} \E^\ast_{p_c,n}\!\left[|S_{n_2}-S_{n_1}|^r
    |S_{n_3}-S_{n_2}|^r\right].
\end{equation}
Write $n_4=n$ and define
\begin{equation}
	\tau_{n_1,n_2,n_3,n_4}(x_1,x_2,x_3,x_4)
	:=\prob_{p_c}(S_{n_1}=x_1,S_{n_2}=x_2,S_{n_3}=x_3, S_{n_4}=x_4),
\end{equation}
so that
\begin{equation}
	\E^\ast_{p_c,n}\!\left[|S_{n_2}-S_{n_1}|^r
    	|S_{n_3}-S_{n_2}|^r\right]=\frac{1}{\hat{\tau}_n(0)}
	\sum_{x_1,x_2,x_3,x_4} |x_2-x_1|^r|x_3-x_2|^r\tau_{n_1,n_2,n_3,n_4}(x_1,x_2,x_3,x_4).
\end{equation}
Since we already know by Corollary~\ref{corr:pivballboundary} that $\hat{\tau}_n(0)\rightarrow A>0$, it suffices to prove that there exist $C<\infty$ and $a>1$ such that
\begin{equation}
	\label{tightness-condition-aim}
	\sum_{x_1,x_2,x_3,x_4} |x_2-x_1|^r|x_3-x_2|^r\tau_{n_1,n_2,n_3,n_4}(x_1,x_2,x_3,x_4)
	\leq C n^{2r/(\am)} |t_3-t_1|^a.
\end{equation}
We will take $r = \tfrac12 \twa +\vep$ for some small $\vep>0$ and $a = 2r/(\am) > 1$, so that \eqref{tightness-condition-aim} follows when there exist $C<\infty$ such that
\begin{equation}
	\label{tightness-condition-aim-2}
	\sum_{x_1,x_2,x_3,x_4} |x_2-x_1|^r|x_3-x_2|^r\tau_{n_1,n_2,n_3,n_4}(x_1,x_2,x_3,x_4)
	\leq C  |n_3-n_1|^a.
\end{equation}
	
\begin{figure}
	\includegraphics[width = .9\textwidth]{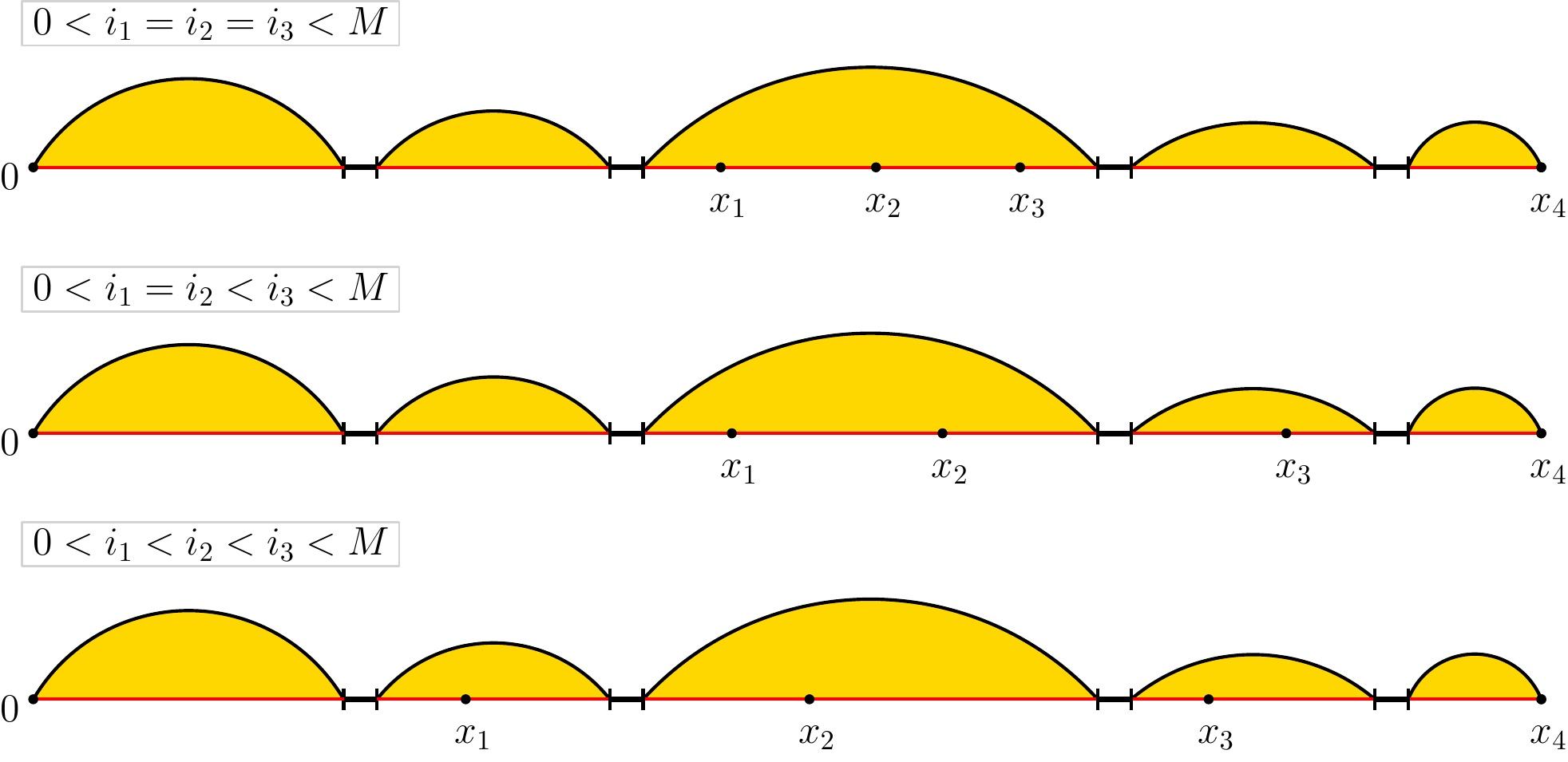}
	\caption{A schematic drawing of three of the cases that occur for the expansion of $\tau_{n_1, n_2, n_3, n_4}(x_1,x_2,x_3,x_4)$. The red line represents the backbone pivotals $\bar{b}_{[n]}$, the black dashes represent the pivotals $b_{I_1},\dots b_{I_4}$. The yellow filled arcs represent the factors $\pi$. Compare with Figure~\ref{fig:pifixed}.\label{fig-piv-intervals-three-arrows}}
\end{figure}

We expand out $\tau_{n_1,n_2,n_3,n_4}(x_1,x_2,x_3,x_4)$ using \eqref{taun-exp-pivs-indef}, and recall that the interval $[0,n]$ is partitioned into the $M+1$ disjoint intervals $[I_{i-1}-1, I_{i}]$. 
For $j\in\{1,2,3\}$, let $i_j$ be the index such that $n_j\in [I_{i_j-1}-1, I_{i_j}]$. This gives rise to different contributions, depending on the cardinality of the set $\{i_1,i_2,i_3\}$, see Figure \ref{fig-piv-intervals-three-arrows}. We will bound these contributions to \eqref{tightness-condition-aim-2} separately. 

Recall the definition of $\pi_{m_1,m_2}(y_1,y_2)$ from \eqref{pnm-def}. Similarly, we define $\pi_{m_1,m_2,m_3,m_4}(y_1,y_2,y_3,y_4)$, starting with \eqref{pimdef} but leaving out the summation over $\olb_{m_1}, \olb_{m_2}, \olb_{m_3},\olb_{m_4}$ and writing $y_1,y_2,y_3,y_4$ for the free variables that denote the positions of $\olb_{m_1}, \olb_{m_2}, \olb_{m_3},\olb_{m_4}$, respectively. The coefficient $\pi_{m_1,m_2,m_3}(y_1,y_2,y_3)$ is defined similarly. At places, where it simplifies notation, we will write these coefficients with an additional argument for the position of the root, i.e., we write $\pi_{m_1,m_2}(w,x,y) := \pi_{m_1,m_2}(x-w,y-w)$.  The following lemma provides the necessary bounds on the lace-expansion coefficients for tightness:

\begin{lemma}[Bounds on lace-expansion coefficients with extra weights and sums]
\label{lem-pi-worst-assump}
Under the assumptions of Theorem \ref{thm-backbone-scal-limit}, for $r = \tfrac12 \twa +\vep$ with $\vep>0$ sufficiently small,
\begin{align}
	\label{pi-worst-assump-rep-1}
	C_{\pi}^{\sss(1)}&:= \sum_{0 \le m_1 < m_2} 
	\sum_{y_1,y_2}(|y_1|^r+|y_2-y_1|^r+1)|\pi_{m_1,m_2}(y_1,y_2)|<\infty,\\
	\label{pi-worst-assump-rep-2}
	C_{\pi}^{\sss(2)}&:= \sum_{0 \le m_1 < m_2 < m_3} 
	\sum_{y_1,y_2,y_3}|y_2-y_1|^r(|y_3-y_2|^r+1)|\pi_{m_1,m_2,m_3}(y_1,y_2,y_3)|
	<\infty,\\	
	\label{pi-worst-assump}
	C_{\pi}^{\sss(3)}&:= \sum_{0 \le m_1 < m_2 < m_3 < m_4} 
	\sum_{y_1,y_2,y_3,y_4} |y_2-y_1|^r  |y_3-y_2|^r |\pi_{m_1,m_2,m_3,m_4}(y_1,y_2,y_3,y_4)|<\infty.
\end{align}
\end{lemma}

\longversion{We prove this lemma in Appendix~\ref{sec:multiweight}.}\shortversion{This lemma is proved in \cite[Appendix~B.3]{HeyHofHulMie17b}.} Let us remark that we need $d > 6\twa$ here and only here. The other proofs in this paper only require that the strong triangle condition is satisfied with some sufficiently small $\beta$. We need a higher dimension here because the bounds above are suboptimal in the sense that they require heavier weights than we believe to be necessary, and those weights make the sums divergent unless we have a sufficiently high dimension. We explain this more precisely in Remark~\ref{rem-dimension-tightness} below.
\medskip

We define
	\eqn{
	\label{tilde-tau-def}
	\tilde\tau_{m}(x,y):=(\tau_{m-1} * J)(x,y), \qquad \text{ and } \qquad \tilde{\tilde\tau}_{m}(x,y):=(J * \tau_{m-2} * J)(x,y).
	}
We further use the bounds from Theorems \ref{thm-endpoint} and \ref{thm-xiR} that uniformly in $m\geq 1, x \in \Zd$ there exists a constant $K$ such that
\begin{equation}
	\label{tau-ass}
	\sum_{y} \tilde\tau_m(x,y)\leq K,
	\qquad \text{ and } \qquad 
	\sum_{y} |y-x|^r \tilde\tau_m(x,y)\leq Km^{r/(\am)},
\end{equation}
and the same bounds hold for $\tilde{\tilde{\tau}}_m$. We now bound all different contributions. The bounds are quite similar to each other, and therefore some explanations will be relatively brief. For the time being, we ignore the boundary cases where $I_{i_1-1}=0$ (which means that the interval that contains $n_1$ also contains 0), and $I_{i_3}=M$ (which means that the interval that contains $n_3$ also contains $n=n_4$). We will deal with these cases (which are very similar) at the end. \bigskip

\begin{figure}
	\includegraphics[width=\textwidth]{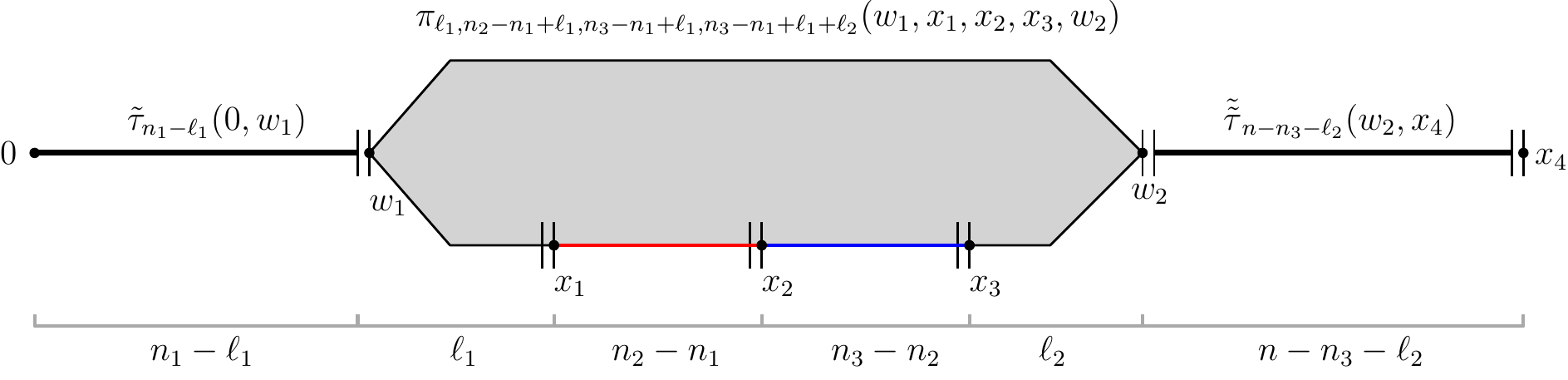}
	\caption{A schematic drawing of \eqref{e:pi3in1}. The $\pi$-factor is indicated by the gray polygon, the two two-point functions by thick lines. The factor $|x_2-x_1|^r$ is indicated in red, the factor $|x_3-x_2|^r$ in blue. \label{fig:pi3in1}}
\end{figure}

\paragraph{{\bf The contribution due to $0< i_1=i_2=i_3 <M$.}} The contribution to the left-hand side of \eqref{tightness-condition-aim-2} due to $i_1=i_2=i_3$ is equal to
\begin{multline}\label{e:pi3in1}
	\sum_{x_1,x_2,x_3,x_4} \sum_{\ell_1 =1}^{n_1-1} \sum_{\ell_2 =n_3 +1}^{n-n_3 - 1} \sum_{w_1, w_2} |x_2 -x_1|^{r} |x_3-x_2|^{r} \ttau_{n_1 - \ell_1}(0,w_1)\\
	\times \pi_{\ell_1, n_2-n_1+\ell_1, n_3-n_1+\ell_1, n_3-n_1+\ell_1+\ell_2}(w_1,x_1, x_2, x_3, w_2) \tilde{\tilde{\tau}}_{n-n_3-\ell_2},(w_2,x_4),
\end{multline}
see Figure \ref{fig:pi3in1}. Note that we have used \eqref{taun-exp-pivs-indef} to isolate the $\pi$-factor that contains $x_1,x_2,$ and $x_3$, and then resummed the remaining terms, again using \eqref{taun-exp-pivs-indef}, to obtain the two $\tau$-factors.
Since the summation indices $\ell_1$ and $\ell_2$ both fall within the $\pi$-factor, we may take the supremum over the indices of the two $\tau$-factors and use the first bound in \eqref{tau-ass} to bound the factors $\ttau_{n_1-\ell_1}$ and $\tilde{\tilde{\tau}}_{n-n_3+\ell_2}$ first (for the former we use translation invariance and symmetry). Then we may extend the sum of the indices $(\ell_1, n_2-n_1+\ell_1, n_3-n_1+\ell_1, n_3-n_1+\ell_1+\ell_2)$ to all possible values the indices may take, so that we may apply \eqref{pi-worst-assump} (and relabel the summation indices and arguments) to bound this from above by
\begin{equation}\label{e:3in1pibd}
	K^2 \sum_{y_1,y_2, y_3, y_4} \sum_{0\le m_1 < m_2 < m_3 < m_4} |y_2-y_1|^r |y_3-y_2|^r |\pi_{m_1,m_2,m_3,m_4} (y_1,y_2,y_3,y_4)| \le K^2 C_\pi^{\sss (3)}.
\end{equation}	
Since $|n_3-n_1|\geq 1$, this satisfies the required bound in \eqref{tightness-condition-aim-2}.
\medskip

\begin{figure}[tb]
	\includegraphics[width=.98\textwidth]{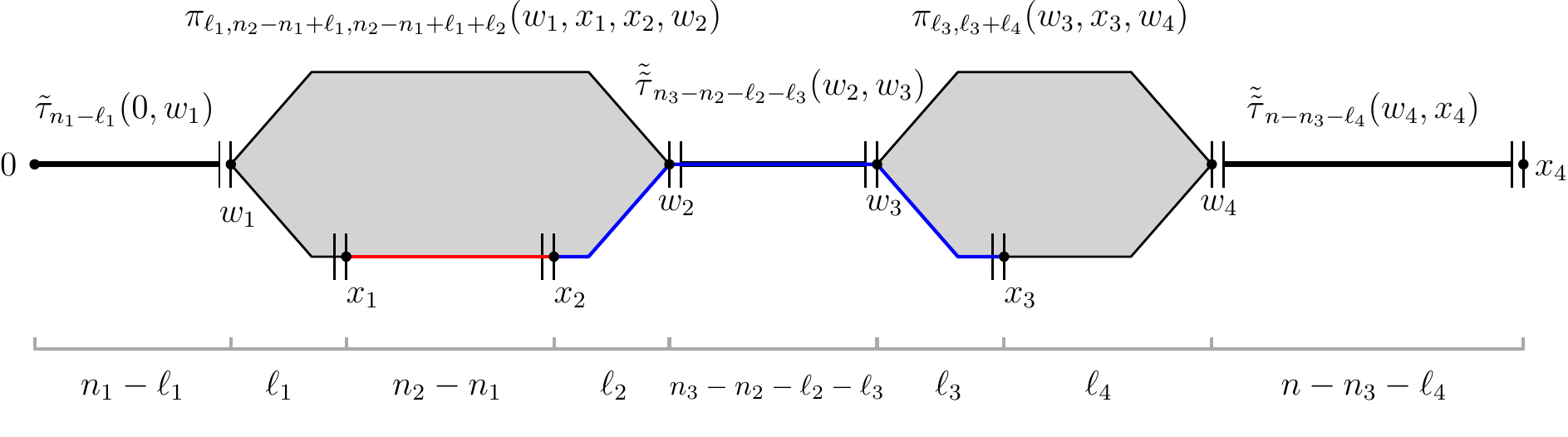}
	\caption{A schematic drawing of \eqref{e:2in1pi}. The factor $|x_2-x_1|^r$ is indicated in red, the factor $|x_3-x_2|^r$ in blue. \label{fig:2in1pi}}
\end{figure}

\paragraph{{\bf The contributions due to $0<i_1=i_2<i_3 <M$ and $0<i_1< i_2=i_3 <M$.}} The contribution to the left-hand side of \eqref{tightness-condition-aim-2} due to $i_1=i_2< i_3$ is equal to
\begin{multline}\label{e:2in1pi}
	\sum_{x_1,x_2,x_3,x_4} \sum_{\ell_1=1}^{n_1-1} \sum_{\ell_2 =n_2-n_1+\ell_1+1}^{n_3-2} \sum_{\ell_3=1}^{n_3-n_2-\ell_2-1} \sum_{\ell_4 = n_3-\ell_3+1}^{n-n_3-\ell_3 -1} \sum_{w_1,w_2, w_3, w_4} |x_2-x_1|^r |x_3-x_2|^r \ttau_{\ell_1}(0,w_1) \\
	\times \pi_{\ell_1,n_2-n_1+\ell_1,n_2-n_1+\ell_1+\ell_2}(w_1, x_1, x_2, w_2) \tilde{\tilde{\tau}}_{n_3-n_2-\ell_2-\ell_3}(w_2, w_3) \\
	\times\pi_{\ell_3, \ell_4}(w_3, x_3, w_4) 
	\tilde{\tilde{\tau}}_{n-n_3-\ell_4}(w_4,x_4),
\end{multline}
see Figure~\ref{fig:2in1pi}.

For all $v,w,x,y \in \Zd$ and $r \ge 0$ we have
\begin{equation}\label{e:distribute}
	|x-y|^r \le 3^r (|x-u|^r + |v-u|^r + |y-v|^r).
\end{equation}
We use this to ``distribute'' the factor $|x_3-x_2|^r$ in \eqref{e:2in1pi} over the three components that are involved. From here on, the steps needed to bound the sum are straightforward but lengthy, so we leave it at a brief description of the computations: 

Note that again the summation indices $\ell_1,\dots,\ell_4$ all fall within $\pi$-factors, so that we may again apply the first bound in \eqref{tau-ass} to the first and last $\tau$-factors. Then we can extend the summation over $(n_3 - \ell_3, \ell_4)$ to all possible values these indices may take, similar to \eqref{e:3in1pibd}, to bound the second $\pi$-factor. Then, we take the supremum over $\ell_3$ in the middle $\tau$-factor (after summing over its displacement) for an upper bound. Lastly, we extend the summation over $(n_1 - \ell_1, n_2 -\ell_2, \ell_2)$ to all possible values these indices may take.

Using the bounds \eqref{pi-worst-assump-rep-1} and \eqref{pi-worst-assump-rep-2} for the $\pi$-factors and the bounds \eqref{tau-ass} for the middle $\tilde{\tilde{\tau}}$-factor, we thus obtain the upper bound
\begin{equation}
	 3^{1+r} K^3 C_\pi^{\sss (2)} C_\pi^{\sss (1)} |n_3-n_2|^{r/\twa}.
\end{equation}
Since $|n_3-n_2| \le |n_3-n_1|$, the bound \eqref{tightness-condition-aim-2} is satisfied for this contribution also.

The contribution due to $0< i_1< i_2=i_3 < M$ is analogous.
\medskip

\begin{figure}[tb]
	\includegraphics[width=\textwidth]{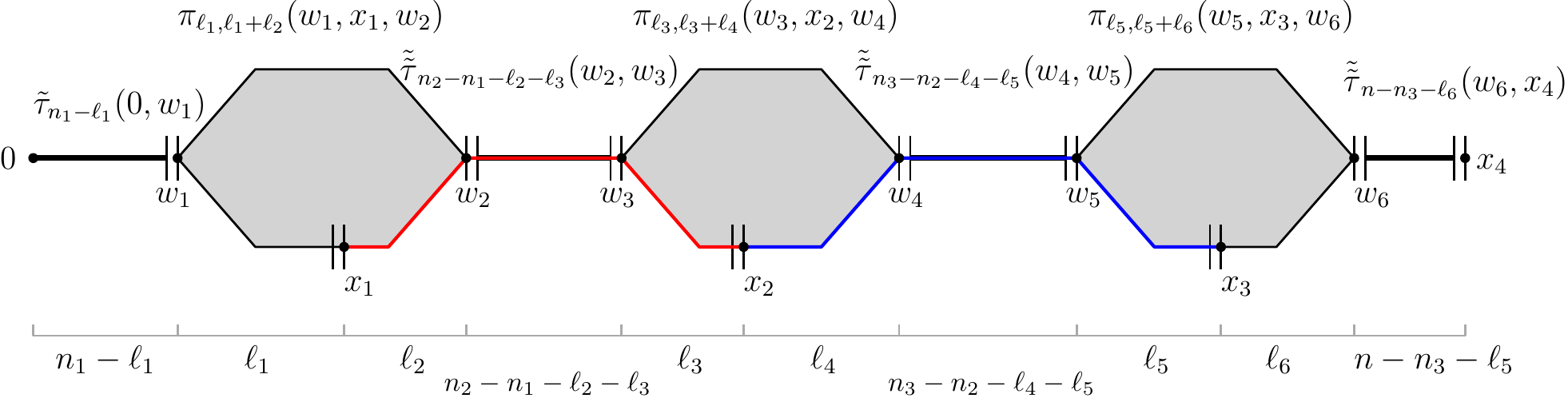}
	\caption{A schematic drawing of \eqref{e:1in1pi}. The factor $|x_2-x_1|^r$ is indicated in red, the factor $|x_3-x_2|^r$ in blue. \label{fig:1in1pi}}
\end{figure}

\paragraph{{\bf The contribution due to $0< i_1< i_2< i_3 < M$.}} We again follow the same steps as in the previous case. The contribution to the left-hand side of \eqref{tightness-condition-aim-2} due to $i_1< i_2< i_3$ is equal to
\begin{multline}\label{e:1in1pi}
	\sum_{x_1,x_2,x_3,x_4} \sum_{\ell_1 =1}^{n_1-1} \sum_{\ell_2 =1}^{n_2-n_1-\ell_1-2} \sum_{\ell_3=1}^{n_2 - n_1-\ell_1-\ell_2 -1} \sum_{\ell_4=1}^{n_3-n_2-2} \sum_{\ell_5 = 1}^{n_3-n_2-\ell_4-1} \sum_{\ell_6 = 1}^{n-n_3-1} \sum_{w_1,\dots,w_6} |x_2-x_1|^r |x_3-x_2|^r \\
	\times \ttau_{n_1-\ell_1}(0,w_1) \pi_{\ell_1, \ell_1+\ell_2}(w_1, x_1, w_2) \tilde{\tilde{\tau}}_{n_2-n_1-\ell_2-\ell_3}(w_2,w_3) \pi_{\ell_3,\ell_3+ \ell_4}(w_3,x_2,w_4) \\
	\times \tilde{\tilde{\tau}}_{n_3-n_2-\ell_4-\ell_5}(w_4,w_5) \pi_{\ell_5, \ell_5+ \ell_6}(w_5,x_3, w_6) \tilde{\tilde{\tau}}_{n-n_3-\ell_6}(w_6,x_4),
\end{multline}
see Figure~\ref{fig:1in1pi}.

As before, we use the first bound in \eqref{tau-ass} to bound the first and the last $\tau$-factors by $K^2$, and we use \eqref{e:distribute} to distribute the weights $|x_2-x_1|^r$ and $|x_3-x_2|^r$ over the coefficients. We then bound the sum from right to left, extending the sums over the indices $(\ell_1, \ell_1+ \ell_2)$, $(\ell_3, \ell_3+ \ell_4)$, and $(\ell_5, \ell_5+ \ell_6)$ to all possible values, and taking the suprema over indices of the two intermediate $\tau$-factors (after summing them over their respective displacements), subject to the constraints that they cannot contain more than $n_2-n_1$ and $n_3-n_2$ pivotals, respectively. 

For each of the nine terms that arise this way, we use \eqref{pi-worst-assump-rep-1} to bound the $\pi$-factors and the second bound in \eqref{tau-ass} to bound the factors $\tilde{\tilde{\tau}}$. This obtains the upper bound
\begin{equation}
	3^{2r+2} K^4 (C_\pi^{\sss (1)})^3 |n_2 - n_1|^{r/\twa} |n_3-n_2|^{r /\twa} \le C |n_3 - n_1|^{2r/\twa},
\end{equation}
so this term also satisfies \eqref{tightness-condition-aim-2}.
\medskip

\paragraph{{\bf The boundary cases involving $i_1=0$ and $i_3=M$.}} The boundary cases where $I_{i_1-1}=0$ (which means that the interval that contains $n_1$ also contains 0), and $I_{i_3}=M$ (which means that the interval that contains $n_3$ also contains $n_4$) are identical. In fact, in these cases, either the factor $\tilde\tau_{n_1-\ell_1}(x_1)$ vanishes and $\ell_1$ is set to $n_1$, or the final factor $\tilde{\tilde{\tau}}_{n-n_3-\ell_i}$ with $i =2,4,$ or $6$ vanishes, respectively, and $\ell_i$ is set to $n-n_3$. This actually makes these contributions smaller, so that the bounds are simpler. See \eqref{eqFinDimProof2} for an example of how such boundary terms arise. We omit further details. This completes the proof for~$Y_n$.
\medskip

\paragraph{{\bf Tightness of $(X_n)$.}}
The proof for $X_n$ is similar. We replace $\E_{p_c}^*$ by $\Eiic$ in \eqref{tightness-condition} above and define $\rho_{n_1,n_2,n_3,n_4}(x_1,x_2,x_3,x_4)$ similarly to $\tau_{n_1,n_2,n_3,n_4}(x_1,x_2,x_3,x_4)$, and use the expansion \eqref{rhon-exp-pivs-indef} to isolate the $\pi$ and $\psi$-terms that contain the points $x_1,x_2,x_3$ and $x_4$. Since the expansions for $\rho$ and $\tau$ only differ in the final coefficient, with $\psi$ replacing the final $\pi$ in the expansion of $\rho$, the cases where $i_3 < M$ are the same as before. 

The case where $i_3 =M$ requires a bit more care now, because the $\psi$-factor cannot be bounded exactly as a $\pi$-factor. This means that for a full proof we would need the analogue of Lemma~\ref{lem-pi-worst-assump} for these $\psi$-factors. We will not state such a lemma, or prove it, because it would be too repetitive. 
\qed
\medskip

We conclude this section by describing what we would need in order to obtain tightness in dimensions $d>3(\am)$, which we currently are not able to do:

\begin{rk}[Why the critical dimension in our tightness proof is too high]
\label{rem-dimension-tightness}
\rm
The following is a somewhat esoteric explanation of the requirement that $d > 6\twa$ is needed in the proof of tightness of $Y_n$ and $X_n$ above. A reader who is not very familiar with diagrammatic estimates may find it difficult to follow, and is advised to study the proofs in the \longversion{appendices first, especially those in Appendices \ref{sec-LEcoefficients}, \ref{sec:spatfrac} and \ref{sec:multiweight},}\shortversion{extended version of this paper first (especially the proofs in \cite[Appendices~A, B.2, and B.3]{HeyHofHulMie17b}),} if they are determined to find out the reasons here.

The proofs of tightness above are suboptimal because in \eqref{pi-worst-assump-rep-1}--\eqref{pi-worst-assump} we sum over $m_i$. If we inspect the contribution due to $i_1=i_2=i_3$ in more detail, we see that it would suffice to prove that
\begin{equation}
	\label{pi-better-assump}
	\sum_{m_1,m_4} \sup_{\substack{m_2,m_3\colon \\m_2-m_1\leq m_3-m_1\leq n_3-n_1}}
	\sum_{y_1,y_2,y_3,y_4} |y_2-y_1|^r(|y_3-y_2|^r+1)|\pi_{m_1,m_2,m_3,m_4}(y_1,y_2,y_3,y_4)|
\end{equation}
is bounded by $C_{\pi}^{\sss(3)}|n_3-n_1|^{2r/(\am)}$ for some $C_{\pi}^{\sss(3)}<\infty$. We conjecture that this is true whenever $d>3(\am)$. 

The reason that we cannot prove this bound on \eqref{pi-better-assump} is that our bounds on $\pi_{m_1,m_2,m_3,m_4}$ are phrased in terms of quantities like $\tau_m^{A}(x)=\prob_{p_c}(\{0\conn x, |\Piv(0,x)|=m\}\off A)$ (see e.g.\ \eqref{e:offexample}) and we actually do not know that $\sum_x \tau_m^{A}(x)\leq K$ uniformly in $m$ and $A\subseteq \Z^d$, nor that $\sum_x |x|^r \tau_m^{A}(x)\leq Km^{r/(\am)}$ uniformly in $m$ and $A$. We resolve this by summing over $m$, as in Lemma~\ref{lem-bd-E}, and obtain a bound in terms of $\tau^{A}(x)=\prob_{p_c}(\{0\conn x\}\off A)$. Here the dependence on $A$ is monotone, so that now we {\em can} remove the complicating $A$ dependence for the upper bound $\tau^{A}(x) \leq \tau(x)$, but this comes at the cost of a few dimensions:

Replacing the suprema over $m_2$ and $m_3$ by sums has the effect of adding two extra vertices to the diagrams bounding \eqref{pi-better-assump} (see e.g.\ Lemma~\ref{lem:fixedbond}), so that we need not just the triangle condition, which is satisfied when $d > 3\twa$, but also a similarly defined pentagon condition for $\sum_{m_1,m_2,m_3,m_4} \sum_{y_1,y_2,y_3,y_4} |\pi_{m_1,m_2,m_3,m_4}(y_1,y_2,y_3,y_4)|$ to be finite. Furthermore, the factors $|y_2-y_1|^r$ and $(|y_3-y_2|^r+1)$ together have a similar effect on the divergence that adding an extra vertex would have (see the proof of Lemma~\ref{lem-pi-worst-assump} in \longversion{Appendix~\ref{sec:multiweight}}\shortversion{\cite[Appendix~B.3]{HeyHofHulMie17b}} and also the proof of Proposition~\ref{prop-LEcoefficients2}(ii) in \longversion{Appendix~\ref{sec:spatfrac}, and \eqref{e:hybridbd} in particular,)}\shortversion{\cite[Appendix~B.2, (B.10)]{HeyHofHulMie17b}),} so that we actually need a \emph{hexagon condition,} which is satisfied when $d > 6(\am)$.  
As we explained above, tightness for $X_n$ also requires the same bounds on a weighted $\psi$-factor. It can be shown that the diagrams required here are also bounded when $d > 6 \twa$.
This explains the restrictions on the dimension in Theorem~\ref{thm-backbone-scal-limit}.
\end{rk}


\section{Convergence of the backbone as a set}
\label{sec-set-conv}

In this section we prove Theorem \ref{sec:main-results-1} and Proposition \ref{sec:main-results} under Hypothesis \ref{hyp:H}. We only prove the convergence of the processes restricted to the time-interval $[0,1]$, for the reasons that we have given at the beginning of Section \ref{sec-overview}.

\subsection{Hausdorff convergence subject to Hypothesis~\ref{hyp:H}}\label{sec:hausdorff}

Our main goal is to show that under Hypothesis \ref{hyp:H} the ``sausages'' $(\Scal_i)_{0\leq i\leq n}$ are all small compared to the scale of the pivotal walk $S_0,S_1,\ldots,S_n$. This is formalized as follows:

\begin{lemma}\label{sec:conv-backb-as}
Under Hypothesis \ref{hyp:H}, as $n\to\infty$,
\begin{equation}
    	f_\alpha(n)\max_{0\leq i\leq n}\diam(\Scal_i)\rightarrow 0
\end{equation}
in probability under $\Piic$, where for $A\subset \R^d$, $\diam(A)$ denotes the diameter of $A$.
\end{lemma}

We Prove Lemma~\ref{sec:conv-backb-as} below. 

\proof[Proof of Proposition~\ref{sec:main-results} subject to Lemma~\ref{sec:conv-backb-as}] 
Note that $\R^d$ is isometrically embedded in $(\Kcal,d_H)$ by the mapping $x\mapsto \{x\}$. Because of this embedding, the convergence in distribution of the process $(X_n(t),0\leq t\leq 1)$ in the space $\mathbb{D}([0,1],\R^d)$ implies the convergence of $(\{X_n(t)\},0\leq t\leq 1)$ to $(\{B^{\sss(\am)}_t\},0\leq t\leq 1)$ in the space $\mathbb{D}([0,1],\Kcal)$.

Next, Lemma \ref{sec:conv-backb-as} implies that
\begin{equation}
    	f_\alpha(n)\sup_{0\leq t\leq 1}d_H \big(\{S_{\lfloor nt\rfloor}\},\Scal_{\lfloor nt\rfloor}\big)
	\leq f_\alpha(n) \max_{0\leq i\leq n}\diam(\Scal_i) \longrightarrow 0
\end{equation}
in probability w.r.t.\ $\Piic$ as $n\to\infty$, because $S_i\in \Scal_i$ by definition. Since the latter uniform estimate dominates the Skorokhod distance, we get that $(f_\alpha(n)\Scal_{\lfloor nt\rfloor})_{0\leq t\leq 1}$ converges in distribution in $\mathbb{D}([0,1],\Kcal)$ to $(\{B^{\sss(\am)}_t\},0\leq t\leq 1)$. This implies the first statement of Proposition \ref{sec:main-results}.

Now, if $g$ is a c\`adl\`ag function from an interval $I$ to $\Kcal$, its \emph{historical path} is the function $\tilde{g}\colon I\to \Kcal$ defined by
\begin{equation}
    	\tilde{g}(t):=\Big(\bigcup_{s\in I,0\le s\leq t}g(s)\Big)^{\mathrm{cl}}\, .
\end{equation}
Indeed, $\tilde{g}$ takes its values in $\Kcal$, because the function $t\mapsto \diam(g(t))$ is right-continuous with left limits, and therefore it is bounded. And the same goes for the function $t\mapsto d(0,g(t))$, where by definition $d(x,A)=\inf_{y\in A}|x-y|$. This comes directly from the fact that $A\mapsto \diam(A)$ and $A\mapsto d(0,A)$ are continuous functions on $\Kcal$. Noting that $(\Bcal_{\lfloor nt\rfloor},t\geq 0)$ is the historical path associated with $(\Scal_{\lfloor nt\rfloor},t\geq 0)$, we see that the second statement of Proposition \ref{sec:main-results} is an immediate consequence of the first statement, the Continuous Mapping Theorem, and the following lemma:
\begin{lemma}
  \label{sec:conv-backb-as-1}
Let $(g_n)_{n\geq 1}$ be a sequence of functions converging in
$\mathbb{D}([0,1],\Kcal)$ to a limit $g$. Then the historical
paths $\tilde{g}_n$ converge to $\tilde{g}$ in
$\mathbb{D}([0,1],\Kcal)$ as well.
\end{lemma}
\proof
The Skorokhod convergence of $g_n$ to $g$ means that there exists a sequence of time-changes $\lambda_n,n\geq 1$, i.e., a sequence of increasing continuous functions from $[0,1]$ onto $[0,1]$, such that $\lambda_n$ converges uniformly to the identity, and such that
\begin{equation}
    	\eps_n=\sup_{0\leq t\leq 1}d_H \big(g_n\circ\lambda_n(t),g(t)\big)\longrightarrow 0\, .
\end{equation}
Now we have, for every $t\in [0,1]$,
\begin{equation}
    	d_H\big(\tilde{g}_n\circ\lambda_n(t),\tilde{g}(t)\big) 
	=d_H\Big(\bigcup_{0\leq s\leq \lambda_n(t)}g_n(s),\bigcup_{0\leq s\leq t}g(s)\Big) 
	=d_H\Big(\bigcup_{0\leq s\leq t}g_n\circ\lambda_n(s),\bigcup_{0\leq s\leq t}g(s)\Big).
\end{equation}
This equality holds because $d_H(A,B)=d_H(A^{\mathrm{cl}},B^{\mathrm{cl}})$ for any two subsets $A,B\subseteq \R^d$. By definition, if $x\in \bigcup_{0\leq s\leq t}g_n\circ\lambda_n(s)$ then it is in $g_n\circ\lambda_n(s)$ for some $s\in [0,t]$, and thus we can find some $y$ in $g(s)$ at distance at most $\eps_n$ from $x$. It follows that the converse also holds when we exchange the roles of $g_n\circ\lambda_n$ and $g$. This shows that $\sup_{0\leq t\leq 1}d_H\big(\tilde{g}_n\circ\lambda_n(t),\tilde{g}(t)\big)\leq \eps_n\longrightarrow 0$ as $n\to\infty$, as desired.
\qed
\medskip

\proof[Proof of Theorem \ref{sec:main-results-1} subject to Lemma~\ref{sec:conv-backb-as}]
Recall that the process $B^{\sss(\am)}$ is almost surely continuous at time~$1$. It immediately follows that the historical process $(\{B^{\sss(\am)}_s:0\leq s\leq t\},0\leq t\leq 1)$ is almost surely continuous at time $1$ as well. From this, we deduce that the projection $g\mapsto g(1)$ from $\mathbb{D}([0,1],\Kcal)$ to $\Kcal$ is almost everywhere continuous with respect to the law of the limiting process of $(f_\alpha(n)\Bcal_{\lfloor nt\rfloor},0 \le t\le 1)$. By standard properties of weak convergence of probability measures, we conclude that $(f_\alpha(n)\Bcal_{\lfloor nt \rfloor}, 0 \le t \le 1)$ converges to $\{B^{\sss(\am)}_s:0\leq s\leq 1\}$. The convergence of $f_\alpha(n)\Bcal_{\lfloor nT\rfloor}$ for a general $T>0$ follows by a scaling argument.\qed
\medskip

It remains to prove Lemma~\ref{sec:conv-backb-as}.

\proof[Proof of Lemma \ref{sec:conv-backb-as}.] For $\vep >0$ we may apply the union bound,
\begin{equation}
	\Piic\Big(f_\alpha(n)\max_{0\leq i\leq n}
	\diam(\Scal_i)\geq
	\eps\Big)\leq n\max_{0\leq i\leq n}\Piic \big(f_\alpha(n)\diam(\Scal_i)\geq \eps \big).
\end{equation}
By the Backbone Limit Reversal Lemma \cite[Lemma 4.2]{HeyHofHul14a}, we can write
\begin{align}
	\label{e:Qiicbb}
	\Piic(f_\alpha(n)\diam(\Scal_i)\geq \eps) &=\lim_{p\ua p_c}\frac{1}{\chi(p)} \sum_{x} \P_p\left(\{f_\alpha(n)\diam(\Scal_i^x)\geq \eps\}\cap \{0\conn x\}\right),
\end{align}
where $\Scal_i^x$ is the $i$th sausage along the path $0\conn x$ (but recall that $\Scal_i^x = \varnothing$ when there are fewer than $i$ pivotals). 
\bigskip

To investigate the latter probability, we observe that
\begin{align}
	&\{f_\alpha(n)\diam(\Scal_i^x)\geq \eps\}\cap \{0\conn x\}\nn\\
	&\quad \subseteq
	\bigcup_{e,b} \{|\Piv(0,\underline{e})|=i-1\}\circ \{e\text{ occupied}\}
	\circ \big\{\exists y\colon |y|>\eps f_\alpha (n)^{-1}\colon
	\{\bar{e} \conn y\}\circ \{\bar{e} \conn \bb\}\circ \{y\conn \bb\}\big\}\nn\\
	& \qquad \qquad \circ \{b \text{ occupied}\}\circ \{x-\tb\}.\label{BKR-set-up}
\end{align}
We see that $\{|\Piv(0,\underline{e})|=i-1\}$ is the only non-increasing event on the right-hand side of \eqref{BKR-set-up}. To ensure that the disjoint occurrence above is valid, we need to identify the right sets of witnesses. For the increasing events, we simply take bond-disjoint occupied paths that realize the connections they describe (in the case of $e$ or $b$ occupied, it is simply the bond $e$ or $b$). The witnesses for $\{|\Piv(0,\underline{e})|=i-1\}$ consist of all the bonds in all the paths between $0$ and $\underline{e}$, together with {\em all} the closed bonds in the graph. By definition, these sets of witnesses are all disjoint. This means that we can apply the BKR-inequality \eqref{BKR-ineq}, to obtain
\begin{align}
	\label{BKR-set-up-b}
	&\P_p\left(\{f_\alpha(n)\diam(\Scal_i^x)\geq \eps\}\cap \{0\conn x\}\right)\\
	&\quad \leq 
	\sum_{e,b} \tau_{i-1}(\underline{e})
	\P_p \big(\{e\text{ occ.}\}\circ \big\{\exists y\colon |y|>\eps f_\alpha (n)^{-1}\colon
	\{\bar{e} \conn y\}\circ \{\bar{e} \conn \bb\}\circ \{y\conn \bb\}\big\}
	\circ \{b \text{ occ.}\}\circ \{\tb \conn x\}\big).\nn
\end{align}

Applying the BK-inequality \eqref{BKR-ineq} to \eqref{BKR-set-up-b}, performing the summation over $x$ and letting $p\ua p_c$, it follows from \eqref{e:Qiicbb} that
\begin{multline}
	\Piic\big(f_\alpha(n)\diam(\Scal_i)\geq \eps\big)\\
    	\leq p_c \sum_{u} \tilde\tau_{i}(u)
	\sum_v\P_{p_c}\big(\exists y\colon |y|>\eps f_\alpha (n)^{-1}\colon
	\{0\conn y\}\circ \{0\conn v\}\circ \{y\conn v\}\big).
\end{multline}
We use that by Corollary~\ref{corr:pivballboundary}, $\sum_u \tau_{i-1}(u)\leq C$ for some $C$ to arrive at
\begin{equation}
    \begin{split}
	\Piic \big(f_\alpha(n)\diam(\Scal_i)\geq \eps \big) 
	& \leq C\sum_v\P_{p_c} \big(\exists y\colon |y|>\eps f_\alpha (n)^{-1}\colon
	\{0\conn y\}\circ \{0\conn v\}\circ \{y\conn v\} \big)\\
    	& = C \sum_{v} \Ppc\big(\{0 \conn v\} \circ \{\exists y\colon |y|>\eps f_\alpha (n)^{-1}\colon 
	\{0\conn y\}\circ \{y\conn v\}\}\big),
    \end{split}
\end{equation}
so that
\begin{equation}
	   \Piic\Big(f_\alpha(n)\max_{0\leq i\leq n} \diam(\Scal_i)\geq \eps\Big) 
	   \leq  C n \sum_{v} \Ppc\left(\{0 \conn v\} \circ \{\exists y\colon |y|>\eps f_\alpha (n)^{-1}
	   \colon \{0\conn y\}\circ \{y\conn v\}\}\right).
\end{equation}
Thus, we are left to show that
\begin{equation}
	\sum_z\P_{p_c}\left(\{0\conn v\}\circ \{\exists y\colon |y|>\eps f_\alpha (n)^{-1}\colon
	\{0\conn y\}\circ  \{y\conn v\}\}\right)=o(1/n).
\end{equation}
By the BK-inequality, we have the upper bound
\begin{equation}
	 C\sum_v\tau(v)\, \P_{p_c}\left(\exists y\colon |y|>\eps f_\alpha (n)^{-1} \colon
	\{0\conn y\}\circ\{y\conn v\}\right).
\end{equation}
Let $a_n = \frac{\vep}{2} (n/ \log n)^{1/(\am)}$. We split the sum over $v$ into $|v|\leq a_n$ and $|v|>a_n$. For $|v|> a_n$, we bound the sum by
\begin{equation}
	\sum_{|v|> a_n}\tau(v)^2
	\leq a_n^{-(\am)-\delta}
	\sum_{|v|> a_n}|v|^{(\am)+\delta}\tau(v)^2
	=O(1) a_n^{-(\am)-\delta}  =o(1/n),
\end{equation}
where, for the second inequality we have used \eqref{eqSpatDerivative} and the fact that $\tau(v)^2$ is the upper bound on $\sum_{n} \pi^{\sss(0)}_n(v)$ that is used in the proof of \eqref{eqSpatDerivative} (see \cite[Section 7]{HeyHofHul14a}).

Uniformly in $|v|\leq a_n$, we use that under Hypothesis \ref{hyp:H},
\begin{equation}
	\P_{p_c} \big(\exists y\colon |y|>\eps f_\alpha (n)^{-1} \colon
	\{0\conn y\}\circ \{y\conn v\} \big)\leq C(f_\alpha (n)/\vep)^{2 (\am)},
\end{equation}
and \cite[Theorem 1.5]{HeyHofHul14a} to bound
\begin{equation}
	\begin{split}
		\sum_{|v|\leq a_n} \tau(v) \P_{p_c} \big(\exists y\colon |y|>\eps f_\alpha (n)^{-1}
		\colon \{0\conn y\}\circ\{y\conn v\} \big) 
		&\leq C(f_\alpha (n)/\vep)^{2 (\am)} \sum_{|v|\leq a_n} \tau(v)\\
		&\leq C(f_\alpha (n)/\vep)^{2 (\am)} a_n^{(\am)}=o(1/n).
	\end{split}
\end{equation}
This completes the proof of Lemma \ref{sec:conv-backb-as}.
\qed
\medskip

\subsection{Verification of Hypothesis~\ref{hyp:H}}

\proof[Proof of Proposition \ref{prop-assum-H}.]
For the settings (i) and (iii) of the proposition, namely {finite}-{range} percolation under the strong triangle condition, and long-range percolation with $d>6$, $\alpha>4$ and satisfying \eqref{e:twoptasymp} we use the results from \cite{KozNac11} and \cite{Huls14}, respectively, that the extrinsic one-arm probability
is bounded by $C/r^2$, i.e.,
\begin{equation}
	\label{Gady-Asaf}
	\P_{p_c}\big(0\conn Q_r^c \big)\leq C/r^2,
\end{equation}
where $Q_r$ is the Euclidean ball of radius $r$ and $Q_r^c$ is its complement. For $|x|\leq m$, we can apply the BK-inequality and \eqref{Gady-Asaf} to bound
\begin{equation}
	\begin{split}
		\P_{p_c}\big(\exists y\in \Z^d:|y|>2m,\{0\conn y\}\circ\{x\conn y\}\big) 
		& \leq \P_{p_c}\big(\{0\conn Q_m^c\}\circ\{x\conn Q_m(x)^c\}\big)\\
		& \leq \P_{p_c}\big(0\conn Q_m^c\big)^2\leq C^2/m^4,
	\end{split}
\end{equation}
where $Q_m(x)$ is the Euclidean ball of radius $m$ around $x$.
This proves the claim in the {finite}-{range} case.

For setting (ii), namely long-range percolation with $d > 4\twa$, we bound
\begin{align}
	\P_{p_c}\big(\exists y\in \Z^d:|y|>2m,\{0\conn y\}\circ\{x\conn y\}\big)
	&\leq \sum_{|y|>2m} \tau(y)\tau(y-x)\\
	&\leq m^{-2(2\wedge \alpha)} \sum_{|y|>2m} |y|^{(2\wedge \alpha)}
	\tau(y)|y-x|^{(2\wedge \alpha)}\tau(y-x)\nn\\
	&\leq m^{-2(2\wedge \alpha)} \sup_x \sum_{y} |y|^{(2\wedge \alpha)}
	\tau(y)|y-x|^{(2\wedge \alpha)}\tau(y-x).\nn
\end{align}
We claim that for $d>4(2\wedge \alpha)$,
\begin{equation}
	\label{weighted-square}
	\sup_x \sum_{y} |y|^{(2\wedge \alpha)}
	\tau(y)|y-x|^{(2\wedge \alpha)}\tau(y-x)<\infty.
\end{equation}
This bound completes the proof. In the course of the proof of  \cite[Proposition 2.5]{HeyHofHul14a} it is proved for $d>3(2\wedge \alpha)$ that
\begin{equation}
	\label{weighted-bubble}
	\sup_x \sum_{y} |y|^{(2\wedge \alpha)+\delta}
	\tau(y)\tau(y-x)<\infty.
\end{equation}
The proof of \eqref{weighted-square} is very similar to this bound and is omitted here (see also \longversion{Appendix~\ref{sec:spatfrac}}\shortversion{\cite[Appendix~B.2]{HeyHofHulMie17b}} for a similar proof).
\qed

\vspace{0.3cm}

\noindent
{\bf Acknowledgements.}
This work was initiated during the Trimester Program \emph{Statistical Physics, Combinatorics and Probability: from discrete to continuous models} at Institut Henri Poincar\'e.
The work of MH, TH and RvdH is supported by the Netherlands Organisation for Scientific Research (NWO), through VICI grant 639.033.806, the Gravitation {\sc Networks} grant 024.002.003, and VENI grant 639.031.035. GM thanks the hospitality of CNRS/PIMS at University of British Columbia, and acknowledges support of the grant ANR-08 BLAN-0190. We thank a referee for pointing out a serious error in a previous version of the manuscript.

\longversion{\color{black}
\newpage

\appendix

\section{Weak bounds on the lace-expansion coefficients}
\label{sec-LEcoefficients}
In this appendix we prove Proposition \ref{prop:smallpi} and Lemmas \ref{lem:rhon-exp-lemma} and \ref{lem:rhon-exp-lemma-2}. These results pave the way for the lace-expansion analysis for $\hat{\tau}_n(k)$ and $\hat{\rho}_n(k)$, and are the necessary ingredients for the infrared bounds in Proposition \ref{prop:positivity}. These bounds only rely on the \emph{infrared bound} 
\begin{equation}\label{e:irbd}
	0 \le \hat{\tau}^{p}(k)\leq \frac{C}{1-\hat{D}(k)} \qquad \text{ when } \quad p \le p_c,
\end{equation}
as proved for models with $d > 3\twa$ that satisfy Assumptions \ref{ass:A} and \ref{ass:D} in \cite{HarSla90a, HeyHofSak08}. The most important tool that we use here is the Fourier transformation of entire diagrams. This tool was introduced in \cite{HeyHofHul14a}.


\subsection{Estimates on basic diagrams}
In the proofs that follow, we assume for brevity that any summation is over $\Zd$ unless otherwise specified, and we will often leave the $z$-dependence implicit. Besides writing $c$ and $C$ for generic constants, we will also sometimes write $C_1, C_2$ etc.\ for generic positive constants. The values of these numbered constants may thus also change from line to line. In this section we will also use the convention that the empty product equals 1.

Before starting with the proof, recall that $\Td = [-\pi,\pi)^d $ and define
\begin{eqnarray}\label{eqtridef}
     	\tri_p    &:= &  \int\limits_{\Td} \frac{\d^d k}{(2\pi)^d}   \htau^p(k)^3,
	\label{eqtridef}\\
    	\btri_p   & := &  \sup_{a,b,c \in \Td}  \int\limits_{\Td} \frac{\d^d k}{(2\pi)^d} 
	|\hat{D}(k+a)| \htau^p(k+a) \htau^p (k+b) \htau^p (k+c). \label{eqtridef2}
\end{eqnarray}
Note that  $\htau^p(k) \ge 0$ (see e.g.\ \cite[Chapter~10.3]{Grim99}) and 
\begin{equation}
	\max_x (\tau^p * \tau^p * \tau^p)(x) = \max_x \int\limits_{\Td} \frac{\dd k}{(2 \pi)^d} 
	\e^{i k \cdot x} \htau^p(k)^3 \le \tri_p,
\end{equation}
and similarly,
\begin{equation}\label{e:opentri}
	\max_x (\tau^p*\tau^p*\tau^p*D)(x) \le \btri_p.
\end{equation}
It is well known in the percolation literature (cf.\ \cite{BorChaHofSlaSpe05b},\cite{HeyHofSak08}) that there exists a $C >0$ such that
\begin{equation}\label{tribd}
    	\tri_p \le 1 + C \beta \qquad \text{ when } \quad p \le p_c,
\end{equation}
where $\beta =1/d$ or $\beta=1/L^d$, depending on whether the model is nearest-neighbor or {spread}-{out}. This is known as the \emph{strong triangle condition.} Moreover, by an application of H\"older's inequality (see e.g.\cite[Lemma 7.3]{HeyHofHul14a}), $|\hat D(k)| \le 1$ (by Parseval's Theorem and the fact that $D$ is a probability distribution), and the fact that the Fourier transform is periodic with period $2\pi$ along each axis by Assumption~\ref{ass:A},
\begin{equation}\label{tribd2}
		\btri_p \le \left(\,\int\limits_{\Td} \frac{\d^d k}{(2\pi)^d} |\hat{D}(k)|^2 
		\htau^p(k)^3 \right)^{1/3}  \left(\,\int\limits_{\Td} \frac{\d^d k}{(2\pi)^d} 
		\htau^p(k)^3 \right)^{2/3}
		\le C' \beta^{1/3} \qquad \text{ when } \quad p \le p_c,
\end{equation}
where the final inequality is due to \eqref{tribd} and \cite[Proposition 2.2]{HeyHofSak08}. We apply these bounds often in the proofs that follow, and sometimes we will do so without explicitly mentioning \eqref{tribd} or \eqref{tribd2}. Moreover, we will from here on no longer write the parameter $p$. Instead, we prove everything under the assumption that the strong triangle condition \eqref{tribd} holds. In our setting, this implies that the results hold for $p \le p_c$.

The following lemma states bounds on unweighted diagrams:
\begin{lemma}[Unweighted diagrams]\label{lem:basic1} Whenever the strong triangle condition \eqref{tribd} holds,
    \begin{eqnarray}
        \label{eqbasic1a} \sum_{x} \sum_{m \ge 0} \bar \pi_m^{\sss (0)}(x) &\le& 1+ C_1 \beta,\\
        \label{eqbasic1b} \sum_{x } \sum_{m \ge 1} \bar \pi_{m}^{\sss (1)} (x) &\le& C_2 \beta^{1/3},\\
        \label{eqbasic1c} \sum_{s,t} \Phl^{\sss (j)}(0, s, t) &\le& C_3 (C_4 \beta^{1/6})^{j} \qquad \text{ for } j \ge 2,\\
             \label{eqbasic2c} \sum_{x,v}   \Phr^{\sss (N-i-1)} (u,v,x) &\le& C_5(C_6 \beta^{1/6})^{N-j-1} \qquad \text{ for } 0 \le j \le N-1,\\
        \label{eqbasic1d}   \max_t  \sum_{x} \Upsilon(s,t,x) &\le& C_7 \beta^{1/3}.
    \end{eqnarray}
\end{lemma}
We prove Lemma~\ref{lem:basic1} below, but first we introduce a standard method in the lace-expansion literature: \emph{graphical bounds.}

We present the proofs of the main estimates on our lace-expansion coefficients in Lemmas \ref{lem:basic1}, \ref{lem:basic2}, and \ref{lem:basic3} largely in terms of diagrams, for reasons of legibility and brevity. There is a one-to-one relation between the diagrammatic proof and a traditional ``written'' proof in terms of formulae. The diagrams fall into two categories: $x$-space diagrams and Fourier-space diagrams. Fourier-space diagrams, and their relation to $x$-space diagrams are discussed in detail in \cite[Section 7.2]{HeyHofHul14a}. We give a brief overview here.

To start, we will henceforth refer to diagrams as defined in Section~\ref{sec-bounds-lace-exp} as $x$-space diagrams. Such diagrams can be represented graphically in an unambiguous way using the following transscription rules, which are an extension of the original lace-expansion transscription rules of \cite{HarSla90a}:
For any vertex of the diagram
\begin{itemize}
	\item  that is not summed over: draw either a small solid circle and the vertex's label 
	$\big(\stackrel{\bullet}{x}\big)$, or a solid black square ($\blacksquare$);
	
	\item  that is summed over: draw a small solid circle ($\bullet$);
	
	\item  that is maximized over: draw a cross ($\times$),
\end{itemize}
and for any two vertices that appear in the diagram as the arguments of
\begin{itemize}
	\item a two-point function $\tau$ (without a backbone label): connect their vertices with 
	a solid black line (---);
	
	\item a backbone two-point function $\tau^{\mathsf{b}}$: connect their vertices with 
	a solid red line (\color{red}---\color{black});
	
	\item a pivotal generating function $\Tau_z$: connect their vertices with a thick solid blue line
	 (\color{blue}{\bf ---}\color{black});
	
	\item a convolution with any of the three above functions with a factor $J$: append to 
	the line two vertical black dashes ($| \: |$).
\end{itemize}
We also draw a dashed gray line (\color{gray} - - -\color{black}) between a vertex that is maximized over and a reference vertex on the same diagram. This is not necessary for an unambiguous graphical representation, but helps to clarify the way in which the Fourier transform is applied to the diagram.

Recall the definition of the (discrete) Fourier transform in \eqref{e:Fourierdef}. In \cite[Section 7.2]{HeyHofHul14a} it is proved that there is a consistent way of taking the Fourier transform of a diagram in such a way that this yields an upper bound on the $x$-space diagram. For instance, applying this method, it is easy to show that $\sum_x \Phi^{\sss (1)}(x) $ can be bounded from above by
\begin{equation}\label{e:Fourierexample}
	  \int\limits_{\Td} \frac{\dd p_a}{(2 \pi)^d} \int\limits_{\Td} \frac{\dd p_b}{(2 \pi)^d}
	  \int\limits_{\Td} \frac{\dd p_c}{(2 \pi)^d} \htau(p_a)^2 \htau(p_b - p_a) |\hat D(p_b)| 
	  \htau(p_b) \htaub(p_b) \htau(p_c-p_b) \htau(p_c) \htaub(p_c).
\end{equation}
The precise application of this method is rather involved, so instead of repeating it here, we refer to reader to \cite{HeyHofHul14a}, and here restrict ourselves to a brief outline of the method, which is quite simple. The core idea of Fourier space representation of diagrams is this: 

To each diagram in $x$-space we can associate a \emph{dual diagram} of the Fourier space coordinates. If the $x$-space diagram can be graphically represented as a planar graph (as is the case for all possible diagrams discussed in this paper), then the Fourier transform of this diagram can be represented as the (planar) dual of this graph. More preciesly, given an $x$-space diagram, such a dual graph can be constructed as follows:

To start, construct a graph by identifying the vertices at both ends of any dashed gray line. Label the faces of the diagram. Then draw a solid black circle ($\bullet$) for each internal face of this graph and give it the same label as the face. A face vertex labelled $a$, say, represents an integral $\int_{\Td} \frac{\dd p_a}{(2 \pi)^d}$ in the Fourier dual diagram. Draw a circle ($\circ$) for the external face of the diagram, and label it $0$. This face vertex does not correspond to an integral. For any line in the $x$-space diagram, draw a similar line connecting the two faces on opposite sides of the diagram. This line represents the Fourier transform of the $x$-space line, so e.g.\ if the line in the $x$-space diagram is blue, representing a factor $\Tau_z(u-v)$, which separates faces $a$ and $b$, then the line in the Fourier diagram is also blue and represents a factor $\hat \Tau_z(p_a - p_b)$ (or equivalently, $\hat \Tau_z(p_b-p_a)$), where for the external vertex we set $p_0 \equiv 0$. Finally, for every factor $\hat D(p_b -p_a)$, we draw an arrowhead ($\blacktriangleright$) on a line with the same argument (there is always such a line in the diagrams that we study here). The result of this exercise is the graphical representation of the Fourier space diagram (and its formulaic equivalent). See for an example Figure \ref{fig:Fourierexample}, where the graphical representation of \eqref{e:Fourierexample} is given.
\medskip

\begin{figure}
	\includegraphics[width = \textwidth]{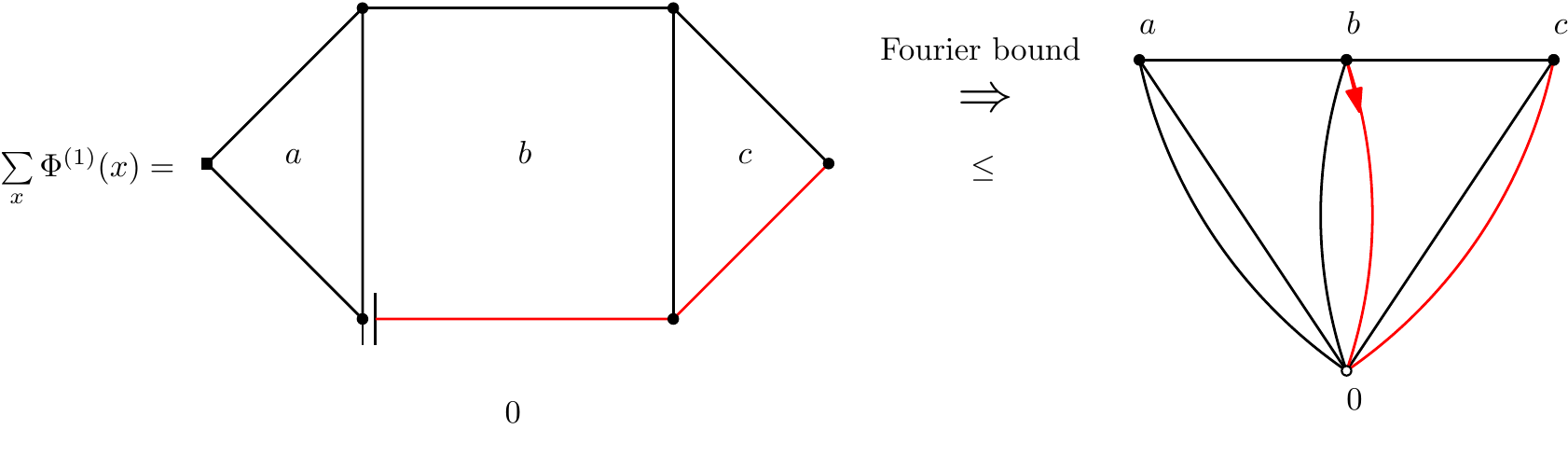}
	\caption{\label{fig:Fourierexample} An illustration of the graphical representation of the diagram $\sum_x \Phi^{\sss (1)}(x)$ and the Fourier dual diagram that bounds it. Compare the structure of these diagrams with \eqref{piexpansion} and \eqref{e:Fourierexample}.}
\end{figure}

In the course of the proofs below we use a number of simple identities, bounds, and operations, which we identify here with symbols. We define these symbols here:

\begin{defn}[Shorthand for standard identities, bounds, and operations]\color{white} . \color{black}

\begin{itemize}
	\item[$\lbrack\Acal \rbrack$] This can denote either of the two following bounds:
	\begin{itemize}
	 	
		\item If we are considering an $x$-space diagram $F$, then $\Acal$ denotes an
		 application of the open triangle bound \eqref{eqtridef}. In particular, if $F$ is of 
		 the form $F = \sum_{x,y} G(x,y) (\tau * \tau * \tau)(x,y)$, then we bound
		\begin{equation}
			\sum_{x,y} G(x,y) (\tau * \tau * \tau)(x,y) \le \sum_{x,y} G(x,y) 
			\max_{y'} (\tau * \tau * \tau)(x,y') \le \tri \sum_{x,y} G(x,y).
		\end{equation}
	
		\item If we are considering a Fourier-space diagram $\hat F$, then $\Acal$ denotes 
		an application of the Fourier-space triangle bound in \eqref{eqtridef}. In particular,
		if $\hat F$ is of the form 
		\begin{equation}
			\hat F =  \hat G \int\limits_{\Td} \frac{\dd k}{(2 \pi)^d} 
			\htau(k+a) \htau(k+b) \htau(k+c)
		\end{equation}
	 	for $a, b, c \in \Td$, such that there are no $k$-dependent terms in $\hat G$ (but it may contain $a,b,c$ dependent terms),
		then we bound $\hat F \le \tri \hat G$.
	\end{itemize}
	We write $\Acal^2$ to indicate that $\Acal$ is applied twice, etc.
	
	\item[$\lbrack \Bcal \rbrack$] This can denote either of the following two bounds:
	\begin{itemize}
		\item If we are considering an $x$-space diagram $F$, then $\Bcal$ denotes an
		 application of the open triangle bound \eqref{eqtridef2}. In particular, if $F$ is of 
		 the form $F = \sum_{x,y} G(x,y) (D * \tau * \tau * \tau)(x,y)$ for $n\ge 1$, then we 
		 bound $F \le \btri \sum_{x,y} G(x,y)$.

		\item If we are considering a Fourier-space diagram $\hat F$, then $\Bcal$ denotes 
		an application of the Fourier-space triangle bound \eqref{eqtridef2}. In particular, 
		if $\hat F$ is of the form 
		\begin{equation}
			\hat F = \hat G   \int\limits_{\Td} \frac{\dd k}{(2 \pi)^d} |\hat D(k+a)| 
			\htau(k+a) |\htau(k+b) \htau(k+c)
		\end{equation}
		then we bound $\hat F \le \btri \hat G$.
	\end{itemize}
	
	\item[$\lbrack \Ccal \Scal \rbrack$] Apply the Cauchy-Schwarz inequality: In diagrammatic 
	form this is expressed as taking the product of two copies of the diagram, both raised to the 
	power $1/2$, and in the first instance doubling a set of lines and removing another set of 
	lines, and in the second instance removing those lines that were doubled in the first instance, 
	and doubling the set of lines that were removed in the first instance.
	
	\item[$\lbrack \Dcal \rbrack$] We bound $|\hat D(k)| \le 1$.
	
	\item[$\lbrack \Fcal \rbrack$]  Take a Fourier transform of the diagram. This yields a 
	Fourier dual diagram as an upper bound.
	
	\item [$\lbrack \Scal \rbrack$] If we are considering a Fourier-space diagram $\hat F$, 
	then $\Scal$ denotes a bound by the Fourier-space weighted square diagram 
	\eqref{e:squzdef}. In particular, if $\hat F$ is of the form 
	\begin{equation}
		\hat F = \hat G  \int\limits_{\Td} \frac{\dd k}{(2 \pi)^d} 
		\htau(k+a) \htau(k+b) \htau(k+c) |\hat D(k+d)| z |\hat \Tau_z(k+d)|
	\end{equation}
	with $a,b,c,d \in \Td$, then we bound $\hat F \le \square_z \hat G$.
	
	\item [$\lbrack \Xcal N \rbrack$] Apply a bound obtained in the course of the proof of the 
	bound for the diagram labelled $N$.
\end{itemize}
\end{defn}

One final remark before we start with the proof of Lemma \ref{lem:basic1}: Here and below it will often be unimportant whether a line has a backbone label $\mathsf{b}$ or not, so where it does not matter, we will ignore this information and simply draw a black line.

\proof[Proof of Lemma \ref{lem:basic1}]
We start with \eqref{eqbasic1a}.
By \eqref{e:pinzerodef}, the BK-inequality and the definition \eqref{eqtridef}, we can bound
\begin{equation}\label{e:pizerotriangle}
    	\sum_{x} \sum_{m\ge 0} |\pi_m^{\sss (0)}(x)|  = \sum_x \P (0 \dbc x) 
	\le  \sum_x \tau(x)^2 \le \triangle.
\end{equation}

To further illustrate the graphical method described above, we perform the bound \eqref{eqbasic1b} both in written form and in graphical notation. In written form, the bound reads
\begin{equation*}
    \begin{split}
        \mathrm{A.} \qquad \sum_{x } \sum_{m \ge 1} |\pi_{m}^{\sss (1)} (x)| 
        \le & \sum_{s,t,u,v,x} A_3(0,s,t) B_1(s,t,u,v) A_3(x,u,v)\\
        & = \sum_{s,t,u,v,x} \tau(s) \tau(t) \tau(t-s) \ttau(u-s) \tau(v-t) \tau(u-v) \tau(x-u) \tau(x-v)\\
        & \le \sum_{s,t,u,v}  \tau(s) \tau(t) \tau(t-s) \ttau(u-s) \tau(v-t) \tau(u-v) 
        \max_{v'} \sum_x \tau(x-u) \tau(v'-x)\\
        & \stackrel{[\Acal]}{\le} \triangle \sum_{s,t,u,v}  \tau(s) \tau(t) \tau(t-s) \ttau(u-s) 
        \tau(v-t) \tau(u-v)\\
        & \le \triangle \sum_{s,t} \tau(s) \tau(t) \tau(s-t) \max_{t'} \sum_{u,v} \ttau(u-s) 
        \tau(v-u) \tau(t'-v)\\
        &\stackrel{[\Bcal]}{\le} \triangle \bar \triangle  \sum_{s,t} \tau(s) \tau(t) \tau(s-t)\\
        &\stackrel{[\Acal]}{\le} \triangle^2 \bar \triangle,
      \end{split}
\end{equation*}
where, for the fifth inequality we also used a symmetry argument. In graphical notation this becomes
\begin{figure}[h!]
  \includegraphics[scale=0.85, left]{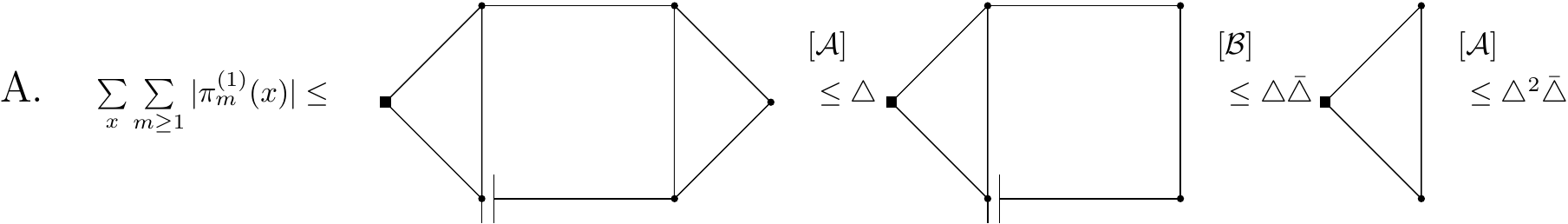}\\
\end{figure}

Now we establish a bound for \eqref{eqbasic1c}. Recall the definition of $\Phl^{\sss (i)}$ in \eqref{eqphildef}. We bound the terms of $\max_t \sum_{u,v,x,y} B_1(s,t,u,v) B_2^{\sss (i)}(u,v,x,y)$ for $i =0,1,2,3$ separately. We then iterate these bounds.

To start we bound
\[
\begin{split}
	\mathrm{B.} \qquad \max_t \sum_{u,v,x,y} & B_1(s,t,u,v) B_2^{\sss (0)}(u,v,x,y) \\
	& = \max_t \sum_{s,u,v,x,y,a,b} \ttau(u-s) \tau(x-u) \tau(a-u) \tau(a-x) \ttau(b-a) 
	\tau(y-b) \tau(v-b) \tau(v-y) \tau(t-v)\\
	& \stackrel{[\Fcal,\Dcal]}{\le}  \int\frac{\d^d p_1}{(2 \pi)^d}  \int\frac{\d^d p_2}{(2 \pi)^d}  
	\int\frac{\d^d p_3}{(2 \pi)^d} \htau(p_1)^2 \htau(p_2-p_1) |\hat D(p_2)| 
	\htau(p_2)^3 \htau(p_3-p_1) \htau(p_3)^2\\
	&\stackrel{[\Acal^2]}{\le} \triangle^2  \int\frac{\d^d p_1}{(2 \pi)^d} |\hat D(p_1)| \htau(p_1)^3 \\
	& \le \triangle^2 \bar \triangle.
\end{split}
\]
(Note that this diagram could also have been bounded without relying on the Fourier transform.)
In graphical notation this becomes 
\begin{figure}[h!]
  \includegraphics[scale=0.85, left]{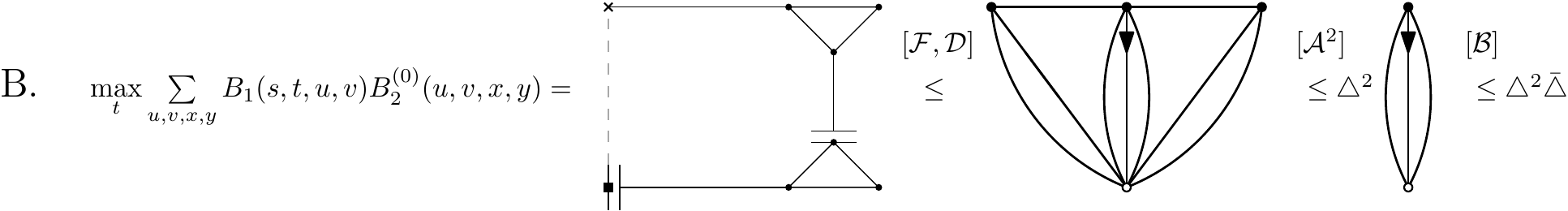}\\
\end{figure}

We proceed,
\[
\begin{split}
	\mathrm{C.} \qquad \max_t \sum_{u,v,x,y} B_1(s,t,u,v) B_2^{\sss (1)}(u,v,x,y) 
	& = \max_t \sum_{s,u,v,x,y} \ttau(u-s)  \tau(v-u) \tau(t-v) \tau(x-u) \tau(y-x)  \tau(v-y)\\
	& \stackrel{[\Acal]}{\le} \max_t \sum_{u,v} \ttau(u-s)  \tau(v-u) \tau(t-v) 
	\stackrel{[\Bcal]}{\le} \triangle \bar \triangle.
\end{split}
\]

By the same computation,
\[
	\mathrm{D.} \qquad \max_t \sum_{u,v,x,y} B_1(s,t,u,v) B_2^{\sss (2)}(u,v,x,y) 
	\stackrel{[\Xcal C]}{\le} \triangle \bar \triangle.
\]

The fourth case is most involved. Here we crucially rely on the Fourier transform and the Cauchy-Schwarz inequality:
\[
\begin{split}
	\mathrm{E.} \qquad \max_t \sum_{u,v,x,y} & B_1(s,t,u,v) B_2^{\sss (3)}(u,v,x,y) \\
	& = \max_t \sum_{s,u,v,x,y} \ttau(u-s)  \tau(t-v) \tau(x-u) \tau(v-x) \tau(y-u) \tau(v-y)\\
 	&\stackrel{[\Fcal]}{\le} \int\frac{\d^d p_1}{(2 \pi)^d}  \int\frac{\d^d p_2}{(2 \pi)^d} |\hat D(p_1)| 
	\htau(p_1)^2 \htau(p_2 -p_1)^2 \htau(p_2)^2 \\
	& \stackrel{[\Ccal \Scal]}{\le} \left(\int\frac{\d^d p_1}{(2 \pi)^d}  \int\frac{\d^d p_2}{(2 \pi)^d} 
	|\hat D(p_1)|^2 \htau(p_1)^3 \htau(p_2 -p_1)^2 \htau(p_2)\right)^{1/2} \\
	& \quad \times \left(\int\frac{\d^d p_1}{(2 \pi)^d}  \int\frac{\d^d p_2}{(2 \pi)^d}  \htau(p_1) 
	\htau(p_2 -p_1)^2 \htau(p_2)^3 \right)^{1/2}\\
	& \stackrel{[\Acal, \Bcal]}{\le} \sqrt{\triangle \bar \triangle} \left(\int\frac{\d^d p_1}{(2 \pi)^d} 
	 \int\frac{\d^d p_2}{(2 \pi)^d} \htau(p_1) \htau(p_2 -p_1)^2 \htau(p_2)^3\right)^{1/2}\\
	& \stackrel{[\Acal^2]}{\le}  \sqrt{\triangle \bar \triangle} \triangle.
\end{split}
\]
In graphical notation this becomes 
\begin{figure}[h!]
  \includegraphics[scale=0.85, left]{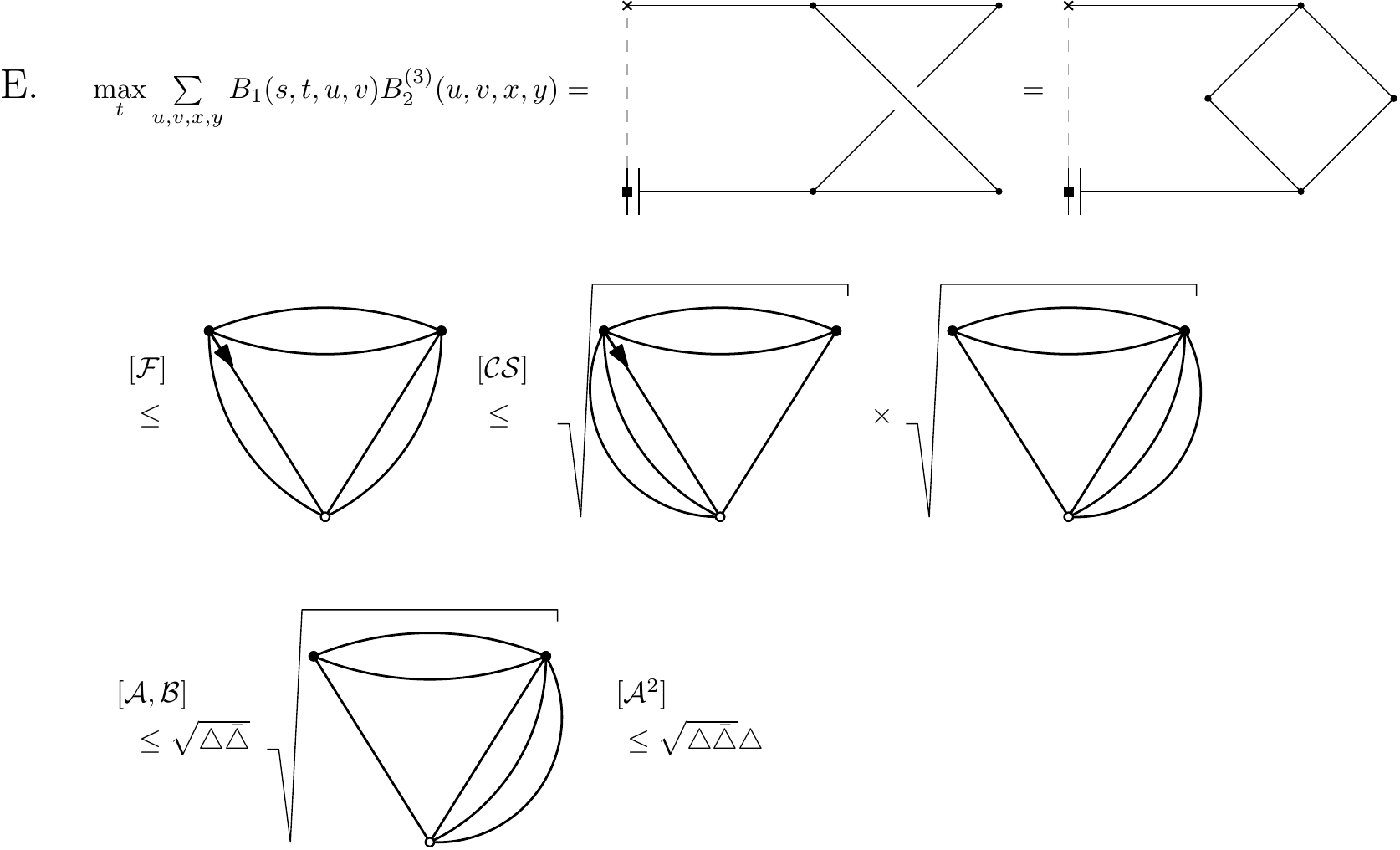}\\
\end{figure}

We use the bounds (B)--(E) iteratively to bound
\begin{equation}
    \begin{split}
        \sum_{s_j, t_j} \Phl^{\sss (j)} (s_j, t_j) =& \sum_{s_1, \dotso, s_{j}}\sum_{t_1, \dotso, t_{j}}
        \sum_{u_1, \dotso, u_{j-1}}\sum_{v_1, \dotso, v_{j-1}} A_3 (0, s_1, t_1)\\
        &\qquad \times \prod_{k=1}^{j-1} \left[B_1 (s_k, t_k, u_k, v_k) 
        B_2 (u_k, v_k, s_{k+1}, t_{k+1}) \right]\\
       \le &  \sum_{s_1, \dotso, s_{j-1}}\sum_{t_1, \dotso, t_{j-1}}\sum_{u_1, \dotso, u_{j-2}}
       \sum_{v_1, \dotso, v_{j-2}} A_3(0,s_1,t_1)\\
        &\qquad \times \prod_{k=1}^{j-1} \left[B_1 (s_k, t_k, u_k, v_k) 
        B_2 (u_k, v_k, s_{k+1}, t_{k+1}) \right]\\
        & \qquad \quad \times \sum_{i=0}^3 \max_{t'} \sum_{u_{j-1}, v_{j-1}, s_j, t_j} 
        B_1(s_{j-1}, t', u_{j-1}, v_{j-1}) B_{2}^{\sss (i)}(u_{j-1}, v_{j-1}, s_j, t_j)\\
        \stackrel{[\Xcal B, \dotso, \Xcal E]}{\le}&  ( \tri^2 \btri + 2 \tri \btri +  \tri \sqrt{\tri \btri}) 
        \sum_{s_{j-1}, t_{j-1}} \Phl^{\sss (j-1)}(s_{j-1}, t_{j-1}) \\
        \le & (  \tri^2 \btri +  2 \tri \btri + \tri \sqrt{\tri \btri})^j \sum_{s_1,t_1} A_3(0,s_1,t_1)\\
        \stackrel{[\Acal]}{\le} & \tri ( \tri^2 \btri +  2 \tri \btri + \tri \sqrt{\tri \btri})^j 
        \le  C_3 (C_4 \beta^{1/6})^j,
    \end{split}
\end{equation}
where the final inequality follows from \eqref{tribd} and \eqref{tribd2}.
\medskip

The proof of \eqref{eqbasic2c} is analogous, so we omit it here.
\medskip

Finally, we prove \eqref{eqbasic1d}:
\begin{equation}
\begin{split}
	 \max_t \sum_{x} \Upsilon(s,t,x) &= \max_t \sum_{u,v,x} B_1(s,t,u,v) A_3(x,u,v) \\
	 &\stackrel{[\Acal]}{\le} \tri \max_t \sum_{u,v} \ttau(u-s) \tau(v-u) \tau(t-v) 
	 \stackrel{[\Bcal]}{\le} \tri \btri \le C_5 \beta^{1/3},
	\end{split}
\end{equation}
where the final inequality follows from \eqref{tribd} and \eqref{tribd2}.\qed
\medskip


\subsection{Singly-weighted diagrams}

The following lemma states the neccesary bounds on diagrams with a single weight:
\bl[Singly weighted diagrams]\label{lem:basic2} Whenever the strong triangle condition holds,
    \begin{eqnarray}
        \label{eqbasic2b} \sum_x \sum_{m \ge 1} m |\pi_m^{\sss (1)}(x)| &\le& C_8 \beta^{1/3},\\
        \label{eqbasic2e} \max_t \sum_{x} \Upsilon(s,t,x;1_\mathsf{b}) &\le& C_9 \beta^{1/3},\\
        \label{eqbasic2d} \max_{t}  \sum_{v} \max_{u-v} \Omega(s,t,u,v;1_\mathsf{b}) &\le& C_{10} \beta^{1/6}.
    \end{eqnarray}
\el

\proof
To bound \eqref{eqbasic2b} we use Lemma \ref{lem:singleLE} to bound
\begin{figure}[h!]
  \includegraphics[scale=0.85, left]{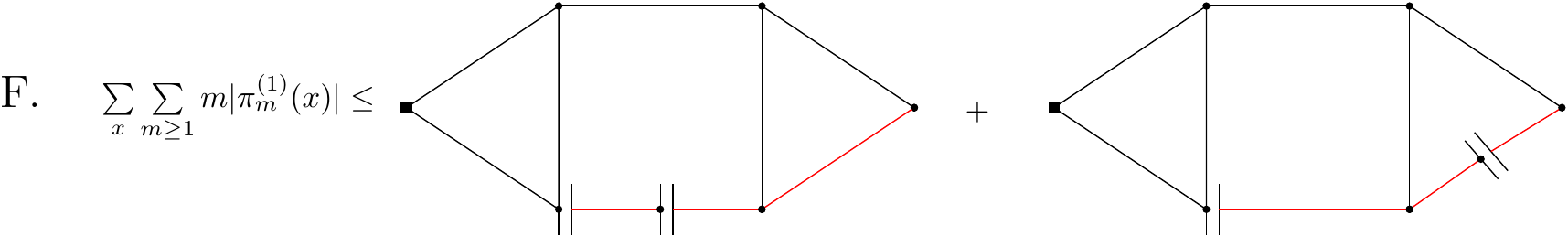}\\
\end{figure}

\noindent
The second term can be bounded by $\tri \btri^2$ by applying the bound $\Bcal$ twice and $\Acal$ once. To bound the first term we first apply a Fourier transform before we apply the bound $\Acal$ twice and then $\Bcal$, to obtain
\begin{figure}[h!]
  \includegraphics[scale=0.85, left]{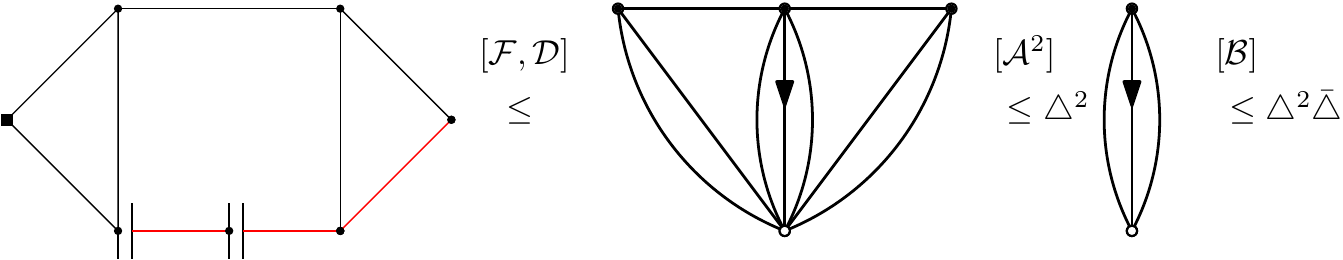}\\
\end{figure}

\noindent Combining these bounds with \eqref{tribd} and \eqref{tribd2} proves \eqref{eqbasic2b}.
\medskip

Next we prove \eqref{eqbasic2e}. This is the sum of two diagrams:
\begin{figure}[h!]
  \includegraphics[scale=0.85, left]{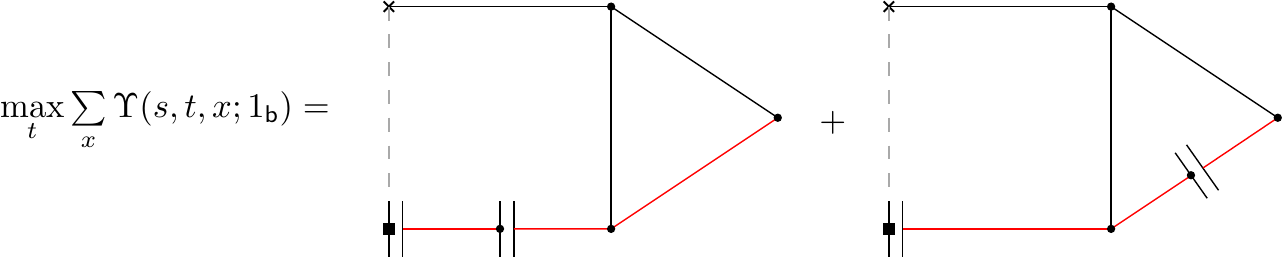}\\
\end{figure}

\noindent These terms are easily bounded as
\begin{figure}[h!]
  \includegraphics[scale=0.85, left]{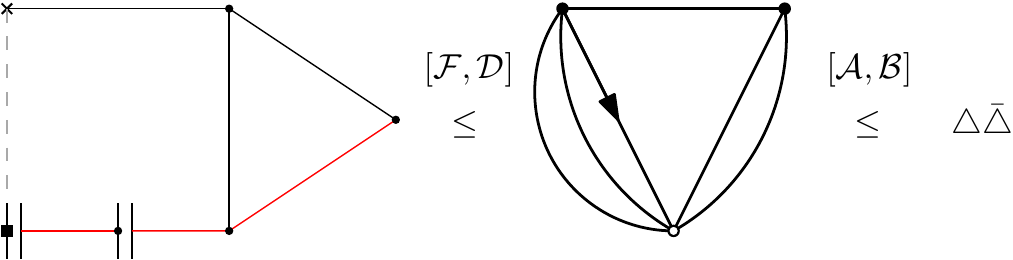}\\
\end{figure}

\noindent
and
\begin{figure}[h!]
  \includegraphics[scale=0.85, left]{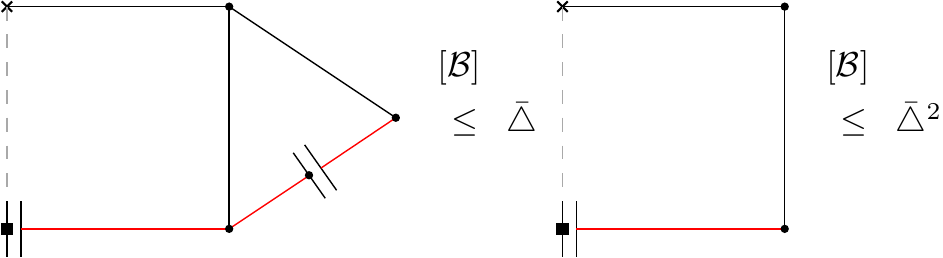}\\
\end{figure}

\noindent 
Combining these two bounds with \eqref{tribd} and \eqref{tribd2} we obtain \eqref{eqbasic2e}.

We proceed with the proof of \eqref{eqbasic2d}. 
The quantity $\Omega(s,t,u,v; 1_{\mathsf{b}})$ can be stated in graphical terms as the sum of fifteen diagrams, bounded in (G)--(U) below.

Although each of the fifteen diagrams is easy to bound using the same methods that we have used for bounds (A)--(F), we explicitly compute the bounds here, because they will be useful in the proof of Lemma \ref{lem:basic3} below.

\begin{figure}[h!]
  \includegraphics[scale=0.85, left]{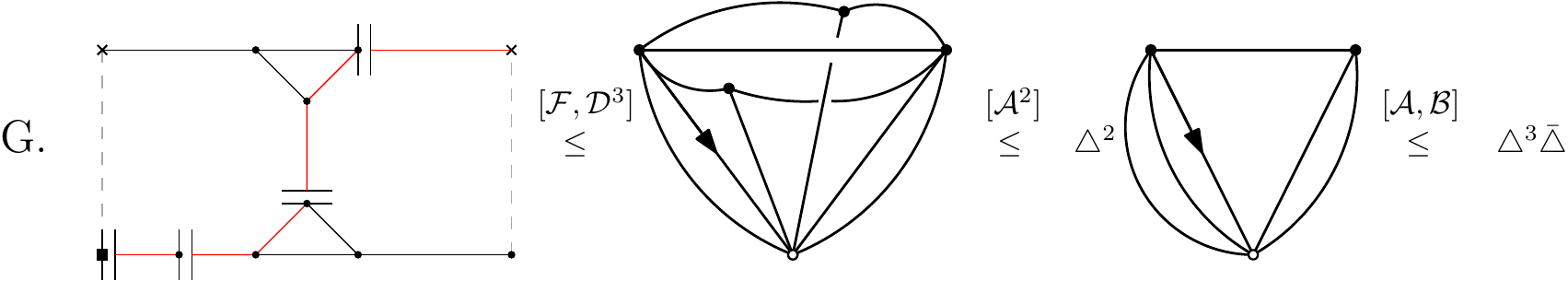}\\
\end{figure}
\begin{figure}[h!]
  \includegraphics[scale=0.85, left]{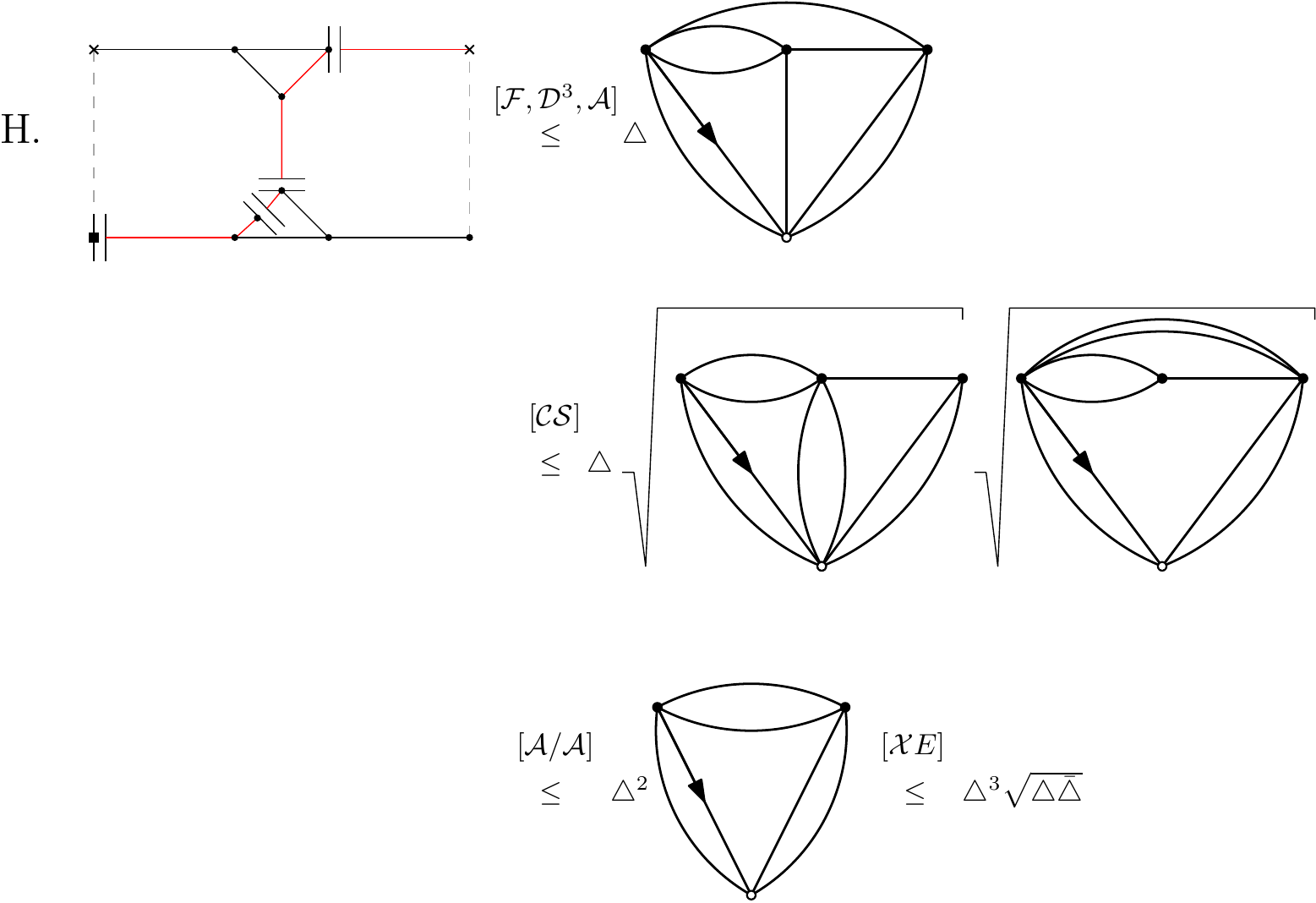}\\
  \includegraphics[scale=0.85, left]{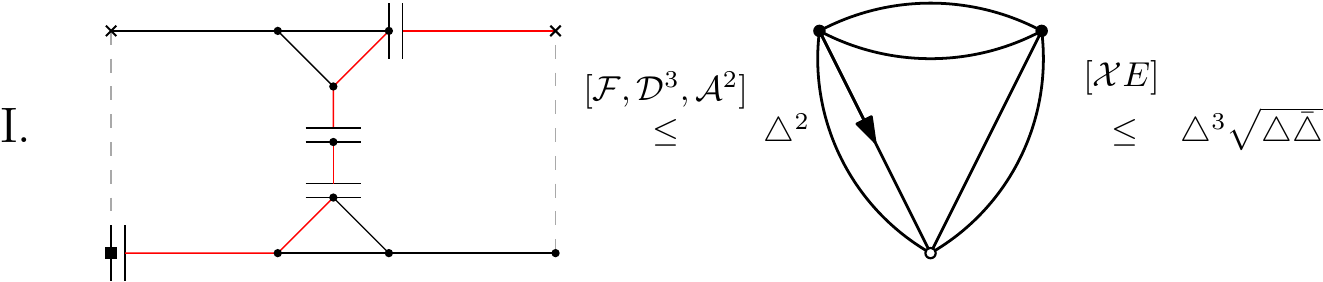}\\
  \includegraphics[scale=0.85, left]{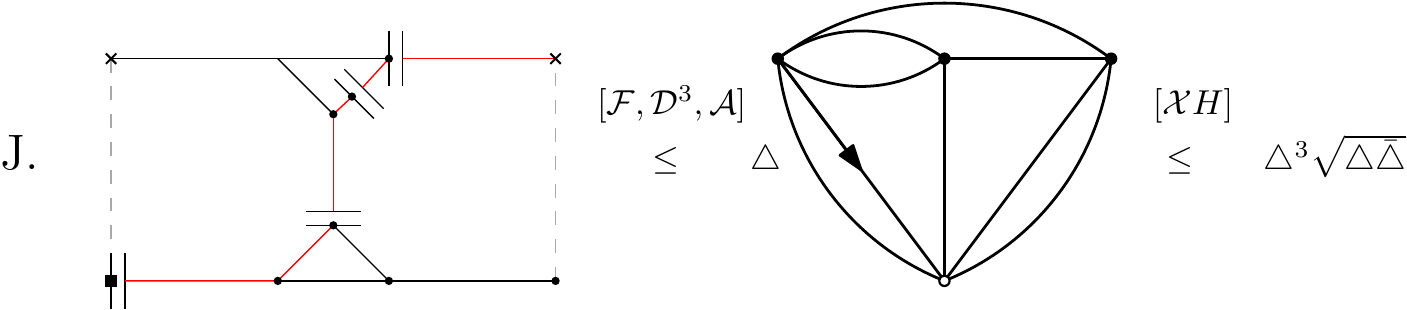}\\
  \includegraphics[scale=0.85, left]{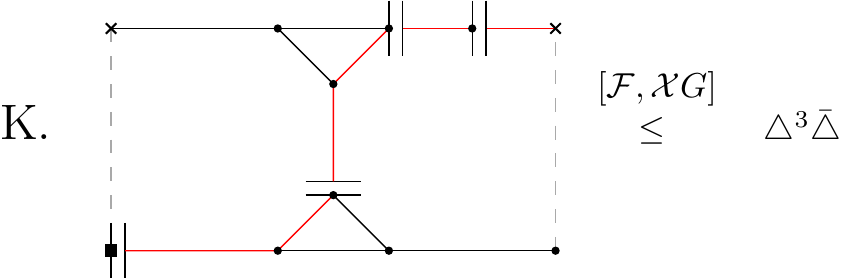}\\
  \includegraphics[scale=0.85, left]{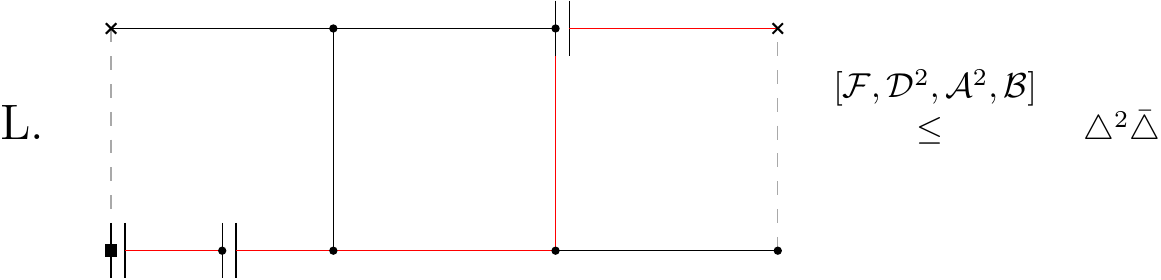}\\
  \includegraphics[scale=0.85, left]{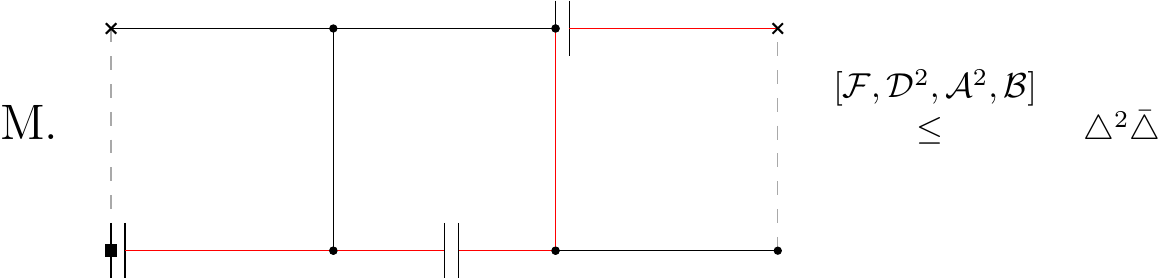}\\
\end{figure}
\clearpage
\begin{figure}[h!]
  \includegraphics[scale=0.85, left]{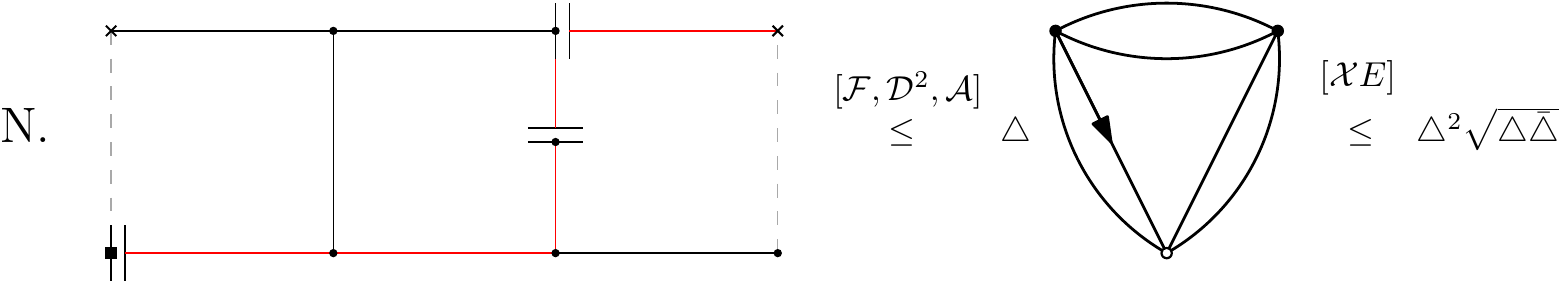}\\
  \includegraphics[scale=0.85, left]{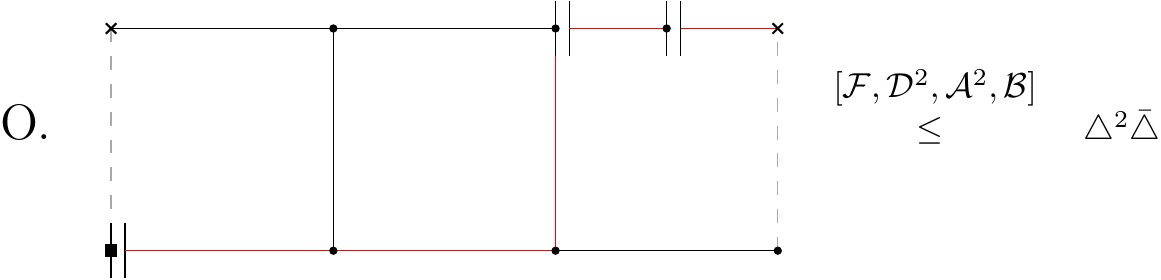}\\
  \includegraphics[scale=0.85, left]{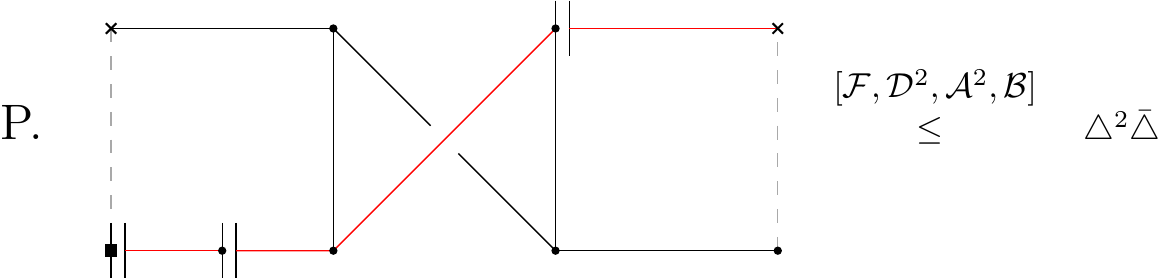}\\
  \includegraphics[scale=0.85, left]{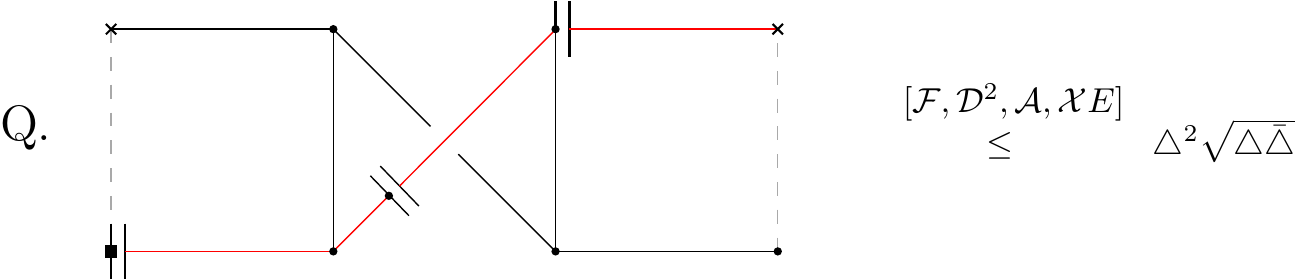}\\
  \includegraphics[scale=0.85, left]{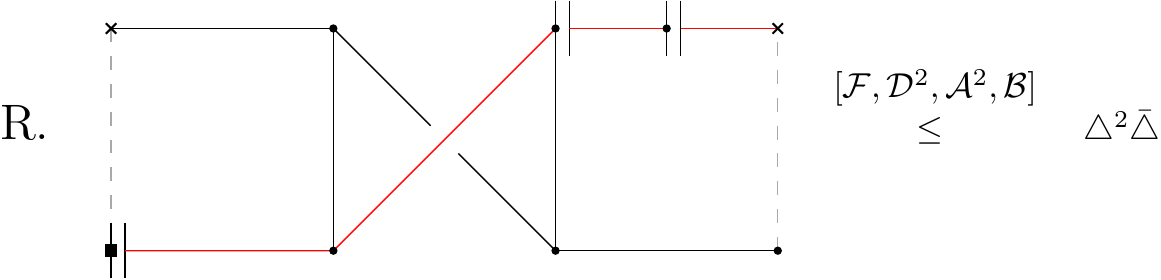}\\
  \includegraphics[scale=0.85, left]{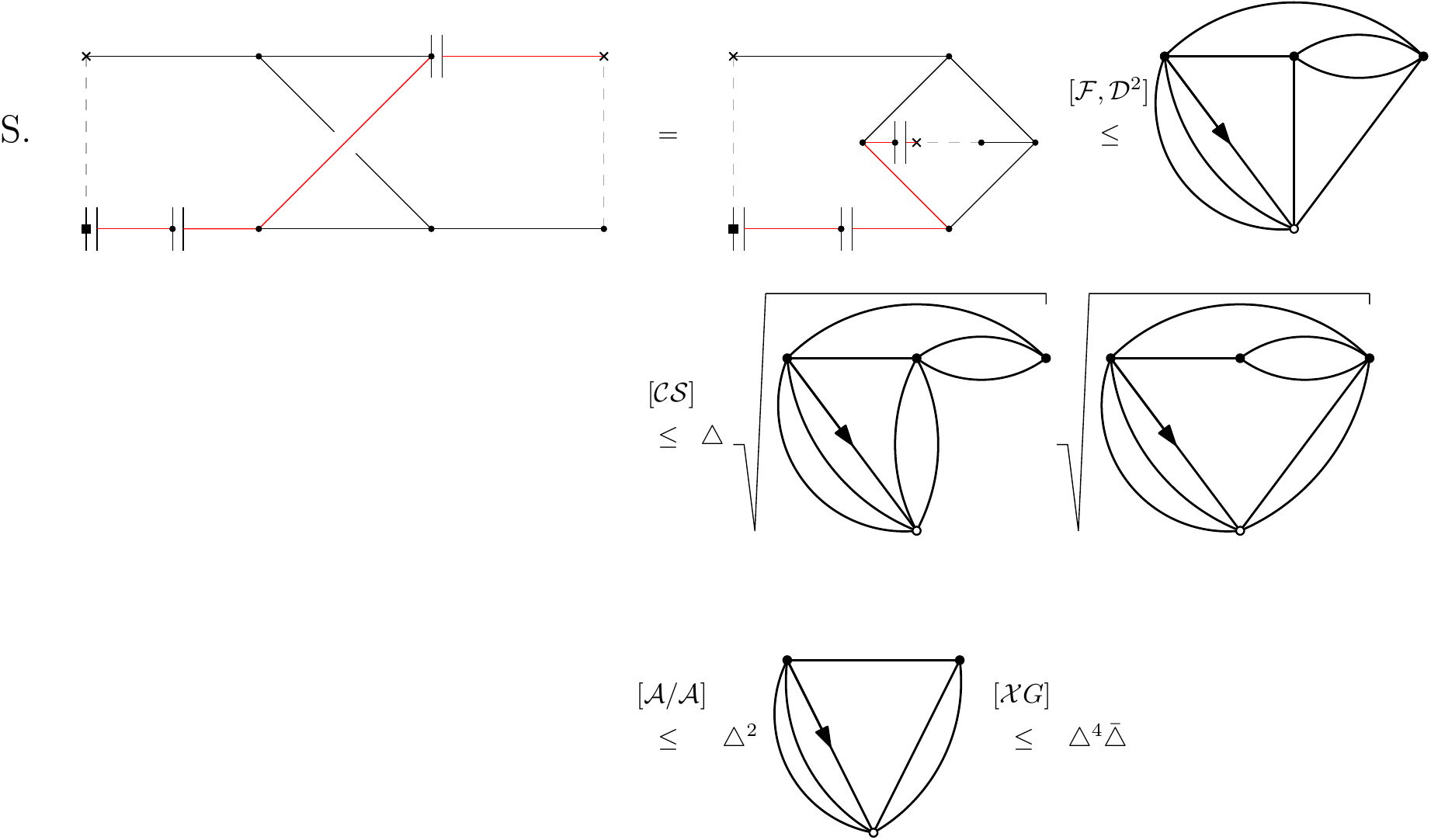}\\
\end{figure}
\clearpage
\begin{figure}[h!]
  \includegraphics[scale=0.85, left]{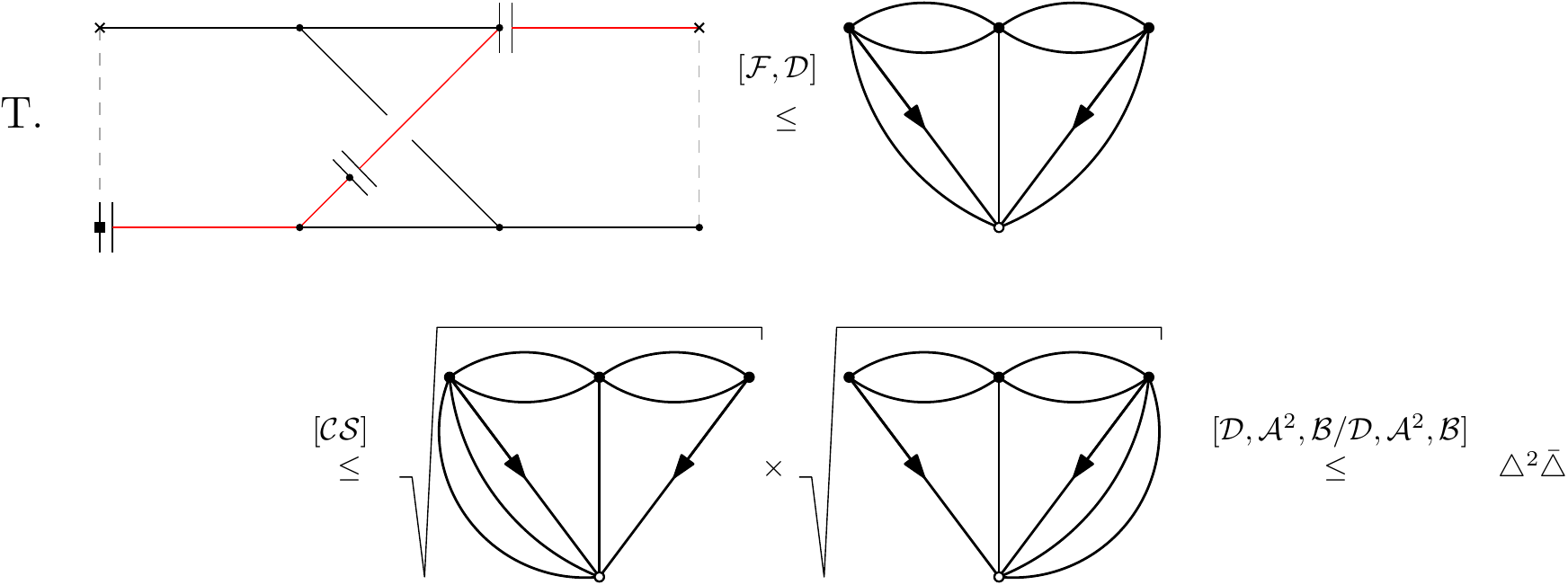}\\
\end{figure}
\begin{figure}[h!]
  \includegraphics[scale=0.85, left]{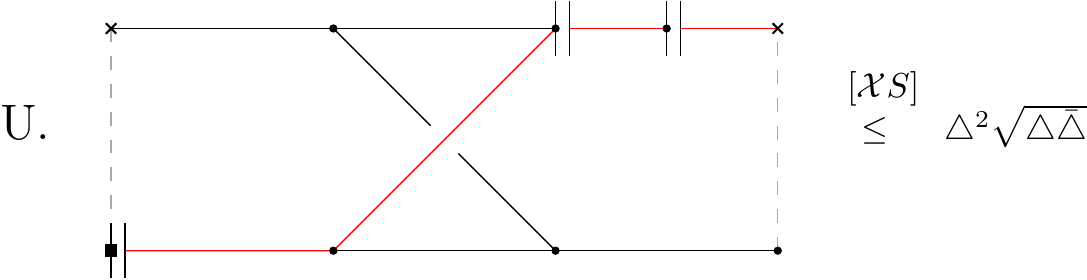}\\
\end{figure}

Adding the bounds (G)--(U) together, and applying \eqref{tribd} and \eqref{tribd2} we obtain
\begin{equation}
 \max_{t}  \sum_{v} \max_{u-v} \Omega(s,t,u,v;1_\mathsf{b})  \le \tri^4 \btri + 2 \tri^3 \btri + 2 \tri^3 \sqrt{\tri \btri} + 4 \tri^2 \sqrt{\tri \btri} + 6 \tri^2 \btri \le C_{10} \beta^{1/6},
\end{equation}
as claimed.\qed
\medskip


\subsection{Proof of Proposition \ref{prop:smallpi}}
We next complete the proof of Proposition \ref{prop:smallpi}:

\proof[Proof of Proposition \ref{prop:smallpi}]
Using \eqref{eqphildef} and \eqref{equpsilondef}, and applying the bounds in Lemma \ref{lem:basic1} to \eqref{piexpansion} we get the upper bound for $N\ge 2$,
\begin{equation}\label{piexpansion2}
	  \sum_{m \ge 0} \bar\pi_m^{\sss (N)} (x) \le \sum_{s,t} \Phi_A^{\sss (N)}(0,s,t)   
	  \max_{t'} \sum_x \Upsilon(s,t',x)  \le C_1 \beta^{1/3} \left(C_2 \beta^{1/6}\right)^{N-1}.
\end{equation}
Combining this with the bounds for $N=0$, $N=1$, we get
\begin{equation}
    	\sum_{x} \sum_{m \ge 0}  \bar\pi_m (x) \le 1 + C \beta 
	+ C_1 \beta^{1/3} \sum_{N \ge 2} (C_2 \beta^{1/6})^{N-1} 
	= 1 + O(\beta^{1/2}),
\end{equation}
when $\beta$ is small enough.
\medskip

Now we prove \eqref{eqsmallpi2}. 
It follows directly from Definitions \eqref{eqphildef}--\eqref{equpsilondef} that we can rewrite Lemma~\ref{lem:singleLE} for $N \ge 2$ as
\begin{equation}\label{piexpansion3}
\begin{split}
    	\sum_x \sum_{m \ge 1} m  \bar\pi_m^{\sss (N)}(x) \le \sum_x  \Phi^{\sss (N)}(x; 1_\mathsf{b}) 
	&= \sum_{i=1}^{N-1} \sum_{s,t,u,v,x} \Phl^{\sss (i)}(s,t) \Omega (s,t,u,v; 1_\mathsf{b}) 
	\Phr^{\sss (N-i-1)}(u,v,x)\\
    	&\quad + \sum_{s,t,x} \Phl^{\sss (N)}(s,t)\Upsilon(s,t,x;1_\mathsf{b}).
   \end{split}
\end{equation}
Applying the bounds of Lemmas \ref{lem:basic1} and \ref{lem:basic2} to \eqref{piexpansion3} we get for $N \ge 2$
\begin{equation}
    	\sum_x \sum_{m \ge 1} m \bar\pi_m^{\sss (N)}(x)
	\le C_4 \beta^{1/6} C_2 (C_3 \beta^{1/6})^{N-2}.
\end{equation}
Finally, combining the bounds for $N=0$, $N=1$ from Lemma \ref{lem:basic2} with the above bound, we get
\begin{equation}
    	\sum_x \sum_{m \ge 1} m \bar\pi_m(x) \le C_1 \beta^{1/3}  
	+ C_4 \beta^{1/6}  \sum_{N \ge 2} C_2 (C_3 \beta^{1/6})^{N-2} 
	= O(\beta^{1/6}),
\end{equation}
when $\beta$ is small enough. This completes the proof of Proposition \ref{prop:smallpi}. \qed

\subsection{Ingredients for the backbone lace expansion: Proof of Lemmas \ref{lem:rhon-exp-lemma} and \ref{lem:rhon-exp-lemma-2}}\label{app:cont}

\proof[Proof of Lemma~\ref{lem:rhon-exp-lemma}(a)]
Lemma~\ref{lem:fixedbond} states the bound
\begin{equation}\label{e:unifpmnbd}
	\sup_{p \le p_c} \sum_{w,y} \sum_{m \ge n} |\pi_{n,m}^p (y,w)| \le  \sum_{N \ge 0} \sum_{w} \Phi^{\sss (N)}(w ; 1_{\mathsf{b}}).
\end{equation}

We bound the $N=0$ term on the right-hand side using \eqref{e:opentri}. A comparison with Lemma~\ref{lem:singleLE} shows that the remaining terms can be bounded in the same manner as $\sum_{x} \sum_{m \ge 1} m \bar \pi_m (x)$ above. We conclude that the right-hand side of \eqref{e:unifpmnbd} is finite and independent of $n$ when $\beta$ is sufficiently small. This concludes the proof of Lemma~\ref{lem:rhon-exp-lemma}(a). \qed	
\medskip

\proof[Proof of Lemma~\ref{lem:rhon-exp-lemma-2}(a)] This is an immediate consequence of Lemma~\ref{lem:rhon-exp-lemma}(a) and the definition of $\pi_{n,m}^p$ in \eqref{pnm-def}. \qed
\medskip
	
\proof[Proof of Lemma~\ref{lem:rhon-exp-lemma-2}(b)]
We follow the proof of Hara in \cite[Appendix A]{Hara08}.
Let $\{\varnothing, < , >\}^N$ denote the set of all strings of length $N$ with alphabet $\varnothing, <,$ and $>$. 
Given an event $A$ measurable with respect to $\Fcal$, the $\sigma$-algebra of $\Pp$, we define the event $A^{\sss [0,N]} = A_0 \times \dotsm \times A_N$ with $A_i = A$, which is measurable with respect to $\Fcal^{\sss [0,N]}$, the $\sigma$-algebra of the measure $\P_{\sss (0)} \times \P_{\sss (1)} \times \dotsm \times \P_{\sss (N)}$ that arises in the lace expansion of Section~\ref{sec-le}. Now we define the functions
\begin{equation}
\begin{split}
	\label{pimNdef-pivs-2}
	\pi_{n,l}^{\sss(N)}(y,\vec{b}_{[l]}, u | A^{\sss [0,N]}) := &   \sum_{\vec{\bullet} \in \{\varnothing, >, <\}^N} (-1)^{\sum_{i=1}^N \delta_{\bullet_i, <}} \sum_{\substack{m_1, \dots, m_N\colon \\ m_1 
	+ \dotsm + m_N = l}} 
	\left[\prod_{i=1}^N J^p (b_{s_i}) \right] 
	\E_{\sss (0)} \Big[\indi_{\{y \Conn \ulb_1\} \cap A_0} \times \\
	& \qquad \times \E_{\sss (1)} \big[\indi_{E^{\bullet_1}_{m_1-1}(\olb_{s_1}, 
	\vec{b}_{[2,s_1-1]}, \ulb_{s_1}; \tCcal^{b_{1}}_0(0)) \cap A_1} \E_{\sss (2)} \big[\indi_{E^{\bullet_2}_{m_2-1}(\olb_{s_1},\vec{b}_{[s_1+1,s_2-1]},
	\ulb_{s_2}; \tCcal_1^{b_{s_1}}(\olb_{1}) ) \cap A_2} \\
	&\qquad \times \dotsm \, \times \E_{\sss (N)}[\indi_{E^{\bullet_N}_{m_N-1} (\olb_{s_{N-1}}, 
	\vec{b}_{[s_{N-1}+1,s_N-1]}, u ; \tCcal^{b_{s_N}}_{N-1}(\olb_{s_{N-1}})) \cap A_N}\big] 
	\dotsm \big]\big]\Big]\\
	& =: \sum_{\vec{\bullet} \in \{\varnothing, >, <\}^N} (-1)^{\sum_{i=1}^N \delta_{\bullet_i, <}} \pi_{n,l}^{\sss (N)}(y, \vec{b}_{[n]}, u ; \vec{\bullet}| A),
\end{split}
\end{equation}
where $\olb_n = y$ and $\bullet_i$ denotes the $i$th entry of the string $\vec{\bullet}$.

Observe that we can rewrite $\pi_{n,l}(y,\vec{b}_{[n]},u)$ as $\pi_{n,l}^{\sss(N)}(y,\vec{b}_{[l]}, u | \Omega^{\sss [0,N]})$, where $\Omega$ is the entire state space.

We will show that for any $\vec{\bullet} \in \{\varnothing, < , > \}^N$, the function $\pi_{n,l}^{\sss (N)}(y, \vec{b}_{[l]}, u ; \vec{\bullet} \,|\, \Omega^{\sss [0,N]})$ is continuous for $p<p_c$ and left-continuous at $p_c$. 
This follows by showing that $\pi_{n,l}^{\sss (N)}(y, \vec{b}_{[l]}, u ; \vec{\bullet} | \Omega^{\sss [0,N]})$ can be written as an increasing and as a decreasing limit of continuous functions, and is therefore both upper and lower semi-continuous on $p \le p_c$, which implies what we want to show. These functions will be chosen as functions that depend only on a finite number of bonds, which implies continuity.
We will henceforth omit most arguments and simply write $\pi_{n,l}^{\sss (N)}(\vec{\bullet})= \pi_{n+}^{\sss (N)}(y, \vec{b}_{[l]}, u ; \vec{\bullet} \,|\, \Omega^{\sss [0,N]})$.

Define the sets $Q_r := y+ [-r,r]^d$ for $r \ge 1$, and $Q_{[r,s]} := Q_s\setminus Q_r$ for $s > r \ge 1$. We will consider limits along the events $\{y \conn Q_r^c\}^{\sss [0,N]}$ and $\{y \nconn Q_{[r,s]}\}^{\sss [0,N]}$.

Define the following functions:
\begin{equation}
	\begin{split}
	\tilde \pi_{n,l}^{\sss (N)}(\vec{\bullet}; r) := \pi_{n,l}^{\sss (N)}(\vec{\bullet} \,|\,\{ y \conn Q_{r}^c\}^{\sss [0,N]}), & \qquad \qquad
 	\tilde{\pi}_{n,l}^{\sss (N)}(\vec{\bullet}; r,s): = \pi_{n,l}^{\sss (N)}(\vec{\bullet} \,|\,\{ y \conn Q_{[r,s]}\}^{\sss [0,N]}),\\
	\tilde{\tilde{ \pi}}_{n,l}^{\sss (N)}(\vec{\bullet}; r):= \pi_{n,l}^{\sss (N)}(\vec{\bullet} \,|\, \{ y \nconn Q_{r}^c\}^{\sss [0,N]}),& \qquad \qquad
	\tilde{\tilde{\pi}}_{n,l}^{\sss (N)}(\vec{\bullet}; r,s): = \pi_{m}^{\sss (N)}(\vec{\bullet} \,|\, \{ y \nconn Q_{[r,s]} \}^{\sss [0,N]}).
	\end{split}
\end{equation}
Then both $\tilde \pi_{n,l}^{\sss (N)}(\vec{\bullet}; r,s)$ and $\tilde{\tilde{\pi}}_{n,l}^{\sss (N)}(\vec{\bullet}; r,s)$ depend only on bonds that have both end-points inside $Q_{s}$, so they are continuous for any finite $r$ and $s$. It is not a-priori clear that the same holds for $\tilde \pi_{n,l}^{\sss (N)}(\vec{\bullet}; r)$ and $\tilde{\tilde{ \pi}}_{n,l}^{\sss (N)}(\vec{\bullet}; r)$ for long-range models.
To see that this is the case also, observe that
\begin{equation}
	\tilde \pi_{n,l}^{\sss (N)}(\vec{\bullet}; r,s) \le  \tilde\pi_{n,l}^{\sss (N)}(\vec{\bullet}; r) \le \tilde \pi_{n,l}^{\sss (N)}(\vec{\bullet}; r,s) + \Pp(\exists v \in Q_r, \exists w \in Q_s^c \text{ s.t. } \{v,w\} \text{ occ.}).
\end{equation}
By Assumption~\ref{ass:A}, the second term on the right-hand side tends to $0$ as $s \to \infty$ uniformly for $p \le p_c$, so we conclude that $\tilde \pi_{n,l}^{\sss (N)}(\vec{\bullet}; r)$ is also a continuous function on $[0,p_c]$ for all $r$. We can similarly reason that
\begin{equation}
	\tilde{\tilde{\pi}}_{n,l}^{\sss (N)}(\vec{\bullet}; r,s) - \Pp(\exists v \in Q_r, \exists w \in Q_s^c \text{ s.t. } \{v,w\} \text{ occ.})\le  \tilde{\tilde{\pi}}_{n,l}^{\sss (N)}(\vec{\bullet}; r) \le \tilde{\tilde{\pi}}_{n,l}^{\sss (N)}(\vec{\bullet}; r,s),
\end{equation}
so $\tilde{\tilde{\pi}}_{n,l}^{\sss (N)}(\vec{\bullet}; r)$ is also continuous on $[0,p_c]$ for all $r$.

Now, since $\tilde \pi_{n,l}^{\sss (N)}(\vec{\bullet}; r)$ is decreasing in $r$ we conclude that the limit as $r \to \infty$, $\pi_{n,l}^{\sss (N)}(\vec{\bullet})$, is upper semicontinuous for $p \in [0,p_c]$.
And $\tilde{\tilde{\pi}}_{n,l}^{\sss (N)}(\vec{\bullet}; r)$ is an increasing function in $r$ when $p \in [0,p_c]$, so now the limit as $r \to \infty$, which is again $\pi_{n,l}^{\sss (N)}(\vec{\bullet})$, is lower semicontinuous for $p \in [0,p_c]$. We conclude that $\pi_{n,l}^{\sss (N)}(\vec{\bullet})$ is continuous for $p \in [0,p_c)$ and left-continuous at $p_c$.

Now, by Lemma~\ref{lem:rhon-exp-lemma-2}(a), $\pi^{\sss (N)}_{n,l} (y, \vec{b}_{[l]}, u)$ is absolutely summable in $y,u, N, l$ and $b_{n+1},\dots,b_{l}$, uniformly in $n$, $\vec{b}_{[n]}$ and $p \le p_c$. Further, by \eqref{pimNdef-pivs-2}, each term $\pi_{n,l}^{\sss (N)}(y, \vec{b}_{[l]}, u ; \vec{\bullet} )$ of $\pi_{n,l} (y, \vec{b}_{[l]}, u ; \vec{\bullet})$ is continuous for $p \in [0,p_c)$ and left-continuous at $p_c$. We conclude that $\sum_{l \ge n} \sum_{b_{n+1},\dots,b_l} \sum_u \pi_{n,l} (y, \vec{b}_{[n]}, u )$ is also continuous for $p \in [0,p_c)$ and left-continuous at $p_c$, which proves the claim. \qed
\medskip

\proof[Proof of Lemma~\ref{lem:rhon-exp-lemma}(b)] This now follows straighforwardly from Lemma~\ref{lem:rhon-exp-lemma}(a) and Lemma~\ref{lem:rhon-exp-lemma-2}(b) by summing over $\vec{b}_{[n]}$. \qed

\section{Improved bounds on lace-expansion coefficients}
\label{sec-squares}
In this section we prove the improved bounds on the lace-expansion coefficients in Proposition~\ref{prop-LEcoefficients2}(ii)--(v) and Lemma \ref{lem:Diagbd}. We only give an outline for the proof of Proposition \ref{prop-LEcoefficients2}(ii)--(v), since many steps are similar to other proofs in the literature.

\subsection{Temporal doubly-weighted diagrams involving $z$-dependent square diagrams}

The main result in this section is the following lemma:

\begin{lemma}[Doubly weighted intermediate bounds]
\label{lem:basic3} 
Whenever the strong triangle condition holds,
\begin{eqnarray}
    \label{pimone} 
    	\sum_{x} \sum_{m \ge 1} m^2 \pi_{m}^{\sss (1)} (x) z^m 
	&\le& C_{11} \squz,\\
    \label{eqbasic3f} 
    	\max_{t} \sum_x \Upsilon (s,t,x; 1_\mathsf{b}, 1^z) 
	&\le& C_{12} \squz,\\
    \label{eqbasic3e} 
    	\max_t \sum_v \max_{u-v} \Omega (s,t,u,v; 1_\mathsf{b}, 1^z) 
	&\le& C_{13} \squz.
\end{eqnarray}
\end{lemma}

\proof
To prove \eqref{pimone} we first apply Lemma \ref{lem:doubleLE}, which yields (using the associativity of the discrete convolution to reduce the number of distinct diagrams from six to four),
\begin{figure}[h!]
  \includegraphics[scale=0.85, left]{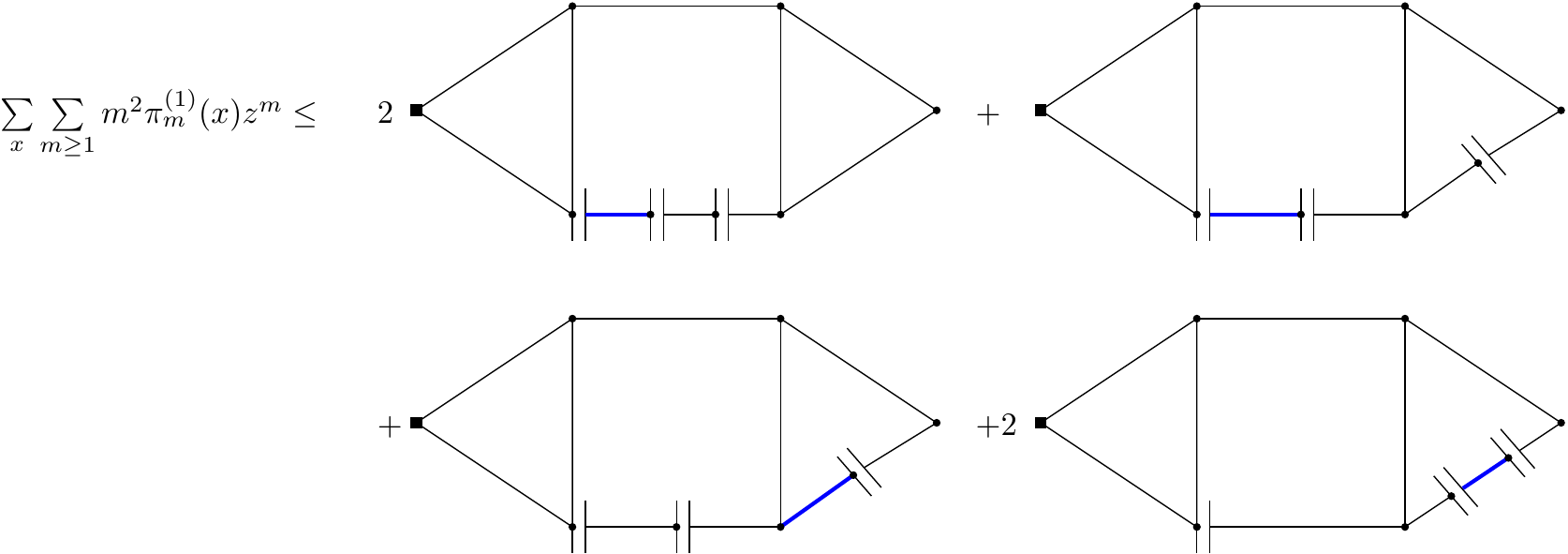}\\
\end{figure}

\noindent
We bound these four diagrams separately:

\begin{figure}[h!]
  \includegraphics[scale=0.85, left]{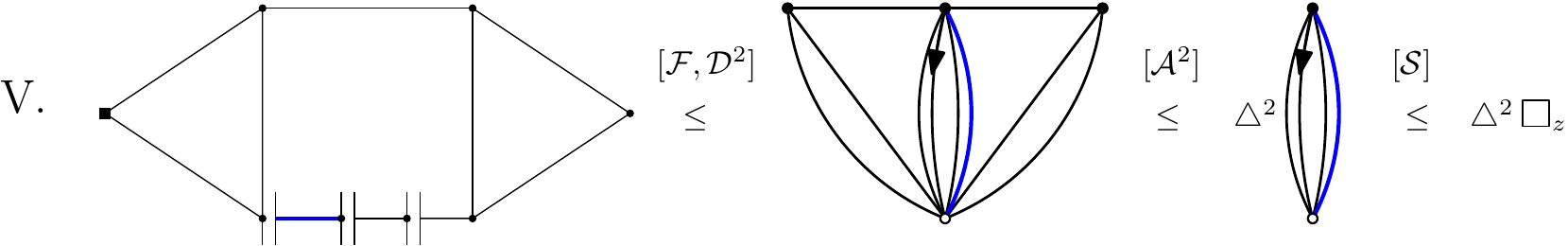}\\
\end{figure}
\begin{figure}[h!]
  \includegraphics[scale=0.85, left]{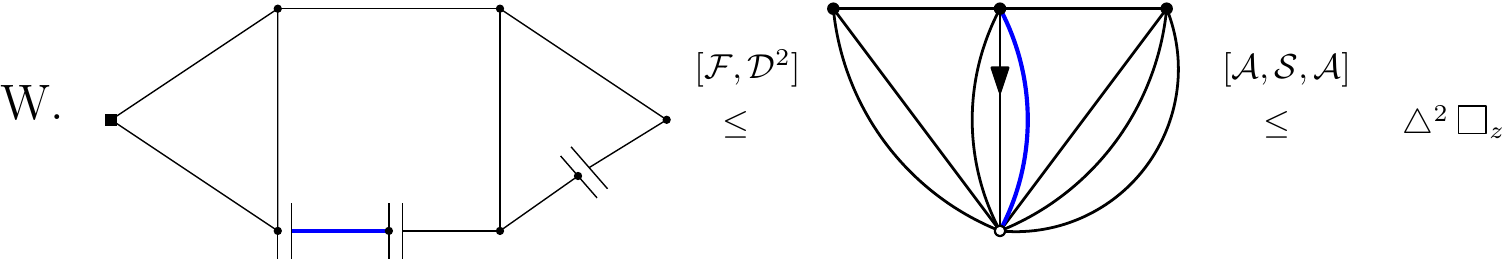}\\
  \includegraphics[scale=0.85, left]{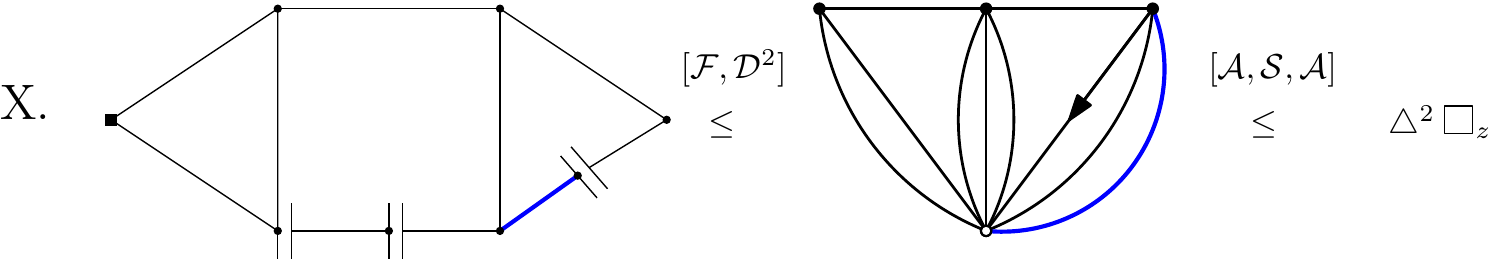}\\
  \includegraphics[scale=0.85, left]{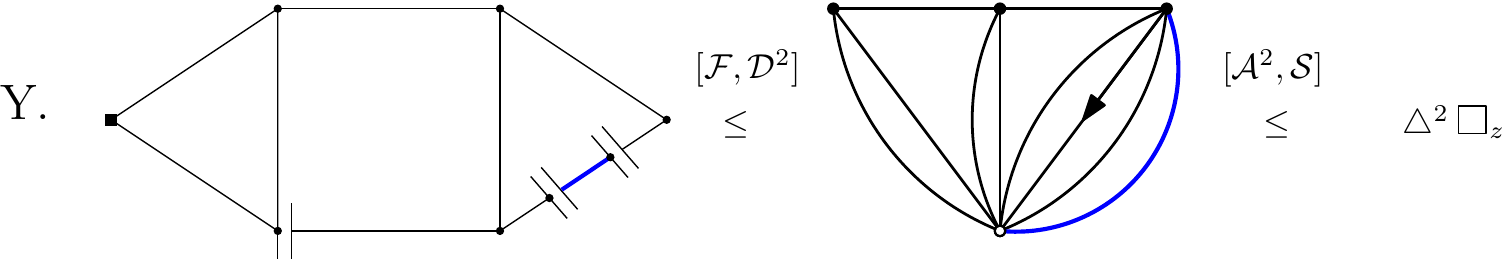}\\
\end{figure}

Adding the bounds (V)--(Y), we obtain
\begin{equation}
	\sum_x \sum_{m \ge 1} m^2 \pi_m^{\sss (1)}(x) z^m \le 6 \triangle^2 \square_z.
\end{equation}
\medskip

We now present a different way of deriving \eqref{pimone} that will be allow us to conclude \eqref{eqbasic3f} and \eqref{eqbasic3e}:

Observe that the diagrams bounding \eqref{pimone} can be derived from \eqref{eqbasic2b} by taking as many copies of each diagram of \eqref{eqbasic2b} as it contains red lines, and then replacing in each of these diagrams one uniquely chosen red line with a convolution $(z \Tau_z * J * \tau)$. Moreover, the bounds (V)--(Y) can easily be derived from the bounds of (F) by noting that simply one application of the open triangle bound $\Bcal$ should be replaced with an application of the open square bound. Likewise, the bounds on \eqref{eqbasic3f} and \eqref{eqbasic3e} can be derived by taking as many copies of a diagram of \eqref{eqbasic2e} and \eqref{eqbasic2d}, respectively, replacing for each copy one red line with a convolution $(z \Tau_z * J * \tau)$, following the steps of the associated bound, and simply replacing one application of $\Bcal$ with $\Scal$ (in some cases this will involve applying the bound $\Dcal$ to a different factor $|\hat D(k)|$, but this is always possible, and it is easy to see how). The result is that we obtain the bounds \eqref{eqbasic3f} and \eqref{eqbasic3e}, which completes the proof. \qed

\medskip

We now prove Lemma \ref{lem:Diagbd}:

\proof[Proof of Lemma \ref{lem:Diagbd}]
By Lemma \ref{lem:doubleLE},
\begin{equation}\label{piexpansion4}
\begin{split}
    \sum_x \sum_{m \ge 1} m^2  \bar\pi_m^{\sss (N)}(x) z^m & \le N^3 \Phi^{\sss (N)}(x) + N \Phi^{\sss (N)}(x; 1_\mathsf{b})\\
&\quad + \sum_{i=1}^{N-1} \sum_{s,t,u,v,x} \Phl^{\sss (i)}(s,t) \Omega (s,t,u,v; 1_\mathsf{b},1^z) \Phr^{\sss (N-i-1)}(u,v,x)\\
    &\quad + \sum_{s,t,x} \Phl^{\sss (N)}(s,t)\Upsilon(s,t,x;1_\mathsf{b},1^z),
   \end{split}
\end{equation}
for $N \ge 1$.
Applying the bounds from Lemmas \ref{lem:basic1}, \ref{lem:basic2}, and \ref{lem:basic3} to the right-hand side completes the proof of Lemma~\ref{lem:Diagbd}.\qed

\subsection{Spatially weighted diagrams: proof of Proposition~\ref{prop-LEcoefficients2}}\label{sec:spatfrac}
In this section we prove bounds on spatio-temporal weighted diagrams of the form $\sum_x\sum_{m \ge 1} |x|^{\delta} m \bar \pi_m^{\sss(N)}(x)$ and $\sum_x \sum_{m \ge 1} |x|^\delta \bar \psi_m^{\sss (N)}(x)$.

\proof[Proof of Proposition~\ref{prop-LEcoefficients2}(ii)] 
We give a very brief proof, since many of the details are identical to what is already in the literature. For a detailed proof of a similar bound we refer the reader to \cite[Section~7]{HeyHofHul14a}.

Lemma~\ref{lem:timespaceLE} establishes that
\begin{equation}
	\sum_x \sum_{m \ge 1} |x|^\delta m \bar \pi_m^{\sss (N)}(x) \le \sum_x \Phi^{\sss(N)}(x; 1_{\mathsf{b}}, 1^x_{\delta}).
\end{equation}
To prove that the right-hand side is finite for all $N\ge 0$ we use the same strategy as in the proof of Proposition~\ref{prop:smallpi} in Appendix~\ref{sec-LEcoefficients} to deal with the effect of the construction $1_{\mathsf{b}}$, and we use the same strategy as in \cite[Section~7]{HeyHofHul14a} to deal with the effect of the construction $1_{\delta}^x$. We know from the proof of Proposition~\ref{prop:smallpi} what the effect is of the construction $1_\mathsf{b}$.

The effect of the construction $1_\delta^x$ can be computed as follows. Recall that construction $1_\delta^x$ entails that for some path of lines in the diagram from $0$ to $x$ we take one line, $\tau(u,v)$ say, and multiply it with $|v-u|^\delta$. Then we sum over all lines along the path. We are free to choose the path, so we choose the path that starts by taking the path in a diagram of $\Phi^{\sss (N)}(x ; 1_{\mathsf{b}})$ that starts with the line $\tau(u_1)$ in $A_3(0,u_1,v_1)$, and then follows the backbone lines from $u_1$ to $x$. 

Now we express the factor $|v-u|^\delta$ in a convenient way by again using the identity
\begin{equation}\label{e:tdeltaid}
	t^\delta = c_\delta \int\limits_0^\infty \frac{1-\cos(st)}{s^{1+\delta}} \mathrm{d} s,
\end{equation}
where $c_\delta$ is a positive finite constant depending on $\delta$.
Using this identity and some simple geometric arguments (see e.g. \cite[Section~4.2.3]{BorChaHofSlaSpe05b}) we can give a bound that involves the factor $1-\cos(\vec{s} \cdot(v-u))$ with $\vec{s} = (s,0,\dots,0)$, and this is convenient because this factor arises naturally when taking the discrete Laplacian of the Fourier transform of $\tau$, namely, defining 
\begin{equation}
	\tau_{[s]}(x) := [1-\cos(\vec{s} \cdot x)]\tau(x),
\end{equation}
it is easy to show that
\begin{equation}
	\htau_{[s]}(k) = \htau(k) - \tfrac12 \htau(k-\vec{s}) - \tfrac12 \htau(k +\vec{s}).
\end{equation}

The part of the integral over $s$ from $1$ to $\infty$ is easy, because here we may use that $1-\cos(t) \le 2$ always to obtain an upper bound in terms of the diagrams for $\sum_x \sum_{m \ge 1} m \bar \pi_m (x)$, which is bounded by Proposition~\ref{prop:smallpi}.

To deal with the integral over $s$ from $0$ to $1$, we use that it is shown in \cite[(7.71)]{HeyHofHul14a} that for any $\theta \in [0,1]$ and $k \in \R^d$,
\begin{multline}\label{e:hybridbd}
		|\htau_{[s]}(k)| \le C \left(\frac{1}{1-\hat D(k-\vec{s})} + \frac{1}{1-\hat D(k)} + \frac{1}{1- \hat D(k+\vec{s})} \right)^{1-\theta} [1-\hat D(\vec{s})]^\theta\\
		   \times  \left(\frac{1}{1-\hat D(k-\vec{s})}\frac{1}{1-\hat D(k)} + \frac{1}{1- \hat D(k)} \frac{1}{1-\hat D(k+\vec{s})} + \frac{1}{1-\hat D(k-\vec{s})} \frac{1}{1-\hat D(k+\vec{s})} \right)^\theta.
	\end{multline}
Following the same bounding steps as in the proof of Proposition~\ref{prop:smallpi} as given in Appendix~\ref{sec-LEcoefficients}, we at some point encounter a Fourier space triangle with one factor $\htau(k)$ (say) replaced with a factor $\htau_{[s]}(k)$, to which we may apply the above bound. 
Roughly speaking, the effect of the above bound is that we replace one factor $[1-\hat D(k)]^{-1}$ by a factor $[1-\hat D(\vec{s})]^\theta [1-\hat D(k)]^{-1-\theta}$. (The shifts $-\vec{s}$ and $\vec{s}$ above are unimportant and can be dealt with using an appropriate H\"older inequality, similar to \eqref{tribd2} above, see also \cite[(7.75)]{HeyHofHul14a}.) So choosing $\theta$ sufficiently small so that $3 \twa + \theta < d$, for fixed $s$ this modified Fourier triangle diagram remains convergent when we integrate over $k$. 

What remains is the factor $[1-\hat D(\vec{s})]^\theta s^{-1-\delta}$, which we still need to integrate from $0$ to $1$. For this we simply use Assumption~\ref{ass:D}, which states that $1-\hat D(\vec{s}) \le w s^\twa$ for some constant $w$ when $s$ is sufficiently small. So choosing $\delta$ sufficiently small such that $-1-\delta + \theta \twa > -1$, we ensure that the integral from $0$ to $1$ is also convergent.

The remaining bounding steps can then be carried out precisely as in the proof of Proposition~\ref{prop:smallpi} again, and thus, the claim of Proposition~\ref{prop-LEcoefficients}(ii) follows. \qed 
\medskip

\proof[Proof of Proposition~\ref{prop-LEcoefficients2}(iii)]
Lemma~\ref{lem:psispatial} establishes that
\begin{equation}
	\sum_x \sum_{m \ge 1} |x|^\delta  \bar \psi_m^{\sss (N)}(x) \le \sum_{v,w,x} \Phi^{\sss(N)}(w; 1_{\mathsf{b}}(v,x), 1^x_{\delta}).
\end{equation}
The proof that the right-hand side is bounded goes along the exact same lines as the proof of Proposition~\ref{prop-LEcoefficients}(ii). The only difference is that now the point $x$ is located somewhere else along the backbone lines in the diagram. This simply means that there are \emph{fewer} terms in the bound on $\sum_x \sum_{m \ge 1} |x|^\delta  \bar \psi_m^{\sss (N)}(x)$ than there are in the bound on $\sum_x \sum_{m \ge 1} |x|^\delta m \bar \pi_m^{\sss (N)}(x)$, and therefore convergence of the former immediately follows from the latter. \qed

\proof[Proof of Proposition~\ref{prop-LEcoefficients2}(iv)]
By Lemma~\ref{lem:pidouble} this proof is again similar to the proof of Proposition~\ref{prop-LEcoefficients}(ii), except now an additional point $y$ and a small spatial weight $|y|^{\delta_4}$ also need to be taken into account. The additional point has been dealt with before in the proof of \eqref{eqsmallpi2} in Appendix~\ref{sec-LEcoefficients} above. Moreover, since $y$ is associated to a point on the path of backbone lines, the same strategy as was used for the weight $|x|^\delta$ in the proof above can also be applied here to the weight $|y|^{\delta_4}$. The bound follows.
\qed
\medskip

\proof[Proof of Proposition~\ref{prop-LEcoefficients2}(v)]
The claim follows from \cite[Proposition~5.2 and Lemma~5.4 and their proofs]{BorChaHofSlaSpe05b} (with $\Pi^{\sss (N)} (x)$ in that proof replaced by $\sum_{m \ge 0} \bar \pi_n(x)$) if we can prove that for all $N \ge 0$,
\begin{equation}\label{e:cosweightbd}
	\sum_x \sum_{m \ge 0} [1-\cos(k \cdot x)] \bar \pi^{\sss (N)}_m (x) \le C \beta^{(N-1 \wedge 0)/6} \sup_y \sum_x [1-\cos(k \cdot x)]\ttau(0,x) \tau(x,y).
\end{equation}
The proof of this bound is again similar to various of the proofs already given, and in particular it is highly analogous to the proof of \cite[Proposition~4.1]{BorChaHofSlaSpe05b}, so we give a brief, three-step outline only.

The first step to proving this is to use the bound \cite[(4.51)]{BorChaHofSlaSpe05b} that for $t = t_1 + \dotsm + t_N$ with $t,t_1,\dots,t_N \in \R$,
\begin{equation}
	1 - \cos t \le (2N+1) \sum_{j=1}^N [1-\cos t_j]
\end{equation}
to distribute the weight $[1-\cos(k \cdot x)]$ along a path from $0$ to $x$ in the diagram. The path we choose is that of the backbone lines together with the line from $0$ leading into this path. Using these new weights we can define a new construction that applies the weight $[1-\cos(k \cdot (u-v)]$ to a line of $\Phi^{\sss(N)}$ with displacement $u-v$, analogous to Construction~$1^{x,y}_\delta$ in Definition~\ref{def-con1delta}. This is the second step. The third step is then to bound these new diagrams using the same methods as in Appendix~\ref{sec-LEcoefficients}. Inspecting the diagrammatic bounds F--U in Appendix~\ref{sec-LEcoefficients}, it is clear that the effect of this new construction is diagrammatically almost identical to an application of Construction~$1_{\mathsf{b}}$, in that instead of associating an additional edge to a backbone line, we may now associate a weight $[1-\cos t_j]$ to such a line. In the Fourier space diagrams, the weight is then associated to \emph{a triplet of lines:} the line with the arrowhead and two parallel lines. Following the same bounding steps as in F--U, each in each of these diagrams except the second term in F, the final weighted term is the last one to be bounded. So for these diagrams we may follow the same bounding steps until we arrive at the Fourier transform of $\sup_y \sum_x [1-\cos(k \cdot x)]\ttau(0,x) \tau(x,y)$. Taking the inverse Fourier transform we arrive at the desired bound.  The other factors are bounded by factors $\tri$ and $\btri$, giving the factor $C\beta^{(N-1 \wedge 0)/6}$, similar to the bounds derived in the proof of Lemma~\ref{lem:basic2}. The second term in the bound F does not involve a Fourier space bound and can therefore be bounded directly to give the desired bound. This establishes \eqref{e:cosweightbd} and thus concludes the proof of Propostion~\ref{prop-LEcoefficients2}(v). \qed

\subsection{Multi-weighted diagrams: Proof of Lemma \ref{lem-pi-worst-assump}}\label{sec:multiweight}

The full proof is rather long due to the large number of different terms, but very much in the same spirit as the proof of Proposition~\ref{prop-LEcoefficients}, so we just indicate the main ideas for the proof of \eqref{pi-worst-assump}, which is the most involved bound to prove.

To start, observe that similar to Lemma~\ref{lem:pidouble}, we can bound \eqref{pi-worst-assump} from above by
\begin{equation}
	\sum_{N\ge0}
	 \sum_{y_1,y_2,y_3,y_4} \sum_{v_1,v_2,v_3,v_4} \Phi^{\sss (N)}(y_4; 1_{\mathsf{b}}(v_1,y_1),1_{\mathsf{b}}(v_2,y_2),1_{\mathsf{b}}(v_3,y_3), 1_{r}^{y_1,y_2}, 1_r^{y_2,y_3})\indi_{(v_1,y_1) \prec (v_2,y_2) \prec (v_3,y_3)},
\end{equation}
where the argument $(v_1,y_1) \prec (v_2,y_2) \prec (v_3,y_3)$ of the indicator means that the order of the edges should be such that $(v_1,y_1)$ comes before $(v_2,y_2)$ and $(v_2,y_2)$ comes before $(v_3,y_3)$ in the order induced by the product in \eqref{piexpansion}, and the paths for the constructions $1_{r}^{y_1,y_2}$ and $1_r^{y_2,y_3}$ are chosen such that they run along the backbone lines between vertices $y_1$ and $y_2$ and between $y_2$ and $y_3$, respectively. Because of the imposed order on the three edges, these paths do not intersect.

Let us first consider the case where $d> 7\twa$, and leave the case $d \in (6\twa, 7\twa]$ for a brief discussion below.

We can bound the effect of the constructions $1_{r}^{y_1,y_2}$ and $1_r^{y_2,y_3}$ in the same manner as was done in the proof of Proposition~\ref{prop-LEcoefficients} above. In particular, for each of these constructions we can go through the steps \eqref{e:tdeltaid}--\eqref{e:hybridbd} (with $\theta=0$) to write the weights as sums over additional lines in the Fourier diagram. In total this procedure adds precisely five lines to the each Fourier diagram of $\Phi^{\sss (N)}$.

Similar to \eqref{eqtridef2}, define the \emph{Fourier $n$-gon diagrams}
\begin{equation}\label{e:ngon}
	\sup_{a_1,\dots,a_n \in \T^d} \int\limits_{\T^d} \frac{\dd k}{(2\pi)^d} \prod_{i=1}^n \frac{1}{1-\hat D(k+a_i)}.
\end{equation}
It follows by an application of H\"older's inequality (as in \cite[Lemma~7.3]{HeyHofHul14a}) and Assumption~\ref{ass:D} that this integral is finite when $d > n\twa$.

The proof of \eqref{pi-worst-assump} now follows the same steps as the proof of Lemma~\ref{lem:basic1}, except that now the Fourier diagrams have more lines. We deal with these additional lines as they arise by bounding these terms using \eqref{e:ngon} with $n=4,5,6,$ or $7$ instead of \eqref{eqtridef2}, whichever applies. This completes the proof of \eqref{pi-worst-assump}. The other bounds in Lemma~\ref{lem-pi-worst-assump} are similar.

It remains to discuss the case $d \in (6\twa,7\twa]$. Here we use that we chose $r = \tfrac12 \twa +\vep$, so that in fact we may choose $\theta = \tfrac12 - \delta$ for some appropriately chosen small $\delta$ in \eqref{e:hybridbd} when writing the weights $|y_2-y_1|^r$ and $|y_3-y_2|^r$ as lines in the Fourier diagrams. Therefore, we can instead use the fact that
\begin{equation}
	\sup_{a_1,\dots,a_n \in \T^d} \int\limits_{\T^d} \frac{\dd k}{(2\pi)^d} \frac{1}{(1-\hat D(k+a_1))^{1-\theta}} \frac{1}{(1-\hat D(k+a_2))^{1-\theta}}\prod_{i=3}^{n} \frac{1}{1-\hat D(k+a_i)}  < \infty
\end{equation}
whenever $d> (n-2\theta)\twa$ to bound the various terms that appear in the bounding of the diagrams. Choosing $\delta$ sufficiently small thus allows us to bound the case $d \in (6\twa,7\twa]$ as well.\qed

}

\begin{small}
\bibliographystyle{abbrv}
\bibliography{TimsBib}
\end{small}

\end{document}